\documentclass[aos]{imsart}

\RequirePackage{amsfonts,amssymb,amsmath,amsthm,dsfont,graphicx,hyperref,mathrsfs,euscript,enumitem,framed,bm,soul}
\RequirePackage[capitalise]{cleveref}
\usepackage{ifthen}
\newboolean{standalone}
\setboolean{standalone}{true}
\RequirePackage[numbers]{natbib}
\RequirePackage{graphicx}
\RequirePackage{stackengine,scalerel}
\RequirePackage{etoolbox}
\RequirePackage{anyfontsize}
\RequirePackage{thmtools}

\startlocaldefs
\theoremstyle{plain}
\newtheorem{assumption}{Assumption}

\theoremstyle{remark}

\newtheorem{lem}{Lemma}
\newtheorem{thm}{Theorem}
\newtheorem{coro}{Corollary}
\newtheorem{prop}{Proposition}
\theoremstyle{definition}             
\newtheorem{remark}{Remark}

\newtheorem{example}{Example}

\makeatother

\usepackage{defs}

\DeclareFontShape{T1}{ptm}{m}{scit}{<-> ssub * ptm/m/it}{}

\endlocaldefs

\begin{document}

\begin{frontmatter}
  \title{Uniform estimation and inference\\ for nonparametric partitioning-based M-estimators}
  \runtitle{Partitioning-Based M-Estimators}

  \begin{aug}
    \author[A]{\fnms{Matias D.}~\snm{Cattaneo}\ead[label=e1]%
    {cattaneo@princeton.edu}},
    \author[B]{\fnms{Yingjie}~\snm{Feng}\ead[label=e2]%
    {fengyj@sem.tsinghua.edu.cn}}
    \and
    \author[C]{\fnms{Boris}~\snm{Shigida}\ead[label=e3]%
    {bs1624@princeton.edu}}
    \address[A]{Department of Economics,
      Princeton University\printead[presep={,\ }]{e1}}
    \address[B]{School of Economics and Management, Tsinghua University\printead[presep={,\ }]{e2}}
    \address[C]{Department of Operations Research and Financial Engineering,
      Princeton University\printead[presep={,\ }]{e3}}
  \end{aug}

  \received{\smonth{9} \syear{2024}}
  \revised{\smonth{8} \syear{2025}}

  \begin{abstract}
    This paper presents uniform estimation and inference theory for a large class of nonparametric partitioning-based M-estimators. The main theoretical results include:
(i) uniform consistency for convex and non-convex objective functions;
(ii) rate-optimal uniform Bahadur representations;
(iii) rate-optimal uniform (and mean square) convergence rates;
(iv) valid strong approximations and feasible uniform inference methods; and
(v) extensions to functional transformations of underlying estimators.
Uniformity is established over both the evaluation point of the nonparametric functional parameter and a Euclidean parameter indexing the class of loss functions. The results also account explicitly for the smoothness degree of the loss function (if any), and allow for a possibly non-identity (inverse) link function. We illustrate the theoretical and methodological results in four examples: quantile regression, distribution regression, $L_p$ regression, and logistic regression. Many other possibly non-smooth, nonlinear, generalized, robust M-estimation settings are covered by our results. We provide detailed comparisons with the existing literature and demonstrate substantive improvements: we achieve the best (in some cases optimal) known results under improved (in some cases minimal) requirements in terms of regularity conditions and side rate restrictions. The supplemental appendix reports complementary technical results that may be of independent interest, including a novel uniform strong approximation result based on Yurinskii's coupling.

  \end{abstract}



\end{frontmatter}



\section{Introduction}\label{sec: Introduction}

Let $(y_1,\bx_1),(y_2,\bx_2),\cdots,(y_n,\bx_n)$ be independent and identically distributed copies of the random vector $(Y,\bX)\in\mathcal{Y}\times\mathcal{X} \subseteq \mathbb{R}\times\mathbb{R}^d$. Given a loss function $\rho\colon \mathcal{Y}\times \mathcal{E} \times\mathcal{Q} \to \mathbb{R}$ with $\mathcal{E}\subseteq\mathbb{R}$ an open connected set and $\mathcal{Q}\subseteq\mathbb{R}$ a connected compact set, and $\eta\colon \mathbb{R} \to \mathcal{E}$ a strictly monotonic transformation function, consider the functional parameter $\mu_{0}\colon \mathcal{X}\times\mathcal{Q} \to \mathbb{R}$ satisfying
\begin{equation}\label{eq: Rho problem}
  \mu_{0}(\cdot,q) \in \argmin_{\mu\in\mathscr{M}} \E\big[\rho\big(y_{i} , \eta\big(\mu(\bx_{i})\big); q\big)\big],
\end{equation}
where the minimization is over the space of measurable functions from $\mathcal{X}$ to $\mathbb{R}$. In particular, we assume that the (local) minimum is achieved, which is true in most cases. This setup covers settings of interest in nonparametric statistics, econometrics, and data science, including generalized linear models, robust nonlinear regression, and generalized conditional quantile regression. In practice, the parameter of interest may be $\mu_{0}$ itself, or otherwise specific transformations thereof such as $\eta(\mu_0(\cdot, \cdot))$ or its partial derivatives. This paper presents uniform over $\mathcal{X}\times\mathcal{Q}$ estimation and inference results for $\mu_{0}$, and transformations thereof, based on nonparametric partitioning-based $M$-estimation.

\subsection{Partitioning-Based Methodology}

The series (or sieve) nonparametric partitioning-based $M$-estimator is
\begin{equation}\label{eq: M estimator}
    \widehat{\mu}(\bx, q)=\bp(\bx)\trans
    \widehat{\bbeta}(q),\qquad
    \widehat\bbeta(q) \in \argmin_{\bb\in\mathcal{B}} \sum_{i=1}^n\rho\big(y_i,\eta(\bp(\bx_i)\trans \bb);q\big),
\end{equation}
where $\mathcal{B} \subseteq \mathbb{R}^K$ is the feasible set of the optimization problem, and $\bx \mapsto \bp(\bx) = \bp(\bx; \Delta, m) = \big( p_1(\bx; \Delta, m), \ldots, p_K(\bx; \Delta, m) \big)\trans$ is a dictionary of $K$ locally supported basis functions of order $m$ based on a quasi-uniform partition $\Delta = \left\{ \delta_l: 1 \leq l \leq \kappa \right\}$ containing a collection of open disjoint polyhedra in $\mathcal{X}$ such that the closure of their union covers $\mathcal{X}$.

\cite{Gyorfi-etal_2002_book} gives a textbook introduction to the partitioning-based estimation literature. In this nonparametric framework, the \textit{Haar} basis
\begin{align*}
    \bp_{\mathtt{H}}(\bx) = \big( \I(\bx\in \delta_1), \ldots, \I(\bx\in \delta_{\kappa}) \big)\trans
\end{align*}
corresponds to the canonical basis with $m=1$ and $K=\kappa$, and is an essential building block for the construction of other basis functions. Since the Haar basis is ``unconnected'' across cells (i.e., each basis function is supported on a single cell), the resulting estimator in \eqref{eq: M estimator} reduces to $K$ separate M-estimators, each only using observations with $\bx_i\in\delta_k$, for $k=1,\ldots,K$. For estimation and inference, ``small'' cells decrease bias but increase variance, while ``large'' cells have the opposite effect.

A natural generalization is the \textit{piecewise polynomial} basis
\begin{align*}
    \bp_{\mathtt{P}}(\bx) = \bp_{\mathtt{H}}(\bx) \otimes \br_m(\bx),
\end{align*}
where the vector $\br_m(\bx)$ contains the unique terms of an $(m-1)$th-degree polynomial expansion based on $\bx$ and $\otimes$ is the Kronecker product operator, and thus $K = \frac{(m+d-1)!}{(m-1)!d!}\kappa$. The resulting piecewise polynomial fit within each cell gives more flexible approximation, thus decreasing bias, but the estimation approach remains unconnected.

Since the piecewise polynomial estimator may be discontinuous over $\mathcal{X}$, it is sometimes preferred to impose smoothness restrictions across cells: for example, if the partition $\Delta$ admits a tensor product representation with equal number of partitions along the $d$ axes, then the \textit{Splines} basis is
\begin{align*}
    \bp_{\mathtt{S}}(\bx)=\otimes_{k=1}^d \bT_s\bp_{\mathtt{P}}(\be\trans_k\bx),
\end{align*}
where $\be_k$ denotes the $k$-th unit vector $(1\leq k\leq d)$, and $\bT_{s}$ denotes a transformation matrix that ensures the estimator $\bx\mapsto\widehat{\mu}(\bx, q)$ is $(s-1)$-times continuously differentiable ($s<m$) over $\mathcal{X}$, and thus
$K=((m-s)\kappa^{1/d}+s)^d$.
Due to the global smoothness restrictions, the spline basis is no longer unconnected, and the resulting estimator in \eqref{eq: M estimator} cannot be reduced to separate local estimators. Other spline constructions on more general partitioning schemes are available, and compactly supported wavelets are yet another example of a local basis constructed recursively out of the Haar basis. See \cite{Belloni-Chernozhukov-Chetverikov-Kato_2015_JoE}, \cite{Cattaneo-Farrell_2013_JoE}, \cite{Cattaneo-Farrell-Feng_2020_AOS}, and \cite{Chen-Christensen_2015_JOE} for more discussion on these and other bases of approximation. Furthermore, partitioning-based estimation naturally arises in the recursive partitioning literature \citep{Devroye-etal2013_book,Zhang-Singer_2010_Book}.

To enable good statistical performance, we need to restrict the partition of $\mathcal{X}$, and the local basis constructed on it. The first assumption concerns the regularity of the cells in the partition. Let $a_n\lesssim b_n$ denote $\limsup_{n\to\infty} |a_n/b_n|<\infty$.

\begin{assumption}[Quasi-uniform partition]\label{Assumption: gl-quasi-uniform-partition}
    The ratio of the sizes of inscribed and circumscribed balls of each $\delta \in \Delta$ is bounded away from zero uniformly in $\delta \in \Delta$, and
    \begin{equation*}
        \frac{\max \{\operatorname{diam}(\delta): \delta \in \Delta\}}{\min \{\operatorname{diam}(\delta): \delta \in \Delta\}} \lesssim 1
    \end{equation*}
    where $\operatorname{diam}(\delta)$ denotes the diameter of $\delta$. Further, for $h = \max \{\operatorname{diam}(\delta): \delta \in \Delta\}$, assume $h = o(1)$ and $\log (1/h)\lesssim \log n$, as $n\to\infty$.
\end{assumption}

Assumption \ref{Assumption: gl-quasi-uniform-partition} requires the partition $\Delta$ be quasi-uniform: the elements in the partition $\Delta$ do not differ too much in size asymptotically. As a consequence, we can use the maximum diameter $h$ as a universal measure of mesh sizes. The next assumption requires the basis be ``locally supported'', non-collinear, and bounded in a proper sense. A function $p(\cdot)$ on $\mathcal{X}$ is said to be \textit{active} on $\delta\in\Delta$ if it is not identically zero on $\delta$; we also employ standard multi-index notation (see Section \ref{sec: Notation} for details).

\begin{assumption}[Local basis]\label{Assumption: gl-local-basis}\leavevmode$ $
    \begin{enumerate}[label=\emph{(\roman*)}]
        \item For each basis function $p_k$, $k=1, \ldots, K$, the union of elements of $\Delta$ on which $p_k$ is active is a connected set, denoted by $\mathcal{H}_k$. For all $k=1, \ldots, K$, both the number of elements of $\mathcal{H}_k$ and the number of basis functions which are active on $\mathcal{H}_k$ are bounded by a constant.

        \item For any $\ba=\left(a_1, \ldots, a_K\right)\trans \in \mathbb{R}^K$, $a_k^2 h^d \lesssim \ba\trans \int_{\mathcal{H}_k} \bp(\bx) \bp(\bx)\trans d \bx \ \ba$ for $k=1, \ldots, K$.

        \item Let $\multindnorm{\bv} < m$. There exists an integer $\varsigma \in[\multindnorm{\bv}, m)$ such that, for all $\boldsymbol{\varsigma},\multindnorm{\boldsymbol{\varsigma}} \leq \varsigma$,
        \begin{equation*}
            h^{-\multindnorm{\boldsymbol{\varsigma}}}
            \lesssim \inf _{\delta \in \Delta} \inf _{\bx \in \closure(\delta)}\big\|\bp^{(\boldsymbol{\varsigma})}(\bx)\big\|
            \leq \sup _{\delta \in \Delta} \sup _{\bx \in \closure(\delta)}\big\|\bp^{(\boldsymbol{\varsigma})}(\bx)\big\| \lesssim h^{-\multindnorm{\boldsymbol{\varsigma}}},
        \end{equation*}
        where $\closure(\delta)$ is the closure of $\delta$.
        \end{enumerate}
\end{assumption}
In Assumption \ref{Assumption: gl-local-basis}, condition (i) implies that each basis function in $\bp(\bx)$ is supported by a region consisting of a finite number of cells in $\Delta$ (independent of $n$). Then, as $\kappa\to\infty$, all basis functions are locally supported relative to the whole support of the data. Condition (ii) can be read as ``non-collinearity'' of the basis functions in $\bp(\bx)$. Since local support condition has been imposed, it suffices to require the basis functions are not too collinear ``locally''. Condition (iii) controls the magnitude of the local basis in a uniform sense.

Assumptions \ref{Assumption: gl-quasi-uniform-partition} and \ref{Assumption: gl-local-basis} implicitly relate the number of approximating series terms, the number of cells in $\Delta$, and the maximum mesh size: $K\asymp \kappa\asymp h^{-d}$, where $a_n\asymp b_n$ means $a_n\lesssim b_n$ and $b_n\lesssim a_n$. Under appropriate assumptions on the statistical model (Assumption \ref{Assumption: dgp-holder} in Section \ref{sec: Assumptions}), the parameter $m$ will control how well $\mu_0$ can be approximated by linear combinations of the local basis (via Assumption \ref{Assumption: gl-approximation-error} in Section \ref{sec: Assumptions}). We consider large sample approximations where $d$ and $m$ are fixed constants, and $\kappa\to\infty$ (and thus $K\asymp h^{-d}\to\infty$) as $n\to\infty$. As a consequence, appropriate choices of $\Delta$ and $\bp(\cdot)$ will enable valid nonparametric approximations of $\mu_{0}$, and transformations thereof, in large samples.

In practice, the parameter \eqref{eq: Rho problem} and its associated plug-in estimator \eqref{eq: M estimator} may not be unique (e.g., when the objective function is not convex and hence several local minima may exist). In such cases the interpretation of the estimator and its probability limit may depend on the specific (algorithmic) implementation used. This paper does not study these additional complications, but rather assumes that the estimator \eqref{eq: M estimator} has been computed, and then relies on the assumptions in Section \ref{sec: Assumptions} concerning the data generating process to study the large sample statistical properties of the partitioning-based $M$-estimator and transformations thereof.

\subsection{Summary of Contributions}

The objective function in \eqref{eq: M estimator} may not be convex. To address this challenge, we first provide primitive conditions for uniform over $\mathcal{X}\times\mathcal{Q}$ (and mean square) consistency of the partitioning-based estimator $\widehat{\mu}(\bx, q)$, taking explicitly into account whether the loss function $\theta\mapsto\rho(y,\eta(\theta);q)$ is convex: setting $\mathcal{B}=\mathbb{R}^K$ if it is convex, or otherwise $\mathcal{B}=\{\bb\in\mathbb{R}^K:\|\bb\|_\infty\leq R\}$ for some large enough fixed constant $R>0$, we establish $\sup_{q \in \mathcal{Q}}\big\| \widehat{\boldsymbol{\beta}}(q) - \bbeta_0(q) \big\|_{\infty} = o_{\P}(1)$, where $\|\cdot\|_{\infty}$ denotes the $\ell^\infty$-norm, and $\bbeta_0\colon \mathcal{Q} \to \mathbb{R}^K$ denotes coefficients such that $\bbeta_0(q)\trans \bp(\bx)$ approximates $\mu_{0}(\bx,q)$ well enough uniformly over $\mathcal{X}\times\mathcal{Q}$ (Assumption \ref{Assumption: gl-approximation-error} in Section \ref{sec: Assumptions}). For the non-convex case, the resulting ``fixed box'' constrained optimization is arguably a mild assumption in practice, and may be justified in theory under different regularity conditions. These results are presented in Section \ref{sec: consistency}.

Taking the uniform consistency of the partitioning-based estimator as given, and hence being agnostic about the shape of the objective function and other optimization-related aspects, we establish three theoretical results for the partitioning-based series $M$-estimator in \eqref{eq: M estimator}:
\begin{enumerate}[label=(\roman*)]
    \item rate-optimal Bahadur representation uniformly over $\mathcal{X}\times\mathcal{Q}$,
    \item rate-optimal convergence rates in mean square and uniformly over $\mathcal{X}\times\mathcal{Q}$, and
    \item valid strong approximation and feasible distribution theory uniformly over $\mathcal{X}\times\mathcal{Q}$.
\end{enumerate}
These results allow for a large class of possibly non-smooth loss functions. In addition, we precisely characterize how the degrees of smoothness of $\rho$ and $\eta$ affect the order of the remainder in the uniform Bahadur representation for $\widehat{\mu}$, its convergence rates, and the validity of the associated uniform inference procedures. Results (i) and (ii) are presented in Section \ref{sec: Bahadur representation}, while result (iii) is presented in Sections \ref{sec: uniform inference} and \ref{sec: Feasible Uniform Inference}.

Section \ref{sec: Motivating Examples} introduces four examples: \textit{Generalized Conditional Quantile Regression}, \textit{Generalized Conditional Distribution Regression}, \textit{Generalized $L_p$ Regression Estimation}, and \textit{Maximum Likelihood Logistic Regression}. These examples are used to both motivate our high-level assumptions and demonstrate the broad applicability of our uniform estimation and inference results. Our most general results cover other applications such as nonparametric partitioning-based (quasi-maximum likelihood) Poisson regression, censored and truncated regression, as well as Tukey and Huber regression. Section \ref{sec: Assumptions} presents the slightly simplified high-level technical assumptions used throughout the paper, but their most general form is given in the supplemental appendix to streamline the presentation. Section \ref{sec: Examples} demonstrates how our general sufficient conditions are verified for each of our motivating examples.

Section \ref{sec: Extensions} discusses how our results can be extended to cover other parameters of interest, while Section \ref{sec: conclusion} concludes. The supplemental appendix reports simulation evidence, collects all the technical proofs, and presents other theoretical results that may be of independent interest. In particular, our more general theoretical results (i) allow for $\mathcal{Q}$ to be a set of vectors rather than scalars, which can be useful in other examples beyond those studied in this paper; and (ii) consider more complex (VC-type) classes of loss and transformation functions, thereby covering a broader class of settings than those studied herein, but at the cost of additional, cumbersome notation and technicalities. In addition, the supplemental appendix presents new strong approximation results for a class of $K$-dimensional linear stochastic processes indexed by $\mathcal{X}\times\mathcal{Q}$ under standard complexity and smoothness conditions, leveraging a conditional Strassen's theorem \citep{chen2020jackknife,monrad1991nearby} and generalizing prior Yurinskii's coupling results in the literature \citep{Belloni-Chernozhukov-Chetverikov-FernandezVal_2019_JoE,yurinskii1978error}.

\subsection{Related Literature}

Our paper contributes to the literature on nonparametric curve estimation and inference, focusing in particular on series (or sieve) partitioning-based methods. See, for example, \cite{Eggermont-LaRiccia_2009_Book} and \cite{Gyorfi-etal_2002_book} for textbook introductions. This literature is mature and well-developed for the special case of a square loss function $\rho(y,\eta(\theta);q) = (y - \theta)^2$ with identity transformation $\eta(u)=u$, and hence not a function of $q\in\mathcal{Q}$. See, for example, \cite{Belloni-Chernozhukov-Chetverikov-Kato_2015_JoE}, \cite{Cattaneo-Farrell_2013_JoE}, \cite{Cattaneo-Farrell-Feng_2020_AOS}, \cite{Cattaneo-Crump-Farrell-Feng_2024_AER}, \cite{Chen-Christensen_2015_JOE}, \cite{Huang_2003_AoS}, \cite{Zhou-Shen-Wolfe_1998_AoS} for pointwise and uniform over $\mathcal{X}$ estimation and inference results at different levels of generality, and with increasingly weaker technical conditions. This strand of the literature explicitly exploits the special structure, which leads to a closed-form solution of the estimator in \eqref{eq: M estimator}, and hence results are often obtained under minimal assumptions and technical regularity conditions. To be more precise, up to $\polylog(n)$ terms and mild regularity conditions, \cite{Cattaneo-Farrell-Feng_2020_AOS} shows that the minimal requirement $K/n\to0$ is (necessary and) sufficient for rate-optimal convergence rates for any $d\geq1$, and for strong approximations uniformly over $\mathcal{X}$ when $d=1$. They also establish valid strong approximations uniformly over $\mathcal{X}$ for $d>1$ under the requirement $K^3/n\to0$, up to $\polylog(n)$ terms and mild regularity conditions.

Despite aiming for generality, that is, allowing for a large class of loss functions with different levels of smoothness and a non-identity transformation function, this paper establishes rate-optimal uniform over \textit{both} $\mathcal{X}$ and $\mathcal{Q}$ estimation results under the same minimal assumption $K/n\to0$ for unconnected bases, and under the slightly stronger assumption $K^2/n\to0$ for general partitioning-based estimators. Furthermore, we establish valid uniform over \textit{both} $\mathcal{X}$ and $\mathcal{Q}$ inference under the same condition $K^3/n\to0$, leveraging a new strong approximation result given in the supplemental appendix. Compared to the prior literature focusing on the special case of square loss and identity transformation, we are able to achieve the same best known (in some cases rate-optimal) estimation and inference results, under the same (in some cases minimal) side rate restrictions and conditions on the partitioning-based method (Assumptions \ref{Assumption: gl-quasi-uniform-partition} and \ref{Assumption: gl-local-basis}). Furthermore, our results on uniform consistency disentangling convex and non-convex loss functions (Section \ref{sec: consistency}), rate-optimal uniform Bahadur representation and convergence capturing explicitly the smoothness degree of the loss function (Section \ref{sec: Bahadur representation}), and uniform feasible inference validity (Sections \ref{sec: uniform inference} and \ref{sec: Feasible Uniform Inference}), are necessarily new relative to prior work studying the special case of least squares partitioning-based methods.

Going beyond square loss and identity transformation, there are only a handful of results available in the literature. The closest antecedent is \cite{Belloni-Chernozhukov-Chetverikov-FernandezVal_2019_JoE}, who consider nonparametric conditional quantile series regression estimation and inference uniformly over $\mathcal{X}\times\mathcal{Q}$ with $\eta(u)=u$, and under the side rate restriction $K^4/n\to0$, up to $\polylog(n)$ terms, and other regularity conditions. As a comparison, for the special case of nonparametric quantile regression (Example \ref{example: Generalized Conditional Quantile Regression} below), this paper allows for a non-identity (inverse) link function $\eta(\cdot)$, and establishes convergence rates under the minimal condition $K/n\to0$ for piecewise polynomials, and the improved condition $K^2/n\to0$ for connected bases, while for uniform inference we require the weaker condition $K^3/n\to0$, in all cases up to $\polylog(n)$ terms. We also weaken other assumptions imposed in \cite{Belloni-Chernozhukov-Chetverikov-FernandezVal_2019_JoE}: see Example \ref{example: Generalized Conditional Quantile Regression} in Sections \ref{sec: Motivating Examples} and \ref{sec: quantile reg}, and Sections \ref{sec: Comparison with Existing Results -- Bahadur} and \ref{sec: Comparison with Existing Results -- SA}. On the other hand, it is worth noting that \cite{Belloni-Chernozhukov-Chetverikov-FernandezVal_2019_JoE} also considers generic dimension-increasing covariates, while our paper focuses exclusively on partitioning-based local basis.

Our contributions can also be compared to recent work on nonparametric M-estimation employing other smoothing techniques. For example, \cite{Kong-Linton-Xia_2013_ET} considers local polynomial methods, and \cite{Shang-Cheng_2013_AOS} considers smoothing spline methods. Sections \ref{sec: Comparison with Existing Results -- Bahadur} and \ref{sec: Comparison with Existing Results -- SA} give a more detailed comparison with prior literature, and explain precisely how our general results are either on par with or improve upon prior work. In a nutshell, we present estimation and inference results for partitioning-based $M$-estimators that (i) allow for a large class of possibly non-smooth loss functions, (ii) are uniformly valid over both $\mathcal{X}$ and $\mathcal{Q}$, (iii) achieve the best known (in some cases rate-optimal) convergence rates, and (iv) require substantially weaker (in some cases minimal) side rate restrictions and regularity conditions. Our results, for example, permit the use of Haar $M$-estimators for estimation and inference, which were ruled out by prior work.

\subsection{Notation}\label{sec: Notation}

We employ standard notation in probability, statistics and empirical process theory \citep{bhatia2013matrix,dudley2014uniform,Kallenberg2021,vaart96_weak_conver_empir_proces}.
For any vector $\ba=(a_1, \cdots, a_M)\in\mathbb{R}^M$, we write $\|\ba\|=(\sum_{j=1}^Ma_j^2)^{1/2}$ and $\|\ba\|_\infty=\max_{1\leq j\leq M} |a_j|$.
For any real function $f$ depending on $d$ variables $(t_1, \ldots, t_d)$ and any vector $\bv=(v_1, \cdots, v_d)$ of nonnegative integers, denote
$f^{(\bv)} = \frac{\partial^{\multindnorm{\bv}}}{\partial t_1^{v_1} \ldots \partial t_d^{v_d}} f$ where $\multindnorm{\bv} = \sum_{k = 1}^d v_j$.
For functions that depend on $(\bx,q)$, the multi-index derivative notation is taken with respect to the first argument $\bx$, unless otherwise noted.
We say a function $f$ is $\alpha$-H\"{o}lder on a set $\mathcal{I}$ if for some constant $C>0$ and $\alpha>0$, $|f(\bx_1)-f(\bx_2)|\leq C\|\bx_1-\bx_2\|^{\alpha}$ for any $\bx_1,\bx_2\in\mathcal{I}$.
For any two numbers $a$ and $b$, $a\vee b=\max\{a,b\}$, and $a\wedge b=\min\{a,b\}$. Let $\E_n[g(x_i)]=\frac{1}{n}\sum_{i=1}^ng(x_i)$ and
$\mathbb{G}_n[g(x_i)]=\frac{1}{\sqrt{n}}\sum_{i=1}^n(g(x_i)-\E[g(x_i)])$.
For sequences, $a_n=O(b_n)$ or $a_n\lesssim b_n$ denotes $\limsup_n |a_n/b_n|$ is finite, $a_n=O_\P(b_n)$ denotes $\limsup_{\epsilon\to\infty} \limsup_{n\to\infty} \P[|a_n/b_n| \geq \epsilon ] = 0$, $a_n = o(b_n)$ denotes $a_n/b_n\to 0$, and $a_n = o_\P(b_n)$ denotes $a_n/b_n\to_\P 0$, where $\to_\P$ is convergence in probability. Limits are taken as $n \to \infty$, and the dependence on $n$ is often suppressed, e.\,g. $K = K_n$. Also, we say a random variable $\xi$ is sub-Gaussian conditional on $\bX$ if for some constant $\sigma^2>0$, $\P(|\xi|\geq t|\bX=\bx)\leq 2\exp(-t^2/\sigma^2)$ for all $t\geq 0$ and $\bx\in\mathcal{X}$.

\section{Examples}\label{sec: Motivating Examples}

We discuss four motivating examples of interest covered by our theoretical results. Section \ref{sec: Examples} demonstrates how our high-level assumptions, introduced in Section \ref{sec: Assumptions}, are verified for these examples in order to obtain uniform estimation and inference results; the supplemental appendix collects omitted details.

Our first example generalizes the work of \cite{Belloni-Chernozhukov-Chetverikov-FernandezVal_2019_JoE}, who studies the large sample properties of nonparametric conditional quantile series regression with $\eta(u)=u$. We allow for non-identity transformation under substantially weaker technical conditions.

\begin{example}[Generalized Conditional Quantile Regression]\label{example: Generalized Conditional Quantile Regression}
    The quantile regression loss function is
    \[ \rho(y,\eta;q)=(q-\mathds{1}(y<\eta))(y-\eta),\]
    where $q\in\mathcal{Q}=[\varepsilon_0,1-\varepsilon_0]$ denotes the quantile position, with $\varepsilon_0>0$. Then, $\eta(\mu_0(\bx,q))$ is the $q$-th conditional quantile function of $Y$ given $\bX=\bx$, and the partitioning-based quantile regression estimator is $\eta(\widehat{\mu}(\bx,q))$ as defined in \eqref{eq: M estimator}. In the classical case, $\eta(\cdot)$ is the identity function, but our theory accommodates other transformations. Interest lies in the quantile process estimator $(\eta(\widehat{\mu}(\bx,q)):(\bx,q)\in\mathcal{X}\times\mathcal{Q})$, which can be used to characterize heterogeneous effects of covariates on the outcome distribution and to conduct specification testing. See Section \ref{sec: quantile reg} for our main results, and  Section \ref{sec: Extensions} for results on transformations.
\end{example}

\cite{chernozhukov-et-al_2013_ECMA} obtains large sample estimation and inference results for parametric ($K$ fixed) generalized conditional distribution regression, and applies them to counterfactual analysis and causal inference. The following example discusses a nonparametric partitioning-based generalized conditional distribution regression estimator.

\begin{example}[Generalized Conditional Distribution Regression]\label{example: Conditional Distribution Regression}
    Non-linear least squares conditional distribution regression employs
    \[\rho(y,\eta;q)= ( \mathds{1}(y \leq q) - \eta)^2,\]
    where we can, for example, use the complementary log-log link $\eta(\theta)=1-\exp(-\exp(\theta))$. Estimands of interest are  $\eta(\mu_0(\bx,q))$, which corresponds to the conditional distribution function of $Y$ given $\bX=\bx$ (i.e., $F_{Y|\bX}(q|\bx)=\E[\mathds{1}(Y \leq q)|\bX=\bx]$), and derivatives thereof. Uniform estimation and inference results based on $(\eta(\widehat{\mu}(\bx,q)):(\bx,q)\in\mathcal{X}\times\mathcal{Q})$ are useful for a variety of purposes, including heterogeneous treatment effect estimation and specification testing. See Section \ref{sec: distribution reg} for our main results, and Section \ref{sec: Extensions} for transformations.
\end{example}

\cite{Lai-Lee_JASA_2005} studies $L_p$ regression estimation with identity transformation $\eta(\cdot)$ in a parametric setting ($K$ fixed). The next example considers a class of nonparametric partitioning-based generalized $L_p$ regression with $p\in[1,2]$, covering the full interpolation between nonparametric generalized median regression ($p=1$) and nonparametric nonlinear least squares regression ($p=2$), with possibly non-identity $\eta(\cdot)$.

\begin{example}[Generalized \texorpdfstring{$L_p$}{Lp} Regression]\label{example: $L_p$ regression}
    The (possibly nonlinear) $L_p$ regression estimator is defined by taking
    \[\rho(y, \eta)=|y-\eta|^{p},\]
    for a fixed $p>0$. In particular, $p=2$ leads to nonlinear least squares, and $p=1$ leads to generalized least absolute deviations, when a non-identity transformation function $\eta(\cdot)$ is used. The estimand of interest is usually the transformed regression function $\eta(\mu_0(\bx))$, which needs to be interpreted in context. Our general theory applies to any choice $p\in[1,2]$, delivering uniform estimation and inference methods based on $(\widehat{\mu}(\bx):\bx\in\mathcal{X})$, and transformations thereof. See Section \ref{sec: lp reg} for our main results, and Section \ref{sec: Extensions} for transformations.
\end{example}

The final example considers (nonparametric) Generalized Linear Models \citep{Mccullagh_2019_book}. For specificity, we focus on (quasi-)maximum likelihood logistic regression, but our results cover many other examples within this class such as regression models with limited dependent variables (e.g., Poisson, fractional, censored and truncated regression).

\begin{example}[Generalized Linear Models]\label{example: Generalized linear models}
    The classical logistic regression model, or binary classification with sigmoid (inverse) link, employs
    \[\rho(y,\eta)=-y\log \eta-(1-y)\log (1-\eta),\qquad \eta(\theta)=\exp(\theta)/(1+\exp(\theta)),\]
    with $\mathcal{Y}=\{0,1\}$. The estimand $\eta(\mu_0(\bx))$ characterizes the conditional probability of $Y=1$ given $\bX=\bx$. See Section \ref{sec: logit est} for uniform estimation and inference methods based on $(\widehat{\mu}(\bx):\bx\in\mathcal{X})$, and Section \ref{sec: Extensions} for transformations thereof. Furthermore, our results also cover other related quasi-maximum likelihood (and non-linear least squares) problems such as fractional regression where $\mathcal{Y}=[0,1]$.
\end{example}

The four examples illustrate distinct settings from a technical perspective. In Example \ref{example: Generalized Conditional Quantile Regression} uniformity over $\mathcal{X}\times\mathcal{Q}$ is of interest, and the loss function is non-smooth as a function of $\bx\in\mathcal{X}$ but smooth as a function of $q\in\mathcal{Q}$. Example \ref{example: Conditional Distribution Regression} is the antithesis of Example \ref{example: Generalized Conditional Quantile Regression} because now the loss function is smooth as a function of $\bx\in\mathcal{X}$ and non-smooth as a function of $q\in\mathcal{Q}$, while uniformity over $\mathcal{X}\times\mathcal{Q}$ is still of interest. In Example \ref{example: $L_p$ regression} only uniformity over $\mathcal{X}$ is of interest because $q\in\mathcal{Q}$ is not present in the loss function, but its smoothness depends on $p\in[1,2]$; the a.\,e. derivative of $\eta\mapsto\rho(y,\eta)$ ranges from discontinuous ($p=1$), to H\"{o}lder continuous ($p\in(1,2)$), to linear ($p=2$). Likewise, Example \ref{example: Generalized linear models} only involves uniformity over $\mathcal{X}$ because $q\in\mathcal{Q}$ is not present in the loss function, but now the loss function is smooth and well-behaved; this last example serves as a benchmark for our theoretical development. All of the examples above have a convex loss function when $\eta(u)=u$, but can be non-convex when $\eta(\cdot)$ is not the identity function.

Our theoretical results cover other examples. For instance, Tukey and Huber regressions are popular methods in robust statistics, and our theory allows for their generalizations to nonparametric partitioning-based uniform estimation and inference. Specifically, Tukey regression employs the loss function $\rho(y,\eta;q) = q^2 (1 - [1 - (y-\eta)^2/q^2]^3) \I(|y-\eta|\leq q) + q^2 \I(|y-\eta|>q)$, while Huber regression uses the loss function $\rho(y,\eta;q) = (y-\eta)^2\I(|y-\eta|\leq q)+ q(2|y-\eta|-q)\I(|y-\eta|>q)$, where $q$ is treated as a tuning parameter that balances the robustness and the bias of the estimation. We do not discuss these and other examples for brevity.

\section{Assumptions}\label{sec: Assumptions}

Our theoretical work proceeds under Assumptions \ref{Assumption: gl-quasi-uniform-partition} and \ref{Assumption: gl-local-basis} on the partitioning-based estimation framework, three assumptions concerning the data generating process and the loss function, and a final assumption linking the statistical model and partition-based approximation.

\begin{assumption}[Data Generating Process]

\label{Assumption: dgp-holder}\leavevmode
  \begin{enumerate}[label=\emph{(\roman*)}]
  \item $((y_i, \bx_i) \colon 1 \leq i \leq n)$ is a random sample satisfying~\eqref{eq: Rho problem}.

  \item The distribution of $\bx_i$ admits a Lebesgue density $f_{X}(\cdot)$ which is continuous and bounded away from zero on support $\mathcal{X} \subset \mathbb{R}^d$, where $\mathcal{X}$ is the closure of an open, connected and bounded set.

  \item The conditional distribution of $y_{i}$ given $\bx_i$ admits a conditional density $f_{Y|X}(y | \bx)$ with support $\mathcal{Y}_{\bx}$ with respect to some (sigma-finite) measure $\mathfrak{M}$, and $\underset{\bx \in \mathcal{X}}{\sup}\, \underset{y \in \mathcal{Y}_{\bx}}{\sup}\, f_{Y|X}(y \mid \bx) < \infty$.

  \item $\bx\mapsto\mu_{0}(\bx,q)$ is $m \geq 1$ times continuously differentiable for every $q\in\mathcal{Q}$, $\bx\mapsto\mu_{0}(\bx,q)$ and its derivatives of order no greater than $m$ are bounded uniformly over $(\bx,q) \in \mathcal{X}\times\mathcal{Q}$, and
  \begin{align*}
      \sup_{\bx\in\mathcal{X}}\sup_{q_1\neq q_2} \frac{|\mu_{0}(\bx,q_1) - \mu_{0}(\bx,q_2)|}{|q_1 - q_2|}
  \lesssim 1.
  \end{align*}
  \end{enumerate}
\end{assumption}

Assumption \ref{Assumption: dgp-holder} imposes standard conditions from the nonparametric regression literature, including basic support and smoothness restrictions. Minimal additional regularity is imposed to accommodate uniformity over $q\in\mathcal{Q}$, and different types of conditional distributions of $Y|\bX$ (e.\,g., absolutely continuous, discrete or mixed) are allowed. We consider continuously distributed covariates for simplicity, but with additional notation, and by appropriate modification of our assumptions and proof, it is possible to accommodate $\bx_i$ with continuous and discrete components.

The next assumption requires regularity conditions on the loss and transformation functions. Define $B_{q}(\bx)=\{\zeta: |\zeta-\mu_0(\bx,q)|\leq r\}$ for some fixed (small enough) constant $r>0$, which is a ``ball'' around the true value $\mu_0(\bx,q)$ with radius $r$.

\begin{assumption}[Loss Function]\label{Assumption: loss function}\leavevmode
    \begin{enumerate}[label=\emph{(\roman*)}]
    \item Let $\mathcal{Q}\subseteq\mathbb{R}$ be a connected compact set. For each $q \in \mathcal{Q}$, $y \in \mathcal{Y}$, and some open connected subset $\mathcal{E}$ of $\mathbb{R}$ not depending on $y$, $\eta\mapsto\rho(y,\eta;q)$ is absolutely continuous on closed bounded intervals within $\mathcal{E}$, and admits an a.\,e. derivative $\psi(y,\eta; q)$.

    \item The first-order optimality condition $\E[\psi(y_i, \eta(\mu_{0}(\bx_i,q));q)| \bx_i] =0$ holds; the function $\E[\psi(y_i, \eta(\mu_{0}(\bx_i,q));q)^{2} | \bx_i = \bx]$ is continuous in both arguments $(\bx,q) \in \mathcal{X}\times\mathcal{Q}$, bounded away from zero, and Lipschitz in $q$ uniformly in $\bx$; there is a positive measurable envelope function $\overline{\psi}(\bx_i, y_i)$ such that $\sup_{q\in\mathcal{Q}} |\psi(y_i, \eta(\mu_{0}(\bx_i,q));q)| \leq \overline{\psi}(\bx_i, y_i)$ with $\sup_{\bx \in \mathcal{X}}\E[\overline{\psi}(\bx_i, y_i)^{\nu} | \bx_i = \bx] < \infty$ for some $\nu > 2$.\sloppy

    \item\label{shape of derivative} $\eta(\cdot)$ is strictly monotonic and twice continuously differentiable. Furthermore, for some fixed constant $\alpha \in (0, 1]$, for any $(\bx,q) \in \mathcal{X}\times\mathcal{Q}$, and a pair of points $\zeta_1, \zeta_2 \in B_{q}(\bx)$, $\psi(\cdot)$ satisfies the following (constants hidden in $\lesssim$ do not depend on $\bx$, $q$, $\zeta_1$, $\zeta_2$):
    \begin{itemize}
        \item if $\mathfrak{M}$ is the Lebesgue measure, then
        \begin{align*}
            &\sup_{\lambda \in [0, 1]} \sup_{y \not\in [\eta(\zeta_1)\wedge \eta(\zeta_2), \eta(\zeta_1)\vee \eta(\zeta_2)]}\; | \psi ( y, \eta(\zeta_1 + \lambda(\zeta_2 - \zeta_1));q) -
            \psi ( y, \eta(\zeta_2);q) | \lesssim | \zeta_1 - \zeta_2 |^{\alpha},\\
            &\sup_{\lambda \in [0, 1]} \sup_{y \in [\eta(\zeta_1)\wedge \eta(\zeta_2), \eta(\zeta_1)\vee \eta(\zeta_2)]}\, | \psi ( y, \eta(\zeta_1 + \lambda(\zeta_2 - \zeta_1));q ) - \psi ( y, \eta(\zeta_2); q) | \lesssim 1;
        \end{align*}
        \item if $\mathfrak{M}$ is \emph{not} the Lebesgue measure, then
        \begin{align*}
           \sup_{\lambda \in [0, 1]} \sup_{y\in\mathcal{Y}}| \psi(y,\eta(\zeta_1 + \lambda(\zeta_2 - \zeta_1));q) - \psi(y,\eta(\zeta_2);q)|\lesssim |\zeta_1-\zeta_2|^\alpha.
        \end{align*}
      \end{itemize}

    \item $\Psi(\bx,\eta;q)= \E[\psi(y_{i}, \eta;q) | \bx_i=\bx]$ is twice continuously differentiable with respect to $\eta$,
    \begin{align*}
        \sup_{\bx \in \mathcal{X}, q \in \mathcal{Q}} \sup_{\zeta \in B_{q}(\bx)} |\Psi_k(\bx, \eta(\zeta);q)| < \infty, \quad \Psi_k(\bx, \eta; q)= \frac{\partial^k}{\partial \eta^k} \Psi(\bx, \eta; q),\quad k=1,2,
    \end{align*}
    and
    \begin{align*}
        &\inf_{\bx \in \mathcal{X}, q \in \mathcal{Q}} \inf_{\zeta \in B_{q}(\bx)} \Psi_1(\bx, \eta(\zeta);q)(\eta^{(1)}(\zeta))^{2} > 0.
    \end{align*}

    \end{enumerate}
\end{assumption}

Assumption \ref{Assumption: loss function} is carefully crafted to accommodate all the examples discussed in Section \ref{sec: Motivating Examples}, and many others. Part (i) allows for different degrees of smoothness in the loss function, assuming only absolute continuity (with respect to $\eta$). It also makes clear that $q$ is scalar, which is assumed only to simplify the notation; see the supplemental appendix for the general case $\mathcal{Q}\subseteq\mathbb{R}^{d_\mathcal{Q}}$ with $d_\mathcal{Q}\geq1$. Part (ii) formalizes the idea that $\mu_0(\bx,q)$ may not be a unique (global) minimizer in \eqref{eq: Rho problem}, and consequently it is only required to be a root of the (conditional) first-order condition; the rest of the assumptions in that part are mild regularity conditions. In some applications, $\mu_0(\bx,q)$ can be the unique minimizer; see, for example, \cite{Lin-Kulasekera_2007_Biometrika}, \cite{matzkin2007nonparametric}, and references therein.

Part (iii) of Assumption \ref{Assumption: loss function} imposes additional structure on the a.\,e. first derivative of the loss function, allowing for all types of outcome data (discrete, mixed, and continuous) and rescalings emerging in some of the motivating examples. Importantly, this part characterizes precisely the role of (H\"older) smoothness, which is controlled by the parameter $\alpha\in(0,1]$. We illustrate the full power of this general assumption in Section \ref{sec: Examples}, where $\alpha=1$ in Examples \ref{example: Generalized Conditional Quantile Regression} and \ref{example: Conditional Distribution Regression}, $\alpha=p-1$ in Example \ref{example: $L_p$ regression} when $p>1$, and $\alpha=1$ in Example \ref{example: Generalized linear models}. The strict monotonicity condition on $\eta(\cdot)$ is satisfied by usual (inverse) link functions used in generalized linear models. Finally, part (iv) of Assumption \ref{Assumption: loss function} collects mild regularity conditions on the smoothed-out a.\,e. derivative of the loss function.

Assumptions \ref{Assumption: dgp-holder} and \ref{Assumption: loss function} have restricted basic aspects of the statistical model, imposing standard support, moment, and smoothness conditions, in addition to other minimal structure required on the loss and transformation functions. These conditions are sufficient for pointwise estimation and inference, but more is needed for uniform over $\mathcal{X}\times\mathcal{Q}$ results. In the supplemental appendix, our theoretical results are established under one more condition that governs the complexity of the loss function and related function classes. To avoid a long list of complexity bounds, we present a more restrictive but simpler assumption motivated by the examples discussed in Section \ref{sec: Motivating Examples}: we consider a loss function $\rho(y, \eta; q)$ that can be expressed as a linear combination of certain ``simple'' functions. See Section \saref{sec:simpler-generalization-of-examples} for omitted details.

\begin{assumption}[Simplified Setup]\label{Assumption: simple setup} $ $

    \begin{enumerate}[label=\emph{(\roman*)}]
        \item $q\mapsto\mu_0(\bx, q)$ is non-decreasing.
        \item $\rho(y,\eta;q)=\sum_{j=1}^4\rho_j(y,\eta;q)$, where the functions $\rho_j(\cdot)$ are of the following types:
        \begin{itemize}
        	\item Type I:  $\rho_1(y,\eta;q)=(f_1(y)+D_1\eta)\I(y\leq \eta)$,
        	\item Type II: $\rho_2(y,\eta;q)=(f_2(y)+D_2\eta)\I(y\leq q)$,
        	\item Type III: $\rho_3(y,\eta;q)=(f_3(y)+D_3\eta)q$,
        	\item Type IV: $\rho_4(y,\eta;q)=\mathcal{T}(y,\eta)$,
        \end{itemize}
        where $f_{j}(\cdot)$ are fixed continuous functions, $D_j$ are universal constants, and $\mathcal{T}(\cdot)$ is a fixed continuous function.

        \item $\eta\mapsto\E[\tau(y_i,\eta)|\bx_i=\bx]$ is differentiable, where $\tau(y,\eta) = \frac{\partial}{\partial \eta}\mathcal{T}(y,\eta)$, $\tau(y,\eta)$ and $\frac{\partial}{\partial \eta}\E[\tau(y_i,\eta)|\bx_i=\bx]$ are continuous in their arguments and $\alpha$-H\"{o}lder continuous ($\alpha\in(0,1]$) in $\eta$ for $\eta$ in any fixed compact subset of $\mathcal{E}$ with the H\"{o}lder constants independent of $(y,\bx)$.

        \item $\sup_{q\in\mathcal{Q}} |\tau(y,\eta(\mu_{0}(\bx,q)))|\leq \bar\tau(\bx,y)$ with $\sup_{\bx \in \mathcal{X}} \E[\bar{\tau}(\bx_i, y_i)^{\nu} | \bx_i = \bx] < \infty$ for some $\nu > 2$.

        \item If $D_1 \neq 0$, then $F_{Y|X}$ is differentiable with a Lebesgue density $f_{Y|X}$, and $f_{Y|X}$ is continuous in both arguments and $y\mapsto f_{Y|X}(y|\bx)$ is continuously differentiable.
    \end{enumerate}
\end{assumption}

This third assumption imposes an additional weak monotonicity condition on $\mu_0(\bx, q)$ as a function of $q$, which is compatible with all the examples in Section \ref{sec: Motivating Examples}. Finally, the key restriction emerging from Assumption \ref{Assumption: simple setup} is on the structure of the loss function, which allows for linear combinations of smooth loss functions of $y$ and $\eta$, and non-smooth loss functions involving indicator functions of either $y$ and $\eta$, or $y$ and $q$. These restrictions are still general enough to cover the four motivating examples: the loss function in Example \ref{example: Generalized Conditional Quantile Regression} is a combination of Type I and Type III functions with $f_1$ and $f_3$ being linear functions of $y$; the loss function in Example \ref{example: Conditional Distribution Regression} is a combination of Type II and Type IV functions; the loss function in Example \ref{example: $L_p$ regression} is of Type IV for $p>1$, and of the same type as Example \ref{example: Generalized Conditional Quantile Regression} when $p=1$ (median regression); and the loss function in Example \ref{example: Generalized linear models} is usually a Type IV function. See Section \ref{sec: Examples} for details.

Our final assumption concerns the approximation power of the basis $\bp(\cdot)$ in connection with the underlying functional parameter.

\begin{assumption}[Approximation Error]\label{Assumption: gl-approximation-error}
  There exists a vector of coefficients $\bbeta_0(q) \in \mathbb{R}^K$ such that for all $\boldsymbol{\varsigma}$ satisfying $\multindnorm{\boldsymbol{\varsigma}} \leq \varsigma$ in Assumption~\ref{Assumption: gl-local-basis},
    \begin{equation*}
        \sup_{q\in\mathcal{Q}}\sup_{\bx \in \mathcal{X}}
        \big| \mu_0^{(\boldsymbol{\varsigma})}(\bx,q) -
        \bbeta_0(q) \trans \bp^{(\boldsymbol{\varsigma})}(\bx) \big|
        \lesssim h^{m - \multindnorm{\boldsymbol{\varsigma}}}.
    \end{equation*}
\end{assumption}
The vector $\bbeta_0(q)$ can be viewed as a pseudo-true value, and does not have to be unique. The existence of such $\bbeta_0(q)$ can be established using approximation theory or related methods, and necessarily depends on the specific underlying structure of the statistical model (determining $\mu_0(\bx,q)$) and the partitioning-based method (determining $\bp^{(\boldsymbol{\varsigma})}(\bx)$). For standard local bases, the assumption can be verified by imposing smoothness conditions on $(\bx,q)\mapsto\mu_0^{(\boldsymbol{\varsigma})}(\bx,q)$ (see Assumption \ref{Assumption: dgp-holder}(iv)). For more discussion, see \cite{Belloni-Chernozhukov-Chetverikov-Kato_2015_JoE}, \cite{Belloni-Chernozhukov-Chetverikov-FernandezVal_2019_JoE},  \cite{Cattaneo-Farrell_2013_JoE}, \cite{Cattaneo-Farrell-Feng_2020_AOS}, \cite{Chen-Christensen_2015_JOE}, \cite{Huang_2003_AoS}, and references therein.

\section{Consistency}\label{sec: consistency}

We show that the partitioning-based $M$-estimator is consistent, which is the starting point for establishing its main point estimation and inference asymptotic properties. We endeavor to impose the weakest possible conditions, which requires careful consideration of the specific shape of the loss function in \eqref{eq: M estimator}: we consider two cases, either the loss function $\rho(y,\eta(\theta);q)$ is convex with respect to $\theta$ or not; in the latter case, we will have to restrict the feasibility region $\mathcal{B}$.

\subsection{Convex Loss Function}

For the case of convex $\theta\mapsto\rho(y,\eta(\theta);q)$, consistency can be established for general unconstrained estimators ($\mathcal{B}=\mathbb{R}^K$ in \eqref{eq: M estimator}) under mild conditions. The proof is deferred to the supplemental appendix (Lemma~\saref{lem:holder-consistency}).

\begin{lem}[Consistency, convex case]\label{lem:consistency, convex}
  Suppose that Assumptions \ref{Assumption: gl-quasi-uniform-partition}--\ref{Assumption: gl-approximation-error} hold,
  $\rho(y, \eta(\theta); q)$ is convex with respect to $\theta$ with left or right derivative $\psi(y, \eta(\theta); q) \eta^{(1)}(\theta)$, $\mathcal{B}=\mathbb{R}^K$ in \eqref{eq: M estimator}, and $m > d / 2$. Furthermore, assume that one of the following two conditions holds:
  \begin{enumerate}[label=\emph{(\roman*)}]
    \item $\frac{( \log n )^{\frac{\nu}{\nu - 1}}}{n h^{\frac{2 \nu}{\nu - 1} d}} = o(1)$, or
    \item $\frac{(\log n)^{3/2}}{nh^{2d}}=o(1)$ and $\overline\psi(\bx_i,y_i)$ is sub-Gaussian conditional on $\bx_i$.
  \end{enumerate}

  Then
  \begin{align}
    \sup_{q \in \mathcal{Q}}\; \big\| \widehat{\bbeta}(q) - \bbeta_0(q) \big\| &= o_{\P}(1),\label{eq: consistency, beta}\\
    \sup_{\bx \in \mathcal{X}} \sup_{q \in \mathcal{Q}}\; \big|\widehat{\mu}^{(\bv)}(\bx,q) - \mu_{0}^{(\bv)}(\bx,q) \big| &= o_{\P}(h^{- \multindnorm{\bv}}),\label{eq: consistency, sup}\\
    \sup_{q \in \mathcal{Q}}\; \int \big(\widehat\mu^{(\bv)}(\bx,q) - \mu_{0}^{(\bv)}(\bx,q) \big)^2 f_X(\bx)\, \mathrm{d} \bx &= o_{\P} ( h^{d - 2\multindnorm{\bv}} )\label{eq: consistency, l2}.
  \end{align}
\end{lem}

Lemma~\ref{lem:consistency, convex} shows that the function estimator $\widehat\mu$ is uniform-in-$q$ consistent for the true value $\mu_0$ in both $L_2$-norm and $\sup$-norm over $\mathcal{X}$, whereas for the derivative estimator $\widehat\mu^{(\bv)}$ (with $|\bv|>0$) the lemma only provides a bound on its deviation from the estimand $\mu_0^{(\bv)}$. Technically, all we need from this lemma to establish the Bahadur representation later is the uniform-in-$q$ consistency of the coefficients estimator $\widehat\bbeta(q)$ for the pseudo-true coefficients $\bbeta_0(q)$ in $\sup$-norm, i.e., $\|\widehat\bbeta(q)-\bbeta_0(q)\|_\infty=o_\P(1)$, which is immediate from \eqref{eq: consistency, beta}, the uniform-in-$q$ consistency in the Euclidean norm.

Two kinds of rate restrictions are imposed in Lemma \ref{lem:consistency, convex}, depending on the moment condition assumed for the generalized residual $\psi(y_i,\eta(\mu_0(\bx_i,q));q)$. In the best case when the residual has a sub-Gaussian envelope, we need $1/(nh^{2d}) \asymp K^2/n = o(1)$, up to $\polylog(n)$ terms, while in the worst case when the envelope of $\psi(y_i,\eta(\mu_0(\bx_i,q));q)$ has a bounded $\nu$-th moment with $\nu$ close to $2$, we roughly need $1/(nh^{4d}) \asymp K^4/n = o(1)$, up to $\polylog(n)$ terms.

A feature of Lemma \ref{lem:consistency, convex} is that \textit{no} constraints are imposed on the coefficients in the optimization procedure, which allows the estimation space to be, for example, piecewise polynomials. In contrast, many studies of series (or sieve) methods restrict the functions in the estimation space to satisfy certain smoothness conditions, e.g., Lipschitz continuity, to derive the uniform consistency \citep[e.g.,][]{chernozhukov-imbens-newey_2007_JoE}.

\subsection{Non-Convex Loss Function}

Consider the case when the loss $\rho(y,\eta(\theta); q)$ is possibly non-convex with respect to $\theta$. This setting naturally arises, for example, in nonlinear regression when $\rho(y,\eta(\theta);q) = (y - \eta(\theta))^2$ with $\eta(\cdot)$ non-identity: while $\eta\mapsto\rho(y,\eta;q)$ is a square loss function, hence convex, introducing a transformation function $\eta$ such as the (inverse) logistic link will often make $\theta\mapsto\rho(y,\eta(\theta);q)$ non-convex.

A proof of consistency for the unconstrained estimator in \eqref{eq: M estimator} with a non-convex loss function is not available, but we are able to establish consistency of a regularized $M$-estimator. Specifically, we add a \textit{fixed} ``box'' constraint: for some fixed constant $R>0$,
\begin{equation*}
    \widehat{\bbeta}(q) \in \argmin_{\| \bb \|_{\infty} \leq R}\; \sum_{i=1}^n
    \rho ( y_i, \eta(\bp(\bx_i) \trans \bb);q ).
\end{equation*}
In the supplemental appendix we show that the pseudo-true coefficients $\bbeta_0(q)$ from Assumption \ref{Assumption: gl-approximation-error} are bounded in $\sup$-norm by a universal constant: $\sup_{q\in\mathcal{Q}}\|\bbeta_0(q)\|_\infty \lesssim 1$ (because $\bp(\bx)\trans \bbeta_0(q)$ has to be close to $\mu_0(\bx, q)$ which is uniformly bounded). Therefore, we can always choose a sufficiently large constant $R$ in the optimization procedure, making the box constraint set contain $\bbeta_0(q)$ as an interior point. The following lemma, proven in the supplemental appendix (Lemma~\saref{lem:consistency-nonconvex}), establishes consistency of the constrained estimator.

\begin{lem}[Consistency, non-convex case]\label{lem:consistency, nonconvex}
  Suppose that Assumptions \ref{Assumption: gl-quasi-uniform-partition}--\ref{Assumption: gl-approximation-error} hold, $\mathcal{B}=\{\bb\in\mathbb{R}^K: \|\bb\|_\infty\leq R\}$ with $R\geq 2\sup_{q\in\mathcal{Q}}\|\bbeta_0(q)\|_\infty$ in \eqref{eq: M estimator},
  $m > d / 2$, and that there exists some constant $c>0$ such that $\inf \Psi_{1}(\bx,\zeta;q)>c$, where the infimum is over $\bx \in \mathcal{X}$, $q \in \mathcal{Q}$, $\zeta$ between $\eta(\bp(\bx)\trans\bbeta)$ and $\eta( \mu_0(\bx, q))$,
  and $\bbeta \in \mathcal{B}$.
  Furthermore, assume one of the following two conditions holds:
  \begin{enumerate}[label=\emph{(\roman*)}]
      \item $\frac{\left( \log n \right)^{\frac{\nu}{\nu - 1}}}{n h^{\frac{2 \nu}{\nu - 1} d}} = o(1)$, or
      \item $\frac{(\log n)^{3/2}}{nh^{2d}}=o(1)$ and $\overline\psi(\bx_i,y_i)$ is sub-Gaussian conditional on $\bx_i$.
  \end{enumerate}

  Then \eqref{eq: consistency, beta}, \eqref{eq: consistency, sup}, and \eqref{eq: consistency, l2} hold.
\end{lem}

Compared to Lemma \ref{lem:consistency, convex}, two additional restrictions are imposed in this lemma. The first one, $R\geq 2\sup_{q\in\mathcal{Q}}\|\bbeta_0(q)\|_\infty$, can be theoretically justified by Lemma~\saref{lem:assumptions-of-nonconvex-consistency-are-reasonable} in the supplemental appendix, and in practice a large enough $R$ is recommended. The other restriction concerns a lower bound for $\Psi_1$, which implies that the (population) loss function is strongly convex in a neighborhood of the true value $\eta(\mu_0(\bx,q))$, making the (constrained) minimizer well defined. (This condition does \textit{not\/} break because of the shape of $\eta(\cdot)$, in contrast with the convexity of $\rho(y, \eta(\theta); q)$ in $\theta$.) The other conditions in this lemma are the same as those in the convex case, and thus all improvements discussed before also apply to this second consistency result.

\subsection{Weaker Conditions for Special Cases}\label{sec: Weaker Conditions for Special Cases}

In the supplemental appendix we provide additional consistency results for two special cases:
\begin{itemize}
    \item the loss function is strongly convex and smooth (i.e.,  the second ``derivative'' of $\rho(y, \eta(\cdot); q)$ is bounded and bounded away from zero), or
    \item an unconnected basis (i.e., each basis function is supported on a single cell of $\Delta$) is employed.
\end{itemize}

The first case covers the usual square loss function with identity transformation. The second case covers partitioning-based $M$-estimation using the (Haar and) piecewise polynomial basis. Notably, in these two special cases, the consistency result $\|\widehat\bbeta(q)-\bbeta_0(q)\|_\infty=o_\P(1)$ is established for any $m$ and $d$, so the requirement $m>d/2$ imposed in Lemmas \ref{lem:consistency, convex} and \ref{lem:consistency, nonconvex} is not needed. Furthermore, in these two cases, it only requires the minimal side rate restrictions $1/(nh^d) \asymp K/n = o(1)$ in the sub-Gaussian case, and $1/(nh^{\frac{\nu}{\nu-1}d}) \asymp K^{\frac{\nu}{\nu-1}}/n = o(1)$ in the bounded $\nu$-th moment case, up to $\polylog(n)$ terms. See Section~\saref{sec:consistency} in the supplemental appendix for more details.

\subsection{Comparison with Existing Results}\label{sec: Comparison with Existing Results -- Consistency}

These restrictions imposed in Lemmas \ref{lem:consistency, convex} and \ref{lem:consistency, nonconvex}, and the associated results in the supplemental appendix for special cases, are either comparable to or improve upon the existing literature. In the special case of the square loss function and the identity transformation, uniform (over $\bx\in\mathcal{X}$) consistency of the partitioning-based estimator is essentially automatic due to the intrinsic closed-form and linearity of the estimator. Nevertheless, compared to the best known result in that strand of the literature \citep{Cattaneo-Farrell_2013_JoE,Belloni-Chernozhukov-Chetverikov-Kato_2015_JoE,Cattaneo-Farrell-Feng_2020_AOS}, our general results are essentially on par in terms of side rate restrictions and regularity conditions. For example, in the sub-Gaussian case, and up to $\polylog(n)$ terms, the best side rate restriction in that literature requires $1/(nh^{d}) \asymp K/n = o(1)$, while Lemmas \ref{lem:consistency, convex} and \ref{lem:consistency, nonconvex} require $1/(nh^{2d}) \asymp K^2/n = o(1)$, and our improved results in the supplemental appendix for unconnected bases require $1/(nh^{d}) \asymp K/n = o(1)$. Therefore, our results are on par with the best available results for series-based least squares regression \citep{Cattaneo-Farrell_2013_JoE,Belloni-Chernozhukov-Chetverikov-Kato_2015_JoE,Cattaneo-Farrell-Feng_2020_AOS}, despite being able to cover a large class of $M$-estimation settings such as piecewise-polynomial-based quantile, nonlinear, or robust regression.

In the case of quantile regression with tensor-product $B$-splines, Corollary 1 of \cite{Belloni-Chernozhukov-Chetverikov-FernandezVal_2019_JoE} implies that the $L_2$-consistency \eqref{eq: consistency, l2} can be obtained under $1/(nh^{2d}) \asymp K^2/n = o(1)$, and their Corollary 2 implies that the uniform consistency \eqref{eq: consistency, sup} can be obtained under $1/(nh^{4d}) \asymp K^4/n = o(1)$. In contrast, since the generalized residual from quantile regression has a sub-Gaussian envelope, we only require $1/(nh^{2d}) \asymp K^2/n = o(1)$ to establish both kinds of consistency. Moreover, when an unconnected basis is used, or when the loss function is strongly convex and smooth (e.g., the square loss for mean regression), we establish consistency under the weakest possible restriction: $1/(nh^d) \asymp K/n = o(1)$, up to $\polylog(n)$ terms. It remains an open question whether it is possible to establish consistency under the weakest condition $1/(nh^d) \asymp K/n = o(1)$ for general partitioning-based $M$-estimators.

\section{Bahadur Representation and Convergence Rates}\label{sec: Bahadur representation}

The Bahadur representation is
\begin{align*}
      \mathsf{L}^{(\bv)}(\bx,q)
      = -\bp^{(\bv)}(\bx)\trans \bQ_{0,q}^{-1}
         \E_{n}\big[\bp(\bx_i) \eta^{(1)}(\mu_0(\bx_i,q)) \psi(y_i,\eta(\mu_{0}(\bx_i,q)); q)\big]
\end{align*}
where
\begin{align*}
    \bQ_{0,q} = \E\big[\bp(\bx_i)\bp(\bx_i)\trans \Psi_1(\bx_i,\eta(\mu_0(\bx_i,q));q)[\eta^{(1)}(\mu_0(\bx_i,q))]^2\big].
\end{align*}

The following theorem takes the sup-norm consistency of the coefficient estimators $\widehat\bbeta(q)$ as a high-level assumption, and thus avoids imposing any of the specific sufficient conditions discussed in Section \ref{sec: consistency}. The proof is provided in the supplemental appendix (Theorem~\saref{th:bahadur-repres}).

\begin{thm}[Bahadur representation]\label{theorem: bahadur-repres}

  Suppose that Assumptions \ref{Assumption: gl-quasi-uniform-partition}--\ref{Assumption: gl-approximation-error} hold.
  Furthermore, assume the following four conditions:
  \begin{enumerate}[label=\emph{(\roman*)}]
      \item $\sup_{q \in \mathcal{Q}}\big\| \widehat{\boldsymbol{\beta}}(q) - \bbeta_0(q) \big\|_{\infty} = o_{\P}(1)$;
      \item there exists a fixed constant $c>0$ such that
      $\{\bb\in\mathbb{R}^K:\|\bb-\bbeta_0(q)\|_\infty\leq c, q\in\mathcal{Q}\}\subseteq\mathcal{B}$;
      \item $\frac{(\log n)^{d + 2}}{nh^d}=o(1)$;
      \item either $\frac{(\log n)^d}{n^{1-2/\nu} h^d} = o(1)$ or
  $\overline{\psi}(\bx_i,y_i)$ is sub-Gaussian conditional on $\bx_i$.
  \end{enumerate}

  Then
    \begin{small}
    \begin{align}\label{eq: Bahadur rep, disc}
      \sup_{q\in\mathcal{Q}}\sup_{\bx\in\mathcal{X}}
      \big|\widehat{\mu}^{(\bv)}(\bx,q) - \mu_{0}^{(\bv)}(\bx,q) - \mathsf{L}^{(\bv)}(\bx,q) \big|
      \lesssim_{\mathbb{P}}
      h^{-\multindnorm{\bv}}
      \Big(\frac{(\log n)^d}{nh^d}\Big)^{\frac{1}{2}+(\frac{\alpha}{2}\wedge \frac{1}{4})} \log n +
      h^{m - \multindnorm{\bv}}.
    \end{align}
    \end{small}

  If, in addition, $\sup_{y\in\mathcal{Y},q\in\mathcal{Q}} | \psi(y,\eta(\zeta_1);q) - \psi(y,\eta(\zeta_2);q)| \lesssim |\zeta_1-\zeta_2|^\alpha$ without any restrictions on $y$ in Assumption~\ref{Assumption: loss function}\ref{shape of derivative}, then
    \begin{small}
      \begin{align}\label{eq: Bahadur rep, cont}
      \sup_{q\in\mathcal{Q}}\sup_{\bx\in\mathcal{X}}
      \big|\widehat\mu^{(\bv)}(\bx,q) - \mu_{0}^{(\bv)}(\bx,q) - \mathsf{L}^{(\bv)}(\bx,q) \big|
      \lesssim_{\mathbb{P}} h^{- \multindnorm{\bv}} \Big(\frac{(\log n)^d}{n h^d}\Big)^{\frac{1 + \alpha}{2}} \log n +
      h^{m - \multindnorm{\bv}}.
    \end{align}
    \end{small}
\end{thm}

The Bahadur representation \eqref{eq: Bahadur rep, disc} applies to the case where the ``derivative'' $\psi(\cdot,\cdot;q)$ of the loss function may be discontinuous. One typical example is quantile regression (Example \ref{example: Generalized Conditional Quantile Regression}), where the ``derivative'' $\psi(y,\eta;q)=\I(y-\eta<0)-q$, as a function of $(y-\eta)$, is piecewise constant with a jump at zero. In this case we can let $\alpha=1$, and \eqref{eq: Bahadur rep, disc} implies that the order of the remainder in the Bahadur representation for partitioning-based quantile regression is $O(h^{-|\bv|}(nh^{d})^{-3/4}+h^{m-|\bv|})$, up to $\polylog(n)$ terms. Another example is $L_p$ regression with $p\in(1, 2)$, Example \ref{example: $L_p$ regression}, where the derivative of the loss function is $\psi(y,\eta)\equiv \psi(y-\eta)=p|y-\eta|^{p-1}\sgn(\eta-y)$ with $\sgn(\cdot)$ denoting the sign function. As a function of $(y-\eta)$, $\psi(\cdot)$ is $\alpha$-H\"{o}lder on $[0,\infty)$ or $(-\infty,0]$ for all $\alpha\in(0, p-1]$ but not for $\alpha>p-1$. Thus, \eqref{eq: Bahadur rep, disc} applies with the order of the remainder depending on $p$, which is the same as that for quantile regression when $p\geq 3/2$.

On the other hand, the Bahadur representation \eqref{eq: Bahadur rep, cont} applies to the case where the ``derivative'' of the loss is a continuous function of $(y-\eta)$. Nonlinear least squares regression (Example \ref{example: Conditional Distribution Regression}) and quasi-maximum likelihood estimation of generalized linear models (Example \ref{example: Generalized linear models}) fall into this category with the H\"{o}lder parameter $\alpha=1$. In such cases, \eqref{eq: Bahadur rep, cont} implies that the order of the remainder in the Bahadur representation is $O(h^{-|\bv|}(nh^{d})^{-1}+h^{m-|\bv|})$, up to $\polylog(n)$ terms, which is a tighter upper bound than that implied by \eqref{eq: Bahadur rep, disc}. See Section \ref{sec: Examples} and the supplemental appendix for more details.

In both cases, the remainder of the Bahadur representation consists of two terms. The last term $h^{m-|\bv|}$ corresponds to the error from approximating the function $\mu_0$ using the partitioning basis (cf. Assumption \ref{Assumption: gl-approximation-error}), whereas the first term arises from the (potential) nonlinearity underlying the $M$-estimation, and reflects explicitly the role of non-smoothness of the loss function. When the ``derivative'' of the loss function has discontinuity points,
the order of the remainder in \eqref{eq: Bahadur rep, disc} is greater than that in the continuous case \eqref{eq: Bahadur rep, cont}; with a smaller H\"{o}lder parameter $\alpha$, the order of the remainder in both cases could increase.

\subsection{Rates of Convergence}\label{sec: rates of convergence}

The uniform Bahadur representations (Theorem \ref{theorem: bahadur-repres}) can be used to establish convergence rates for the general partitioning-based $M$-estimators. We focus first on uniform convergence over $\bx\in\mathcal{X}$ and $q\in\mathcal{Q}$.

\begin{coro}[Uniform Rate of Convergence]\label{coro: convergence}

    Suppose that Assumptions \ref{Assumption: gl-quasi-uniform-partition}--\ref{Assumption: gl-approximation-error} and the four conditions (i)--(iv) in Theorem \ref{theorem: bahadur-repres} hold. Furthermore, assume one of the following two conditions holds:
    \begin{enumerate}[label=\emph{(\roman*)}]
        \item $\frac{(\log n)^{d+\frac{d+1}{\alpha\wedge 0.5}}}{nh^d}=O(1)$, and $h^{(\alpha\wedge 0.5)m}(\log n)^{0.5 d}=O(1)$, or

        \item the additional condition for \eqref{eq: Bahadur rep, cont} holds, $\frac{(\log n)^{d+\frac{d+1}{\alpha}}}{nh^d}=O(1)$, and $h^{\alpha m}(\log n)^{0.5d}=O(1)$.
    \end{enumerate}

    Then
    \begin{equation}\label{eq: uniform convergence}
        \sup_{q\in\mathcal{Q}}\sup_{\bx\in\mathcal{X}}
        \big|\widehat{\mu}^{(\bv)}(\bx,q)-\mu_0^{(\bv)}(\bx,q)\big|
        \lesssim_\P
        h^{-|\bv|}\sqrt{\frac{\log n}{nh^d}}+h^{m-|\bv|}.
    \end{equation}

\end{coro}

By setting $h\asymp \big(\frac{\log n}{n}\big)^{\frac{1}{2m+d}}$, Corollary \ref{coro: convergence} implies that the partitioning-based $M$-estimator can achieve the uniform convergence rate $\big(\frac{\log n}{n}\big)^{\frac{m}{2m+d}}$. This matches the optimal rate of convergence in $\sup$-norm for nonparametric estimators of the conditional mean \citep{Stone_1982_AoS} and conditional quantiles \citep{chaudhuri1991global}. In this sense, the rate of convergence in Corollary \ref{coro: convergence} is optimal and cannot be further improved at our level of generality.

Theorem \ref{theorem: bahadur-repres} can also be used to obtain the mean square convergence rate of the partitioning-based $M$-estimator uniformly-in-$q$.

\begin{coro}[Mean Square Rate of Convergence]\label{coro: Mean Square Rate of Convergence}

    Suppose that Assumptions \ref{Assumption: gl-quasi-uniform-partition}--\ref{Assumption: gl-approximation-error} and the four conditions (i)--(iv) in Theorem \ref{theorem: bahadur-repres} hold. Furthermore, assume one of the following two conditions holds:
    \begin{enumerate}[label=\emph{(\roman*)}]
        \item $\frac{(\log n)^{d+\frac{d+2}{\alpha\wedge 0.5}}}{nh^d}=o(1)$ and $h^{(\alpha\wedge 0.5)m}(\log n)^{\frac{d+1}{2}}=o(1)$, or

        \item the additional condition for \eqref{eq: Bahadur rep, cont} holds, $\frac{(\log n)^{d+\frac{d+2}{\alpha}}}{nh^d}=o(1)$, and $h^{\alpha m}(\log n)^{\frac{d+1}{2}}=o(1)$.
    \end{enumerate}

    Then
    \begin{equation}\label{eq: l2 convergence}
        \sup_{q\in\mathcal{Q}}\int_{\mathcal{X}} \big(\widehat{\mu}^{(\bv)}(\bx,q)-\mu_0^{(\bv)}(\bx,q)\big)^2 f_X(\bx)d\bx
        \lesssim_\P \frac{1}{nh^{d+2|\bv|}}+h^{2(m-|\bv|)}.
    \end{equation}
\end{coro}

By setting $h\asymp n^{-\frac{1}{2m+d}}$, Corollary \ref{coro: Mean Square Rate of Convergence} implies that the partitioning-based $M$-estimator can also achieve the $L_2$ convergence rate $n^{-\frac{m}{2m+d}}$, uniformly over $\mathcal{Q}$, thereby matching the optimal rate of convergence in $L_2$-norm for nonparametric estimators of conditional means \citep{Stone_1980_AoS} and conditional quantiles \citep{chaudhuri1991global}.

The convergence rates in \eqref{eq: uniform convergence} and \eqref{eq: l2 convergence} capture two contributions: the first term reflects the variance of the estimator, while the second term arises from the error of approximating the unknown $\mu_0$ by the partitioning basis. In the case of Corollary \ref{coro: Mean Square Rate of Convergence}, it is possible to further leverage Theorem \ref{theorem: bahadur-repres} to obtain a precise first-order asymptotic approximation for the integrated mean square error of the partitioning-based $M$-estimator, uniformly over $\mathcal{Q}$, which in turn could be used to develop plug-in asymptotically optimal rules for selecting $K\asymp h^{-d}$. See, for example, Theorem 4.2 in \cite{Cattaneo-Farrell-Feng_2020_AOS} for a similar result in the special case of the square loss and identity transformation functions. We do not pursue this result here for brevity.

\subsection{Comparison with Existing Results}\label{sec: Comparison with Existing Results -- Bahadur}

The results in this section establish uniformly valid Bahadur representations for partitioning-based $M$-estimators at a high level of generality (Theorem \ref{theorem: bahadur-repres}), which imply the convergence rates in Corollaries \ref{coro: convergence} and \ref{coro: Mean Square Rate of Convergence}. The restriction on the tuning parameter $h$ required by the theorem is seemingly minimal: when the envelope of the generalized residual $\psi(y_i,\eta(\mu_0(\bx_i,q));q)$ is sub-Gaussian (or its $\nu$-th moment is bounded with a large $\nu$), we roughly only need $1/(nh^d) \asymp K/n = o(1)$, up to $\polylog(n)$ terms. Having noted this, verification of the high-level consistency assumption $\|\widehat\bbeta(q)-\bbeta_0(q)\|_\infty=o_\P(1)$ in the sub-Gaussian case may require a more stringent condition on $h$, as discussed in Section \ref{sec: consistency}. In the best scenario (e.g., an unconnected basis is used), the minimal restriction $1/(nh^d) \asymp K/n = o(1)$ suffices, while in the worst scenario we need at most $1/(nh^{2d}) \asymp K^2/n = o(1)$, up to $\polylog(n)$ terms.

The rest of this section discusses precisely how our results improve on the prior literature.

\subsubsection*{Mean Regression}

The usual mean regression is a special case of our general setup where $\rho(\cdot,\cdot)$ is the square loss, $\eta(\cdot)$ is the identity link, and $\mathcal{Q}$ is a singleton. Bahadur representations for this special case were established by \cite{Belloni-Chernozhukov-Chetverikov-Kato_2015_JoE} and \cite{Cattaneo-Farrell-Feng_2020_AOS}. Since the derivative of the square loss for mean regression is linear, the first term in \eqref{eq: Bahadur rep, disc} or \eqref{eq: Bahadur rep, cont} does not show up in the uniform linearization of least squares series estimators. See, for example, Lemma SA-4.2 of \cite{Cattaneo-Farrell-Feng_2020_AOS}; $R_{1n,q}$ defined therein has been implicitly included in the leading variance term in \eqref{eq: Bahadur rep, cont} above. Theorem \ref{theorem: bahadur-repres} substantially extends these prior results to other nonlinear settings, under minimal additional conditions.

Finally, Corollaries \ref{coro: convergence} and \ref{coro: Mean Square Rate of Convergence} demonstrate the convergence rate optimality of general partitioning-based series $M$-estimation, recovering in particular known results for mean regression \citep{Belloni-Chernozhukov-Chetverikov-Kato_2015_JoE,Cattaneo-Farrell-Feng_2020_AOS} under essentially the same minimal conditions.

\subsubsection*{Quantile Regression}

Theorem \ref{theorem: bahadur-repres} improves upon prior theoretical results for nonparametric series quantile regression estimators. The most recent advance in this literature is due to \cite{Belloni-Chernozhukov-Chetverikov-FernandezVal_2019_JoE}, which establishes a uniform linear approximation for general series-based quantile regression estimators. In comparison, we exploit the ``local support'' feature of the partitioning basis, and make improvements in (at least) four aspects. To summarize these improvements without additional cumbersome notation, we set $\bv=\bm{0}$ and ignore the smoothing bias $h^m$ in the Bahadur approximation remainders.

First, \cite{Belloni-Chernozhukov-Chetverikov-FernandezVal_2019_JoE} shows that the order of the remainder in the Bahadur representation is $O((nh^d)^{-3/4} h^{-d/2})$, up to $\polylog(n)$ terms (see proofs of Theorem 2 and Corollary 2 therein for details). In contrast, Theorem \ref{theorem: bahadur-repres} implies that the remainder in the Bahadur representation for partitioning-based quantile regression estimators is $O((nh^d)^{-3/4})$, up to $\polylog(n)$ terms, which is not only a much tighter bound but also matches the optimal parametric bound when taking $nh^d$ as the effective sample size.\sloppy

Second, the rate restriction $1/(nh^{4d}) \asymp K^4/n = o(1)$ is required for $B$-spline-based estimators in \cite{Belloni-Chernozhukov-Chetverikov-FernandezVal_2019_JoE}. In contrast, the restriction on $h$ in Theorem \ref{theorem: bahadur-repres} depends on the tail behavior of the generalized residuals and becomes weaker as $\nu$ gets larger. In the best case (the residuals have a sub-Gaussian envelope) we only need the seemingly weakest restriction $1/(nh^d) \asymp K/n = o(1)$, up to $\polylog(n)$ terms, along with the consistency condition for $\widehat\bbeta(q)$. Recall that in the sub-Gaussian scenario we need at worst $1/(nh^{2d}) \asymp K^2/n =o(1)$, up to $\polylog(n)$ terms, to satisfy the consistency requirement.

Third, the restriction $h^{m-d}=o(n^{-\varepsilon})$ for some $\varepsilon>0$ in \cite{Belloni-Chernozhukov-Chetverikov-FernandezVal_2019_JoE} implicitly requires the smoothness $m$ of the conditional quantile function be greater than the dimensionality $d$ of the covariates. In contrast, the proof of Theorem \ref{theorem: bahadur-repres} does not need such a restriction, though a weaker condition $m>d/2$ might be needed to verify the consistency condition on $\widehat{\bbeta}(q)$; see Lemmas \ref{lem:consistency, convex} and \ref{lem:consistency, nonconvex}. Furthermore, when an unconnected basis (e.g., piecewise polynomials) is used for approximation, the condition $m>d/2$ is unnecessary for consistency, and thus we have \textit{no} constraint on the relation between smoothness $m$ and dimensionality $d$; see Section \ref{sec: Weaker Conditions for Special Cases}.

Fourth, compared to \cite{Belloni-Chernozhukov-Chetverikov-FernandezVal_2019_JoE}, we allow for a possibly non-identity link. Introducing a link function may lead to non-convexity of the loss $\rho(y,\eta(\theta);q)$ with respect to $\theta$, making the usual proof strategies for consistency and Bahadur representation under convexity inapplicable. For example, non-convex quantile regression is covered in Theorem \ref{theorem: bahadur-repres} by virtue of our general consistency results in Lemma \ref{lem:consistency, nonconvex}.

All of the aforementioned improvements are practically relevant. For example, they accommodate univariate quantile regression using the piecewise constant basis with the IMSE-optimal choice of the mesh size $h$ (in this case $h\asymp n^{-1/3}$ and $m=d$), which was theoretically excluded in the prior literature.

Finally, Corollaries \ref{coro: convergence} and \ref{coro: Mean Square Rate of Convergence} establish the optimal rate of convergence for general partitioning-based series $M$-estimators, which substantially improve on prior work on quantile series regression in particular. More specifically, the conditions on the mesh size $h$, the smoothness $m$, and the dimensionality $d$ in both corollaries are weaker than in prior work. In the best case (e.g., an unconnected basis is used and a sub-Gaussian envelope for residuals exists), we only require the seemingly minimal restriction $1/(nh^d) \asymp K/n =o(1)$, up to $\polylog(n)$ terms, and an arbitrary relation between $m$ and $d$ is permitted. For comparison, in the special case of quantile regression, \cite{Belloni-Chernozhukov-Chetverikov-FernandezVal_2019_JoE} shows that series estimators can achieve the fastest possible uniform-in-$q$ rate in both $L_2$-norm and $\sup$-norm (see Comments 3 and 4 therein), but under more stringent conditions: $\eta$ is the identity function, $m>d$, and $1/(nh^{4d}) \asymp K^4/n = o(n^{-\varepsilon})$ for some $\varepsilon>0$ (see their Corollary 2). Such conditions exclude, e.g., the IMSE-optimal choice of $h$ for the Haar basis when $d=1$ (since $h\asymp n^{-1/3}$ and $m=d=1$), or piecewise linear fit when $d=2$ (since $m=d=2$).

\subsubsection*{Other Nonparametric Smoothing Methods}

\cite{Kong-Linton-Xia_2013_ET} establishes a similar Bahadur representation for kernel-based $M$-estimators using weakly stationary time series data. They consider a special case of our setup in Assumption \ref{Assumption: loss function}: their loss function class $\mathcal{Q}$ is a singleton, $\eta$ is an identity function, and the ``derivative'' of the loss can be written as $\psi(y,\eta)\equiv\psi(y-\eta)$ and is assumed to be piecewise Lipschitz continuous. In a comparable cross-sectional context with $\alpha=1$ and $\bv=\bm{0}$, the order of the remainder in our Bahadur representation \eqref{eq: Bahadur rep, disc} is $O((nh^d)^{-3/4})$, up to $\polylog(n)$ terms, and thus Theorem \ref{theorem: bahadur-repres} matches their approximation error up to a minor difference in $\log n$ terms. Taking $nh^d$ to be the effective sample size, the approximation rate can not be further improved at this level of generality, and hence Theorem \ref{theorem: bahadur-repres} establishes that the partitioning-based series $M$-estimator in \eqref{eq: M estimator} can achieve the same best Bahadur approximation as local polynomial kernel methods, up to $\polylog(n)$ terms.

Compared to \cite{Kong-Linton-Xia_2013_ET} or other similar contributions in the literature, Theorem \ref{theorem: bahadur-repres} exhibits (at least) two novel features. First, the Bahadur representations \eqref{eq: Bahadur rep, disc} and \eqref{eq: Bahadur rep, cont} hold uniformly not only over the evaluation point $x\in\mathcal{X}$, but also over the loss function index $q\in\mathcal{Q}$, which may be important, for example, to study simultaneous quantile regression where the \textit{entire} conditional quantile process may be of interest. Second, Theorem \ref{theorem: bahadur-repres} also covers the more general setup where the ``derivative'' function may exhibit different degrees of smoothness, reflected by discontinuity points and/or the H\"{o}lder parameter $\alpha$, or admits a more complex structure so that $\psi(y,\eta;q)$ cannot be written as $\psi(y-\eta;q)$. Thus, we cover more examples such as distribution regression (Example \ref{example: Conditional Distribution Regression}) and $L_p$ regression with $p\in(1,2)$ (Example \ref{example: $L_p$ regression}). Finally, \cite{Kong-Linton-Xia_2013_ET} does not discuss convergence rates as we do in Corollaries \ref{coro: convergence} and \ref{coro: Mean Square Rate of Convergence}.

In the context of nonparametric penalized smoothing spline methods, \cite{Shang-Cheng_2013_AOS} also establishes a uniform Bahadur representation (and other results) that can be compared to Theorem \ref{theorem: bahadur-repres}. However, their paper imposes more stringent assumptions and hence covers a smaller class of settings: using our notation, they assume that (i) $\mathcal{Q}$ is a singleton so their uniformity is only over $\mathcal{X}$; (ii) $d=1$ so they consider only scalar covariate $\bx_i$; and (iii) $\rho(\cdot,\eta(\cdot))$ is smooth so they rule out many important examples such as quantile regression, and Tukey and Huber regression. Furthermore, their results do not take explicit advantage of specific moment and boundedness conditions, or the structure of the nonparametric estimator, and instead impose the generic side condition $nh^2\to\infty$, which is comparable to our condition $K^2/n\to\infty$, up to $\polylog(n)$ terms. Most importantly, in the closest comparable case ($d=1$, $\alpha=1$, and $\bv=\bm{0}$), and only focusing on the variance component for simplicity, the order of the remainder in their uniform Bahadur representation (a combination of Theorem 3.4 and Lemma 3.1 in \cite{Shang-Cheng_2013_AOS}) is $O((nh)^{-1} h^{-(6m-1)/(4m)})$, while \eqref{eq: Bahadur rep, cont} in Theorem \ref{theorem: bahadur-repres} gives the optimal result $O((nh)^{-1})$, thereby demonstrating a substantial improvement over their result. Finally, as for convergence rates, Proposition 3.3 in \cite{Shang-Cheng_2013_AOS} and our Corollary \ref{coro: Mean Square Rate of Convergence} are essentially equivalent, both delivering optimal mean square convergence. They do not explicitly discuss uniform convergence rates as we do in Corollary \ref{coro: convergence}.

\section{Strong Approximation}\label{sec: uniform inference}

The uniform Bahadur representations in Theorem \ref{theorem: bahadur-repres} can also be leveraged to establish uniform distribution theory for $\widehat\mu^{(\bv)}$. The infeasible conditional variance of the estimator can be written as
\begin{align*}
    \bar\Omega_{\bv}(\bx,q) = \bp^{(\bv)}(\bx)\trans\bar\bQ_q^{-1}\bar\bSigma_q\bar\bQ_q^{-1}\bp^{(\bv)}(\bx),
\end{align*}
where
\begin{align*}
    \bar\bQ_{q} &= \E_n\big[\bp(\bx_i)\bp(\bx_i)\trans \Psi_1(\bx_i,\eta(\mu_0(\bx_i,q));q)[\eta^{(1)}(\mu_0(\bx_i,q))]^2\big], \qquad \text{and}\\
    \bar\bSigma_q &= \E_n\big[\bp(\bx_i)\bp(\bx_i)\trans\E[\psi(y_i,\eta(\mu_0(\bx_i,q));q)^2|\bx_i]
[\eta^{(1)}(\mu_0(\bx_i,q))]^2\big]
\end{align*}
Accordingly, we define a feasible variance estimator as
\begin{align*}
    \widehat{\Omega}_{\bv}(\bx,q) = \bp^{(\bv)}(\bx)\trans\widehat\bQ_q^{-1}\widehat\bSigma_q\widehat\bQ_q^{-1}\bp^{(\bv)}(\bx),
\end{align*}
where $\widehat\bQ_q$ and $\widehat\bSigma_q$ are some estimators of $\bar\bQ_q$ and $\bar\bSigma_q$, respectively, which are consistent in a sense described below. Therefore, $\widehat\Omega_{\bv}(\bx,q)$ is an estimator of the infeasible conditional variance $\bar\Omega_{\bv}(\bx,q)$.

Statistical inference on $\mu_0^{(\bv)}$ usually relies on the following $t$-statistic process:
\begin{align*}
    T(\bx,q) = \frac{\widehat\mu^{(\bv)}(\bx,q)-\mu_0^{(\bv)}(\bx,q)}
                    {\sqrt{\widehat{\Omega}_{\bv}(\bx,q)/n}}, \qquad (\bx,q)\in\mathcal{X}\times\mathcal{Q},
\end{align*}
where we drop the dependence of $T(\cdot)$ on $\bv$ for simplicity.

Employing Theorem \ref{theorem: bahadur-repres}, or more precise arguments under slightly weaker conditions, it is easy to show that $T(\bx,q)$ converges in distribution to $\mathsf{N}(0,1)$ for each $(\bx,q)\in\mathcal{X}\times\mathcal{Q}$. However, the stochastic process $(T(\bx,q):(\bx,q)\in\mathcal{X}\times\mathcal{Q})$ is generally not asymptotically tight and, therefore, does not converge weakly in $\ell^\infty(\mathcal{X}\times\mathcal{Q})$, where $\ell^\infty(\mathcal{X}\times\mathcal{Q})$ denotes the set of all (uniformly) bounded real functions on $\mathcal{X} \times \mathcal{Q}$ equipped with uniform norm \citep{vaart96_weak_conver_empir_proces}. Nevertheless, we can construct a Gaussian process, in a possibly enlarged probability space, that approximates the entire process $T(\cdot)$ sufficiently fast, which can then be used to approximate the finite sample distribution of the function $M$-estimator $\widehat{\mu}^{(\bv)}(\cdot)$.

More precisely, under some mild consistency conditions on $\widehat\Omega_{\bv}(\bx,q)$, our Theorem \ref{theorem: bahadur-repres} guarantees that $\sup_{q\in\mathcal{Q}} \sup_{\bx\in\mathcal{X}} \big| T(\bx,q) - t(\bx,q) \big| \to_\P 0$ sufficiently fast, where
\begin{align*}
    t(\bx,q) = -\frac{\bp^{(\bv)}(\bx)\trans\bar\bQ_q^{-1}}{\sqrt{\bar\Omega_{\bv}(\bx,q)}}
                \mathbb{G}_n\big[\bp(\bx_i)\eta^{(1)}(\mu_0(\bx_i,q)) \psi(y_i,\eta(\mu_0(\bx_i,q));q)\big].
\end{align*}
It follows that, conditional on $\{\bx_i\}_{i=1}^n$, the randomness of $t(\bx,q)$ comes exclusively from the $K$-dimensional vector $\mathbb{G}_n[\bp(\bx_i)\eta^{(1)}(\mu_0(\bx_i,q)) \psi(y_i,\eta(\mu_0(\bx_i,q));q)]$. Thus, our proof strategy is to further ``discretize'' this vector with respect to $q\in\mathcal{Q}$, and then apply Yurinskii's coupling \citep{yurinskii1978error} to construct a (conditional) Gaussian process that is close to the original $t$-statistic process $T(\bx,q)$ uniformly over both $\bx\in\mathcal{X}$ and $q\in\mathcal{Q}$. Our construction leverages a conditional Strassen's theorem \citep[Theorem B.2]{chen2020jackknife} to generalize prior coupling results \citep[Lemma 36]{Belloni-Chernozhukov-Chetverikov-FernandezVal_2019_JoE}. See Section~\saref{sec:strong-approximation} in the supplemental appendix for details.

Our strong approximation approach is formalized in the next theorem. We employ high-level conditions to ease the exposition, but those conditions can be verified using Corollaries \ref{coro: convergence} and \ref{coro: Mean Square Rate of Convergence}, and Theorem \ref{theorem: bahadur-repres}, as well as using the more general results in the supplemental appendix. Let $r_{\tt UC}$, $r_{\tt BR}$, $r_{\tt VC}$, and $r_{\tt SA}$ be positive non-random sequences as $n\to\infty$. The proof is available in the supplemental appendix (Theorem~\saref{th:strong-approximation-binscatter-yurinski}).

\begin{thm}[Strong approximation]\label{thm: strong approximation}
    Suppose that Assumptions \ref{Assumption: gl-quasi-uniform-partition}--\ref{Assumption: gl-approximation-error} with $\nu\geq3$ hold. Furthermore, assume the following four conditions hold:
    \begin{enumerate}[label=\emph{(\roman*)}]

        \item $\sup_{q\in\mathcal{Q}}\sup_{\bx\in\mathcal{X}} |\widehat\mu^{(\bv)}(\bx,q)-\mu_0^{(\bv)}(\bx,q)| \lesssim_\P h^{-|\bv|}r_{\tt UC}$.

        \item $\sup_{q\in\mathcal{Q}}\sup_{\bx\in\mathcal{X}}
              \big|\widehat{\mu}^{(\bv)}(\bx,q) - \mu_0^{(\bv)}(\bx,q) - \mathsf{L}^{(\bv)}(\bx,q) \big|
              \lesssim_\P h^{-|\bv|}r_{\tt BR}$, with $\frac{\log n}{nh^d} \lesssim r_{\tt BR}$.

        \item $\sup_{q\in\mathcal{Q}}\sup_{\bx\in\mathcal{X}}
              \big|\widehat\Omega_{\bv}(\bx,q)-\bar\Omega_{\bv}(\bx,q)\big|
              \lesssim_\P h^{-2|\bv|-d}r_{\tt VC}$, with $r_{\tt VC} = o(1)$.

        \item \label{eq: strong approx,  lipschitz condition} $\E\big[\big|
                \psi(y_i,\eta(\mu_0(\bx_i,q));q) \eta^{(1)}(\mu_{0}(\bx_i,q))$
        \item[] $\qquad\qquad
                - \psi(y_i, \eta(\mu_{0}(\bx_i,\tilde{q})); \tilde{q}) \eta^{(1)}(\mu_{0}(\bx_i,\tilde{q}))
                \big|^2 \,\big|\, \bx_i \big]
                \lesssim | q - \tilde{q} |$,
              for all $q, \tilde{q} \in \mathcal{Q}$.

    \end{enumerate}

    Then (provided the probability space is rich enough) there exists a stochastic process $(Z(\bx,q):(\bx,q)\in\mathcal{X}\times\mathcal{Q})$ such that, conditional on $\bX_n=(\bx_1,\cdots,\bx_n)$, $Z$ is a mean-zero Gaussian process with $\E[Z(\bx,q)Z(\tilde\bx,\tilde{q})|\bX_n]=\E[t(\bx,q)t(\tilde\bx,\tilde{q})|\bX_n]$ for all $(\bx,q),(\tilde{\bx},\tilde{q})\in\mathcal{X}\times\mathcal{Q}$, and
    \begin{align*}
        \sup_{q\in\mathcal{Q}}\sup_{\bx\in\mathcal{X}}
        \big|T(\bx,q)-Z(\bx,q)\big| \lesssim_\P r_{\tt SA} + \sqrt{nh^d}(r_{\tt UC} \; r_{\tt VC} + r_{\tt BR}).
    \end{align*}
    where $r_{\tt SA}=o(1)$ is any positive sequence satisfying
    \[\Big(\frac{1}{nh^{3d}}\Big)^{\frac{1}{10}}\sqrt{\log n} + \frac{\log n}{\sqrt{n^{1-2/\nu}h^{d}}} = o(r_{\tt SA}).\]

    Furthermore, if $\overline\psi(\bx_i,y_i)$ is sub-Gaussian conditional on $\bx_i$, then the same result holds with any positive sequence $r_{\tt SA}=o(1)$ satisfying
    \[\Big(\frac{1}{nh^{3d}}\Big)^{\frac{1}{10}}\sqrt{\log n} + \frac{(\log n)^{3/2}}{\sqrt{nh^d}} = o(r_{\tt SA}).\]
\end{thm}

The speed of strong approximation in Theorem \ref{thm: strong approximation} is determined by four factors:
the uniform convergence rate $r_{\tt UC}$,
the order of the remainder in the Bahadur representation $r_{\tt BR}$,
the convergence rate $r_{\tt VC}$ of the variance estimator $\widehat\Omega_{\bv}$, and
the strong approximation rate $r_{\tt SA}$.
Therefore, our strong approximation results are established at a high level of generality, building on our prior theoretical results: Corollary \ref{coro: convergence} for $r_{\tt UC}$, and Theorem \ref{theorem: bahadur-repres} for $r_{\tt BR}$, while $r_{\tt VC}$ is a high-level condition that needs to be verified on a case-by-case basis. See Sections \ref{sec: Feasible Uniform Inference} and \ref{sec: Examples} for more discussion.

With respect to the strong approximation rate, Theorem \ref{thm: strong approximation} lays down two versions of lower bounds on $r_{\tt SA}$, depending on the tail behavior of the generalized residuals. Such restrictions may not be optimal, but are still weak enough to cover almost all partition size choices commonly used in practice. In particular, the restriction on $r_{\tt SA}$ in Theorem \ref{thm: strong approximation} allows the use of the MSE-optimal choice $h\asymp n^{-\frac{1}{2m+d}}$ in all cases except when $m\leq d$ (e.g., unidimensional Haar basis approximation); there is also room for undersmoothing in order to make the smoothing bias negligible in all cases but $m\leq d$. The strong approximation for one-dimensional partitioning-based series estimators in the special case of square loss and identity transformation functions was studied in \cite{Cattaneo-Farrell-Feng_2020_AOS,Cattaneo-Crump-Farrell-Feng_2024_AER} via a different coupling strategy, which delivered tighter approximation results allowing for an MSE-optimal choice of $h$. We conjecture those techniques could be adapted to cover the case $m=d=1$ for the general partitioning-based $M$-estimator in \eqref{eq: M estimator}, but we do not pursue this line of research here because it would require a different theoretical treatment.

\subsection{Comparison with Existing Results}\label{sec: Comparison with Existing Results -- SA}

Theorem \ref{thm: strong approximation} establishes strong approximation results for general partitioning-based $M$-estimators. In the prior literature, similar results are usually available only in specific scenarios such as least squares regression or quantile regression. To be more precise, in the least squares context ($\mathcal{Q}$ is a singleton), \cite{Cattaneo-Farrell-Feng_2020_AOS} establishes uniform inference theory for univariate regression ($d=1$) and multivariate regression ($d>1$) separately via different strong approximation methods. In particular, when $d>1$, the same Yurinskii coupling technique is employed to obtain strong approximation for $t$-statistic processes, leading to similar rate restrictions on $h$. Theorem \ref{thm: strong approximation} is a substantial generalization of results therein, not just covering other loss functions, but also providing distributional approximation uniformly over the loss function index $q\in\mathcal{Q}$.

In the quantile regression context, \cite{Belloni-Chernozhukov-Chetverikov-FernandezVal_2019_JoE} provides two strong approximations for general series-based estimators. When $B$-splines are used, their first strategy relies on a pivotal coupling, imposing $1/(nh^{10 d}) = o(n^{-\varepsilon})$ and $h^{m-d}=o(n^{-\varepsilon})$ for some constant $\varepsilon>0$ (see Theorem 11 therein), while the second strategy uses a Gaussian coupling (as in this paper), imposing $1/(nh^{4d\vee(2+3d)}) = o(n^{-\varepsilon})$ and $h^{m-d}=o(n^{-\varepsilon})$ (see Theorem 12 and Comment 13 therein). In comparison, our Theorem \ref{thm: strong approximation} requires weaker conditions on the tuning parameter $h$ and the relation between the smoothness $m$ and the dimensionality $d$. Specifically, we assume $1/(nh^{3d})=o(1)$ up to $\polylog(n)$ terms for a valid approximation. This improvement is practically relevant: for example, it allows for Gaussian approximation of linear-spline-based univariate quantile regression estimators with the MSE-optimal mesh size $h\asymp n^{-1/5}$. In addition, our general strategy to verify the consistency condition on $\widehat\bbeta(q)$ in Theorem \ref{theorem: bahadur-repres} only requires $m>d/2$ (not required for the special case of unconnected basis), which is weaker than $m>d$ as implicitly assumed in \cite{Belloni-Chernozhukov-Chetverikov-FernandezVal_2019_JoE}. In practice, this improvement can accommodate, for example, the use of cubic splines for trivariate quantile regression.

Finally, \cite{Shang-Cheng_2013_AOS} establishes uniform inference results for nonparametric penalized smoothing spline M-estimators. As mentioned before, their work is more specialized because they assume that $d=1$, $\mathcal{Q}$ is a singleton, and $\rho(\cdot,\eta(\cdot))$ is smooth. Furthermore, their approach to constructing valid confidence bands and related uniform inference methods relies on approximating the suprema of the stochastic process directly via extreme value theory \citep[Theorem 5.1]{Shang-Cheng_2013_AOS}, which leads to substantially slower approximation rates and requires stronger assumptions and side rate restrictions; see \cite{Belloni-Chernozhukov-Chetverikov-Kato_2015_JoE} and \cite{Cattaneo-Farrell-Feng_2020_AOS} for more discussion in the context of nonparametric least squares series estimation. In contrast, Theorem \ref{thm: strong approximation}, and our related uniform inference methods, provide a pre-asymptotic approximation with better finite sample properties, faster approximation rates, and weaker regularity conditions.

\section{Feasible Plug-in Uniform Inference}\label{sec: Feasible Uniform Inference}

The (conditional) Gaussian process $(Z(\bx,q):(\bx,q)\in\mathcal{X}\times\mathcal{Q})$ in Theorem \ref{thm: strong approximation} is still infeasible since its covariance structure contains unknowns. This section establishes the validity of a generic \textit{plug-in method} to construct a feasible version of $Z(\bx,q)$. We later employ this result in Section \ref{sec: Examples} to develop feasible uniform inference in the context of our motivating examples.

The core idea behind the plug-in method is to estimate the covariance structure of $Z(\bx,q)$ and then simulate its feasible version $\widehat{Z}(\bx,q)$, a Gaussian process conditional on the data. If the covariance estimate converges to the true covariance sufficiently fast, $\widehat{Z}(\bx,q)$ will be ``close'' to a copy of $Z(\bx,q)$. The covariance structure of the process $Z(\bx,q)$ in Theorem \ref{thm: strong approximation} is
\begin{align*}
    \E[Z(\bx,q)Z(\tilde\bx,\tilde{q})|\bX_n]
    = \frac{\bp^{(\bv)}(\bx)\trans\bar{\bQ}_q^{-1} \bar\bSigma_{q,\tilde{q}}\bar{\bQ}_{\tilde{q}}^{-1}\bp^{(\bv)}(\tilde\bx)}
           {\sqrt{\bar{\Omega}_{\bv}(\bx,q)\bar{\Omega}_{\bv}(\tilde\bx,\tilde{q})}},
\end{align*}
for all $(\bx,q),(\tilde{\bx},\tilde{q})\in\mathcal{X}\times\mathcal{Q}$, where
\begin{align*}
    \bar\bSigma_{q,\tilde{q}}
    = \E_n\big[S_{q, \tilde{q}}(\bx_i)  \eta^{(1)}(\mu_0(\bx_i, q)) \eta^{(1)}(\mu_0(\bx_i, \tilde{q})) \bp(\bx_i) \bp(\bx_i)\trans]
\end{align*}
with
\begin{align*}
    S_{q, \tilde{q}}(\bx)
    = \E\big[\psi(y_i, \eta(\mu_0(\bx_i, q)); q) \psi(y_i, \eta(\mu_0(\bx_i, \tilde{q})); \tilde{q}) \big| \bx_i = \bx \big].
\end{align*}
Given context-specific estimates $\widehat{\bQ}_q$ and $\bx \mapsto \widehat{S}_{q, \tilde{q}}(\bx)$, we can put
\begin{align*}
  \widehat\bSigma_{q,\tilde{q}}
  = \E_n\big[\widehat{S}_{q, \tilde{q}}(\bx_i)  \eta^{(1)}(\widehat{\mu}(\bx_i, q)) \eta^{(1)}(\widehat{\mu}(\bx_i, \tilde{q})) \bp(\bx_i) \bp(\bx_i)\trans\big],
\end{align*}
and $\widehat{\Omega}_{\bv}(\bx,q) = \bp^{(\bv)}(\bx)\trans\widehat\bQ_q^{-1}\widehat\bSigma_{q, q} \widehat\bQ_q^{-1}\bp^{(\bv)}(\bx)$ as above. Section \ref{sec: Examples} illustrates how the estimates $\widehat{\bQ}_q$ and $\bx \mapsto \widehat{S}_{q, \tilde{q}}(\bx)$ can be constructed in specific examples. Then, a feasible Gaussian approximation $\widehat{Z}(\bx,q)$ can be constructed as a mean-zero Gaussian process conditional on the data $\bD_n = ((y_1,\bx_1),\cdots,(y_n,\bx_n))$ with conditional covariance structure
\begin{align}\label{eq: feasible cov}
    \E[\widehat{Z}(\bx,q)\widehat{Z}(\tilde{\bx},\tilde{q})|\bD_n]
    = \frac{\bp^{(\bv)}(\bx)\trans \widehat\bQ_q^{-1} \widehat\bSigma_{q,\tilde{q}}\widehat\bQ_{\tilde{q}}^{-1}\bp^{(\bv)}(\tilde\bx)}
           {\sqrt{\widehat\Omega_{\bv}(\bx,q)\widehat\Omega_{\bv}(\tilde\bx,\tilde{q})}},
\end{align}
for all $(\bx,q),(\tilde{\bx},\tilde{q})\in\mathcal{X}\times\mathcal{Q}$.

The following theorem establishes the validity of the plug-in approach.

\begin{thm}[Feasible Plug-in Strong Approximation]\label{thm: Feasible Plug-in Strong Approximation}
    Suppose that the assumptions in Theorem \ref{thm: strong approximation} hold. Furthermore, assume the following three conditions hold:
    \begin{enumerate}[label=\emph{(\roman*)}]

        \item $\|\widehat{\bQ}_{q} - \bar{\bQ}_{q}\|_{\infty} \lesssim_{\P} h^d r_{\mathtt{Q}}$ and $\|\widehat{\bQ}^{-1}_{q}\|_{\infty} \lesssim_{\P} h^{-d}$, with $r_{\mathtt{Q}} = o(1)$.

        \item $S_{q,\tilde{q}}(\bx)$ is continuous for all $q, \tilde{q}$, and $\sup_{\bx, q, q_1 \neq q_2} \frac{|S_{q,q_1}(\bx) - S_{q,q_2}(\bx)|}{|q_1 - q_2|} \lesssim 1$.

        \item $\sup_{q, \tilde{q}, \bx} |\widehat{S}_{q,\tilde{q}}(\bx) - S_{q,\tilde{q}}(\bx)| \lesssim_{\P} r_{\mathtt{S}}$, with $r_{\mathtt{S}} = o(1)$.

    \end{enumerate}

    Then (provided the probability space is rich enough) there exists a mean-zero Gaussian process, conditional on $\bD_n$, $(\widehat{Z}(\bx,q):(\bx,q)\in\mathcal{X}\times\mathcal{Q})$ satisfying \eqref{eq: feasible cov} and
    \begin{align*}
        \sup_{q\in\mathcal{Q}}\sup_{\bx\in\mathcal{X}}
        \big| \widehat{Z}(\bx,q) - Z^{\star}(\bx,q) \big| \lesssim_\P
        [(r_{\mathtt{UC}} + r_{\mathtt{S}})^{1 / 4} + r_{\mathtt{Q}} + r_{\mathtt{VC}}] \sqrt{\log n},
    \end{align*}
    where $(Z^{\star}(\bx,q):(\bx,q)\in\mathcal{X}\times\mathcal{Q})$ is a process such that, conditional on $\bX_n$, $Z(\cdot)$ and $Z^{\star}(\cdot)$ have the same (conditional) distribution, and $Z^{\star}(\cdot)$ is (conditionally) independent of $(y_1,\ldots,y_n)$.
\end{thm}

Once we have a feasible process $\widehat{Z}(\bx,q)$ that is ``close'' to a copy of $Z(\bx,q)$ uniformly over $\mathcal{X}\times\mathcal{Q}$, conditional on the data, then $\widehat{Z}(\bx,q)$ can be used to conduct inference on the entire function $\mu_0(\bx,q)$, and functionals thereof. For example, our strong approximation results can be converted to convergence of the Kolmogorov distance between the distributions of $\sup_{q\in\mathcal{Q}}\sup_{\bx\in\mathcal{X}}|T(\bx,q)|$ and its feasible Gaussian approximation $\sup_{q\in\mathcal{Q}}\sup_{\bx\in\mathcal{X}}|\widehat{Z}(\bx,q)|$. See Theorem \saref{th:k-s-distance} in the supplemental appendix for the formal result.

Furthermore, Theorem \saref{th:confidence-bands} in the supplemental appendix establishes the asymptotic validity of the uniform confidence band for $\mu_0^{(\bv)}$ given by
\begin{align}\label{eq: confidence band}
    \mathsf{CB}_{1-\alpha}(\bx,q) = \Big[ \widehat{\mu}^{(\bv)}(\bx,q) \pm \mathfrak{c}_{1-\alpha}\sqrt{\widehat\Omega_{\bv}(\bx,q)/n} \Big]
\end{align}
with $\mathfrak{c}_{1-\alpha}$ satisfying
\[\P\Big(\sup_{q\in\mathcal{Q}}\sup_{\bx\in\mathcal{X}}|\widehat{Z}(\bx,q)|\leq \mathfrak{c}_{1-\alpha} \Big| \bD_n \Big)
  =1-\alpha+o_\P(1),
\]
provided the smoothing (or misspecification) bias relative to the standard error of the estimator is small, which could be achieved by undersmoothing, bias correction \citep{hall1992effect}, simply ignoring the bias \citep{hall2001bootstrapping}, robust bias correction
\citep{Calonico-Cattaneo-Farrell_2018_JASA,
Calonico-Cattaneo-Farrell_2022_Bernoulli},
or the Lepskii's method \citep{lepskii1992asymptotically,birge2001alternative},
among other possibilities. Thus, under regularity conditions, the confidence band \eqref{eq: confidence band} covers $\mu_0^{(\bv)}$ with probability approximately $1-\alpha$ in large samples, that is,
\begin{align*}
    \lim_{n\to\infty} \P \big( \mu_0^{(\bv)}(\bx,q) \in \mathsf{CB}_{1-\alpha}(\bx,q), \;\;\text{for all } (\bx,q)\in\mathcal{X}\times\mathcal{Q} \big) = 1-\alpha.
\end{align*}

\section{Verification of Assumptions in Examples}\label{sec: Examples}

We verify assumptions in our four motivating examples (Section \ref{sec: Motivating Examples}). Assumptions \ref{Assumption: gl-quasi-uniform-partition} and \ref{Assumption: gl-local-basis} concern the partitioning-based methodology itself, Assumption \ref{Assumption: dgp-holder} are already primitive conditions on the data generating process, and Assumption \ref{Assumption: gl-approximation-error} can be verified for usual local bases (e.g., piecewise polynomials and splines) when we assume the functional parameter $\mu_0$, such as the conditional quantile function in Example 1 or conditional distribution function in Example 2, is smooth enough (see Assumption \ref{Assumption: dgp-holder}(iv) for details).
Thus, we focus attention on two major issues that remain: (i) how the high-level conditions imposed in Assumptions \ref{Assumption: loss function} and \ref{Assumption: simple setup}, and Condition \emph{\ref{eq: strong approx,  lipschitz condition}} in Theorem \ref{thm: strong approximation}, can be verified under intuitive primitive assumptions; and (ii) how to implement uniform inference based on our theory in Section \ref{sec: Feasible Uniform Inference}.

\subsection{Example \ref{example: Generalized Conditional Quantile Regression}: Generalized Conditional Quantile Regression}\label{sec: quantile reg}

This example considers generalized conditional quantile regression with a possibly non-identity link: $\rho(y,\eta;q)=(q-\I(y<\eta))(y-\eta)$, where $q\in\mathcal{Q}$ denotes the quantile position. Thus, let $\eta(\mu_0(\bx, q))$ be the conditional $q$-quantile of $Y$ given $\bX = \bx$; we verify in the supplemental appendix that such $\mu_0$ solves \eqref{eq: Rho problem}. For this example, the following simple proposition, proven in the supplemental appendix (Proposition~\saref{prop:gl-holder-quantile-regression}), gives sufficient conditions to verify the general Assumptions \ref{Assumption: loss function} and \ref{Assumption: simple setup}, and Condition \emph{\ref{eq: strong approx,  lipschitz condition}} in Theorem \ref{thm: strong approximation}.

\begin{prop}[Quantile Regression]\label{lem: quantile reg}
  Suppose Assumption~\ref{Assumption: dgp-holder} holds with $\mathcal{Q}=[\varepsilon_0, 1-\varepsilon_0]$ for some $\varepsilon_0\in(0,0.5)$, the loss is $\rho(y,\eta;q) = (q-\I(y<\eta))(y-\eta)$, the first moment of $Y$ is finite $\E[|Y|] < \infty$. Assume further that $\eta(\cdot)\colon \mathbb{R} \to \mathcal{E}$ is strictly monotonic and twice continuously differentiable with $\mathcal{E}$ an open connected subset of $\mathbb{R}$ containing the conditional $q$-quantile of $Y|\bX=\bx$, given by $\eta(\mu_{0}(\bx, q))$ for all $(\bx,q)$, and that $|\eta^{(1)}(\zeta)|$ is bounded away from zero uniformly over $\bx \in \mathcal{X}$, $q \in \mathcal{Q}$, and $\zeta$ satisfying $|\zeta-\mu_0(\bx,q)|\leq r$ for some $r>0$; $f_{Y|X}(\eta(\mu_{0}(\bx,q)) | \bx)$ is bounded away from zero uniformly over $q \in \mathcal{Q}$ and $\bx \in \mathcal{X}$, and the derivative of $y \mapsto f_{Y|X}(y | \bx)$ is continuous and bounded in absolute value from above uniformly over $y \in \mathcal{Y}_{\bx}$ and $\bx \in \mathcal{X}$.
    Then Assumptions \ref{Assumption: loss function}--\ref{Assumption: simple setup} and Condition \emph{\ref{eq: strong approx,  lipschitz condition}} in Theorem \ref{thm: strong approximation} hold.
  \end{prop}

The additional conditions in this proposition are primitive and easy-to-interpret, only restricting the conditional density of $Y$ given $\bX$ to be bounded and smooth in a mild sense. Our assumptions are on par with or are weaker than those imposed in \cite{Belloni-Chernozhukov-Chetverikov-FernandezVal_2019_JoE}, despite the high level of generality of our theoretical results.

We can implement uniform inference following the plug-in method described in Section~\ref{sec: Feasible Uniform Inference}. In this context $S_{q, \tilde{q}}(\bx) = q \wedge \tilde{q} - q \tilde{q}$ is known and constant in $\bx$, so a natural plug-in estimator of $\bar\bSigma_{q, \tilde{q}}$ is
\begin{align*}
\widehat\bSigma_{q,\tilde{q}} = (q \wedge \tilde{q} - q \tilde{q}) \E_n\big[\eta^{(1)}(\widehat{\mu}(\bx_i, q)) \eta^{(1)}(\widehat{\mu}(\bx_i, \tilde{q})) \bp(\bx_i) \bp(\bx_i)\trans\big].
\end{align*}
On the other hand, the matrix $\bar\bQ_q$
\begin{align*}
\bar\bQ_q = \E_n\big[\bp(\bx_i)\bp(\bx_i)\trans f_{Y|X}(\eta(\mu_0(\bx_i,q))|\bx_i) [\eta^{(1)}(\mu_0(\bx_i,q))]^2\big]
\end{align*}
depends on the unknown conditional density $f_{Y|X}$, and a plug-in estimator is not immediately available. However, many estimation strategies have been proposed in the literature \citep{Koenker_2005_book}. We do not recommend a particular choice, but rather any estimator satisfying the mild convergence rate requirement in Condition (iii) of Theorem \ref{thm: strong approximation} may be used.

\subsection{Example \ref{example: Conditional Distribution Regression}: Generalized Conditional Distribution Regression}\label{sec: distribution reg}

The loss function is $\rho(y, \eta;q) = ( \I(y \leq q) - \eta )^2$ with a possibly non-identity inverse link function $\eta(\cdot)$. The derivative function is $\psi(y,\eta;q)=-2(\I(y\leq q)-\eta)$. The following proposition, proven in the supplemental appendix (Proposition~\saref{prop:gl-holder-distribution-regression}), verifies our high-level assumptions under mild regularity conditions on the conditional distribution function of $Y$ given $\bX$.

\begin{prop}[Distribution Regression]\label{lem: distribution reg}
  Let $\mathcal{Q} = [-A, A]$ for some $A > 0 $. Suppose that Assumption \ref{Assumption: dgp-holder} holds, the loss is $\rho(y, \eta;q) = \left( \I(y \leq q) - \eta \right)^2$, $\eta(\cdot)\colon \mathbb{R} \to (0,1)$ is strictly monotonic and twice continuously differentiable, $|\eta^{(1)}(\zeta)|$ is bounded away from zero uniformly over $\bx \in \mathcal{X}$, $q \in \mathcal{Q}$, and $\zeta$ satisfying $|\zeta-\mu_0(\bx,q)|\leq r$ for some $r>0$, $\bx \mapsto F_{Y|X}(q|\bx)$ is a continuous function, and $F_{Y|X}(q|\bx) = \eta(\mu_0(\bx, q))$ lies in a compact subset of $(0,1)$ for all $q\in\mathcal{Q}$ and $\bx\in\mathcal{X}$ (this subset does not depend on $q$ and $\bx$).
  Then Assumptions \ref{Assumption: loss function}--\ref{Assumption: simple setup} and Condition \emph{\ref{eq: strong approx,  lipschitz condition}} in Theorem \ref{thm: strong approximation} hold.
\end{prop}

The implementation of uniform inference follows the plug-in method described in Section~\ref{sec: Feasible Uniform Inference}. To construct the prerequisite estimators, in this case $S_{q, \tilde{q}}(\bx_i) = 4 F_{Y|X}(q \wedge \tilde{q} | \bx_i) \big(1 - F_{Y|X}(q \vee \tilde{q} | \bx_i)\big)$. Therefore, a simple plug-in estimator of $\bar\bSigma_{q,\tilde{q}}$ is
\begin{align*}
    \widehat\bSigma_{q,\tilde{q}}
    = 4\E_n\big[\bp(\bx_i)\bp(\bx_i)\trans\eta(\widehat{\mu}(\bx_i,q\wedge\tilde{q}))(1-\eta(\widehat\mu(\bx_i,q\vee\tilde{q})))\eta^{(1)}(\widehat\mu(\bx_i,q))\eta^{(1)}(\widehat\mu(\bx_i,\tilde{q}))\big].
\end{align*}
In addition, a plug-in estimator of the matrix $\bar\bQ_q$ is $\widehat\bQ_q=2\E_n[(\eta^{(1)}(\widehat\mu(\bx_i,q)))^2 \bp(\bx_i)\bp(\bx_i)\trans]$.

\subsection{Example \ref{example: $L_p$ regression}: Generalized \texorpdfstring{$L_p$}{Lp} Regression}\label{sec: lp reg}

The loss function is $\rho(y, \eta) = | y - \eta |^p$, $p \in(1,2]$ with a possibly non-identity link. The case $p=1$ is equivalent to quantile (median) regression discussed previously. The derivative function is $\psi(y,\eta)\equiv\psi(y-\eta)=p|y-\eta|^{p-1}\sgn(\eta-y)$. In this example the family $\mathcal{Q}$ of the loss functions is a singleton, and hence the dependence on the index $q$ can be dropped to simplify notation.

The following proposition, proven in the supplemental appendix (Proposition~\saref{prop:lp-regression}), provides a set of simple regularity conditions that ensure our general theory can be applied to study generalized $L_p$ regression estimation and inference.

\begin{prop}[\texorpdfstring{$L_p$}{Lp} Regression]\label{lem: lp reg}
  Suppose that Assumption \ref{Assumption: dgp-holder} holds with the loss function $\rho(y, \eta) = \left| y - \eta \right|^p$, $p \in(1,2]$, and $\eta(\cdot)\colon \mathbb{R} \to \mathcal{E}$ is strictly monotonic and twice continuously differentiable with $\mathcal{E}$ an open connected subset of $\mathbb{R}$, and that $|\eta^{(1)}(\zeta)|$ is bounded away from zero uniformly over $\bx \in \mathcal{X}$ and $\zeta \in B(\bx)$, where $B(\bx)=\{\zeta:|\zeta-\mu_0(\bx)|\leq r\}$ for some $r>0$. Denoting by $a_l$ and $a_r$ the left and right ends of $\mathcal{E}$ respectively (possibly $\pm \infty$), assume that $\int_{\mathbb{R}} \psi(y;a_l) f_{Y |X}(y|x)dy < 0$ if $a_l$ is finite, and $\int_{\mathbb{R}} \psi(y;a_r)f_{Y |X}(y|x)dy > 0$ if $a_r$ is finite. Also assume that $\E[|Y|^{\nu (p - 1)}]<\infty$ for some $\nu > 2$, and that $\bx \mapsto f_{Y|X}(y | \bx)$ is continuous for any $y \in \mathcal{Y}$. In addition, assume that for each $\bx$, the map $u \mapsto \int_{\mathbb{R}} |u - y |^{p - 1} \sgn(u - y) f_{Y|X}(y|\bx)\,\mathrm{d}y$ is twice continuously differentiable with derivatives $\frac{\mathrm{d}^j}{\mathrm{d} u^j} \int_{\mathbb{R}} |u - y |^{p - 1} \sgn(u - y) f_{Y|X}(y|\bx)\,\mathrm{d}y = \int_{\mathbb{R}} |u - y |^{p - 1} \sgn(u - y) \frac{\partial^j}{\partial y^j} f_{Y|X}(y|\bx)\,\mathrm{d}y$ for $j \in \{1, 2\}$.
  Moreover, the function $\int_{\mathbb{R}} |\eta(\zeta) - y |^{p - 1} \sgn(\eta(\zeta) - y) \frac{\partial}{\partial y}f_{Y|X}(y|\bx)\,\mathrm{d}y$ is bounded and bounded away from zero uniformly over $\bx \in \mathcal{X}$ and $\zeta \in B(\bx)$, and the function $\int_{\mathbb{R}} | \eta(\zeta) - y |^{p - 1} \sgn(\eta(\zeta) - y) \frac{\partial^2}{\partial y^2}f_{Y|X}(y|\bx)\,\mathrm{d}y$ is bounded in absolute value uniformly over $\bx \in \mathcal{X}$ and $\zeta \in B(\bx)$.
  Then Assumptions \ref{Assumption: loss function}--\ref{Assumption: simple setup} and Condition \emph{\ref{eq: strong approx,  lipschitz condition}} in Theorem \ref{thm: strong approximation} hold.\sloppy
\end{prop}

For implementation, we follow the plug-in method in Section~\ref{sec: Feasible Uniform Inference}. Since $\mathcal{Q}$ is a singleton, dependence on $q$ can be dropped. Direct plug-in choices for estimating the prerequisite matrices take the form
\[
\widehat\bQ=\E_n[\bp(\bx_i)\bp(\bx_i)\trans\widehat\Psi_{1, i}
[\eta^{(1)}(\widehat\mu(\bx_i))]^2]\quad \text{and}\quad
\widehat\bSigma=\E_n[\bp(\bx_i)\bp(\bx_i)\trans\psi(\widehat\epsilon_i)^2[\eta^{(1)}(\widehat\mu(\bx_i))]^2],
\]
where $\widehat\epsilon_i=y_i-\eta(\widehat\mu(\bx_i))$ and
$\widehat\Psi_{1, i}$ is some estimator of the function $\Psi_1(\bx_i,\eta(\mu_0(\bx_i)))$. In $L_p$ regression with $p\in(1,2]$,
$\Psi_1(\bx,\eta)=p(p-1)\E[|Y-\eta|^{p-2}\sgn(\eta-Y)|\bX=\bx]$, and therefore a simple plug-in choice is $\widehat\Psi_{1, i} =p(p-1)|y_i-\eta(\widehat\mu(\bx_i))|^{p-2}\sgn(\eta(\widehat\mu(\bx_i))-y_i)$. As an alternative, bootstrap-based inference could be used.

\subsection{Example \ref{example: Generalized linear models}: Logistic Regression}\label{sec: logit est}

For this final example, the loss function is $\rho(y,\eta)=-y\log \eta-(1-y)\log (1-\eta)$, the inverse link function is $\eta(\theta)=1/(1+e^{-\theta})$, and the derivative function is $\psi(y,\eta)=- y / \eta + (1-y) / (1-\eta)$, and the loss function does not depend on $q\in\mathcal{Q}$. The following proposition, proven in the supplemental appendix (Proposition~\saref{prop:gl-holder-logistic-regression}), gives simple primitive conditions verifying the high-level assumptions for our general theoretical results.

\begin{prop}[Logistic Regression]\label{lem: logit est}
    Suppose that Assumption \ref{Assumption: dgp-holder} holds with the loss function $\rho(y,\eta)=-y\log \eta-(1-y)\log (1-\eta)$ and the inverse link $\eta(\theta)=1/(1+e^{-\theta})$; $\mathcal{Y}=\{0,1\}$; $\P(Y = 1 | \bX = \bx)$ is continuous and lies in the interval $(0, 1)$
    for all $\bx \in \mathcal{X}$.
    Then Assumptions \ref{Assumption: loss function}--\ref{Assumption: simple setup} and Condition \emph{\ref{eq: strong approx,  lipschitz condition}} in Theorem \ref{thm: strong approximation} hold.
\end{prop}

It is easy to construct a feasible Gaussian process $\widehat{Z}(\bx)$ conditional on the data $\bD_n$ with covariance structure \eqref{eq: feasible cov}. Standard choices are
\begin{align*}
    \widehat\bQ = \E_n[\bp(\bx_i)\bp(\bx_i)\trans\widehat{\eta}_i(1-\widehat{\eta}_i)]
    \qquad\text{and}\qquad
    \widehat\bSigma = \E_n[\bp(\bx_i)\bp(\bx_i)\trans\widehat\epsilon_i^2],
\end{align*}
where $\widehat{\eta}_i=\eta(\widehat\mu(\bx_i))$ and $\widehat\epsilon_i=y_i-\widehat\eta_i$. See Section \ref{sec: Feasible Uniform Inference} for more discussion.


\section{Other Parameters of Interest}\label{sec: Extensions}

We focused on uniform estimation and inference for the unknown function $\mu_0$ and derivatives thereof. However, the parameter of interest may be other linear or nonlinear transformations of $\mu_0$. For example, in generalized linear models usually the goal is to estimate the function $\eta(\mu_0(\bx, q))$, or the marginal effect of a covariate on that function $\frac{\partial}{\partial x_k}\eta(\mu_0(\bx,q))=\eta^{(1)}(\mu_0(\bx,q))\mu_0^{(\be_k)}(\bx,q)$. Furthermore, in treatment effect and causal inference settings \citep{Abadie-Cattaneo_2018_ARE}, and references therein, interest often lies in differences of such estimands across two or more subgroups: for two treatment levels $j=1,2$, $\eta(\mu_2(\bx,q))-\eta(\mu_1(\bx,q))$ can be interpreted as a mean, quantile, or other conditional (on $(\bx,q)$) treatment effect, where $\mu_j(\bx,q)$ is estimated using separately the subsample of, say, control ($j=1$) and treated ($j=2$) units. Our results can be applied to all these cases of practical interest with minimal additional effort.

We showcase the generality of our theory by briefly discussing uniform inference on the transformed function
$\eta(\mu_0(\bx,q))$, its first derivative, and differences thereof across subgroups. Given the partitioning-based $M$-estimators $\widehat\mu(\bx,q)$ and $\widehat\mu_j(\bx,q)$, $j=1,2$, where $\widehat\mu_j$ is constructed using only data from the subsample $j$ of the full sample, we can immediately plug in to form the desired estimators.
\begin{itemize}
    \item \textit{Level Estimator}: $\eta(\widehat\mu(\bx,q))$.
    \item \textit{Marginal Effect Estimator}: $\eta^{(1)}(\widehat\mu(\bx,q))\widehat\mu^{(\be_k)}(\bx,q)$.
    \item \textit{Conditional Treatment Effect Estimator}: $\eta(\widehat\mu_1(\bx,q))-\eta(\widehat\mu_2(\bx,q))$.
\end{itemize}
Uniform consistency of the three estimators follows from uniform consistency of $\widehat\mu(\bx,q)$ (Corollary \ref{coro: convergence}) because the transformation function $\eta$ is twice continuously differentiable. A Bahadur representation for each of the transformation estimators can be established via Theorem \ref{theorem: bahadur-repres} and a Taylor expansion. For example, for the level estimator,
\begin{align*}
    \sup_{q\in\mathcal{Q}}\sup_{\bx\in\mathcal{X}}
    \big| \eta(\widehat{\mu}(\bx,q))-\eta(\mu_0(\bx,q)) - \mathsf{L}_{\tt LE}(\bx,q) \big|
    \lesssim_\P r_{\tt LE}
\end{align*}
with
\begin{align*}
    \mathsf{L}_{\tt LE}(\bx,q)
    = -\eta^{(1)}(\mu_0(\bx, q)) \bp(\bx)\trans \bQ_{0,q}^{-1}
    \E_{n}\big[\bp(\bx_i) \eta^{(1)}(\mu_{0}(\bx_i,q)) \psi(y_i, \eta(\mu_{0}(\bx_i,q));q)\big],
\end{align*}
and for the marginal effect of the $k$th covariate,
\begin{align*}
    \sup_{q\in\mathcal{Q}}\sup_{\bx\in\mathcal{X}}
    \big| \eta^{(1)}(\widehat{\mu}(\bx,q))\widehat{\mu}^{(\be_k)}(\bx,q)
         - \eta^{(1)}(\mu_0(\bx,q))\mu_0^{(\be_k)}(\bx,q)
         - \mathsf{L}_{\tt ME}(\bx,q) \big|
    \lesssim_\P r_{\tt ME}
\end{align*}
with
\begin{align*}
    \mathsf{L}_{\tt ME}(\bx,q)
    = -\eta^{(1)}(\mu_0(\bx, q))
    \bp^{(\be_k)}(\bx)\trans \bQ_{0,q}^{-1}
    \E_{n}\big[\bp(\bx_i) \eta^{(1)}(\mu_{0}(\bx_i,q)) \psi(y_i, \eta(\mu_{0}(\bx_i,q));q)\big],
\end{align*}
where the approximation remainders from the Taylor expansion, and their uniform rates $r_{\tt LE}$ and $r_{\tt ME}$, are precisely characterized in the supplemental appendix (Theorem~\saref{th:other-parameters}). The conditional treatment effect estimator is simply a difference of two level estimators, each employing a disjoint sub-sample, and therefore it follows directly that
\begin{align*}
    \sup_{q\in\mathcal{Q}}\sup_{\bx\in\mathcal{X}}
    \big| (\eta(\widehat\mu_2(\bx,q))-\eta(\widehat\mu_1(\bx,q)))
         - (\eta(\mu_2(\bx,q))-\eta(\mu_1(\bx,q)))
    - \mathsf{L}_{\tt CTE}(\bx,q) \big|
    \lesssim_\P r_{\tt LE}
\end{align*}
with $\mathsf{L}_{\tt CTE}(\bx,q) = \mathsf{L}_{{\tt LE},2}(\bx,q) - \mathsf{L}_{{\tt LE},1}(\bx,q)$ with $\mathsf{L}_{{\tt LE},j}(\bx,q)$ denoting the Bahadur approximation $\mathsf{L}_{\tt LE}(\bx,q)$ but when only using the sub-sample $j$.

Given the uniform Bahadur representations for each of the transformation estimators, strong approximations of their corresponding $t$-statistic processes can be constructed as in Section \ref{sec: uniform inference}. For example, conditional on $\bX_n$, the stochastic process $(\mathsf{L}_{\tt LE}(\bx,q) : (\bx,q)\in\mathcal{X}\times\mathcal{Q})$ has mean zero and variance $|\eta^{(1)}(\mu_0(\bx, q))|^2\bar\Omega_{\bm{0}}(\bx,q)/n$. Then, applying our strong approximation strategy, we can construct a conditional Gaussian process $Z_{\tt LE}(\bx, q)$ that approximates the $t$-statistic process of $\eta(\widehat{\mu}(\bx,q))$:
\begin{align*}
    \sup_{q\in\mathcal{Q}}\sup_{\bx\in\mathcal{X}}
    \Big| \frac{\eta(\widehat{\mu}(\bx,q))-\eta(\mu_0(\bx,q))}
               {|\eta^{(1)}(\mu_0(\bx,q))|\sqrt{\bar\Omega_{\bm{0}}(\bx,q)/n}}
          - Z_{\tt LE}(\bx, q) \Big|
    \lesssim_\P r_{\tt SALE}
\end{align*}
with strong approximation rate $r_{\tt SALE}$ as in Theorem \ref{thm: strong approximation}. Similarly, we can also construct a conditional Gaussian process $Z_{\tt ME}(\bx,q)$ that approximates the $t$-statistic process of the marginal effect estimator $\frac{\partial}{\partial x_k}\eta(\widehat\mu(\bx,q))$:
\begin{align*}
    \sup_{q\in\mathcal{Q}}\sup_{\bx\in\mathcal{X}}
    \Big| \frac{\eta^{(1)}(\widehat{\mu}(\bx,q)) \widehat{\mu}^{(\be_k)}(\bx,q)
                 - \eta^{(1)}(\mu_0(\bx,q))\mu_0^{(\be_k)}(\bx, q)}
               {|\eta^{(1)}(\mu_0(\bx,q))|\sqrt{\bar\Omega_{\be_k}(\bx,q)/n}}
          - Z_{\tt ME}(\bx, q) \Big|
    \lesssim_\P r_{\tt SAME}
\end{align*}
with strong approximation rate $r_{\tt SAME}$ as in Theorem \ref{thm: strong approximation}. These results are formalized in the supplemental appendix (Theorem~\saref{th:other-parameters}). An analogous result holds for the conditional treatment effect estimator.

Finally, for implementation we can construct feasible processes to approximate $Z_{\tt LE}(\bx, q)$ and $Z_{\tt ME}(\bx, q)$ via the plug-in method discussed in Section \ref{sec: Feasible Uniform Inference}, and illustrated in Section \ref{sec: Examples}, which then can be employed to approximate the distributions of the entire level process $(\eta(\widehat\mu(\bx,q)):(\bx,q)\in\mathcal{X}\times\mathcal{Q})$, marginal effect process $(\frac{\partial}{\partial x_k}\eta(\widehat\mu(\bx,q)):(\bx,q)\in\mathcal{X}\times\mathcal{Q})$, and conditional treatment effect process $(\eta(\widehat\mu_2(\bx,q)) - \eta(\widehat\mu_1(\bx,q)) :(\bx,q)\in\mathcal{X}\times\mathcal{Q})$.

\section{Conclusion}\label{sec: conclusion}

This paper investigated the asymptotic properties of a large class of nonparametric partitioning-based M-estimators, allowing for different degrees of non-smoothness in the loss function and a possibly non-identity monotonic transformation function. Our main theoretical results include uniform consistency for convex and non-convex objective functions, uniform Bahadur representations with optimal remainder under appropriate conditions, uniform and mean square convergence rates achieving optimal approximation under appropriate conditions, uniform strong approximation methods under general conditions, and uniform inference methods via plug-in approximations. We illustrated our general theory with four examples, and demonstrated how our results improve on the prior literature, in many cases requiring minimal side rate restrictions on tuning parameters and achieving rate-optimal approximation rates. The supplemental appendix collects further theoretical results and generalizations that may be of independent interest. In future work, we plan to investigate optimal tuning parameter selection, including random partitioning schemes, and the validity of bootstrap-based approximations.

\begin{acks}[Acknowledgments]
We thank Richard Crump, Max Farrell, Will Underwood, and Rae Yu for helpful comments and discussions. We also thank the Co-Editor, Associate Editor, and two reviewers for their comments.
\end{acks}

\begin{funding}
Cattaneo gratefully acknowledges financial support from the National Science Foundation through grants DMS-2210561, and SES-2241575, and SES-2342226, the National Institute for Food and Agriculture (NIFA) through grant 2024-67023-42704, and the John Simon Guggenheim Memorial Foundation through a 2026 Guggenheim Fellowship.

Feng gratefully acknowledges financial support from the National Natural Science Foundation of China (NSFC) through grants 72522001 and 72203122.
\end{funding}



%

\begin{supplement}
  \stitle{Proofs and technical results}
  \slink[doi]{10.1214/26-AOS2623SUPP}
  \sdatatype{.pdf}
  \sfilename{aos2623supp.pdf}
  \sdescription{The supplementary material
    \citep{Cattaneo-Feng-Shigida_2026_AOS--SA}
    contains more general results encompassing those reported in the paper and all proofs, omitted technical and methodological details, other theoretical results that may be of independent interest, and simulation evidence.}
\end{supplement}


\bibliographystyle{imsart-nameyear}
\bibliography{bib.bib}       


\end{document}


\title{Supplement to ``Uniform Estimation and Inference\\
       for Nonparametric Partitioning-Based M-Estimators''
\bigskip }
\author{Matias D. Cattaneo\thanks{Department of Economics, Princeton University.} \and
	    Yingjie Feng\thanks{School of Economics and Management, Tsinghua University.} \and
        Boris Shigida\thanks{Department of Operations Research and Financial Engineering, Princeton University.}
        }
\maketitle

\begin{abstract}
    This supplement is self-contained. It presents more general theoretical results than those reported in the paper, as well as their proofs. In particular, a larger class of loss functions is allowed for, and weaker regularity conditions are employed when possible. Additional results not reported in the paper and their proofs are given, and some of the technical results presented herein may be of independent theoretical interest. Omitted details on the motivating examples are also provided. Finally, simulation evidence is reported.
\end{abstract}

\thispagestyle{empty}
\clearpage

\pagenumbering{roman}
\maketitle

\setcounter{tocdepth}{2}

\tableofcontents

\clearpage

\listoftheorems[
title={List of important statements},
ignore={lemma,remark,remark*},
numwidth=4em,
]

\clearpage

\listoftheorems[
title={List of lemmas},
ignoreall,
show={lemma},
numwidth=4em,
]

\clearpage

\pagenumbering{arabic}

\ifthenelse{\boolean{standalone}}{
  \section{Introduction}
}{
  \section{Introduction to supplemental material}
}

\begin{pycontext}
import os
import sys
os.chdir("code/")
sys.path.insert(0, os.getcwd())
from report import *
\end{pycontext}

Let $\mathcal{Q} \subset \mathbb{R}^{d_{\mathcal{Q}}}$, $\mathcal{X} \subset \mathbb{R}^d$ be fixed compact sets, where $d_{\mathcal{Q}}$ and $d$ are positive integers. (In the paper, $d_{\mathcal{Q}}=1$ was set only for simplicity.) Suppose that $\prn[\big]{(y_{i}, \vx_{i})}_{1 \leq i \leq n}$ is a random sample, where $y_i \in \mathcal{Y} \subset \mathbb{R}$ is a scalar response variable, $\vx_i$ is a $d$-dimensional covariate with values in $\mathcal{X}$. Let $\riskfunc[\vq]{\cdot}{\cdot}$ be a loss function parametrized by $\vq \in \mathcal{Q}$ (Borel-measurable in all three arguments), and let $\eta(\cdot)\colon \mathbb{R} \to \rangeeta$ be a strictly monotonic continuously differentiable transformation function. (More detailed assumptions on $\riskfunc[\vq]{\cdot}{\cdot}$ and $\eta(\cdot)$ will be provided below.) We fix a function $\muzero[\cdot]{\cdot}\colon \mathcal{X} \times \mathcal{Q} \to \mathbb{R}$, Borel-measurable in both arguments and such that
\begin{equation}\label{eq:holder-rho-problem}
  \muzero[\vq]{\cdot} \in \argmin_{\mu} \mathbb{E}[\riskfunc[\vq]{y_1}{\eta(\mu(\vx_1))}],
\end{equation}
where the argmin is over the space of Borel functions $\mathcal{X} \to \mathbb{R}$. In particular, we assume that the minimum is finite, and such a minimizer exists.

Our main goal is to conduct uniform (over $\mathcal{X}\times\mathcal{Q}$) estimation and inference for $\mu_{0}$, and transformations thereof, employing the partitioning-based series $M$-estimator
\begin{equation}\label{eq:m-estimator}
    \muhat[\vq]{\vx} = \bp(\vx)\trans
    \betahat[\vq],\qquad
    \betahat[\vq] \in \argmin_{\bb\in\mathcal{B}} \sum_{i=1}^n \riskfunc[\vq]{y_i}{\eta(\bp(\vx_i)\trans \bb)},
\end{equation}
where $\mathcal{B} \subseteq \mathbb{R}^K$ is the feasible set of the optimization problem, and
\begin{align*}
    \vx \mapsto \bp(\vx) = \bp(\vx; \Delta, m) = \big( p_1(\vx; \Delta, m), \ldots, p_K(\vx; \Delta, m) \big)\trans
\end{align*} is a dictionary of $K$ locally supported basis functions of order $m$ based on a quasi-uniform partition $\Delta = \crl{\delta_l: 1 \leq l \leq \kappa}$ containing a collection of open disjoint polyhedra in $\mathcal{X}$ such that the closure of their union covers $\mathcal{X}$. The $m$ parameter controls how well $\mu_0$ can be approximated by linear combinations of the basis (\cref{ass:gl-approximation-error} below); the partition being quasi-uniform intuitively means that the largest size of a cell cannot get asymptotically bigger than the smallest one (Assumption~\ref{ass:gl-quasi-uniform-partition} below). We consider large sample approximations where $d$ and $m$ are fixed constants, and $\kappa \to \infty$ (and thus $K\to\infty$) as $n\to\infty$. Prior literature is discussed in the paper.

\subsection{Notation}

For any real function $f$ depending on $d$ variables $(t_1, \ldots, t_d)$ and any vector $\vv$ of nonnegative integers, denote
\begin{equation*}
  f^{(\vv)} := \partial^{\vv} f := \frac{\partial^{\multindnorm{\vv}}}{\partial t_1^{v_1} \ldots \partial t_d^{v_d}} f.
\end{equation*}
the multi-indexed partial derivative of $f$, where $\multindnorm{\vv} = \sum_{k = 1}^d v_k$. A derivative of order zero is the function itself, so if $v_i = 0$, the $i$th partial differentiation is ignored. For functions that depend on $(\vx, \vq)$, the multi-index derivative notation is taken with respect to the first argument $\vx$, unless otherwise noted.

We will denote $N(\mathcal{F}, \rho, \varepsilon)$ the $\varepsilon$-covering number of a class $\mathcal{F}$ with respect to a semi-metric $\rho$ defined on it.

For a function $f\colon S \to \mathbb{R}$ the set $\{(x, t) \in S \times \mathbb{R}: t < f(x)\}$ is called the \textit{subgraph} of $f$. A class $\mathcal{F}$ of measurable functions from $\mathcal{S}$ to $\mathbb{R}$ is called a \textit{VC-subgraph class} or \textit{VC-class} if the collection of all subgraphs of functions in $\mathcal{F}$ is a VC-class of sets in $S \times \mathbb{R}$, which means that for some finite $m$ no set of size $m$ is shattered by it. In this case, the smallest such $m$ is called the VC-index of $\mathcal{F}$. See \parencite{vaart96_weak_conver_empir_proces} for details.

For a measurable function $f\colon S \to \mathbb{R}$ on a measurable space $(S, \mathcal{S})$, a probability measure $\mathbb{Q}$ on this space and some $q \geq 1$, define the \textit{$(\mathbb{Q}, q)$-norm} of $f$ as $\|f\|_{\mathbb{Q}, q}^q= \mathbb{E}_{X \sim \mathbb{Q}}[f(X)^q]$.

We will say that a class of measurable functions $\mathcal{F}$ from any set $S$ to $\mathbb{R}$ has a measurable \textit{envelope} $F$ if $F\colon S \to \mathbb{R}$ is such a measurable function that $|f(s)| \leq F(s)$ for all $s \in S$ and all $f \in \mathcal{F}$. We will say that this class satisfies the \textit{uniform entropy bound} with envelope $F$ and real constants $A \geq e$ and $V \geq 1$ if
\begin{equation}\label{eq:uniform-entropy-bound}
  \sup_{\mathbb{Q}} N(\mathcal{F}, \norm{\cdot}_{\mathbb{Q}, 2}, \varepsilon \norm{F}_{\mathbb{Q}, 2}) \leq \prn[\bigg]{\frac{A}{\varepsilon}}^{V}
\end{equation}
for all $0 < \varepsilon \leq 1$, where the supremum is taken over all finite discrete probability measures $\mathbb{Q}$ with $\norm{F}_{\mathbb{Q}, 2} > 0$, $\norm{\cdot}_{\mathbb{Q}, 2}$ denotes the $(\mathbb{Q}, 2)$-norm.

We will say that an $\mathbb{R}$-valued random variable $\xi$ is \textit{$\sigma^2$-sub-Gaussian}, where $\sigma^2 > 0$, if
\begin{equation}\label{eq:def-subgaussian}
  \Prob*{\abs{\xi} \geq t} \leq 2 \exp\crl{-t^2 / \sigma^2}\quad\text{for all $t \geq 0$}.
\end{equation}

We will denote by $\newtarget{def:data}{\dataNoLink}$ the random vector of all the data $\crl{\vx_i, y_i}_{i = 1}^n$, $\mX_n$ for $\crl{\vx_i}_{i = 1}^n$ and $\vy_n$ for $\crl{y_i}_{i = 1}^n$.

For two random elements $Z_1$ and $Z_2$, the notation $Z_1 \independent Z_2$ means they are independent, and the notation $Z_1 \underset{\xi}{\independent} Z_2$ means they are independent given a third random element $\xi$.

The notation $Z_1 \overset{d}{=} Z_2$ means the two random elements $Z_1$ and $Z_2$ have the same laws.

If we say the probability space is \textit{rich enough}, it means that it admits a randomization variable, that is, a random variable distributed uniformly on $[0, 1]$ independent of the data.
Since this one random variable can be replicated (see, e.\,g., Lemma 4.21 in \cite{Kallenberg2021}), we can find a random variable distributed uniformly on $[0, 1]$ independent of the data and such random variables previously used, whenever the argument requires.

Finally, we will use the following notations:
\begin{align}
  &\qbar[\vq] := \Expectlet_n \brk{\vp(\vx_i) \vp(\vx_i) \trans \Psione[\vq]{\vx_i}{\eta(\muzero[\vq]{\vx_i})} \eta^{(1)}(\muzero[\vq]{\vx_i})^{2}}\label{eq:qbar-def}\\
  &\qzero[\vq] := \Expect{\vp(\vx_i) \vp(\vx_i) \trans \Psione[\vq]{\vx_i}{\eta(\muzero[\vq]{\vx_i})} \eta^{(1)}(\muzero[\vq]{\vx_i})^{2}}\label{eq:qzero-def}\\
  &\sigmazero[\vq] := \Expect{\vp(\vx_i) \vp(\vx_i) \trans \riskder[\vq]{y_i}{\eta(\muzero[\vq]{\vx_i})}^2 \eta^{(1)}(\muzero[\vq]{\vx_i})^{2}}\label{eq:sigma-zero}\\
  &\sigmabar[\vq] := \Expectlet_n \brk{\vp(\vx_i) \vp(\vx_i) \trans \sigma^2_{\vq}(\vx_i) \eta^{(1)}(\muzero[\vq]{\vx_i})^{2}}\label{eq:sigma-bar}\\
  &\bigomega[\vq] := \vp^{(\vv)}(\vx)\trans \qzero[\vq][-1] \sigmazero[\vq] \qzero[\vq][-1] \vp^{(\vv)}(\vx)\label{eq:bigomega-def}\\
  &\bigomegabar[\vq] := \vp^{(\vv)}(\vx)\trans \qbar[\vq][-1] \sigmabar[\vq] \qbar[\vq][-1] \vp^{(\vv)}(\vx)\label{eq:bigomegabar-def},
\end{align}
where $\riskder[\cdot]{\cdot}{\cdot}$, $\Psione[\cdot]{\cdot}{\cdot}$, $\sigma^2_{\cdot}(\cdot)$ are defined in \cref{ass:gl-holder}.

Many mathematical notations are clickable and link to their definitions; for example, clicking on $\sigmabar[\vq]$ should lead to \eqref{eq:sigma-bar}.

\subsection{Organization}\label{sec:organization}

\Cref{sec:setup} collects the general assumptions used through this supplemental appendix. In subsequent sections, we will list which assumptions are required for each lemma, proposition, or theorem. Note that \cref{ass:gl-holder} is weaker than the assumptions on the loss function described in Section~\mainref{sec: Assumptions} of the main paper. In particular, a complexity condition on $\{\riskder[\vq]{\cdot}{\cdot}\}$ is absent in Section~\mainref{sec: Assumptions} because the main paper contains a stronger but simpler Assumption~\mainref{Assumption: simple setup}. Thus, the setup in this supplement is more general than in the main paper; this fact is formally shown by \cref{prop:simpler-generalization}. The correspondence of assumptions given in the paper and the versions given in the supplement is as follows:
\begin{align*}
  \text{Assumption~\mainref{Assumption: gl-quasi-uniform-partition}} &\Leftrightarrow \text{\cref{ass:gl-quasi-uniform-partition}},\\
  \text{Assumption~\mainref{Assumption: gl-local-basis}} &\Leftrightarrow \text{\cref{ass:gl-local-basis}},\\
  \text{Assumption~\mainref{Assumption: dgp-holder}} &\Leftrightarrow \text{\cref{ass:dgp-holder}},\\
  \text{Assumption~\mainref{Assumption: loss function}} &\Leftrightarrow \text{\cref{ass:gl-holder}},\\
  \text{Assumption~\mainref{Assumption: simple setup}} &\Rightarrow \text{\cref{ass:complexity,ass:complexity-addition}} \quad  \text{(by \cref{prop:simpler-generalization}),}\\
\text{Assumption~\mainref{Assumption: gl-approximation-error}} &\Leftrightarrow \text{\cref{ass:gl-approximation-error}}.
\end{align*}

\Cref{sec:preliminaries} records some well-known results and tools in probability and stochastic process theory that will be useful for our theoretical analysis, and also presents a collection of lemmas that are used repeatedly in many subsequent arguments throughout this supplement.

The rest of the supplement presents our main theoretical results, which encompass and generalize the simplified results presented in the main paper (to improve exposition). In the rest of this subsection, we explain how the two sets of theoretical results relate to each other. All equivalences below are given up to the fact that the supplement contains weaker complexity assumptions \cref{ass:complexity,ass:complexity-addition}, which are not present in the paper because the simplified setup (Assumption~\mainref{Assumption: simple setup}) is used.

\Cref{sec:consistency} presents all our consistency results,
categorized based on whether the objective function is convex or not. Additional results of theoretical and methodological interest, such as the consistency of partitioning-based M-estimators in special cases (e.g., an unconnected basis or a strongly convex and smooth loss function), which were not formally reported in the paper, are also presented. The matching of consistency results given in the paper with the versions given in the supplement is as follows:
\begin{align*}
  \text{Lemma~\mainref{lem:consistency, convex}} &\Leftrightarrow \text{\cref{lem:holder-consistency}},\\
  \text{Lemma~\mainref{lem:consistency, nonconvex}} &\Leftrightarrow \text{\cref{lem:consistency-nonconvex}}.
\end{align*}

\Cref{sec:bahadur-representation} proves our main Bahadur representation result and its corollaries, matched as follows to the results in the paper:
\begin{align*}
  \text{Theorem~\mainref{theorem: bahadur-repres}} &\Leftrightarrow \text{\cref{th:bahadur-repres}},\\
  \text{Corollary~\mainref{coro: convergence}} &\Leftarrow \text{\cref{cor:uniform-rate-of-convergence}},\\
  \text{Corollary~\mainref{coro: Mean Square Rate of Convergence}} &\Leftrightarrow \text{\cref{cor:msq-rate-of-convergence}}.
\end{align*}

\Cref{sec:strong-approximation} develops strong approximation results using a generalized conditional Yurinskii's coupling approach. First, \cref{sec:conditional-yurinskii-coupling} gives general results that may be of independent theoretical interest: they provide generalizations of, and in some cases complement, prior coupling results established in \parencite{Belloni-Chernozhukov-Chetverikov-FernandezVal_2019_JoE}, \parencite{Cattaneo-Feng-Underwood_2024_JASA}, \parencite{cattaneo2024yurinskii}, and references therein. Second, \cref{sec:SA-main-result} deploys those results to the setting of interest in our paper to verify our main result \cref{th:strong-approximation-binscatter-yurinski} that confirms Theorem~\mainref{thm: strong approximation} in the paper:
\begin{align*}
  \text{Theorem~\mainref{thm: strong approximation}} &\Leftarrow \text{\cref{th:strong-approximation-binscatter-yurinski} (see also \cref{rem:strong-approx-correspondence})}.
\end{align*}

\Cref{sec:uniform-inference} discusses results related to the implementation of uniform inference, formally showing the validity of the plug-in approximation method and confidence bands described in the paper. In particular, the following correspondence holds:
\begin{align*}
\text{Theorem~\mainref{thm: Feasible Plug-in Strong Approximation}} &\Leftarrow \text{\cref{th:plug-in} (see also \cref{rem:plug-in-correspondence})}.
\end{align*}

\Cref{sec:examples} discusses in detail the verification of our high-level assumptions for each of the four motivating examples in the paper, proving Propositions 1 through 4 in the paper along with verifying more general complexity assumptions discussed above. The correspondence between the results in the paper and the results in this section is as follows:
\begin{align*}
  \text{Proposition~\mainref{lem: quantile reg}} &\Leftarrow \text{\cref{prop:gl-holder-quantile-regression}},\\
  \text{Proposition~\mainref{lem: distribution reg}} &\Leftarrow \text{\cref{prop:gl-holder-distribution-regression}},\\
  \text{Proposition~\mainref{lem: lp reg}} &\Leftarrow \text{\cref{prop:lp-regression}},\\
  \text{Proposition~\mainref{lem: logit est}} &\Leftarrow \text{\cref{prop:gl-holder-logistic-regression}}.
\end{align*}

\Cref{sec:other-parameters-of-interest} discusses other parameters of potential interest such as the level and marginal effect functions, which formalizes the claims made in Section~\mainref{sec: Extensions} of the paper.

\Cref{sec:simulation-evidence} is devoted to simulation runs with synthetic data, covering pointwise and uniform inference.

\subsection{Comparison to prior literature}

Our general $M$-estimation theoretical results are on par with or improve upon the prior literature on partitioning-based nonparametric least squares regression \citep{Belloni-Chernozhukov-Chetverikov-Kato_2015_JoE,Cattaneo-Farrell_2013_JoE,Cattaneo-Farrell-Feng_2020_AOS,Cattaneo-Crump-Farrell-Feng_2024_AER,Chen-Christensen_2015_JOE,Huang_2003_AoS,Zhou-Shen-Wolfe_1998_AoS}, and on series nonparametric quantile regression \citep{Belloni-Chernozhukov-Chetverikov-FernandezVal_2019_JoE}. Both strands of the literature assume an identity transformation function $\eta(\cdot)$. In addition, we also compare to the closely related setting of nonparametric smoothing spline nonlinear regression \citep{Shang-Cheng_2013_AOS}.

The main paper gives a detailed discussion of our improvements when discussing the theoretical results for consistency (Section 4), Bahadur representation (Section 5), and strong approximation (Section 6). To complement that discussion, we present three tables that summarize the comparison between prior literature and this paper considering only the sub-exponential moment assumption for simplicity.

Table \ref{tab:CompareLit-Consistency} compares the best known results for mean-square ($L_2$) and uniform ($L_\infty$) consistency of series and partitioning-based estimators. In each case considered in this table, rate-optimal $L_2$ and $L_\infty$ convergence rates are achievable, but with different side rate conditions on $h$, $m$, and $d$. The table shows that our results either achieve the best (minimal) restrictions in the literature for the special case of square loss and identity transformation functions, or improve upon prior results in the general $M$-estimation case.

\renewcommand{\arraystretch}{1.4}

\begin{table}[!h]
\centering
\caption{Comparison to prior literature: Consistency (Section 4)}\label{tab:CompareLit-Consistency}
\resizebox{\textwidth}{!}{
\begin{tabular}{lccc}
\hline\hline
                                                  & Side rate restriction  & $m$ vs. $d$ restriction & Uniformity\\ \hline
Mean regression                                   & $nh^{d}\to\infty$      & None                     & $\bx$ only\\
Quantile regression, $L_2$-consistency            & $nh^{2d}\to\infty$     & $m>d$                    & $\bx$, and $q$ (scalar)\\
Quantile regression, $L_\infty$-consistency       & $nh^{4d}\to\infty$     & $m>d$                    & $\bx$, and $q$ (scalar)\\
\hline
This paper, $L_2$- and $L_\infty$-consistency     &   &  \\
$\,\,\,$ General partitioning basis               & $nh^{2d}\to\infty$     & $m>d/2$                  & $\bx$, and $\bq$\\
$\,\,\,$ Strongly convex and smooth loss function & $nh^{d}\to\infty$      & None                     & $\bx$, and $\bq$\\
$\,\,\,$ Unconnected partitioning basis           & $nh^{d}\to\infty$      & None                     & $\bx$, and $\bq$\\
\hline
\end{tabular}
}
\scriptsize{
    \begin{flushleft}Notes: (i) $m$ is the smoothness of $\bx\mapsto\mu_0(\bx,\bq)$ and order of the basis; (ii) $d$ is the dimension of $\bx_i$; and (iii) poly-log-$n$ and similar additional terms are not reported to simplify the exposition.
    \end{flushleft}
}
\end{table}

Table \ref{tab:CompareLit-Bahadur} compares the best known results for uniform Bahadur representation for series and partitioning-based estimators. It again demonstrates that the results in our paper either achieve the best (minimal) restrictions in the literature for the special case of square loss and identity transformation functions, or improve upon prior literature in the general $M$-estimation case.

\begin{table}[!h]
\centering
\caption{Comparison to prior literature: Uniform Bahadur representation (Section 5)}\label{tab:CompareLit-Bahadur}
\resizebox{\textwidth}{!}{
\begin{tabular}{lcccc}
\hline\hline
                                  & Side rate restriction & $m$ vs. $d$ restriction  & remainder                       & Uniformity\\ \hline
Mean regression                   & $nh^{d}\to\infty$     & None                      & $\frac{1}{nh^{d}}$ (optimal)    & $\bx$ only\\
Quantile regression               & $nh^{4d}\to\infty$    & $m>d$                     & $\frac{1}{(nh^d)^{3/4}h^{d/2}}$ & $\bx$, and $q$ (scalar)
\smallskip \\
Smooth loss + smoothing spline    & $nh^{2d}\to\infty$    & $d=1$                     & $\frac{1}{nh^{10m-1}}$          & $x$, and $q$ (scalar)\\
\hline
This paper (assuming consistency)                      &                      &           &  & \\
$\,\,\,$ Smooth weak derivative ($\alpha=1$)           & $nh^{d}\to\infty$    & None      & $\frac{1}{(nh^d)^{3/4}}$ (optimal) & $\bx$, and $\bq$\\
$\,\,\,$ Non-smooth weak derivative ($\alpha\in(0,1)$) & $nh^{d}\to\infty$    & None      & See Theorem 1                      & $\bx$, and $\bq$\\
\hline
\end{tabular}
}
\scriptsize{
    \begin{flushleft}Notes: (i) $m$ is the smoothness of $\bx\mapsto\mu_0(\bx,\bq)$ and order of the basis; (ii) $d$ is the dimension of $\bx_i$; and (iii) poly-log-$n$ and similar additional terms are not reported to simplify the exposition.
    \end{flushleft}
}
\end{table}

Finally, Table \ref{tab:CompareLit-StrongApprox} compares the best known results for uniform strong approximation for series and partitioning-based estimators. Our results achieve the best (minimal) restrictions in the literature for the special case of square loss and identity transformation functions when $d>1$, and improve upon prior literature in the general $M$-estimation case.

\begin{table}[!h]
\centering
\caption{Comparison to prior literature: Strong approximation (Section 6)}\label{tab:CompareLit-StrongApprox}
\resizebox{\textwidth}{!}{
\begin{tabular}{lccc}
\hline\hline
                               & Side rate restriction        & $m$ vs. $d$ restriction  & Uniformity \\ \hline
Mean regression                &                              &                           &            \\
$\,\,\,$ $d=1$                 & $nh^{d}\to\infty$            & None                      & $\bx$ only \\
$\,\,\,$ $d\geq1$              & $nh^{3d}\to\infty$           & None                      & $\bx$ only \\
Quantile regression ($d\geq1$) & $nh^{4d\vee(2+3d)}\to\infty$ & $m>d$                     & $\bx$, and $q$ (scalar)\\
\hline
This paper (assuming consistency, $d\geq1$)  & $nh^{3d}\to\infty$    & None           & $\bx$, and $\bq$\\
\hline
\end{tabular}
}
\scriptsize{
    \begin{flushleft}Notes: (i) $m$ is the smoothness of $\bx\mapsto\mu_0(\bx,\bq)$ and order of the basis; (ii) $d$ is the dimension of $\bx_i$; and (iii) poly-log-$n$ and similar additional terms are not reported to simplify the exposition.
    \end{flushleft}
}
\end{table}

\ifthenelse{\boolean{standalone}}{
  \section{Assumptions}\label{sec:setup}
}{
  \section{Full assumptions}\label{sec:setup}
}

This section collects the assumptions used throughout the supplemental appendix. These assumptions are weaker than (i.e., implied by) the assumptions imposed in the paper.

\begin{fullassumption}[Quasi-uniform partition]\label{ass:gl-quasi-uniform-partition}
  The ratio of the sizes of inscribed and circumscribed balls of each $\delta \in \Delta$ is bounded away from zero uniformly in $\delta \in \Delta$, and
  \begin{equation*}
    \frac{\max \{\operatorname{diam}(\delta): \delta \in \Delta\}}{\min \{\operatorname{diam}(\delta): \delta \in \Delta\}} \lesssim 1
  \end{equation*}
  where $\operatorname{diam}(\delta)$ denotes the diameter of $\delta$. Further, for $h = \max \{\operatorname{diam}(\delta): \delta \in \Delta\}$, assume $h = o(1)$ and $\log(1 / h) \lesssim \log(n)$.
\end{fullassumption}

\begin{fullassumption}[Local basis]\label{ass:gl-local-basis}\leavevmode
  \begin{enumerate}[wide,label=(\roman*),nosep]
  \item\label{ass:localization} For each basis function $p_k$, $k=1, \ldots, K$, the union of elements of $\Delta$ on which $p_k$ is active is a connected set, denoted by $\mathcal{H}_k$. For all $k=1, \ldots, K$, both the number of elements of $\mathcal{H}_k$ and the number of basis functions which are active on $\mathcal{H}_k$ are bounded by a constant.
  \item\label{ass:well-conditioning} For any $\va = (a_1, \ldots, a_K)\trans \in \mathbb{R}^K$
    \begin{equation*}
      \va\trans \int_{\mathcal{H}_k} \vp(\vx ; \Delta, m) \vp(\vx ; \Delta, m)\trans d \vx \ \va \gtrsim a_k^2 h^d, \quad k=1, \ldots, K.
    \end{equation*}

  \item\label{ass:basis-norm-bound} Let $\multindnorm{\vv} < m$. There exists an integer $\varsigma \in[\multindnorm{\vv}, m)$ such that, for all $\boldsymbol{\varsigma},\multindnorm{\boldsymbol{\varsigma}} \leq \varsigma$,
    \begin{equation*}
      h^{-\multindnorm{\boldsymbol{\varsigma}}} \lesssim \inf _{\delta \in \Delta} \inf _{\vx \in \closure(\delta)}\norm{\vp^{(\boldsymbol{\varsigma})}(\vx ; \Delta, m)} \leq \sup _{\delta \in \Delta} \sup _{\vx \in \closure(\delta)}\norm{\vp^{(\boldsymbol{\varsigma})}(\vx ; \Delta, m)} \lesssim h^{-\multindnorm{\boldsymbol{\varsigma}}}
    \end{equation*}
    where $\closure(\delta)$ is the closure of $\delta$.
  \end{enumerate}
\end{fullassumption}

\Cref{ass:gl-local-basis} implicitly relates the number of basis functions and the maximum mesh size: $K \asymp h^{-d} = J^d$.

For the following assumption and throughout the document, when speaking of the conditional distribution of $y_1$ given $\vx_1$, or its functionals (like conditional moments or quantiles), we mean one fixed regular variant of such a distribution satisfying all the assumptions listed.

\begin{fullassumption}[Data generating process]\label{ass:dgp-holder}\leavevmode
  \begin{enumerate}[wide,label=(\roman*),nosep]
  \item\label{ass:random-sample} $\prn[\big]{(y_{i}, \vx_{i})}_{1 \leq i \leq n}$ is a random sample satisfying~\eqref{eq:holder-rho-problem}.

  \item The distribution of $\vx_i$ admits a Lebesgue density $\fxDef{\cdot}$ which is continuous and bounded away from zero on support $\mathcal{X} \subset \mathbb{R}^d$, where $\mathcal{X}$ is the closure of an open, connected and bounded set.

  \item\label{ass:cond-density} The conditional distribution of $y_i$ given $\vx_i$ admits a conditional density $f_{Y|X}(y | \vx)$ with support $\mathcal{Y}_{\vx}$ with respect to some (sigma-finite) measure $\newtarget{def:measureonY}{\measureonYNoLink}$, and $\sup _{\vx \in \mathcal{X}} \sup _{y \in \mathcal{Y}_{\vx}} f_{Y|X}(y| \vx) < \infty$.

  \item\label{ass:mu-bounded} $\vx \mapsto \muzero[\vq]{\vx}$ is $m \geq 1$ times continuously differentiable for every $\vq \in \mathcal{Q}$, $\vx \mapsto \muzero[\vq]{\vx}$ and its derivatives of order no greater than $m$ are bounded uniformly over $(\vx, \vq) \in \mathcal{X} \times \mathcal{Q}$, and
    \begin{equation*}
      \sup_{\vx \in \mathcal{X}} \sup_{\vq_1 \neq \vq_2} \frac{\abs{\muzero[\vq_1]{\vx} - \muzero[\vq_2]{\vx}}}{\norm{\vq_1 - \vq_2}} \lesssim 1.
    \end{equation*}
  \end{enumerate}
\end{fullassumption}

For the following assumption and throughout the document, we fix some (small enough) $r > 0$ and denote
\begin{equation}\label{eq:ballqx-def}
  \ballqxNoLink := \{\zeta \colon |\zeta - \muzero[\vq]{\vx}| \leq r\},
\end{equation}
i.\,e. we will work with a fixed neighborhood of $\muzero[\vq]{\vx}$.

\begin{fullassumption}[Loss function]\label{ass:gl-holder}\leavevmode
  \begin{enumerate}[wide,label=(\roman*),nosep]
  \item\label{ass:abs-cont} Let $\mathcal{Q} \subset \mathbb{R}^{d_{\mathcal{Q}}}$ be a connected compact set. For each $\vq \in \mathcal{Q}$ and $y \in \mathcal{Y}$, and some open connected subset $\rangeeta$ of $\mathbb{R}$ not depending on $y$, $\eta \mapsto \riskfunc[\vq]{y}{\eta}$ is absolutely continuous on closed bounded intervals within $\rangeeta$, and admits an a.\,e. derivative $\newtarget{def:riskder}{\riskder[\vq]{y}{\eta}}$, Borel measurable as a function of $(y, \eta, \vq)$.

\item\label{ass:psi-moments} The first-order optimality condition $\Expect{\riskder[\vq]{y_i}{\eta(\muzero[\vq]{\vx_i})} \given \vx_i} = 0$ holds; the function
    \begin{equation*}
      \sigma_{\vq}^{2}(\vx) := \Expect[\big]{\riskder[\vq]{y_i}{\eta(\muzero[\vq]{\vx_i})}^{2} \given \vx_i = \vx}
    \end{equation*}
    is continuous in both arguments $(\vx, \vq) \in \mathcal{X} \times \mathcal{Q}$, bounded away from zero, and Lipschitz in $\vq$ uniformly in $\vx$; there is a positive measurable envelope function $\newtarget{def:riskderenvel}{\riskderenvelNoLink{\vx_i}{y_i}}$ such that $\sup_{\vq \in \mathcal{Q}} \abs{\riskder[\vq]{y_i}{\eta(\muzero[\vq]{\vx_i})}} \leq \riskderenvel{\vx_i}{y_i}$ with
    \begin{equation*}
      \sup_{\vx \in \mathcal{X}}\Expect[\big]{\riskderenvel{\vx_i}{y_i}^{\nu} \given \vx_i = \vx} < \infty\quad\text{for some $\nu > 2$.}
    \end{equation*}

  \item\label{ass:psi-holder-jumps} The function $\eta(\cdot)\colon \mathbb{R} \to \newtarget{def:rangeeta}{\rangeetaNoLink}$ is strictly monotonic and twice continuously differentiable; for some fixed constant $\alpha \in (0, 1]$, for any $(\vx, \vq) \in \mathcal{X} \times \mathcal{Q}$, and a pair of points $\zeta_1, \zeta_2 \in \ballqx$, $\riskder[\vq]{\cdot}{\cdot}$ satisfies the following (constants hidden in $\lesssim$ do not depend on $\vx, \vq, \zeta_1, \zeta_2$):
    \begin{itemize}
    \item if $\measureonY$ in \cref{ass:dgp-holder}\ref{ass:cond-density} is the Lebesgue measure, then
\begin{equation}\label{eq:holder-assumption-with-cases-nonconvex}
    \begin{aligned}
      &\sup_{\lambda \in [0, 1]} \sup_{y \notin [\eta(\zeta_1) \wedge \eta(\zeta_2), \eta(\zeta_1) \vee \eta(\zeta_{2})]} \abs[\big]{ \riskder[\vq]{y}{\eta(\zeta_1 + \lambda(\zeta_2 - \zeta_1))} - \riskder[\vq]{y}{\eta(\zeta_2)}} \lesssim \abs{\zeta_1 - \zeta_2}^{\alpha},\\
      &\sup_{\lambda \in [0, 1]} \sup_{y \in [\eta(\zeta_1) \wedge \eta(\zeta_2), \eta(\zeta_1) \vee \eta(\zeta_{2})]} \abs[\big]{ \riskder[\vq]{y}{\eta(\zeta_1 + \lambda(\zeta_2 - \zeta_1))} - \riskder[\vq]{y}{\eta(\zeta_2)}} \lesssim 1;
    \end{aligned}
    \end{equation}
    \item if $\measureonY$ is \emph{not} the Lebesgue measure, then
\begin{equation}\label{eq:holder-assumption-without-cases-nonconvex}
\sup_{\lambda \in [0, 1]} \sup_{y \in \mathcal{Y}} \abs[\big]{\riskder[\vq]{y}{\eta(\zeta_1 + \lambda(\zeta_2 - \zeta_1))} - \riskder[\vq]{y}{\eta(\zeta_2)}} \lesssim \abs{\zeta_1 - \zeta_2}^{\alpha}.
\end{equation}
    \end{itemize}

  \item\label{ass:psi-eta}
    $\newtarget{def:expriskder}{\expriskderNoLink[\vq]{\vx}{\eta}} := \Expect*{\riskder[\vq]{y_i}{\eta} \given \vx_i=\vx}$ is twice continuously differentiable with respect to $\eta$,
\begin{equation}\label{eq:PsiOneTwo}
\sup_{\vx \in \mathcal{X}, \vq \in \mathcal{Q}} \sup_{\zeta \in B_{\vq}(\vx)} \abs{\PsikNoLink[\vq]{\vx}{\eta(\zeta)}} < \infty,\quad \PsikNoLink[\vq]{\vx}{\eta} := \frac{\partial^k}{\partial \eta^k} \expriskder[\vq]{\vx}{\eta},\quad k = 1, 2,
\end{equation}
and
\begin{equation*}
  \inf_{\vx \in \mathcal{X}, \vq \in \mathcal{Q}} \inf_{\zeta \in \ballqx} \Psione[\vq]{\vx}{\eta(\zeta)} \eta^{(1)}(\zeta)^{2} > 0.
\end{equation*}
\end{enumerate}
\end{fullassumption}

In addition, we make the following assumption on the complexity of $\crl{\riskder[\vq]{\cdot}{\cdot}}$. It is much more general than Assumption~\mainref{Assumption: simple setup} in the main paper (\cref{prop:simpler-generalization}).

\begin{fullassumption}[Complexity of \texorpdfstring{$\crl{\riskder[\vq]{\cdot}{\cdot}}$}{Complexity of psi}.]\label{ass:complexity}
  For any fixed $r > 0$ and $c > 0$, $l \in \crl{1, \ldots, K}$, the classes of functions
    \begin{align*}
      &\mathcal{G}_1 := \crl[\big]{\mathcal{X} \times \mathcal{Y} \ni (\vx, y) \mapsto \riskder[\vq]{y}{\eta(\vp(\vx)\trans \vbeta)} - \riskder[\vq]{y}{\eta(\muzero[\vq]{\vx})} \colon\\
      &\qquad \norm{\vbeta - \betazero[\vq]}_{\infty} \leq r, \vq \in \mathcal{Q}}\\
      &\mathcal{G}_2 := \crl[\big]{\mathcal{X} \times \mathcal{Y} \ni (\vx, y) \mapsto \riskder[\vq]{y}{\eta(\muzero[\vq]{\vx})} \colon \vq \in \mathcal{Q}},\\
      &\mathcal{G}_3 := \crl[\big]{\mathcal{X} \times \mathcal{Y} \ni (\vx, y) \mapsto\\
      &\qquad [\riskder[\vq]{y}{\eta(\vp(\vx)\trans \vbeta)} - \riskder[\vq]{y}{\eta(\vp(\vx)\trans \betazero[\vq])}] \Indlet(\vx \in \mathcal{N}_{[c \log n]}(\delta)) \colon\\
      &\qquad \norm{\vbeta - \betazero[\vq]}_{\infty} \leq r, \delta \in \Delta, \vq \in \mathcal{Q}},\\
      &\mathcal{G}_4 := \crl[\big]{\mathcal{X} \ni \vx \mapsto \Psione[\vq]{\vx}{\eta(\muzero[\vq]{\vx})}: \vq \in \mathcal{Q}},\\
      &\mathcal{G}_5 := \crl[\big]{\mathcal{X} \times \mathcal{Y} \ni (\vx, y) \mapsto\\
      &\qquad p_l(\vx) [\riskder[\vq]{y}{\eta(\muzero[\vq]{\vx})} - \riskder[\vq]{y}{\eta(\vp(\vx)\trans \betazero[\vq])}]: \vq \in \mathcal{Q}}
    \end{align*}
    satisfy the uniform entropy bound~\eqref{eq:uniform-entropy-bound} with respective envelopes and constants as follows:
    \begin{align*}
      \bar{G}_1 &\lesssim 1, & A_1 &\lesssim 1, & V_1 &\lesssim K \asymp h^{-d};\\
      \bar{G}_2(\vx, y) &\lesssim \riskderenvel{\vx}{y}, & A_2 &\lesssim 1, & V_2 &\lesssim 1;\\
      \bar{G}_3 &\lesssim 1, & A_3 &\lesssim 1, & V_3 &\lesssim \log^d n;\\
      \bar{G}_4 &\lesssim 1, & A_4 &\lesssim 1, & V_4 &\lesssim 1;\\
      \bar{G}_5 &\lesssim 1, & A_5 &\lesssim 1, & V_5 &\lesssim 1,
    \end{align*}
    where for $s \in \mathbb{Z} \cap [0, \infty)$, $\mathcal{N}_s(\delta)$ denotes the $s$-neighborhood of cell $\delta \in \Delta$ which is the union of all cells $\delta' \in \Delta$ reachable from some point of $\delta$ in no more than $s$ steps (following a continuous path).
\end{fullassumption}

\begin{fullassumption}[Approximation error]\label{ass:gl-approximation-error}
  There exists a vector of coefficients $\betazero[\vq] \in \mathbb{R}^K$ such that for all $\boldsymbol{\varsigma}$ satisfying $\multindnorm{\boldsymbol{\varsigma}} \leq \varsigma$ in \cref{ass:gl-local-basis},
  \begin{equation*}
    \sup_{\vq \in \mathcal{Q}, \vx \in \mathcal{X}} \abs[\big]{\muzero[\vq][(\boldsymbol{\varsigma})]{\vx} - \betazerotrans[\vq] \vp^{(\boldsymbol{\varsigma})}(\vx; \Delta, m)} \lesssim h^{m - \multindnorm{\boldsymbol{\varsigma}}}.
  \end{equation*}
\end{fullassumption}

In particular, this requires $\sup_{\vq, \vx} \abs[\big]{\muzero[\vq][(\vv)]{\vx} - \betazerotrans[\vq] \vp^{(\vv)}(\vx; \Delta, m)} \lesssim h^{m - \multindnorm{\vv}}$.

\begin{fullassumption}[Estimator of the Gram matrix]\label{ass:gram-estimate}
  $\qhat[\vq]$ is an estimator of the matrix $\qbar[\vq]$ such that $\norm[\big]{\qhat[\vq] - \qbar[\vq]}_{\infty} \lesssim_{\Problet} h^d \newtarget{def:qhat-rate}{\rqhatNL}$ and $\norm[\big]{\qhat[\vq][-1]}_{\infty} \lesssim_{\Problet} h^{-d}$, where $\rqhat = o(1)$.
\end{fullassumption}


Finally, we state the following technical addition to \cref{ass:complexity} needed for \cref{th:bahadur-repres}. As proven in \cref{prop:simpler-generalization}, it is also automatically true under Assumption~\mainref{Assumption: simple setup} in the main paper.

\begin{fullassumption}[Additional complexity assumption]\label{ass:complexity-addition}
For any fixed $c > 0$, $\gamma > 0$, and a positive sequence $\newtarget{def:epsninassump}{\epsninassumpNoLink} \to 0$ the class of functions
\begin{align*}
  &\crl[\bigg]{
    (\vx, y) \mapsto \int_{-\bp(\vx)\trans\vv}^0 \brk[\big]{\riskder[\vq]{y}{\eta\prn{\bp(\vx)\trans \prn{\betazero[\vq] + \bbeta } + t}} - \riskder[\vq]{y}{\eta\prn{\bp(\vx)\trans \betazero[\vq]}}}\\
  &\quad \phantom{(\vx, y) \mapsto }\times \eta^{(1)} \prn{ \bp(\vx)\trans \prn{ \betazero[\vq] + \bbeta } + t }\,\mathrm{d} t \cdot \Indlet\prn{\vx \in \mathcal{N}_{[c \log n]}(\delta)}:\\
  &\quad \norm{\vbeta - \betazero[\vq]}_{\infty} \leq \gamma, \norm{\vv}_{\infty} \leq \epsninassump, \delta \in \Delta, \vq \in \mathcal{Q}
    }
\end{align*}
with envelope $\epsninassump$ multiplied by a large enough constant (not depending on $n$), satisfies the uniform entropy bound~\eqref{eq:uniform-entropy-bound} with $A \lesssim 1 / \epsninassump$, $V \lesssim \log^d n$, where the constants in $\lesssim$ do not depend on $n$.
\end{fullassumption}

\subsection{Simplifying assumption on loss function}\label{sec:simpler-generalization-of-examples}

As discussed in the paper, it is possible to impose a simplifying assumption on the loss function, which is motivated by the examples considered. More specifically, we aim at simplifying the general complexity \cref{ass:complexity,ass:complexity-addition}.

\begin{proposition}[Simplified setup]\label{prop:simpler-generalization}
  Assume the following conditions.
  \begin{conditions}[wide]
  \item $q \mapsto \muzero[q]{\vx}$ is non-decreasing.
  \item $\riskfunc[q]{y}{\eta} = \sum_{j = 1}^4 \riskfuncj[q]{y}{\eta}$, where the functions $\riskfuncj[q]{y}{\eta}$ are of the following types:
    \begin{itemize}
    \item Type I: $\riskfunci[q]{y}{\eta} = (f_1(y) + D_1 \eta) \Indlet\{y \leq \eta\}$,
    \item Type II: $\riskfuncii[q]{y}{\eta} = (f_2(y) + D_2 \eta) \Indlet\{y \leq q\}$,
    \item Type III: $\riskfunciii[q]{y}{\eta} = (f_3(y) + D_3 \eta) q$,
    \item Type IV: $\riskfunciv[q]{y}{\eta} = \mathcal{T}(y, \eta)$,
    \end{itemize}
    where $f_j(\cdot)$ are fixed continuous functions, $D_j$ are universal constants, and $\mathcal{T}(y; \eta): \mathcal{Y} \times \rangeeta \to \mathbb{R}$ is a measurable function, differentiable in $\eta$ for any fixed $y$ with a derivative $\tau(y, \eta) := \frac{\partial}{\partial \eta} \mathcal{T}(y, \eta)$.
  \item $\Expect*{\tau(y_i,\eta) \given \vx_i = \vx}$ is also differentiable in $\eta$ for any $\vx$.
    The functions
    \begin{align*}
      &(y, \eta) \mapsto \tau(y, \eta),\text{ and}\\
      &(\vx, \eta) \mapsto \frac{\partial}{\partial \eta} \Expect*{\tau(y_i, \eta) \given \vx_i = \vx}
    \end{align*}
    are continuous in their arguments, and $\alpha$-H\"older continuous (with $0 < \alpha \leq 1$) in $\eta$ on $\mathcal{Y} \times \mathcal{K}$ and $\mathcal{X} \times \mathcal{K}$ respectively, where $\mathcal{K}$ is any fixed compact subset of $\rangeeta$ (with the H\"older constant possibly depending on $\mathcal{K}$, but not on $y$ or $\vx$).
      \item \cref{ass:dgp-holder,ass:gl-holder,ass:gl-quasi-uniform-partition,ass:gl-local-basis} hold with $\riskderenvel{\vx}{y} = \overline{\tau}(\vx, y) + |D_1| + |D_2| + |D_3| \max_{q \in \mathcal{Q}}|q|$, where $\overline{\tau}(\vx, y)$ is a measurable envelope of
    $\crl[\big]{ (\vx, y) \mapsto \tau \prn[\big]{ y, \eta \prn{ \muzero[q]{\vx} } } }$.
  \item If $D_1$ is nonzero, $F_{Y|X}$ is differentiable and $f_{Y|X}$ is its derivative (in particular, $\measureonY$ is Lebesgue measure), $(\vx, \eta) \mapsto f_{Y|X}(\eta|\vx)$ is continuous in both arguments and continuously differentiable in $\eta$.
  \end{conditions}

  Then \cref{ass:complexity,ass:complexity-addition} also hold.
\end{proposition}

We will prove this now.

In this case
\begin{equation*}
  \riskder[q]{y}{\eta} = D_1 \Indlet\{y \leq \eta\} + D_2 \Indlet\{y \leq q\} + D_3 q + \tau(y, \eta).
\end{equation*}

\subsubsection{Verifying \cref{ass:complexity}}
\begin{lemma}[Class \texorpdfstring{$G_1$}{G1}]\label{lem:simplified-class-g-one}
  The class $\mathcal{G}_1$ described in \cref{ass:complexity} with a large enough constant envelope satisfies the uniform entropy bound \eqref{eq:uniform-entropy-bound} with $A \lesssim 1$, $V \lesssim K$.
\end{lemma}

\begin{proof}[Proof of \cref{lem:simplified-class-g-one}]
It is shown in the proof of \cref{prop:gl-holder-quantile-regression} (replacing $<$ with $\leq$ does not change the argument) that for any fixed $r > 0$ the class
\begin{equation*}
  \crl[\big]{ (\vx, y) \mapsto \Indlet \crl{ y \leq \eta(\vp(\vx)\trans \vbeta) } - \Indlet \crl{ y \leq \eta\prn{\muzero[q]{\vx}} }: \norm{ \vbeta - \betazero[q] }_{\infty} \leq r, q \in \mathcal{Q}}
\end{equation*}
with envelope 2 satisfies the uniform entropy bound~\eqref{eq:uniform-entropy-bound} with $A \lesssim 1$ and $V \lesssim K$.

Next, assume that the infinity-norms of $\bbeta$ and $\widetilde{\bbeta}$ lie in a fixed bounded interval, and let $q, \tilde{q} \in \mathcal{Q}$. Note that by $\alpha$-H\"older continuity of $\tau(y, \cdot)$ on compacta
\begin{align*}
  &\abs[\big]{ \tau \prn{ y, \eta(\vp(\vx)\trans \vbeta) } - \tau \prn{ y, \eta\prn{\muzero[q]{\vx}} } - \tau \prn{ y, \eta(\vp(\vx)\trans \widetilde{\bbeta}) } + \tau \prn{ y, \eta\prn{\muzero[\tilde{q}]{\vx}} } }\\
  &\quad \lesssim \abs[\big]{ \eta(\vp(\vx)\trans \vbeta) - \eta(\vp(\vx)\trans \widetilde{\bbeta}) }^{\alpha} + \abs[\big]{ \eta\prn{\muzero[q]{\vx}} - \eta\prn{\muzero[\tilde{q}]{\vx}} }^{\alpha}\\
  &\quad \lesssim \norm{ \bbeta - \widetilde{\bbeta} }_{\infty}^{\alpha} + |q - \tilde{q}|^{\alpha},
\end{align*}
where the constants in $\lesssim$ do not depend on $\bbeta$, $\widetilde{\bbeta}$, $q$, $\tilde{q}$, and we used that $\eta(\cdot)$ on a fixed bounded interval is Lipschitz, and $q \mapsto \muzero[q]{\vx}$ is Lipschitz (uniformly in $\vx$). Again by $\alpha$-H\"older continuity for any fixed $r > 0$ the class
\begin{equation*}
  \crl[\big]{ (\vx, y) \mapsto \tau \prn[\big]{ y, \eta(\vp(\vx)\trans \vbeta) } - \tau \prn[\big]{ y, \eta\prn{\muzero[q]{\vx}} }: \norm{ \vbeta - \betazero[q] }_{\infty} \leq r, q \in \mathcal{Q}}
\end{equation*}
has a constant envelope. It follows that it satisfies the uniform entropy bound \eqref{eq:uniform-entropy-bound} with $A \lesssim 1$, $V \lesssim K$.

Combining these results concludes the proof of \cref{lem:simplified-class-g-one} by \cref{lem:permanence-properties}.
\end{proof}

\begin{lemma}[Class \texorpdfstring{$G_2$}{G2}]\label{lem:simplified-class-g-two}
  The class $\mathcal{G}_2$ described in \cref{ass:complexity} with envelope $\riskderenvel{\vx}{y}$ satisfies the uniform entropy bound \eqref{eq:uniform-entropy-bound} with $A \lesssim 1$, $V \lesssim 1$.
\end{lemma}

\begin{proof}[Proof of \cref{lem:simplified-class-g-two}]
The class
\begin{equation*}
  \crl[\big]{\mathcal{X} \times \mathcal{Y} \ni (\vx, y) \mapsto \tau\prn[\big]{y, \eta\prn{\muzero[q]{\vx}}} \colon q \in \mathcal{Q}}
\end{equation*}
with envelope $\overline{\tau}(\vx, y)$ satisfies the uniform entropy bound \eqref{eq:uniform-entropy-bound} with $A \lesssim 1$, $V \lesssim 1$ by $\alpha$-H\"older continuity and since $\eta(\cdot)$ on a fixed bounded interval is Lipschitz, $q \mapsto \muzero[q]{\vx}$ is Lipschitz (uniformly in $\vx$).

The class
\begin{equation*}
  \crl[\big]{\mathcal{X} \times \mathcal{Y} \ni (\vx, y) \mapsto \Indlet \crl{ y \leq \eta\prn{\muzero[q]{\vx}} } \colon q \in \mathcal{Q}}
\end{equation*}
with envelope 1 satisfies the uniform entropy bound \eqref{eq:uniform-entropy-bound} with $A \lesssim 1$, $V \lesssim 1$ because $\muzero[q]{\vx}$ is nondecreasing in $q$, see the proof of \cref{prop:gl-holder-quantile-regression}.

The class
\begin{equation*}
  \crl[\big]{\mathcal{X} \times \mathcal{Y} \ni (\vx, y) \mapsto \Indlet \crl{ y \leq q } \colon q \in \mathcal{Q}}
\end{equation*}
with envelope 1 satisfies the uniform entropy bound \eqref{eq:uniform-entropy-bound} with $A \lesssim 1$, $V \lesssim 1$ because it is VC-subgraph with a constant index, (cf. the proof of \cref{prop:gl-holder-distribution-regression}).

The class
\begin{equation*}
  \crl[\big]{\mathcal{X} \times \mathcal{Y} \ni (\vx, y) \mapsto q \colon q \in \mathcal{Q}}
\end{equation*}
with envelope $\max_{q \in \mathcal{Q}} |q|$ satisfies the uniform entropy bound \eqref{eq:uniform-entropy-bound} with $A \lesssim 1$, $V \lesssim 1$ because it is VC-subgraph with a constant index (as a subclass of a one-dimensional space of functions, namely constants in $\vx, y$).

It is left to apply \cref{lem:permanence-properties}, concluding the proof of \cref{lem:simplified-class-g-two}.
\end{proof}

\begin{lemma}[Class \texorpdfstring{$G_3$}{G3}]\label{lem:simplified-class-g-three}
  The class $\mathcal{G}_3$ described in \cref{ass:complexity} with a large enough constant envelope satisfies the uniform entropy bound \eqref{eq:uniform-entropy-bound} with $A \lesssim 1$, $V \lesssim \log^d n$.
\end{lemma}

\begin{proof}[Proof of \cref{lem:simplified-class-g-three}]
Fix $\delta \in \Delta$ and some large enough $R > 0$. Note that if $\Indlet \bigl\{ \vx \in \mathcal{N}_{[c \log n]}(\delta) \bigr\}$ and $\bp(\vx)\trans \bbeta$ are both nonzero, $\bbeta$ must lie in a vector subspace $\mathcal{B}_{\delta}$ of $\mathbb{R}^K$ of dimension $O(\log^d n)$. For any positive and small enough $\varepsilon$, the class of vectors $\bigl\{\bbeta \in \mathcal{B}_{\delta}, \norm{ \bbeta }_{\infty} < R\bigr\}$ has an infinity-norm $\varepsilon$-net $\bar{\mathcal{B}}_{\delta}^{\varepsilon}$ such that
\begin{equation*}
  \log \abs{ \bar{\mathcal{B}}_{\delta}^{\varepsilon} } \lesssim \log^d n \log(C / \varepsilon),
\end{equation*}
where $C$ is some positive constant.

By $\alpha$-H\"older continuity of $\tau(y, \cdot)$ on compacta and since $\eta(\cdot)$ on a compact is Lipschitz, this means that the class of bounded (by a constant not depending on $n$) functions
\begin{align*}
  \mathcal{G}_{3, \delta} := \crl[\big]{ &(\vx, y) \mapsto \brk[\big]{\tau\prn[\big]{y, \eta(\bp(\vx)\trans \bbeta)} - \tau\prn[\big]{y, \eta(\bp(\vx)\trans \betazero[q])}} \Indlet\crl{\vx \in \mathcal{N}_{[c \log n]}(\delta)}:\\
  &\norm{ \vbeta - \betazero[q] }_{\infty} \leq r, q \in \mathcal{Q}}
\end{align*}
has a covering number bound
\begin{equation*}
  \log N(\mathcal{G}_{3, \delta}, \text{sup-norm}, C' \varepsilon^{\alpha}) \lesssim \log^d n \log(C / \varepsilon),
\end{equation*}
where $C'$ is some other positive constant. This means that also
\begin{equation*}
  \log N(\mathcal{G}_{3, \delta}, \text{sup-norm}, \varepsilon) \lesssim \log^d n \log(C'' / \varepsilon)
\end{equation*}
for some other positive constant $C''$. Finally, from this we can conclude that $\mathcal{G}_{3, \delta}$ with a large enough constant envelope satisfies the uniform entropy bound \eqref{eq:uniform-entropy-bound} with $A \lesssim 1$ and $V \lesssim \log^d n$. Therefore, the union of $O\prn{h^{-d}}$ such classes
\begin{align*}
  \crl[\big]{ &(\vx, y) \mapsto \crl[\big]{\tau\prn[\big]{y, \eta(\bp(\vx)\trans \bbeta)} - \tau\prn[\big]{y, \eta(\bp(\vx)\trans \betazero[q])}} \Indlet\crl{\vx \in \mathcal{N}_{[c \log n]}(\delta)}:\\
  &\norm{ \vbeta - \betazero[q] }_{\infty} \leq r, q \in \mathcal{Q}, \delta \in \Delta}
\end{align*}
with a large enough constant envelope satisfies the uniform entropy bound~\eqref{eq:uniform-entropy-bound} with $A \lesssim 1$ and $V \lesssim \log^d n$, see~\eqref{eq:log-d-entropy-chain}.

Next, since $\eta(\cdot)$ is monotonic, the functions $\vx \mapsto \bp(\vx)\trans\bbeta$ with $\bbeta \in \mathcal{B}_{\delta}$ form a vector space of $O(\log^d n)$ dimension, and $\vx \mapsto \Indlet\crl{\vx \in \mathcal{N}_{[c \log n]}(\delta)}$ is one fixed function, the class
\begin{equation*}
  \crl[\big]{(\vx, y) \mapsto \Indlet\{y \leq \eta(\bp(\vx)\trans\bbeta)\} \Indlet\crl{\vx \in \mathcal{N}_{[c \log n]}(\delta)}: \bbeta \in \mathbb{R}^K}
\end{equation*}
is VC with an index $O(\log^d n)$. Therefore, it satisfies the uniform entropy bound \eqref{eq:uniform-entropy-bound} with $A \lesssim 1$ and $V \lesssim \log^d n$. By \cref{lem:permanence-properties}, a subclass of the difference of two such classes
\begin{align*}
  \crl[\big]{ &(\vx, y) \mapsto \prn[\big]{\Indlet\crl{y \leq \eta(\bp(\vx)\trans \bbeta)} - \Indlet\crl{y \leq \eta(\bp(\vx)\trans \betazero[q])}} \Indlet\crl{\vx \in \mathcal{N}_{[c \log n]}(\delta)}:\\
  &\norm{ \vbeta - \betazero[q] }_{\infty} \leq r, q \in \mathcal{Q}}
\end{align*}
for fixed $\delta \in \Delta$ also satisfies the uniform entropy bound \eqref{eq:uniform-entropy-bound} with $A \lesssim 1$ and $V \lesssim \log^d n$. Therefore, the union of $O\prn{h^{-d}}$ such classes
\begin{align*}
  &\crl[\big]{ (\vx, y) \mapsto \crl{\Indlet\crl{y \leq \eta(\bp(\vx)\trans \bbeta)} - \Indlet\crl{y \leq \eta(\bp(\vx)\trans \betazero[q])}} \Indlet\crl{\vx \in \mathcal{N}_{[c \log n]}(\delta)}:\\
  &\quad \norm{ \vbeta - \betazero[q] }_{\infty} \leq r, \delta \in \Delta, q \in \mathcal{Q}}
\end{align*}
with a large enough constant envelope satisfies the uniform entropy bound~\eqref{eq:uniform-entropy-bound} with $A \lesssim 1$ and $V \lesssim \log^d n$, see~\eqref{eq:log-d-entropy-chain}.

It is left to apply \cref{lem:permanence-properties} once again, concluding the proof of \cref{lem:simplified-class-g-three}.
\end{proof}

\begin{lemma}[Class \texorpdfstring{$G_4$}{G4}]\label{lem:simplified-class-g-four}
  The class $\mathcal{G}_4$ described in \cref{ass:complexity} with a large enough constant envelope satisfies the uniform entropy bound \eqref{eq:uniform-entropy-bound} with $A \lesssim 1$, $V \lesssim 1$.
\end{lemma}

\begin{proof}[Proof of \cref{lem:simplified-class-g-four}]
In this case
\begin{equation*}
  \expriskder[q]{\vx}{\eta} = \Expect*{\tau(y_i, \eta) \given \vx_i = \vx} + D_1 F_{Y|X}(\eta | \vx) + D_2 F_{Y|X}(q | \vx) + D_3 q
\end{equation*}
and
\begin{equation*}
  \Psione[q]{\vx}{\eta} = \frac{\partial}{\partial \eta} \Expect*{\tau(y_i, \eta) \given \vx_i = \vx} + D_1 f_{Y|X}(\eta | \vx).
\end{equation*}
By assumption, if $\eta$ lies in a fixed compact, this function of $\vx, \eta$ is bounded (by continuity). Moreover,
\begin{align*}
  &\abs[\big]{ \Psione[q]{\vx}{\eta(\muzero[q]{\vx})} - \Psione[\tilde{q}]{\vx}{\eta(\muzero[\tilde{q}]{\vx})}}\\
  &\quad \lesssim |\eta(\muzero[q]{\vx}) - \eta(\muzero[\tilde{q}]{\vx})|^{\alpha} + |\eta(\muzero[q]{\vx}) - \eta(\muzero[\tilde{q}]{\vx})| \labrel{1}[\lesssim] |q - \tilde{q}|^{\alpha} + |q - \tilde{q}|,
\end{align*}
where in \labrel{1} we used that $\eta(\cdot)$ on compacta is Lipschitz and $q \mapsto \muzero[q]{\vx}$ is uniformly Lipschitz in $q$. The result of \cref{lem:simplified-class-g-four} follows.
\end{proof}

\begin{lemma}[Class \texorpdfstring{$G_5$}{G5}]\label{lem:simplified-class-g-five}
  The class $\mathcal{G}_5$ described in \cref{ass:complexity} with a large enough constant envelope satisfies the uniform entropy bound \eqref{eq:uniform-entropy-bound} with $A \lesssim 1$, $V \lesssim 1$.
\end{lemma}

\begin{proof}[Proof of \cref{lem:simplified-class-g-five}]
Take $R > 0$ fixed and large enough so that $\norm[\big]{\bbeta_{0, q}}_{\infty} \leq R$ for all $q$ and $n$. Note that for $p_l(\vx)$ and $\bp(\vx)\trans \bbeta$ to be nonzero at the same time, $\bbeta$ must lie in a fixed vector subspace $\mathcal{B}_l$ of $\mathbb{R}^K$ of bounded dimension. For $0 < \varepsilon < 1$, the class of vectors $\bigl\{\bbeta \in \mathcal{B}_l, \norm{ \bbeta }_{\infty} \leq R\bigr\}$ has an infinity-norm $\varepsilon$-net $\bar{\mathcal{B}}_l^{\varepsilon}$ such that
\begin{equation*}
  \log \abs[\big]{ \bar{\mathcal{B}}_l^{\varepsilon} } \lesssim \log(C / \varepsilon),
\end{equation*}
where $C$ is some positive constant.

By $\alpha$-H\"older continuity of $\tau(y, \cdot)$ on compacta, since $\eta(\cdot)$ on a compact is Lipschitz and $\muzero[q]{\cdot}$ is uniformly Lipschitz in $q$, this means that the class of bounded (by a constant not depending on $n$) functions
\begin{equation*}
  \mathcal{G}_5' := \crl[\big]{
    (\vx, y) \mapsto p_l(\vx) \prn[\big]{\tau(y, \eta(\muzero[q]{\vx})) - \tau(y, \eta(\bp(\vx)\trans \betazero[q]))}: q \in \mathcal{Q}
  }
\end{equation*}
has a covering number bound
\begin{equation*}
  \log N(\mathcal{G}_5', \text{sup-norm}, C' \varepsilon^{\alpha}) \lesssim \log(C / \varepsilon),
\end{equation*}
where $C'$ is some other positive constant. It follows that this class with a large enough constant envelope satisfies the uniform entropy bound \eqref{eq:uniform-entropy-bound} with $A \lesssim 1$ and $V \lesssim 1$.

As in the proof of \cref{lem:simplified-class-g-two}, the class
\begin{equation*}
  \crl[\big]{
    \mathcal{X} \times \mathcal{Y} \ni (\vx, y) \mapsto \Indlet \crl{ y \leq \eta\prn{\muzero[q]{\vx}} } \colon q \in \mathcal{Q}
  }
\end{equation*}
with envelope 1 satisfies the uniform entropy bound \eqref{eq:uniform-entropy-bound} with $A \lesssim 1$, $V \lesssim 1$. Since $\eta(\cdot)$ is monotonic, the functions $\vx \mapsto \bbeta\trans\bp(\vx)$ with $\bbeta \in \mathcal{B}_l$ form a vector space, and $p_l(\vx)$ is one fixed function, we have that the class
\begin{equation*}
  \crl[\big]{
    \mathcal{X} \times \mathcal{Y} \ni (\vx, y) \mapsto p_l(\vx) \Indlet \crl{ y \leq \eta\prn{ \bp(\vx)\trans\bbeta } }: \bbeta \in \mathbb{R}^K
  }
\end{equation*}
is VC with a bounded index. Then it also satisfies the uniform entropy bound \eqref{eq:uniform-entropy-bound} with $A \lesssim 1$ and $V \lesssim 1$.

It is left to apply \cref{lem:permanence-properties}, concluding the proof of \cref{lem:simplified-class-g-five}.
\end{proof}

\subsubsection{Verifying \cref{ass:complexity-addition}}\label{sec:simplified-example-verifying-addition-for-nonconvex-bahadur}
Suppose $\theta_1, \theta_2, \theta \in \mathbb{R}$ and $\tilde{\theta}_1, \tilde{\theta}_2, \tilde{\theta} \in \mathbb{R}$ all lie in a fixed compact interval. Then
\begin{align*}
  &\mathcal{T}(y, \eta(\theta_1)) - \mathcal{T}(y, \eta(\theta_2)) - [\eta(\theta_1) - \eta(\theta_2)] \tau(y, \eta(\theta))\\
  &\qquad - \mathcal{T}(y, \eta(\tilde{\theta}_1)) + \mathcal{T}(y, \eta(\tilde{\theta}_2)) + [\eta(\tilde{\theta}_1) - \eta(\tilde{\theta}_2)] \tau(y, \eta(\tilde{\theta}))\\
  &\quad = \tau(y, \zeta_{1, y}) [\eta(\theta_1) - \eta(\tilde{\theta}_1)] - \tau(y, \zeta_{2, y}) [\eta(\theta_2) - \eta(\tilde{\theta}_2)]\\
  &\qquad - [\eta(\theta_1) - \eta(\tilde{\theta}_1) + \eta(\tilde{\theta}_1) - \eta(\theta_2)] \tau(y, \eta(\theta)) + [\eta(\tilde{\theta}_1) - \eta(\theta_2) + \eta(\theta_2) - \eta(\tilde{\theta}_2)] \tau(y, \eta(\tilde{\theta}))\\
  &\quad = [\tau(y, \zeta_{1, y}) - \tau(y, \eta(\theta))] \cdot [\eta(\theta_1) - \eta(\tilde{\theta}_1)] - [\tau(y, \zeta_{2, y}) - \tau(y, \eta(\tilde{\theta}))] \cdot [\eta(\theta_2) - \eta(\tilde{\theta}_2)]\\
  &\qquad - [\eta(\tilde{\theta}_1) - \eta(\theta_2)] \tau(y, \eta(\theta)) + [\eta(\tilde{\theta}_1) - \eta(\theta_2)] \tau(y, \eta(\tilde{\theta}))\\
  &\quad = \underbrace{[\tau(y, \zeta_{1, y}) - \tau(y, \eta(\theta))]}_{\lesssim 1} \cdot \underbrace{[\eta(\theta_1) - \eta(\tilde{\theta}_1)]}_{\lesssim \abs{\theta_1 - \tilde{\theta}_1}} - \underbrace{[\tau(y, \zeta_{2, y}) - \tau(y, \eta(\tilde{\theta}))]}_{\lesssim 1} \cdot \underbrace{[\eta(\theta_2) - \eta(\tilde{\theta}_2)]}_{\lesssim \abs{ \theta_2 - \tilde{\theta}_2 }}\\
  &\qquad + \underbrace{[\eta(\tilde{\theta}_1) - \eta(\theta_2)]}_{\lesssim 1} \underbrace{[\tau(y, \eta(\tilde{\theta})) - \tau(y, \eta(\theta))]}_{\lesssim \abs{ \theta - \tilde{\theta} }^{\alpha}} \lesssim |\theta_1 - \tilde{\theta}_1| + | \theta_2 - \tilde{\theta}_2 | + | \theta - \tilde{\theta} |^{\alpha}
\end{align*}
for some $\zeta_{1, y}$ between $\eta(\theta_1)$ and $\eta(\tilde{\theta}_1)$, $\zeta_{2, y}$ between $\eta(\theta_2)$ and $\eta(\tilde{\theta}_2)$, where we used the $\alpha$-H\"older continuity of $\tau(y, \cdot)$ on a fixed compact and the Lipschitzness of $\eta(\cdot)$ on a compact. This means that the class of functions
\begin{align*}
  &\mathcal{G}' := \crl[\big]{ (\vx, y) \mapsto
    \prn[\big]{
    \mathcal{T}\prn[\big]{y, \eta\prn[\big]{\bp(\vx)\trans(\betazero[q] + \bbeta)}}
    - \mathcal{T}\prn[\big]{y, \eta\prn[\big]{\bp(\vx)\trans(\betazero[q] + \bbeta - \vv)}}\\
  &\qquad - \brk[\big]{ \eta\prn[\big]{\bp(\vx)\trans(\betazero[q] + \bbeta)} - \eta\prn[\big]{\bp(\vx)\trans(\betazero[q] + \bbeta - \vv)} } \tau\prn[\big]{y, \eta\prn{\bp(\vx)\trans \betazero[q]}}
    }\\
  &\quad\qquad \times \Indlet\crl{\vx \in \mathcal{N}_{[c \log n]}(\delta)} : q \in \mathcal{Q},
    \norm{ \vbeta - \betazero[q] }_{\infty} \leq r, \norm{ \vv }_{\infty} \leq \epsninassump
    }
\end{align*}
has a covering number bound
\begin{equation}\label{eq:nonconvex-complexity--verification-smooth-aux}
  \log N(\mathcal{G}', \text{sup-norm}, \varepsilon) \lesssim \log^d n \log\prn[\bigg]{\frac{C_1}{\varepsilon}},
\end{equation}
for all small enough positive $\varepsilon$, and $C_1$ is some positive constant (not depending on $n$), where we used that all $\betazero[q]$, $\betazero[q] + \bbeta$ and $\betazero[q] + \bbeta - \vv$ must lie in a vector space of dimension $O\prn{\log^d n}$ if $\Indlet\crl{\vx \in \mathcal{N}_{[c \log n]}(\delta)}$ is not zero. Applying the mean-value theorem to
\begin{equation*}
  \mathcal{T}\prn[\big]{y, \eta\prn[\big]{\bp(\vx)\trans(\betazero[q] + \bbeta)}} - \mathcal{T}\prn[\big]{y, \eta\prn[\big]{\bp(\vx)\trans(\betazero[q] + \bbeta - \vv)}}
\end{equation*}
and $\alpha$-H\"older continuity again, we see that class $\mathcal{G}'$ has an envelope which is $\epsninassump$ multiplied by a large enough constant $C_2$. Replacing $\varepsilon$ with $C_2 \varepsilon \epsninassump$ (for large enough $n$ this is small enough) in \eqref{eq:nonconvex-complexity--verification-smooth-aux}, we get a covering number bound
\begin{equation*}
  \log N(\mathcal{G}', \text{sup-norm}, C_2 \varepsilon \epsninassump) \lesssim \log^d n \log\prn[\bigg]{\frac{C_1}{C_2 \varepsilon \epsninassump}}.
\end{equation*}
It follows that class $\mathcal{G}'$ satisfies the uniform entropy bound \eqref{eq:uniform-entropy-bound} with $A \lesssim 1 / \epsninassump$ and $V \lesssim \log^d n$.

Therefore, the union of $O\prn{h^{-d}}$ such classes
\begin{align*}
  \mathcal{G}' := \crl[\Big]{
    &\prn[\Big]{
    \mathcal{T}\prn[\big]{y, \eta\prn{\bp(\vx)\trans(\betazero[q] + \bbeta)}}
    - \mathcal{T}\prn[\big]{y, \eta\prn{\bp(\vx)\trans(\betazero[q] + \bbeta - \vv)}}\\
  &- \brk[\big]{ \eta\prn{\bp(\vx)\trans(\betazero[q] + \bbeta)} - \eta\prn{\bp(\vx)\trans(\betazero[q] + \bbeta - \vv)} } \tau\prn{y, \eta\prn{\bp(\vx)\trans \betazero[q]}}
    }\\
  &\times \Indlet\crl{\vx \in \mathcal{N}_{[c \log n]}(\delta)} : q \in \mathcal{Q},
    \norm{ \vbeta - \betazero[q] }_{\infty} \leq r, \norm{ \vv }_{\infty} \leq \epsninassump, \delta \in \Delta
    },
\end{align*}
also with envelope $\epsninassump$ multiplied by a large enough constant, satisfies the uniform entropy bound \eqref{eq:uniform-entropy-bound} with $A \lesssim 1 / \epsninassump$ and $V \lesssim \log^d n$, by the same argument as~\eqref{eq:log-d-entropy-chain}.

Next,
\begin{align*}
  & \brk{ f_1(y) + D_1 \eta(\theta_1) } \Indlet\{y \leq \eta(\theta_1)\} - \brk{ f_1(y) + D_1 \eta(\theta_2) } \Indlet\{y \leq \eta(\theta_2)\}                                                     \\
  & \qquad - \brk{ \eta(\theta_1) - \eta(\theta_2) } D_1 \Indlet\{y \leq \eta(\theta)\}                                                                                                                    \\
  & \quad = f_1(y) \brk{ \Indlet\{y \leq \eta(\theta_1)\} - \Indlet\{y \leq \eta(\theta_2)\} } + D_1 \eta(\theta_1) \Indlet\{y \leq \eta(\theta_1)\} - D_1 \eta(\theta_2) \Indlet\{y \leq \eta(\theta_2)\} \\
  & \qquad - D_1 \brk{ \eta(\theta_1) - \eta(\theta_2) } \Indlet\{y \leq \eta(\theta)\}.
\end{align*}
It is proven by the same argument as in the proof of \cref{prop:gl-holder-quantile-regression} that the class
\begin{align*}
  \crl[\Big]{
  & (\vx, y) \mapsto \prn[\Big]{\brk[\big]{ f_1(y) + D_1 \eta\prn[\big]{\bp(\vx)\trans(\betazero[q] + \bbeta)} } \Indlet\crl[\big]{y \leq \eta\prn[\big]{\bp(\vx)\trans(\betazero[q] + \bbeta)}}\\
  & \qquad - \brk[\big]{ f_1(y) + D_1 \eta\prn[\big]{\bp(\vx)\trans(\betazero[q] + \bbeta - \vv)} } \Indlet\crl{y \leq \eta(\bp(\vx)\trans(\betazero[q] + \bbeta - \vv))}                 \\
  & \qquad - \brk[\big]{ \eta(\bp(\vx)\trans(\betazero[q] + \bbeta)) - \eta(\bp(\vx)\trans(\betazero[q] + \bbeta - \vv)) } D_1 \Indlet\{y \leq \eta(\theta)\}}\\
  & \qquad \times \Indlet\crl{\vx \in \mathcal{N}_{[c \log n]}(\delta)}: q \in \mathcal{Q},
    \norm{ \vbeta - \betazero[q] }_{\infty} \leq r, \norm{ \vv }_{\infty} \leq \epsninassump, \delta \in \Delta
    }
\end{align*}
satisfies the uniform entropy bound \eqref{eq:uniform-entropy-bound} with $A \lesssim 1 / \epsninassump$ and $V \lesssim \log^d n$.

The terms $(f_2(y) + D_2 \eta) \Indlet\{y \leq q\}$ and $(f_3(y) + D_3 \eta) q$ play no role in this verification because they cancel out in the class described in \cref{ass:complexity-addition}.

It is left to apply \cref{lem:permanence-properties}.

The proof of \cref{prop:simpler-generalization} is finished.

\section{Frequently used lemmas}\label{sec:preliminaries}

We collect several lemmas that will be used multiple times throughout this supplemental appendix.
 \Cref{lem:subgaussian-max-bound,lem:tech-stoch-orders-lemma,lem:conv-to-zero-in-cond-p-equiv-to-conv-to-zero-in-p,lem:permanence-properties,lem:maximal-ineq-gaussian-vectors,lem:ueb-maximal-ineq,lem:maximal-gaussian,lem:closeness-in-prob-implies-closeness-of-cond-quantiles,lem:anti-concentration-for-suprema-of-separable-gaussian-processes} are well-known facts, so we provide either brief proofs or references to the literature.

\begin{lemma}[Second moment bound of the max of sub-Gaussian random variables]\label{lem:subgaussian-max-bound}
  Let $n \geq 3$ and $\xi_1, \ldots, \xi_n$ be \subgauss\ random variables (not necessarily independent). Then
  \begin{equation*}
    \Expect[\Big]{\max_{1 \leq i \leq n} \xi_i^2}^{1 / 2} \leq C_U \sigma \sqrt{\log n},
  \end{equation*}
  where $C_U$ is a universal constant.
\end{lemma}

\begin{proof}
If $p$ is an even positive integer, $\Expect*{\xi_i^p} \leq 3 \sigma^{p} p (p / 2)^{p / 2}$. Then
\begin{align*}
  &\Expect[\Big]{\max_{1 \leq i \leq n} \xi_i^2}^{1 / 2} \leq \Expect[\Big]{\max_{1 \leq i \leq n} \xi_i^p}^{1 / p} \leq \prn[\bigg]{\sum_{i = 1}^n \Expect[\big]{\xi_i^{p}}}^{1 / p} \lesssim n^{1 / p} \cdot \sigma \cdot p^{1 / p} \sqrt{p}\\
  &\quad \lesssim \sigma n^{1 / p} \sqrt{p} \tag*{using $p^{1 / p} \leq 2$.}
\end{align*}
It is left to take $p = p_n$ such that $\ln n \leq p \leq 2 \ln n$.
\end{proof}

\begin{lemma}[Boundedness of conditional expectation in probability implies unconditional boundedness in probability]\label{lem:tech-stoch-orders-lemma}
  Let $X_n$ be a sequence of integrable random variables, $\vD_n$ a sequence of random vectors, $r_n$ a sequence of positive numbers. If $\Expect{\abs{X_n} \given \vD_n} \lesssim_{\Problet} r_n$, then $|X_n| \lesssim_{\Problet} r_n$.
\end{lemma}

\begin{proof}
Take any sequence of positive numbers $\gamma_n \to \infty$. By Markov's inequality,
\begin{equation*}
  \Prob{\abs{X_n} > \gamma_n r_n \given \vD_n}
  \leq \frac{\Expect{\abs{X_n} \given \vD_n}}{\gamma_n r_n}
  \lesssim_{\Problet} \frac{1}{\gamma_n} = o(1).
\end{equation*}

In other words, the sequence of random variables $\Prob{\abs{X_n} > \gamma_n r_n \given \vD_n}$ converges to zero in probability. By dominated convergence (in probability), the sequence of numbers $\Prob{\abs{X_n} > \gamma_n r_n}$ converges to zero. Since it is true for any positive sequence $\gamma_n \to \infty$, this means $|X_n| = O_{\Problet}(r_n)$.
\end{proof}

\begin{lemma}[Converging to zero in conditional probability is the same as converging to zero in probability]\label{lem:conv-to-zero-in-cond-p-equiv-to-conv-to-zero-in-p}
  Let $X_n$ be a sequence of random variables, $\vD_n$ a sequence of random vectors. The following are equivalent:
  \begin{enumerate}[label=(\roman*)]
  \item\label{item:conv-to-zero-in-cond-p-equiv-to-conv-to-zero-in-p-1} for any $\varepsilon > 0$, we have $\Prob*{|X_n| > \varepsilon \given \vD_n} = o_{\Problet}(1)$;
  \item\label{item:conv-to-zero-in-cond-p-equiv-to-conv-to-zero-in-p-2} $|X_n| = o_{\Problet}(1)$.
  \end{enumerate}
\end{lemma}

\begin{proof}
The implication \ref{item:conv-to-zero-in-cond-p-equiv-to-conv-to-zero-in-p-1} $\Rightarrow$ \ref{item:conv-to-zero-in-cond-p-equiv-to-conv-to-zero-in-p-2} follows from dominated convergence in probability. To prove the converse, take any $\varepsilon, \gamma > 0$. By Markov's inequality,
\begin{equation*}
  \Prob*{\Prob*{|X_n| > \varepsilon \given \vD_n} > \gamma} \leq \frac{\Prob*{|X_n| > \varepsilon}}{\gamma} \to 0,
\end{equation*}
so by definition $\Prob*{|X_n| > \varepsilon \given \vD_n} = o_{\Problet}(1)$.
\end{proof}

\begin{lemma}[Permanence properties of the uniform entropy bound]\label{lem:permanence-properties}
  Let $\mathcal{F}$ and $\mathcal{G}$ be two classes of measurable functions from $S \to \mathbb{R}$ on a measurable space $(S, \mathcal{S})$ with strictly positive measurable envelopes $F$ and $G$ respectively. Then the uniform entropy numbers of $\mathcal{F} \mathcal{G}=$ $\{f g: f \in \mathcal{F}, g \in \mathcal{G}\}$ satisfy
  \begin{align*}
    &\sup_\mathbb{Q} \log N \prn[\big]{\mathcal{F G}, \norm{\cdot}_{\mathbb{Q}, 2}, \varepsilon \norm{F G}_{\mathbb{Q}, 2}} \\
    &\quad\leq \sup_\mathbb{Q} \log N \prn[\bigg]{\mathcal{F}, \norm{\cdot}_{\mathbb{Q}, 2}, \frac{\varepsilon \norm{F}_{\mathbb{Q}, 2}}{2}}
      + \sup_\mathbb{Q} \log N \prn[\bigg]{\mathcal{G}, \norm{\cdot}_{\mathbb{Q}, 2}, \frac{\varepsilon \norm{G}_{\mathbb{Q}, 2}}{2}}
  \end{align*}
  for all $\varepsilon > 0$. Also, the uniform entropy numbers of $\mathcal{F} + \mathcal{G} = \crl{f + g: f \in \mathcal{F}, g \in \mathcal{G}}$ satisfy
  \begin{align*}
    &\sup_\mathbb{Q} \log N\prn[\big]{\mathcal{F} + \mathcal{G}, \norm{\cdot}_{\mathbb{Q}, 2}, \varepsilon \norm{F + G}_{\mathbb{Q}, 2}} \\
    &\quad\leq \sup_\mathbb{Q} \log N \prn[\bigg]{\mathcal{F}, \norm{\cdot}_{\mathbb{Q}, 2}, \frac{\varepsilon \norm{F}_{\mathbb{Q}, 2}}{2}} + \sup_\mathbb{Q} \log N \prn[\bigg]{\mathcal{G}, \norm{\cdot}_{\mathbb{Q}, 2}, \frac{\varepsilon \norm{G}_{\mathbb{Q}, 2}}{2}}
  \end{align*}
  for all $\varepsilon > 0$. In both cases, $\mathbb{Q}$ ranges over all finitely-discrete probability measures.
\end{lemma}

\begin{proof}
    This lemma is well-known. See, for example, \parencite{vaart96_weak_conver_empir_proces}.
\end{proof}

\begin{lemma}[Maximal inequality for Gaussian vectors]\label{lem:maximal-ineq-gaussian-vectors}
  Take $n \geq 2$. Let $X_i \sim \mathcal{N}(0, \sigma_i^2)$ for $1 \leq i \leq n$ (not necessarily independent), with $\sigma_i^2 \leq \sigma^2$. Then
  \begin{align*}
    &\Expect[\Big]{\max_{1 \leq i \leq n} X_i} \leq \sigma \sqrt{2 \log n}, \\
    &\Expect[\Big]{\max_{1 \leq i \leq n} \abs{X_i}} \leq 2 \sigma \sqrt{\log n}.
  \end{align*}

  If $\mSigma_1$ and $\mSigma_2$ are positive semi-definite $n \times n$ matrices and $\vn \sim \mathcal{N}(0, \mI_n)$, then
  \begin{equation*}
    \Expect{\norm{\mSigma_1^{1 / 2} \vn - \mSigma_2^{1 / 2} \vn}_{\infty}}
    \leq 2 \sqrt{\log n}\norm{\mSigma_1 - \mSigma_2}_2^{1 / 2}.
  \end{equation*}

  If further $\mSigma_1$ is positive definite, then
  \begin{equation*}
    \Expect{\norm{\mSigma_1^{1 / 2} \vn - \mSigma_2^{1 / 2} \vn}_{\infty}} \leq \sqrt{\log n} \lambda_{\min}(\mSigma_1)^{-1 / 2} \norm{\mSigma_1 - \mSigma_2}_2.
  \end{equation*}
\end{lemma}

\begin{proof}
    See Lemma~SA31 in~\parencite{Cattaneo-Feng-Underwood_2024_JASA}.
\end{proof}

\begin{lemma}[A maximal inequality for i.n.i.d. empirical processes]\label{lem:ueb-maximal-ineq}
  Let $X_1, \ldots, X_n$ be independent but not necessarily identically distributed (i.n.i.d.) random variables taking values in a measurable space $(S, \mathcal{S})$. Denote the joint distribution of $X_1, \ldots, X_n$ by $\mathbb{P}$ and the marginal distribution of $X_i$ by $\mathbb{P}_i$, and let $\overline{\mathbb{P}}=n^{-1} \sum_i \mathbb{P}_i$.

  Let $\mathcal{F}$ be a class of Borel measurable functions from $S$ to $\mathbb{R}$ which is pointwise measurable (i.e. it contains a countable subclass which is dense under pointwise convergence), and satisfying the uniform entropy bound~\eqref{eq:uniform-entropy-bound} with parameters $A$ and $V$. Let $F$ be a strictly positive measurable envelope function for $\mathcal{F}$ (i.e. $|f(s)| \leq|F(s)|$ for all $f \in \mathcal{F}$ and $s \in S$). Suppose that $\|F\|_{\overline{\mathbb{P}}, 2}<\infty$. Let $\sigma>0$ satisfy $\sup _{f \in \mathcal{F}}\|f\|_{\overline{\mathbb{P}}, 2} \leq \sigma \leq\|F\|_{\overline{\mathbb{P}}, 2}$ and $M = \max _{1 \leq i \leq n} F(X_i)$.

  For $f \in \mathcal{F}$ define the empirical process
  \begin{equation*}
    G_n(f) = \frac{1}{\sqrt{n}} \sum_{i=1}^n(f(X_i)-\mathbb{E}[f(X_i)]).
  \end{equation*}

  Then we have
  \begin{equation*}
    \Expect[\Big]{\sup _{f \in \mathcal{F}}\abs{G_n(f)}}
    \lesssim \sigma \sqrt{V \log \prn[\big]{A \norm{F}_{\overline{\mathbb{P}}, 2} / \sigma}} + \frac{\norm{M}_{\mathbb{P}, 2} V \log \prn[\big]{A \norm{F}_{\overline{\mathbb{P}}, 2} / \sigma}}{\sqrt{n}},
  \end{equation*}
  where $\lesssim$ is up to a universal constant.
\end{lemma}

\begin{proof}
    See Lemmas SA24 and SA25 in \parencite{Cattaneo-Feng-Underwood_2024_JASA}.
\end{proof}

\begin{lemma}[Maximal inequalities for Gaussian processes]\label{lem:maximal-gaussian}
  Let $Z$ be a separable mean-zero Gaussian process indexed by $x \in \mathcal{X}$. Recall that $Z$ is separable for example if $\mathcal{X}$ is Polish and $Z$ has continuous trajectories. Define its covariance structure on $\mathcal{X} \times \mathcal{X}$ by $\Sigma\prn{x, x'} := \mathbb{E}\brk{Z(x) Z(x')}$, and the corresponding semimetric on $\mathcal{X}$ by
  \begin{equation*}
    \rho(x, x') := \mathbb{E}\brk{\prn{Z(x)-Z(x')}^2}^{1 / 2} = \prn[\big]{\Sigma(x, x)-2 \Sigma\prn{x, x'}+\Sigma\prn{x', x'}}^{1 / 2} .
  \end{equation*}
  Let $N(\mathcal{X}, \rho, \varepsilon)$ denote the $\varepsilon$-covering number of $\mathcal{X}$ with respect to the semimetric $\rho$. Define $\sigma := \sup_x \Sigma(x, x)^{1 / 2}$.

  Then there exists a universal constant $C > 0$ such that for any $\delta > 0$,
  \begin{align*}
    &\mathbb{E}\brk[\Big]{\sup_{x \in \mathcal{X}}|Z(x)|} \leq C \sigma+C \int_0^{2 \sigma} \sqrt{\log N(\mathcal{X}, \rho, \varepsilon)}\, \mathrm{d} \varepsilon, \\
    &\mathbb{E}\brk[\Big]{\sup _{\rho\prn{x, x'} \leq \delta} \abs{Z(x)-Z\prn{x'}}} \leq C \int_0^\delta \sqrt{\log N(\mathcal{X}, \rho, \varepsilon)}\, \mathrm{d} \varepsilon .
  \end{align*}
\end{lemma}

\begin{proof}
    This lemma is well-known. See, for example, \parencite{vaart96_weak_conver_empir_proces}.
\end{proof}

\begin{lemma}[Closeness in probability implies closeness of conditional quantiles]\label{lem:closeness-in-prob-implies-closeness-of-cond-quantiles}
  Let $X_n$ and $Y_n$ be random variables and $\vD_n$ be a random vector. Let $F_{X_n}\prn{x | \vD_n}$ and $F_{Y_n}\prn{x | \vD_n}$ denote the conditional distribution functions, and $F^{-1}_{X_n}\prn{x | \vD_n}$ and $F^{-1}_{Y_n}\prn{x | \vD_n}$ denote the corresponding conditional quantile functions. If $\abs{ X_n - Y_n } = o_{\Problet} \prn{ r_n }$, then there exists a sequence of positive numbers $\nu_n \to 0$, depending on $r_n$, such that w.\,p.\,a.~1
  \begin{equation*}
    F^{-1}_{X_n}\prn{p | \vD_n} \leq F^{-1}_{Y_n}\prn{p + \nu_n | \vD_n} + r_n\quad \text{and}\quad F^{-1}_{Y_n}\prn{p | \vD_n} \leq F^{-1}_{X_n}\prn{p + \nu_n | \vD_n} + r_n
  \end{equation*}
  for all $p \in (\nu_n, 1 - \nu_n)$.
\end{lemma}

\begin{proof}
    See Lemma~13 in~\parencite{Belloni-Chernozhukov-Chetverikov-FernandezVal_2019_JoE}.
\end{proof}

\begin{lemma}[Anti-concentration for suprema of separable Gaussian processes]\label{lem:anti-concentration-for-suprema-of-separable-gaussian-processes}
  Let $X = (X_t)_{t \in T}$ be a mean-zero separable Gaussian process indexed by a semimetric space $T$ such that $\Expect{X_t^2} = 1$ for all $t \in T$. Then for any $\varepsilon > 0$,
  \begin{equation*}
    \sup_{u \in \mathbb{R}} \Prob[\bigg]{\abs[\bigg]{\sup_{t \in T} |X_t| - u} \leq \varepsilon} \leq 4 \varepsilon \prn[\bigg]{\Expect[\bigg]{\sup_{t \in T} |X_t|} + 1}.
  \end{equation*}
\end{lemma}

\begin{proof}
    See Corollary 2.1 in~\parencite{Chernozhukov-Chetverikov-Kato_2014b_AoS}.
\end{proof}

The following lemma appears to be new to the literature at the level of generality considered. It guarantees the existence and gives some properties of the main estimand considered in this paper.

\begin{lemma}[The existence of \texorpdfstring{$\mu(\cdot)$}{mu}]\label{lem:mu-borel-measurable}
  Suppose Assumptions \ref{ass:dgp-holder}\ref{ass:random-sample} and \ref{ass:dgp-holder}\ref{ass:cond-density} hold. We will suppress the dependence of $\riskfunc{\cdot}{\cdot}$ on $\vq$ in this lemma because the result can be applied separately for each $\vq$. Assume $\eta \mapsto \riskfunc{y}{\eta}$ is convex on $\rangeeta$, $\rangeeta$ is an open connected subset of $\mathbb{R}$, $\riskder{y}{\cdot}$ is the left or right derivative of $\riskfunc{y}{\cdot}$ (in particular, it is a subgradient: $\prn{\eta_1 - \eta_0} \riskder{y}{\eta_0} \leq \riskfunc{y}{\eta_1} - \riskfunc{y}{\eta_0}$), and $\riskder{y}{\eta}$ is strictly increasing in $\eta$ for any fixed $y \in \mathcal{Y}$. Assume the real inverse link function $\eta(\cdot)\colon \mathbb{R} \to \rangeeta$ is strictly monotonic and two times continuously differentiable.

  Denoting $a_l$ and $a_r$ the left and right ends of $\rangeeta$ respectively (possibly $\pm \infty$), assume that for each $\vx \in \mathcal{X}$ the expectation $\Expect{\riskder{y_i}{\zeta} \given \vx_i = \vx}$ is negative for (real deterministic) $\zeta$ in a neighborhood of $a_l$, positive for $\zeta$ in a neighborhood of $a_r$, and continuous in $\zeta$ (in particular, it crosses zero).

  Then for each $\vx \in \mathcal{X}$ the number $\min \crl{ \zeta: \Expect{\riskder{y_i}{\zeta} \given \vx_i = \vx} \geq 0 }$ exists and belongs to $\rangeeta$. Moreover,
  \begin{equation}\label{eq:mu-for-lp-reg-definition}
    \begin{aligned}
      &\mu_0(\vx) := \eta^{-1}\prn[\big]{\min \crl[\big]{ \zeta: \Expect{\riskder{y_i}{\zeta} \given \vx_i = \vx} \geq 0}}\\
      &\quad = \eta^{-1}\prn[\big]{\min \crl[\big]{ \zeta: \Expect{\riskder{y_i}{\zeta} \given \vx_i = \vx} = 0}}
    \end{aligned}
  \end{equation}
  defines a Borel-measurable function such that for all $\vx \in \mathcal{X}$
  \begin{equation*}
    \mu_0(\vx) \in \argmin_{\zeta \in \mathbb{R}} \Expect{\riskfunc{y_i}{\eta(\zeta)} \given \vx_i = \vx}.
  \end{equation*}
\end{lemma}

  If $\mathcal{Q}$ is not a singleton, applying this result for each $\vq \in \mathcal{Q}$ gives a function $\muzero[\vq]{\vx}$ which is Borel in $\vx$ for each fixed $\vq$. Measurability in $\vq$ is not asserted by this lemma.

\begin{proof}
The conditions ensure that $\min \crl[\big]{\zeta: \Expect{\riskder{y_i}{\zeta} \given \vx_i = \vx} \geq 0}$ exists and belongs to $\rangeeta$ by continuity.

So the function $\mu_0(\vx)$ is well-defined. It is Borel because $\eta^{-1}$ is continuous and
\begin{equation*}
  \crl{ \vx: \eta\prn{\mu_0(\vx)} > a } = \crl{ \vx: \Expect{\riskder{y_i}{a} \given \vx_i = \vx} < 0 }
\end{equation*}
is a Borel set (the equality of the two sets is true because $\zeta \mapsto \riskder{y}{\zeta}$ is strictly increasing).

For any $\zeta \in \mathbb{R}$, using $\prn{\eta(\zeta) - \eta\prn{\mu_0(\vx)}} \riskder{y}{\eta\prn{\mu_0(\vx)}} \leq \riskfunc{y}{\eta(\zeta)} - \riskfunc{y}{\eta\prn{\mu_0(\vx)}}$, we have
\begin{align*}
  &0 = \prn{\eta(\zeta) - \eta\prn{\mu_0(\vx)}} \left.\Expect[\big]{\riskder{y_i}{a} \given \vx_i = \vx} \right|_{a = \eta(\mu_0(\vx))}\\
  &\quad \leq \Expect{\riskfunc{y}{\eta(\zeta)} \given \vx_i = \vx} - \left.\Expect{\riskfunc{y}{a} \given \vx_i = \vx}\right|_{a = \eta(\mu_0(\vx))},
\end{align*}
so $\mu_0(\vx)$ is indeed the argmin.
\end{proof}

The following lemma establishes basic properties of the ``Gram'' (or Hessian, depending on the perspective) matrix generated by the partitioning-based M-estimator.

\begin{lemma}[Gram matrix]\label{lem:holder-gram}
  Suppose \cref{ass:dgp-holder,ass:gl-quasi-uniform-partition,ass:gl-local-basis,ass:gl-holder,ass:complexity,ass:gl-approximation-error} hold. Then
  \begin{gather}
    h^{d} \lesssim \inf_{\vq} \lambda_{\min}(\qzero[\vq]) \leq \sup_{\vq} \lambda_{\max}(\qzero[\vq]) \lesssim h^{d},\label{eq:q0q-eigenvalues-bound}\\
    \sup_{\vq} \norm[\big]{ \qzero[\vq][-1] }_{\infty} \lesssim h^{-d}.\label{eq:q0q-inv-inf-norm-bound}
  \end{gather}
  If, in addition, $\frac{\log(1 / h)}{n h^{d}} = o(1)$, then uniformly over $\vq \in \mathcal{Q}$
  \begin{gather}
    \sup_{\vq} \crl[\big]{ \norm[\big]{\qbar[\vq] - \qzero[\vq]}_{\infty} \vee \norm[\big]{\qbar[\vq] - \qzero[\vq]}} \lesssim_{\Problet} h^{d} \prn[\bigg]{\frac{\log (1 / h)}{n h^{d}}}^{1 / 2},\label{eq:qbarq-concentration}\\
    h^{d} \lesssim \inf_{\vq} \lambda_{\min}(\qbar[\vq]) \leq \sup_{\vq} \lambda_{\max}(\qbar[\vq]) \lesssim h^{d}\quad \text{w.\,p.\,a.~1},\label{eq:qbarq-eigenvalues-bound}\\
    \sup_{\vq} \norm[\big]{ \qbar[\vq][-1] }_{\infty} \lesssim h^{-d}\quad \text{w.\,p.\,a.~1},\label{eq:qbarq-inv-inf-norm-bound}\\
    \sup_{\vq} \crl[\big]{\norm[\big]{\qbar[\vq][-1]-\qzero[\vq][-1]}_{\infty} \vee \norm[\big]{\qbar[\vq][-1]-\qzero[\vq][-1]}} \lesssim_{\mathbb{P}} h^{-d} \prn[\bigg]{\frac{\log (1 / h)}{n h^{d}}}^{1 / 2}.\label{eq:qbarq-inv-concentration}
  \end{gather}

  For some positive integer $L_n \lesssim 1 / h$, the rows and columns of $\qbar[\vq]$ and its inverse can be numbered by multi-indices $(\boldsymbol{\alpha}, a) = (\alpha_1, \ldots, \alpha_d, a)$ and $(\vbeta, b) = (\beta_1, \ldots, \beta_d, b)$, where
  \begin{equation*}
    \boldsymbol{\alpha}, \vbeta \in \{1, \ldots, L_n\}^d,\quad a \in \crl{ 1, \ldots, T_{n, \boldsymbol{\alpha}} }, b \in \crl{ 1, \ldots, T_{n, \vbeta} },\quad T_{n, \boldsymbol{\alpha}}, T_{n, \vbeta} \lesssim 1,
  \end{equation*}
  in the following way. First, $\qbar[\vq]$ has a multi-banded structure:
  \begin{equation}\label{eq:multiband-property}
    \brk{ \qbar[\vq] }_{(\boldsymbol{\alpha}, a), (\vbeta, b)} = 0\quad\text{if $\norm{ \boldsymbol{\alpha} - \vbeta }_{\infty} > C$},
  \end{equation}
  for some constant $C > 0$ (not depending on $n$). Second, with probability approaching one
  \begin{equation}
    \label{eq:inv-multiband-exp-decay}
    \sup_{\vq} \abs[\big]{ \brk{ \qbar[\vq][-1] }_{(\boldsymbol{\alpha}, a), (\vbeta, b)}} \lesssim h^{-d} \varrho^{\|\boldsymbol{\alpha} - \vbeta\|_{\infty}}
  \end{equation}
  for some constant $\varrho \in (0, 1)$ (not depending on $n$).

  The same results hold with $\qzero[\vq]$ replaced by $\sigmazero[\vq]$ and $\qbar[\vq]$ replaced by $\sigmabar[\vq]$.

  Finally, the same results hold with $\qzero[\vq]$ replaced by $\Expect{\vp(\vx_i) \vp(\vx_i)\trans}$ and $\qbar[\vq]$ replaced by $\Expectlet_n \brk{\vp(\vx_i) \vp(\vx_i)\trans}$.
\end{lemma}

\begin{proof}
The last claim of the lemma, corresponding to the case
\begin{equation*}
  \Psione[\vq]{\vx_i}{\eta( \muzero[\vq]{\vx_i})} \eta^{(1)}\prn{\muzero[\vq]{\vx_i}}^2 \equiv 1,
\end{equation*}
is Lemma~SA-2.1 in~\parencite{Cattaneo-Farrell-Feng_2020_AOS}. The properties~\eqref{eq:multiband-property} and~\eqref{eq:inv-multiband-exp-decay} are not explicitly stated but follow from the proof.

In the general case, by \cref{ass:gl-holder}\ref{ass:psi-eta} $\Psione[\vq]{\vx_i}{\eta(\muzero[\vq]{\vx_i})} \eta^{(1)}\prn{\muzero[\vq]{\vx_i}}^2$ is bounded and bounded away from zero uniformly over $i$, $n$ and $\vq$, so \eqref{eq:q0q-eigenvalues-bound} and \eqref{eq:qbarq-eigenvalues-bound} follow from the previous case. The additional $\Psione[\vq]{\vx_i}{\eta( \muzero[\vq]{\vx_i})} \eta^{(1)}\prn{\muzero[\vq]{\vx_i}}^2$ term does not influence the multi-banded structure of the matrices, so \eqref{eq:q0q-inv-inf-norm-bound}, \eqref{eq:qbarq-inv-inf-norm-bound}, \eqref{eq:multiband-property}, \eqref{eq:inv-multiband-exp-decay} remain true by the same argument as in the previous case. The inequalities
\begin{align*}
  \norm[\big]{\qbar[\vq][-1]-\qzero[\vq][-1]}_{\infty} &\leq \norm[\big]{\qbar[\vq][-1]}_{\infty} \cdot \norm[\big]{\qbar[\vq]-\qzero[\vq]}_{\infty} \cdot
                                                           \norm[\big]{\qzero[\vq][-1]}_{\infty},\\
  \norm[\big]{\qbar[\vq][-1]-\qzero[\vq][-1]} &\leq \norm[\big]{\qbar[\vq][-1]} \cdot \norm[\big]{\qbar[\vq]-\qzero[\vq]} \cdot
                                                  \norm[\big]{\qzero[\vq][-1]}
\end{align*}
show that \eqref{eq:qbarq-inv-concentration} follows from norm bounds \eqref{eq:q0q-eigenvalues-bound}, \eqref{eq:q0q-inv-inf-norm-bound}, \eqref{eq:qbarq-eigenvalues-bound}, \eqref{eq:qbarq-inv-inf-norm-bound} and concentration \eqref{eq:qbarq-concentration}.

Now we prove \cref{eq:qbarq-concentration}.

Define the class of functions
\begin{equation*}
  \mathcal{G} := \crl[\big]{\vx \mapsto p_{k}\prn{\vx} p_{l}\prn{\vx} \Psione[\vq]{\vx}{\eta(\muzero[\vq]{\vx})} \eta^{(1)}\prn{\muzero[\vq]{\vx}}^2\colon 1 \leq k, l \leq K, \vq \in \mathcal{Q}}.
\end{equation*}

We will now prove that the class $\mathcal{G}$ with a large enough constant envelope satisfies the uniform entropy bound~\eqref{eq:uniform-entropy-bound} with $A \lesssim h^{- 2 d}$ and $V \lesssim 1$. By \cref{ass:complexity}, the class
\begin{equation*}
\crl{ \vx \mapsto \Psione[\vq]{\vx}{\eta( \muzero[\vq]{\vx})} \colon \vq \in \mathcal{Q} }
\end{equation*}
with a large enough constant envelope satisfies the uniform entropy bound~\eqref{eq:uniform-entropy-bound} with $A \lesssim 1$ and $V \lesssim 1$. Since it is also true of the class $\crl[\big]{\eta^{(1)}\prn[\big]{\muzero[\vq]{\vx}}^2 \colon \vq \in \mathcal{Q}}$ because $\eta^{(1)}\prn[\big]{\muzero[\vq]{\vx}}^2$ is Lipschitz in $\vq$, by \cref{lem:permanence-properties} the class
\begin{equation*}
  \crl[\big]{ \vx \mapsto \Psione[\vq]{\vx}{\eta( \muzero[\vq]{\vx})} \eta^{(1)}\prn{\muzero[\vq]{\vx}}^2 \colon \vq \in \mathcal{Q} }
\end{equation*}
with a large enough constant envelope satisfies the uniform entropy bound~\eqref{eq:uniform-entropy-bound} with $A \lesssim 1$ and $V \lesssim 1$. The class $\crl{ \vx \mapsto p_{k}\prn{\vx} p_{l}\prn{\vx} \colon 1 \leq k, l \leq K }$ with a large enough constant envelope just contains $K^2$ functions, so it also satisfies the uniform entropy bound~\eqref{eq:uniform-entropy-bound} with $A \lesssim h^{- 2 d}$ and $V = 1$, where we used $K \asymp h^{-d}$. By \cref{lem:permanence-properties}, combining these facts proves the claim about the complexity of $\mathcal{G}$.

Moreover, class $\mathcal{G}$ satisfies the following variance bound:
\begin{equation*}
  \sup_{g \in \mathcal{G}} \Expect*{g(\vx_i)^2} \lesssim h^{d},
\end{equation*}
which follows from the fact that the class is bounded by a large enough constant and the Lebesgue measure of the support of $p_k(\vx) p_l(\vx)$ shrinks (uniformly over $k, l$) at the rate $h^{d}$.

Applying \cref{lem:ueb-maximal-ineq}, we see that
\begin{align*}
  &\sup _{g \in \mathcal{G}}\abs[\bigg]{\frac{1}{n} \sum_{i=1}^{n} g\prn{\vx_{i}}-\mathbb{E}\brk{g\prn{\vx_{i}}}} \lesssim_{\mathbb{P}} h^{d} \prn[\bigg]{\frac{\log (1 / h)}{n h^d}}^{1 / 2} + \frac{\log (1 / h)}{n}\\
  &\quad \lesssim h^{d} \prn[\bigg]{\frac{\log (1 / h)}{n h^d}}^{1 / 2} \tag*{since $\frac{\log (1 / h)}{n h^d} = o(1)$.}
\end{align*}

So we have shown $\max_{k, l} \abs[\big]{(\qbar[\vq] - \qzero[\vq])_{k, l}} \lesssim_{\Problet} h^{d} \prn[\Big]{\frac{\log (1 / h)}{n h^d}}^{1 / 2}$. Since each row and column of $\qbar[\vq] - \qzero[\vq]$ has a bounded number of nonzero entries, this implies
\begin{equation*}
  \norm{\qbar[\vq] - \qzero[\vq]}_{\infty} \asymp \norm{\qbar[\vq] - \qzero[\vq]}_{1} \lesssim_{\Problet} h^{d} \prn[\bigg]{\frac{\log (1 / h)}{n h^d}}^{1 / 2}.
\end{equation*}
To conclude the proof of \cref{eq:qbarq-concentration}, it is left to use the inequality
\begin{equation*}
  \norm{ \qbar[\vq] - \qzero[\vq] } \leq \sqrt{\norm{ \qbar[\vq] - \qzero[\vq] }_1 \norm{ \qbar[\vq] - \qzero[\vq] }_{\infty}}.
\end{equation*}

The claim about $\sigmazero[\vq]$ and $\sigmabar[\vq]$ is proven analogously, using that $\sigma_{\vq}^{2}(\vx)$ is bounded and bounded away from zero and Lipschitz in $\vq$ by \cref{ass:gl-holder}\ref{ass:psi-moments}.

\cref{lem:holder-gram} is proven.
\end{proof}


The following \cref{lem:holder-uniform-convergence-variance,lem:holder-projection-approximation-error}
are needed for the proof of the Bahadur representation (\cref{th:bahadur-repres}), \cref{cor:uniform-rate-of-convergence,cor:msq-rate-of-convergence} and a version of the consistency result (\cref{lem:holder-consistency-strong-shape}).

\begin{lemma}[Uniform convergence: variance]\label{lem:holder-uniform-convergence-variance}
  Suppose \cref{ass:dgp-holder,ass:gl-quasi-uniform-partition,ass:gl-local-basis,ass:gl-holder,ass:complexity,ass:gl-approximation-error} hold. If
  \begin{align}
    &\frac{\brk[\big]{\log (1 / h)}^{\nu / (\nu - 2)}}{n h^{\nu d / (\nu - 2)}} = o(1),\text{ or}\label{eq:uniform-convergence-variance-moment-case}\\
    &\riskderenvel{\vx_i}{y_i} \text{ is \subgauss\ conditionally on $\vx_i$ and }\frac{\log n \log(1 / h)}{n h^d} = o(1),\label{eq:uniform-convergence-variance-subgaussian-case}
  \end{align}
  then
  \begin{equation*}
    \sup _{\vq \in \mathcal{Q}, \vx \in \mathcal{X}}\abs[\big]{\vp^{(\vv)}(\vx)\trans \qbar[\vq][-1] \Expectlet_n \brk[\big]{\vp\prn{\vx_{i}} \eta^{(1)}\prn{\muzero[\vq]{\vx_i}} \riskder[\vq]{y_i}{\eta\prn{\muzero[\vq]{\vx_i}}}}} \lesssim_{\mathbb{P}} h^{- \multindnorm{\vv}}\prn[\bigg]{\frac{\log (1 / h)}{n h^d}}^{1 / 2};
  \end{equation*}
and the same inequality is true with $\qbar[\vq]$ replaced by $\qzero[\vq]$.
\end{lemma}

\begin{proof}
By \cref{ass:gl-local-basis}, $\sup_{\vx \in \mathcal{X}} \norm{ \vp^{(\vv)}(\vx) } \lesssim h^{-\multindnorm{\vv}}$; by \cref{lem:holder-gram}, $\norm{\qbar[\vq][-1] }_{\infty} + \norm{\qzero[\vq][-1] }_{\infty} \lesssim_{\Problet} h^{-d}$. Therefore, it is enough to show
\begin{equation}\label{eq:uniform-conv-variance-enough}
  \sup _{\vq\in\mathcal{Q}}\norm[\big]{\Expectlet_n \brk[\big]{\vp\prn{\vx_{i}} \eta^{(1)}\prn{\muzero[\vq]{\vx_i}} \riskder[\vq]{y_i}{\eta\prn{\muzero[\vq]{\vx_i}}}}}_{\infty} \lesssim_{\mathbb{P}} \sqrt{\frac{h^{d} \log (1 / h)}{n}}.
\end{equation}

Define the function class
\begin{equation*}
  \mathcal{G} := \crl[\big]{\prn{\vx, y} \mapsto p_{l}(\vx) \eta^{(1)}\prn{\muzero[\vq]{\vx}} \riskder[\vq]{y}{\eta\prn{\muzero[\vq]{\vx}}}: 1 \leq l \leq K, \vq \in \mathcal{Q}}.
\end{equation*}

We will now control the complexity of $\mathcal{G}$. Introduce some more classes of functions:
\begin{align*}
  \mathcal{W}_1 &:= \crl{ \prn{\vx, y} \mapsto p_{l}(\vx): 1 \leq l \leq K },\\
  \mathcal{W}_2 &:= \crl{ \prn{\vx, y} \mapsto \eta^{(1)} \prn{ \muzero[\vq]{\vx} }: \vq \in \mathcal{Q} },\\
  \mathcal{W}_3 &:= \crl{ \prn{\vx, y} \mapsto \riskder[\vq]{y}{\eta\prn{\muzero[\vq]{\vx}}}: \vq \in \mathcal{Q} }.
\end{align*}
$\mathcal{W}_1$ with a large enough constant envelope contains $K$ fixed measurable functions, so it satisfies the uniform entropy bound~\eqref{eq:uniform-entropy-bound} with $A \lesssim h^{-d}$ and $V = 1$. $\mathcal{W}_2$ with a large enough constant envelope satisfies the uniform entropy bound~\eqref{eq:uniform-entropy-bound} with $A, V \lesssim 1$ because $\muzero[\vq]{\vx}$ is bounded uniformly over $\vx, \vq$ and Lipschitz in $\vq$, $\eta^{(1)}$ on a fixed bounded interval is Lipschitz. $\mathcal{W}_3$ with envelope $\riskderenvel{\vx}{y}$ satisfies the uniform entropy bound~\eqref{eq:uniform-entropy-bound} with $A, V \lesssim 1$ by \cref{ass:complexity}. By \cref{lem:permanence-properties}, $\mathcal{G}$ with envelope $\riskderenvel{\vx}{y}$ multiplied by a large enough constant satisfies the uniform entropy bound~\eqref{eq:uniform-entropy-bound} with $A \lesssim h^{-d}$ and $V \lesssim 1$.

Moreover, class $\mathcal{G}$ satisfies the following variance bound:
  \begin{equation*}
    \sup_{g \in \mathcal{G}} \Expect*{g(\vx_i, y_i)^2} \lesssim h^{d}.
  \end{equation*}
Indeed, for a fixed $i \in \{1, \ldots, n\}$
\begin{align*}
  \sup_{g \in \mathcal{G}} \Expect*{g(\vx_i, y_i)^2} \lesssim \sup_{l} \Expect*{p_l(\vx_i)^2 \Expect*{\riskderenvel{\vx_i}{y_i}^2 \given \vx_i}} \lesssim \sup_l \Expect*{p_l(\vx_i)^2} \lesssim h^{d}.
\end{align*}

Finally, under~\eqref{eq:uniform-convergence-variance-moment-case}
\begin{align*}
  &\Expect[\Big]{\max_{1 \leq i \leq n} \abs{ \riskderenvel{\vx_i}{y_i} }^2}^{1 / 2} \leq \Expect[\Big]{\max_{1 \leq i \leq n} \abs{\riskderenvel{\vx_i}{y_i}}^\nu}^{1 / \nu} \leq \Expect[\bigg]{\sum_{i = 1}^{n} \abs{ \riskderenvel{\vx_i}{y_i} }^\nu}^{1 / \nu}\\
  &\quad = \prn[\bigg]{ \sum_{i = 1}^n \Expect[\big]{\riskderenvel{\vx_i}{y_i}^{\nu}} }^{1 / \nu} \lesssim \prn[\bigg]{ \sum_{i = 1}^n 1 }^{1 / \nu} = n^{1 / \nu},
\end{align*}
and under~\eqref{eq:uniform-convergence-variance-subgaussian-case}
\begin{equation*}
  \Expect[\Big]{\max_{1 \leq i \leq n} \abs{ \riskderenvel{\vx_i}{y_i} }^2}^{1 / 2} \lesssim \sqrt{\log n}
\end{equation*}
by \cref{lem:subgaussian-max-bound}.

Applying \cref{lem:ueb-maximal-ineq}, we obtain~\eqref{eq:uniform-conv-variance-enough} since
\begin{align*}
  &\frac{n^{1 / \nu} \log(1 / h)}{n} = \sqrt{\frac{h^{d} \log(1 / h)}{n}} \cdot \sqrt{\frac{\log (1 / h)}{n^{1 - 2 / \nu} h^d}} = \sqrt{\frac{h^{d} \log(1 / h)}{n}} \cdot o(1),\text{ and}\\
  &\frac{\sqrt{\log n} \log(1 / h)}{n} = \sqrt{\frac{h^d \log(1 / h)}{n}} \cdot \sqrt{\frac{\log n \log(1 / h)}{n h^d}} = \sqrt{\frac{h^d \log(1 / h)}{n}} \cdot o(1).
\end{align*}

\cref{lem:holder-uniform-convergence-variance} is proven.
\end{proof}

\begin{remark}\label{rem:rates-in-uniform-conv-var}
  Using
\begin{equation}\label{eq:XJQNHz}
\frac{n^{1 / \nu} \log(1 / h)}{n} = o(h^d) \Leftrightarrow \frac{\brk{\log(1 / h)}^{\nu / (\nu - 1)}}{n h^{\nu d / (\nu - 1)}} = o(1)
\end{equation}
instead of
\begin{equation*}
\frac{n^{1 / \nu} \log(1 / h)}{n} = \sqrt{\frac{h^{d} \log(1 / h)}{n}} \cdot o(1)
\end{equation*}
in the argument leads to the bound
\begin{equation*}
  \sup _{\vq\in\mathcal{Q}}\norm[\big]{\Expectlet_n \brk[\big]{\vp\prn{\vx_{i}} \eta^{(1)}\prn{\muzero[\vq]{\vx_i}} \riskder[\vq]{y_i}{\eta\prn{\muzero[\vq]{\vx_i}}}}}_{\infty} = o_{\Problet}(h^d)
\end{equation*}
under just the condition \eqref{eq:XJQNHz} which is slightly weaker than \eqref{eq:uniform-convergence-variance-moment-case}.
\end{remark}

The following lemma gives control on the projection approximation error.

\begin{lemma}[Projection of approximation error]\label{lem:holder-projection-approximation-error}
  Suppose \cref{ass:dgp-holder,ass:gl-quasi-uniform-partition,ass:gl-local-basis,ass:gl-holder,ass:complexity,ass:gl-approximation-error} hold. If
  \begin{align}
    &\frac{\brk[\big]{\log (1 / h)}^{\nu / (\nu - 2)}}{n h^{\nu d / (\nu - 2)}} = o(1),\text{ or}\label{eq:projection-approx-moment-case}\\
    &\riskderenvel{\vx_i}{y_i} \text{ is \subgauss\ conditionally on $\vx_i$ and }\frac{\log n \log(1 / h)}{n h^d} = o(1),\nonumber
  \end{align}
  then
  \begin{align*}
    &\sup _{\vq \in \mathcal{Q}, \vx \in \mathcal{X}}
      \abs[\big]{
      \vp^{(\vv)}(\vx)\trans \qbar[\vq][-1]
      \Expectlet_n \brk[\big]{
      \vp\prn{\vx_{i}}
      \crl[\big]{
      \riskder[\vq]{y_i}{\eta\prn{\muzero[\vq]{\vx_i}}} \eta^{(1)}\prn{\muzero[\vq]{\vx_i}}\\
    &\quad\quad- \riskder[\vq]{y_i}{\eta\prn{\vp(\vx_i)\trans \betazero[\vq]}} \eta^{(1)}\prn{\vp(\vx_i)\trans \betazero[\vq]}
      }
      }
      }\\
    &\quad\lesssim_{\mathbb{P}} h^{m - \multindnorm{\vv}} + h^{\prn{ \alpha \wedge (1 / 2) } m - \multindnorm{\vv}} \prn[\bigg]{\frac{\log (1 / h)}{n h^d}} ^{1 / 2} + \frac{\log (1 / h)}{n h^{\multindnorm{\vv} + d}}.
  \end{align*}
\end{lemma}

\begin{proof}
Denote
\begin{gather*}
  A_{1, \vq}(\vx_i, y_i) :=
  \riskder[\vq]{y_i}{\eta\prn{\muzero[\vq]{\vx_i}}}
  \crl[\big]{
    \eta^{(1)}\prn{\muzero[\vq]{\vx_i}}
    - \eta^{(1)}\prn{\vp(\vx_i)\trans \betazero[\vq]}
  },\\
  A_{2, \vq}(\vx_i, y_i) :=
  \crl[\big]{
    \riskder[\vq]{y_i}{\eta\prn{\muzero[\vq]{\vx_i}}}
    - \riskder[\vq]{y_i}{\eta\prn{\vp(\vx_i)\trans \betazero[\vq]}}
  }
  \eta^{(1)}\prn{\vp(\vx_i)\trans \betazero[\vq]}.
\end{gather*}

By \cref{ass:gl-local-basis}, $\sup_{\vx \in \mathcal{X}} \norm{ \vp^{(\vv)}(\vx) } \lesssim h^{- \multindnorm{\vv}}$; by \cref{lem:holder-gram}, $\norm{\qbar[\vq][-1] }_{\infty} \lesssim_{\Problet} h^{-d}$. Therefore, it is enough to show
\begin{align*}
  &\sup _{\vq \in \mathcal{Q}} \left\|
    \Expectlet_n\left[
    \vp\prn{\vx_{i}}
    \left\{
    \riskder[\vq]{y_i}{\eta\prn{\muzero[\vq]{\vx_i}}} \eta^{(1)}\prn{\muzero[\vq]{\vx_i}}\right.\right.\right.\\
  &\quad\quad\left.\left.\left.- \riskder[\vq]{y_i}{\eta\prn{\vp(\vx_i)\trans \betazero[\vq]}} \eta^{(1)}\prn{\vp(\vx_i)\trans \betazero[\vq]}
    \right\}
    \right]
    \right\|_{\infty}\\
  &\quad\lesssim_{\mathbb{P}} h^{d + m} + h^{\frac{d}{2} + \prn{\alpha \wedge (1 / 2)} m} \prn[\bigg]{\frac{\log (1 / h)}{n}}^{1 / 2} + \frac{\log (1 / h)}{n}.
\end{align*}

We will do this by showing the three bounds
\begin{align}
  &\sup_{\vq \in \mathcal{Q}} \norm[\big]{
    \Expectlet_n \brk{
    \vp(\vx_i) A_{1, \vq}(\vx_i, y_i)
    }
    }_{\infty}
    \lesssim_{\Problet} h^{\frac{d}{2} + m} \prn[\bigg]{\frac{\log (1 / h)}{n}}^{1 / 2},\label{eq:a1q-conc}\\
  &\sup_{\vq \in \mathcal{Q}} \norm[\big]{
    \Expectlet_n \brk[\big]{
    \Expect[\big]{
    \vp(\vx_i) A_{2, \vq}(\vx_i, y_i) \given \vx_i
    }
    }
    }_{\infty}
    \lesssim_{\Problet} h^{d + m},\label{eq:a2condexpbound}\\
  &\sup_{\vq \in \mathcal{Q}} \norm[\big]{
      \Expectlet_n \brk[\big]{
      \vp(\vx_i) \prn{ A_{2, \vq}(\vx_i, y_i)
      - \Expect*{A_{2, \vq}(\vx_i, y_i) \given \vx_i} }
      }
      }_{\infty} \lesssim_{\Problet} h^{\frac{d}{2} + \prn{\alpha \wedge (1 / 2)} m} \sqrt{\frac{\log (1 / h)}{n}} + \frac{\log (1 / h)}{n}.\label{eq:a2diffbound}
\end{align}

To show~\eqref{eq:a1q-conc}, consider the class of functions
\begin{equation*}
  \crl{(\vx, y) \mapsto p_{l}(\vx) A_{1, \vq}(\vx, y): 1 \leq l \leq K, \vq \in \mathcal{Q} }.
\end{equation*}
Note that $\sup_{\vq, \vx}\abs{ \eta^{(1)}\prn{\muzero[\vq]{\vx}} - \eta^{(1)}\prn{\vp(\vx)\trans \betazero[\vq]} } \lesssim h^{m}$ by \cref{ass:gl-approximation-error}. \eqref{eq:a1q-conc} follows by the same concentration argument as in \cref{lem:holder-uniform-convergence-variance}.

To show~\eqref{eq:a2condexpbound}, note that
\begin{align*}
  &\Expect*{\riskder[\vq]{y_i}{\eta\prn{\muzero[\vq]{\vx_i}}}
    - \riskder[\vq]{y_i}{\eta\prn{\vp(\vx_i)\trans \betazero[\vq]}} \given \vx_i}\\
  &\quad= - \expriskder[\vq]{\vx_i}{\eta(\vp(\vx_i)\trans\betazero[\vq])} =\\
  &\quad= \Psione[\vq]{\vx_i}{\zeta} \eta^{(1)}(\tilde{\zeta})\crl{\muzero[\vq]{\vx_i} - \vp(\vx_i)\trans \betazero[\vq]},
\end{align*}
where $\zeta$ is between $\eta\prn{\muzero[\vq]{\vx_i}}$ and $\eta\prn{\vp(\vx_i)\trans \betazero[\vq]}$, $\tilde{\zeta}$ is between $\muzero[\vq]{\vx_i}$ and $\vp(\vx_i)\trans \betazero[\vq]$. By \cref{ass:gl-holder}\ref{ass:psi-eta} and \ref{ass:gl-approximation-error}, it follows that a.\,s.
\begin{equation*}
  \sup_{\vq\in\mathcal{Q}} \abs{ \Expect{\riskder[\vq]{y_i}{\eta\prn{\muzero[\vq]{\vx_i}}}
      - \riskder[\vq]{y_i}{\eta\prn{\vp(\vx_i)\trans \betazero[\vq]}} \given \vx_i} } \lesssim h^{m}.
\end{equation*}
Since $\eta^{(1)}\prn{\vp(\vx_i)\trans \betazero[\vq]}$ is bounded, \eqref{eq:a2condexpbound} follows by applying \cref{lem:ueb-maximal-ineq} to the class $\crl{\vx \mapsto \abs{p_{l}(\vx)}, 1 \leq l \leq K}$.

It is left to show~\eqref{eq:a2diffbound}.

Consider the class of functions
\begin{equation*}
  \mathcal{G} := \crl{ \prn{\vx, y} \mapsto p_{l}(\vx) A_{2, \vq}(\vx, y): 1 \leq l \leq K, \vq \in \mathcal{Q} }.
\end{equation*}

We will  now control the complexity of $\mathcal{G}$. Introduce some more classes of functions:
\begin{align*}
  \mathcal{W}_{1, l} &:= \crl{ \prn{\vx, y} \mapsto p_l(\vx) \brk{\riskder[\vq]{y}{\eta\prn{\muzero[\vq]{\vx}}} - \riskder[\vq]{y}{\eta\prn{\vp(\vx)\trans \betazero[\vq]}}}: \vq \in \mathcal{Q} },\\
  \mathcal{W}_1 &:= \bigcup_{l = 1}^K \mathcal{W}_{1, l},\\
  \mathcal{W}_{2, l} &:= \crl{ \prn{\vx, y} \mapsto \eta^{(1)}\prn{\vp(\vx)\trans \betazero[\vq]} \Indlet\crl{\vx \in \supp p_l}: \vq \in \mathcal{Q} },\\
  \mathcal{W}_2 &:=  \bigcup_{l = 1}^K \mathcal{W}_{2, l} = \crl{ \prn{\vx, y} \mapsto \eta^{(1)}\prn{\vp(\vx)\trans \betazero[\vq]} \Indlet\crl{\vx \in \supp p_l}: 1 \leq l \leq K, \vq \in \mathcal{Q} }.
\end{align*}

By \cref{ass:complexity}, for $l$ fixed $\mathcal{W}_{1, l}$ with a large enough constant envelope (not depending on $l$) satisfies the uniform entropy bound~\eqref{eq:uniform-entropy-bound} with $A, V \lesssim 1$ (not depending on $l$). This immediately implies that $\mathcal{W}_1$ with a large enough constant envelope satisfies the uniform entropy bound~\eqref{eq:uniform-entropy-bound} with $A \lesssim h^{- d}$ and $V \lesssim 1$.

For $l$ fixed, $\mathcal{W}_{2, l}$ is a product of a (bounded) subclass of $\crl{ \eta^{(1)} \prn{ \vp(\vx)\trans \vbeta }: \vbeta \in \mathcal{B}_l }$, where $\mathcal{B}_l$ is a vector space of dimension $O(1)$ (not depending on $l$), and a fixed function. By Lemma~2.6.15 in~\parencite{vaart96_weak_conver_empir_proces}, $\crl{ \vp(\vx)\trans \vbeta: \vbeta \in \mathcal{B}_l }$ is VC with a bounded index. Therefore, since $\eta^{(1)}$ on a bounded interval is Lipschitz, $\mathcal{W}_{2, l}$ with a large enough constant envelope satisfies the uniform entropy bound~\eqref{eq:uniform-entropy-bound} with $A, V \lesssim 1$. This immediately implies that $\mathcal{W}_2$ with a large enough constant envelope satisfies the uniform entropy bound~\eqref{eq:uniform-entropy-bound} with $A \lesssim h^{-d}$ and $V \lesssim 1$.

By \cref{lem:permanence-properties}, it follows from the above that $\mathcal{G}$ with a large enough constant envelope satisfies the uniform entropy bound~\eqref{eq:uniform-entropy-bound} with $A \lesssim h^{-d}$ and $V \lesssim 1$.

Next, we will show that class $\mathcal{G}$ satisfies the following variance bound:
\begin{equation}\label{eq:projection-of-approx-variance}
  \sup_{g \in \mathcal{G}} \Expectlet_n \brk{\Var{g(\vx_i, y_i) \given \vx_i}} \lesssim h^{d + (2 \alpha \wedge 1) m}\quad\text{w.\,p.\,a.~1}.
\end{equation}
We will prove~\eqref{eq:projection-of-approx-variance} under the assumption that $\measureonY$ is Lebesgue measure, so~\eqref{eq:holder-assumption-with-cases-nonconvex} holds; the argument under~\eqref{eq:holder-assumption-without-cases-nonconvex} is similar (and leads to an even stronger variance bound), so it is omitted. For $y$ outside the closed segment between $\eta(\vp(\vx)\trans \betazero[\vq])$ and $\eta(\muzero[\vq]{\vx})$,
\begin{align*}
  &\abs{ \riskder[\vq]{y}{\eta \prn{ \muzero[\vq]{\vx} }} - \riskder[\vq]{y}{\eta \prn{ \vp(\vx)\trans \betazero[\vq] }} }\\
  &\quad \lesssim \prn{ \riskderenvel{\vx}{y} + 1 } \cdot \abs{ \vp(\vx)\trans \betazero[\vq] - \muzero[\vq]{\vx} }^{\alpha}\\
  &\quad \lesssim \prn{ \riskderenvel{\vx}{y} + 1 } h^{\alpha m},\numberthis\label{eq:psi-difference-outside-segment}
\end{align*}
where in~\eqref{eq:psi-difference-outside-segment} we used \cref{ass:gl-approximation-error}. For $y$ between $\eta(\vp(\vx)\trans \betazero[\vq])$ and $\eta(\muzero[\vq]{\vx})$ inclusive,
\begin{align}
  &\abs{ \riskder[\vq]{y}{\eta \prn{ \muzero[\vq]{\vx} }} - \riskder[\vq]{y}{\eta \prn{ \vp(\vx)\trans \betazero[\vq] }} } \lesssim \riskderenvel{\vx}{y} \cdot \abs{ \vp(\vx)\trans \betazero[\vq] - \muzero[\vq]{\vx} } + 1\nonumber\\
  &\quad \lesssim \riskderenvel{\vx}{y} h^m + 1,\label{eq:psi-difference-inside-segment}
\end{align}
where in~\eqref{eq:psi-difference-inside-segment} we again used \cref{ass:gl-approximation-error}.

In the chain below, to avoid cluttering notation we will use $[a, b]$ to denote the closed segment between $a$ and $b$ regardless of their ordering (a more standard notation is $[a \wedge b, a \vee b]$). Using that $\eta^{(1)} \prn{ \vp(\vx)\trans \betazero[\vq] }$ is also bounded uniformly over $\vx \in \mathcal{X}$, we have a.\,s.
\begin{align*}
  &\Expect[\big]{A_{2, \vq}(\vx_i, y_i)^2 \given \vx_i}\\
  &\quad = \Expect[\big]{A_{2, \vq}(\vx_i, y_i)^2 \Indlet \crl[\big]{ y_i \notin \brk[\big]{\eta(\vp(\vx_i)\trans \betazero[\vq]), \eta(\muzero[\vq]{\vx_i})} } \given \vx_i}\\
  &\qquad + \Expect[\big]{A_{2, \vq}(\vx_i, y_i)^2 \Indlet \crl[\big]{ y_i \in \brk[\big]{\eta(\vp(\vx_i)\trans \betazero[\vq]), \eta(\muzero[\vq]{\vx_i})} } \given \vx_i}\\
  &\quad \lesssim h^{2 \alpha m} \Expect[\big]{ \prn{ \riskderenvel{\vx_i}{y_i}^2 + 1 } \Indlet \crl[\big]{ y_i \notin \brk[\big]{\eta(\vp(\vx_i)\trans \betazero[\vq]), \eta(\muzero[\vq]{\vx_i})} } \given \vx_i}\\
  &\qquad + h^{2 m} \Expect[\big]{\riskderenvel{\vx_i}{y_i}^2 \Indlet \crl[\big]{ y_i \in \brk[\big]{\eta(\vp(\vx_i)\trans \betazero[\vq]), \eta(\muzero[\vq]{\vx_i})} } \given \vx_i}\\
  &\qquad + \Prob[\big]{y_i \in \brk[\big]{\eta(\vp(\vx_i)\trans \betazero[\vq]), \eta(\muzero[\vq]{\vx_i})} \given \vx_i}\\
  &\quad \leq h^{2 \alpha m} \Expect*{ \prn{ \riskderenvel{\vx_i}{y_i}^2 + 1 } \given \vx_i} + h^{2 m} \Expect*{\riskderenvel{\vx_i}{y_i}^2 \given \vx_i}\\
  &\qquad + \Prob*{y_i \in \brk[\big]{\eta(\vp(\vx_i)\trans \betazero[\vq]), \eta(\muzero[\vq]{\vx_i})} \given \vx_i}\\
  &\quad \lesssim h^{2 \alpha m} + h^{2 m} + \abs{ \vp(\vx_i)\trans \betazero[\vq] - \muzero[\vq]{\vx_i} }\numberthis\label{eq:proj-of-approx-error-big-chain-restr-lebesgue-meas}\\
  &\quad \lesssim h^{2 \alpha m} + h^{2 m} + h^m\numberthis\label{eq:proj-of-approx-error-big-chain-bounding-length}\\
  &\quad \lesssim h^{(2 \alpha \wedge 1) m},
\end{align*}
where in~\eqref{eq:proj-of-approx-error-big-chain-restr-lebesgue-meas} we used that by \cref{ass:dgp-holder}\ref{ass:cond-density} the conditional density of $y_i | \vx_i$ is bounded and \cref{ass:gl-holder}\ref{ass:psi-moments}, in~\eqref{eq:proj-of-approx-error-big-chain-bounding-length} we used \cref{ass:gl-approximation-error}.

Therefore, uniformly over $l$ and $\vq$
\begin{align*}
  &\Expectlet_n \brk{\Var{g(\vx_i, y_i) \given \vx_i}} \leq \Expectlet_n \brk[\big]{\Expect{g(\vx_i, y_i)^2 \given \vx_i}} = \Expectlet_n \brk[\big]{p_l(\vx_i)^2 \Expect{A_{2, \vq}(\vx_i, y_i)^2 \given \vx_i}}\\
  &\quad \lesssim h^{(2 \alpha \wedge 1) m} \Expectlet_n \brk{p_l(\vx_i)^2}\\
  &\quad \leq h^{(2 \alpha \wedge 1) m} \norm{ \Expectlet_n \brk{\vp(\vx_i) \vp(\vx_i)\trans} }\\
  &\quad \lesssim h^{d + (2 \alpha \wedge 1) m}\quad\text{w.\,p.\,a.~1},\numberthis\label{eq:prof-of-approx-error-final-appl-of-gram-norm-lem}
\end{align*}
where in~\eqref{eq:prof-of-approx-error-final-appl-of-gram-norm-lem} we used \cref{lem:holder-gram}. We have proven~\eqref{eq:projection-of-approx-variance}.

Applying \cref{lem:ueb-maximal-ineq} conditionally on $\crl{\vx_i}_{i = 1}^n$, on an event with probability approaching one, we get~\eqref{eq:a2diffbound}, and the proof of \cref{lem:holder-projection-approximation-error} is finished.
\end{proof}

\begin{remark}\label{rem:rates-in-proj-approx-error}
Inspecting the argument, one can note that instead of \eqref{eq:projection-approx-moment-case} the condition
\begin{equation*}
\frac{\log n}{n h^{d}} = o(1)
\end{equation*}
is enough for \cref{eq:a2condexpbound,eq:a2diffbound}.
\end{remark}

\section{Consistency}\label{sec:consistency}

We first study the convex case, and then move on to the non-convex case.

\subsection{Convex case}

The following lemma gives our most general result for a convex objective function. This is Lemma~\mainref{lem:consistency, convex} in the main paper.

\begin{lemma}[Consistency, convex case]\label{lem:holder-consistency}
  Suppose \cref{ass:dgp-holder,ass:gl-quasi-uniform-partition,ass:gl-local-basis,ass:gl-holder,ass:complexity,ass:gl-approximation-error} hold, $\riskfunc[\vq]{y}{\eta(\theta)}$ is convex with respect to $\theta$ with left or right derivative $\riskder[\vq]{y}{\eta(\theta)} \eta^{(1)}(\theta)$, $\mathcal{B} = \mathbb{R}^K$ in \eqref{eq:m-estimator}, and $m > d / 2$. Furthermore, assume that one of the following two conditions holds:
  \begin{align}
    &\frac{\brk[\big]{ \log (1 / h) }^{\frac{\nu}{\nu - 1}}}{n h^{\frac{2 \nu}{\nu - 1} d}} = o(1),\text{ or}\label{eq:consistency-moment-case}\\
    &\riskderenvel{\vx_i}{y_i} \text{ is \subgauss\ conditionally on $\vx_i$ and }\frac{\sqrt{\log n} \log(1 / h)}{n h^{2 d}} = o(1).\label{eq:consistency-subgaussian-case}
  \end{align}
  Then
  \begin{gather}
    \sup_{\vq \in \mathcal{Q}}\norm[\big]{ \betahat[\vq] - \betazero[\vq] } = o_{\Problet}(1),\label{eq:beta-q-close-to-beta0}\\
    \sup_{\vx \in \mathcal{X}} \sup_{\vq \in \mathcal{Q}} \abs[\big]{ \muhat[\vq][(\vv)]{\vx} - \muzero[\vq][(\vv)]{\vx} } = o_{\Problet}\prn{h^{- \multindnorm{\vv}}},\label{eq:mu-est-close-to-mu0}\\
    \sup_{\vq \in \mathcal{Q}} \prn[\Big]{\int_{\mathcal{X}} \prn[\big]{\muhat[\vq][(\vv)]{\vx} - \muzero[\vq][(\vv)]{\vx} }^2 \fx{\vx}\, \mathrm{d} \vx}^{1 / 2} = o_{\Problet} \prn[\big]{ h^{d / 2 - \multindnorm{\vv}} }.\label{eq:l2-rate}
  \end{gather}
\end{lemma}

\begin{proof}
First, note that \eqref{eq:mu-est-close-to-mu0} follows from \eqref{eq:beta-q-close-to-beta0} since uniformly over $\vx \in \mathcal{X}$ and $\vq \in \mathcal{Q}$
\begin{align*}
  &\abs[\big]{ \betahat[\vq]\trans \vp^{(\vv)}(\vx) - \muzero[\vq][(\vv)]{\vx} }\\
    &\quad \leq
    \abs[\big]{ \betahat[\vq]\trans \vp^{(\vv)}(\vx) - \betazerotrans[\vq] \vp^{(\vv)}(\vx) } + \abs[\big]{ \betazerotrans[\vq] \vp^{(\vv)}(\vx) - \muzero[\vq][(\vv)]{\vx} }\\
  &\quad \lesssim \norm[\big]{ \betahat[\vq] - \betazero[\vq] }_{\infty} \cdot \norm[\big]{ \vp^{(\vv)}(\vx) }_1 + \abs[\big]{ \betazerotrans[\vq] \vp^{(\vv)}(\vx) - \muzero[\vq][(\vv)]{\vx} }\\
  &\quad \lesssim \norm[\big]{ \betahat[\vq] - \betazero[\vq] }_{\infty} \cdot h^{- \multindnorm{\vv}} + \abs[\big]{ \betazerotrans[\vq] \vp^{(\vv)}(\vx) - \muzero[\vq][(\vv)]{\vx} } \tag*{by \cref{ass:gl-local-basis}}\\
  &\quad \lesssim \norm[\big]{ \betahat[\vq] - \betazero[\vq] }_{\infty} \cdot h^{- \multindnorm{\vv}} + h^{m - {\multindnorm{\vv}}} \tag*{by \cref{ass:gl-approximation-error},}
\end{align*}
where we used that only a bounded number of elements in $\vp^{(\vv)}(\vx)$ are nonzero. Similarly, \eqref{eq:l2-rate} follows from \eqref{eq:beta-q-close-to-beta0} since
\begin{align*}
  &\norm[\big]{ \betahat[\vq]\trans \vp^{(\vv)}(\vx) - \muzero[\vq][(\vv)]{\vx} }_{L_2(X)}\\
  &\quad \leq \norm[\big]{ \betahat[\vq]\trans \vp^{(\vv)}(\vx) - \betazerotrans[\vq] \vp^{(\vv)}(\vx) }_{L_2(X)} + \norm[\big]{ \betazerotrans[\vq] \vp^{(\vv)}(\vx) - \muzero[q][(\vv)]{\vx} }_{L_2(X)}\\
  &\quad \leq \norm[\big]{ \betahat[\vq]\trans \vp^{(\vv)}(\vx) - \betazerotrans[\vq] \vp^{(\vv)}(\vx) }_{L_2(X)} + \sup_{\vx \in \mathcal{X}} \abs[\big]{ \betazerotrans[\vq] \vp^{(\vv)}(\vx) - \muzero[\vq][(\vv)]{\vx} }\\
  &\quad \lesssim \norm[\big]{ \betahat[\vq]\trans \vp^{(\vv)}(\vx) - \betazerotrans[\vq] \vp^{(\vv)}(\vx) }_{L_2(X)} + h^{m - \multindnorm{\vv}} \tag*{by \cref{ass:gl-approximation-error}}\\
  &\quad = \prn[\Big]{\prn[\big]{\betahat[\vq] - \betazero[\vq]}\trans \Expect[\big]{\vp^{(\vv)}(\vx_i) \vp^{(\vv)}(\vx_i)\trans} \prn[\big]{\betahat[\vq] - \betazero[\vq]}}^{1 / 2} + h^{m - \multindnorm{\vv}}\\
  &\quad \leq \lambda_{\max}\prn[\big]{ \Expect[\big]{\vp^{(\vv)}(\vx_i) \vp^{(\vv)}(\vx_i)\trans} }^{1 / 2} \norm[\big]{\betahat[\vq] - \betazero[\vq]} + h^{m - \multindnorm{\vv}}\\
  &\quad \leq h^{d / 2 - \multindnorm{\vv}} \cdot \norm[\big]{ \betahat[\vq] - \betazero[\vq] } + h^{m - \multindnorm{\vv}} = o_{\Problet} \prn[\big]{ h^{d / 2 - \multindnorm{\vv}} } + h^{m - \multindnorm{\vv}} = o_{\Problet} \prn[\big]{ h^{d / 2 - \multindnorm{\vv}}  }.
\end{align*}
uniformly over $\vq \in \mathcal{Q}$, where by $\|g(\vx) \|_{L_2(X)}$ we denote $\prn[\big]{\int_{\mathcal{X}} g(\vx)^2 \fx{\vx}\,\mathrm{d}\vx}^{1 / 2}$ for simplicity. In the last equality we used $m > d / 2$ again. We also used that the largest eigenvalue of $\Expect*{\vp^{(\vv)}(\vx_i) \vp^{(\vv)}(\vx_i)\trans}$ is bounded from above by $h^{d - 2 \multindnorm{\vv}}$ up to a multiplicative coefficient, which is proven by the same argument as for $\vv = 0$ in \cref{lem:holder-gram} in combination with \cref{ass:gl-local-basis}.

It is left to prove~\eqref{eq:beta-q-close-to-beta0}. Fix a sufficiently small $\gamma > 0$. Denote for $i \in \{1, \ldots, n\}$ and $\alphavec \in \mathcal{S}^{K - 1}$
\begin{equation*}
  \delta_{\vq, i}(\alphavec) := \alphavec\trans \vp(\vx_i) \riskder[\vq]{y_i}{\eta\prn{\vp(\vx_i)\trans\prn{\betazero[\vq] + \gamma \alphavec } }} \eta^{(1)}\prn{\vp(\vx_i)\trans\prn{\betazero[\vq] + \gamma \alphavec } }.
\end{equation*}

Since $\Expectlet_n \brk[\big]{ \riskder[\vq]{y_i}{\eta\prn{\vp(\vx_i)\trans\vbeta }} \eta^{(1)}\prn{\vp(\vx_i)\trans\vbeta } \vp(\vx_i)}$ is a subgradient of the convex (by \cref{ass:gl-holder}\ref{ass:psi-holder-jumps}) objective function $\Expectlet_n \brk{ \riskfunc[\vq]{y_i}{\eta\prn{\vp(\vx_i)\trans\vbeta }} }$ of $\vbeta$, the strategy is to show that
\begin{equation}\label{eq:holder-consistency-strategy}
  \inf_{\vq, \alphavec} \Expectlet_n \brk{\delta_{\vq, i}(\alphavec)} > 0\quad\text{with probability approaching 1},
\end{equation}
which is enough to prove \cref{lem:holder-consistency} by convexity.

To implement this, we will show
\begin{gather}
  \inf_{\vq, \alphavec} \Expectlet_n \brk{\Expect{\delta_{\vq, i}(\alphavec) \given \vx_i}} \gtrsim \inf_{\alphavec} \Expectlet_n \brk{\alphavec\trans \vp(\vx_i) \vp(\vx_i)\trans \alphavec} + o_{\Problet} \prn{h^{d}}\quad\text{and}\label{eq:holder-consistency-delta-expectation-bound}\\
  \sup_{\vq, \alphavec} \abs[\big]{\Expectlet_n \brk{\delta_{\vq, i}(\alphavec) - \Expect{\delta_{\vq, i}(\alphavec) \given \vx_i}}} \lesssim_{\Problet}
  \begin{cases}
    \sqrt{\frac{\log(1 / h)}{n}} + \frac{\log(1 / h)}{n^{1 - 1 / \nu} h^d} = o \prn{ h^{d} },\\
    \sqrt{\frac{\log(1 / h)}{n}} + \frac{\sqrt{\log n} \log(1 / h)}{n h^d} = o \prn{ h^{d} }
  \end{cases}\label{eq:holder-consistency-delta-concentration}
\end{gather}
under~\eqref{eq:consistency-moment-case} and \eqref{eq:consistency-subgaussian-case} respectively (proof below), and conclude
\begin{align*}
  &\inf_{\vq, \alphavec} \Expectlet_n \brk{\delta_{\vq, i}(\alphavec)} \geq\\
  &\quad\geq \inf_{\vq, \alphavec} \Expectlet_n \brk[\big]{\Expect{\delta_{\vq, i}(\alphavec) \given \mX_n}} - \sup_{\vq, \alphavec} \abs[\big]{\Expectlet_n \brk[\big]{\delta_{\vq, i}(\alphavec) -\Expect{\delta_{\vq, i}(\alphavec) \given \mX_n}}}\\
  &\quad \gtrsim \inf_{\alphavec} \Expectlet_n \brk{\alphavec\trans \vp(\vx_i)\vp(\vx_i)\trans\alphavec} + o_{\Problet}\prn{h^{d}},
\end{align*}
which gives~\eqref{eq:holder-consistency-strategy} by \cref{lem:holder-gram}.

We will now prove~\eqref{eq:holder-consistency-delta-expectation-bound}. By \cref{ass:gl-holder}\ref{ass:psi-eta},
\begin{align*}
  &\Expect*{\delta_{\vq, i}(\alphavec) \given \vx_i} = \alphavec\trans \vp(\vx_i) \expriskder[\vq]{\vx_i}{\eta(\vp(\vx_i)\trans(\betazero[\vq] + \gamma \alphavec))}\\
  &\qquad \times \eta^{(1)}\prn{\vp(\vx_i)\trans\prn{\betazero[\vq] + \gamma \alphavec} }\\
  &\quad = \alphavec\trans \vp(\vx_i) \Psione[\vq]{\vx_i}{\xi_{\vq, i}} \prn{\vp(\vx_i)\trans\prn{\betazero[\vq] + \gamma \alphavec} - \muzero[\vq]{\vx_i}}\\
  &\qquad \times \eta^{(1)}\prn{\zeta_{\vq, i}} \eta^{(1)}\prn{\vp(\vx_i)\trans\prn{\betazero[\vq] + \gamma \alphavec} }\\
  &\quad \gtrsim \alphavec\trans \vp(\vx_i)\vp(\vx_i)\trans\alphavec - C \sup_{\vq \in \mathcal{Q}, \vx \in \mathcal{X}} \abs{\vp(\vx)\trans\betazero[\vq] - \muzero[\vq]{\vx} } \cdot \abs{ \alphavec\trans\vp(\vx_i) }
\end{align*}
almost surely uniformly over $\vq$, where $C$ is some positive constant not depending on $n$ or $i$, $\xi_{\vq, i}$ is between $\eta\prn{\vp(\vx_i)\trans\prn{\betazero[\vq] + \gamma \alphavec}}$ and $\eta\prn{\muzero[\vq]{\vx_i}}$, $\zeta_{\vq, i}$ between $\vp(\vx_i)\trans\prn{\betazero[\vq] + \gamma \alphavec}$ and $\muzero[\vq]{\vx_i}$. We used that $\expriskder[\vq]{\vx_i}{\eta(\muzero[\vq]{\vx_i})} = 0$, $\gamma$ is small enough, $\vp(\vx)\trans\betazero[\vq]$ is (for large enough $n$) uniformly close to $\muzero[\vq]{\vx}$ by \cref{ass:gl-approximation-error} and $\eta(\cdot)$ is strictly monotonic by Assumption~\eqref{ass:gl-holder}\ref{ass:psi-holder-jumps} giving the positivity of the product $\eta^{(1)}\prn{\zeta_{\vq, i}} \eta^{(1)}\prn{\vp(\vx_i)\trans\prn{\betazero[\vq] + \gamma \alphavec} }$.

Again using the uniform approximation bound $\sup_{\vq \in \mathcal{Q}, \vx \in \mathcal{X}} \abs{\vp(\vx)\trans\betazero[\vq] - \muzero[\vq]{\vx}} \lesssim h^{m}$ by \cref{ass:gl-approximation-error}, we obtain
\begin{equation}\label{eq:holder-consistency-main-second-term-estimate}
  \sup_{\vq \in \mathcal{Q}, \vx \in \mathcal{X}} \abs{\vp(\vx)\trans\betazero[\vq] - \muzero[\vq]{\vx} } \cdot \Expectlet_n \brk[\big]{\abs{ \alphavec\trans\vp(\vx_i) }} \lesssim_{\Problet} h^{m + d / 2}
\end{equation}
since $\Expectlet_n \brk[\big]{\abs{ \alphavec\trans\vp(\vx_i) }} \lesssim \Expectlet_n \brk[\big]{\prn{ \alphavec\trans\vp(\vx_i) }^2}^{1 / 2} \lesssim_{\Problet} h^{d / 2}$ by Lyapunov's inequality and \cref{lem:holder-gram}. Note that since $m > d / 2$, $h^{m + d / 2} = o(h^{d})$. \eqref{eq:holder-consistency-delta-expectation-bound} is proven.

We will now prove~\eqref{eq:holder-consistency-delta-concentration}. Define the function class
\begin{equation*}
  \mathcal{G}_1 := \crl{(\vx_i, y_i) \mapsto \delta_{\vq, i}(\alphavec)\colon \norm{\alphavec} = 1, \vq \in \mathcal{Q}}.
\end{equation*}
By \cref{ass:gl-holder}\ref{ass:psi-holder-jumps}, \ref{ass:gl-local-basis}\ref{ass:basis-norm-bound} and \cref{ass:gl-approximation-error}, for $\gamma$ small enough
\begin{equation*}
  \abs[\big]{ \riskder[\vq]{y_i}{\eta\prn{\vp(\vx_i)\trans\prn{\betazero[\vq] + \gamma \alphavec } }} - \riskder[\vq]{y_i}{\eta\prn{ \muzero[\vq]{\vx_i} }} } \lesssim 1 + \riskderenvel{\vx_i}{y_i}.
\end{equation*}
Recalling the envelope condition in \cref{ass:gl-holder}\ref{ass:psi-moments} and that $\muzero[\vq]{\cdot}$ is bounded by \cref{ass:dgp-holder}\ref{ass:mu-bounded}, we see that $\sup_{g \in \mathcal{G}_1} \abs{g} \lesssim 1 + \riskderenvel{\vx_i}{y_i}$, which means that under~\eqref{eq:consistency-moment-case}
\begin{align*}
  &\Expect[\Big]{\max_{1 \leq i \leq n} \abs{\delta_{\vq, i}(\alphavec)}^2 \given \mX_n}^{1 / 2} \leq \Expect[\Big]{\max_{1 \leq i \leq n} \abs{\delta_{\vq, i}(\alphavec)}^\nu \given \mX_n}^{1 / \nu}\\
  &\quad \leq \Expect[\bigg]{\sum_{i = 1}^{n} \abs{\delta_{\vq, i}(\alphavec)}^\nu \given \mX_n}^{1 / \nu} \lesssim \prn[\bigg]{ \sum_{i = 1}^n \Expect[\big]{\prn[\big]{1 + \riskderenvel{\vx_i}{y_i}}^{\nu} \given \mX_n} }^{1 / \nu}\\
  &\quad \lesssim \prn[\bigg]{ \sum_{i = 1}^n 1 }^{1 / \nu} = n^{1 / \nu}\quad\text{a.\,s.}
\end{align*}
and under~\eqref{eq:consistency-subgaussian-case} by \cref{lem:subgaussian-max-bound}
\begin{align*}
  \Expect[\Big]{\max_{1 \leq i \leq n} \abs{\delta_{\vq, i}(\alphavec)}^2 \given \mX_n}^{1 / 2} \lesssim \sqrt{\log n}\quad\text{a.\,s.}
\end{align*}
By similar considerations $\sup_{g \in \mathcal{G}_1} \Expectlet_n \brk[\big]{\Expect{g^2 \given \vx_i}} \lesssim \Expectlet_n \brk{\alphavec\trans \vp(\vx_i) \vp(\vx_i)\trans \alphavec} \lesssim h^{d}$ w.\,p.\,a.~1, where the last inequality holds by \cref{lem:holder-gram}.

By \cref{ass:complexity}, the class
\begin{equation*}
  \crl{(\vx_i, y_i) \mapsto \riskder[\vq]{y_i}{\eta\prn{\vp(\vx_i)\trans\vbeta }} \colon \norm{ \vbeta - \betazero[\vq] } \leq \gamma, \vq \in \mathcal{Q} }
\end{equation*}
with envelope $1 + \riskderenvel{\vx_i}{y_i}$ multiplied by a constant has a uniform entropy bound~\eqref{eq:uniform-entropy-bound} with $A \lesssim 1$, $V \lesssim K \asymp h^{-d}$. Moreover, the class
\begin{equation*}
  \crl{(\vx_i, y_i) \mapsto \boldsymbol{\alpha}\trans \vp(\vx_i) \colon \norm{ \boldsymbol{\alpha} } = 1 }
\end{equation*}
has a constant envelope and is VC with index no more than $K + 2$ by Lemma 2.6.15 in~\parencite{vaart96_weak_conver_empir_proces}, which means it satisfies the uniform entropy bound~\eqref{eq:uniform-entropy-bound} with $A \lesssim 1$, $V \lesssim K \asymp h^{-d}$. Similarly, the same is true of
\begin{equation*}
  \crl{(\vx_i, y_i) \mapsto \vp(\vx_i)\trans \vbeta \colon \norm{ \vbeta - \betazero[\vq] } \leq \gamma, \vq \in \mathcal{Q}}
\end{equation*}
and therefore of
\begin{equation*}
  \crl{(\vx_i, y_i) \mapsto \eta^{(1)}\prn{ \vp(\vx_i)\trans \vbeta } \colon \norm{ \vbeta - \betazero[\vq] } \leq \gamma, \vq \in \mathcal{Q}}
\end{equation*}
since $\eta^{(1)}(\cdot)$ on a bounded interval is Lipschitz. By \cref{lem:permanence-properties}, we conclude that $\mathcal{G}_1$ satisfies the uniform entropy bound~\eqref{eq:uniform-entropy-bound} with envelope $1 + \riskderenvel{\vx_i}{y_i}$ multiplied by a constant, $A \lesssim 1$ and $V \lesssim K \asymp h^{-d}$.

Applying the maximal inequality \cref{lem:ueb-maximal-ineq}, we obtain~\eqref{eq:holder-consistency-delta-concentration}.
\end{proof}

\subsection{Nonconvex case}

Our next goal is to prove the consistency result \cref{lem:consistency-nonconvex} for the nonconvex case. We will need the following lemma.

\begin{lemma}[Preparation for consistency in the nonconvex case]\label{lem:assumptions-of-nonconvex-consistency-are-reasonable}
  Suppose \cref{ass:dgp-holder,ass:gl-quasi-uniform-partition,ass:gl-local-basis,ass:gl-holder} hold. Then the infinity norm of $\betazero[\vq]$ is bounded:
  \begin{equation}\label{eq:betazeroq-is-bounded-in-inf-norm}
    \sup_{\vq \in \mathcal{Q}} \norm{ \betazero[\vq] }_{\infty} \lesssim 1.
  \end{equation}

  Moreover, for any $R > 0$, there is a positive constant $\nC{deltaenvelopeboundauxtwo} = \oC{deltaenvelopeboundauxtwo}(R)$ depending only on $R$ such that for any $\vx \in \mathcal{X}$
  \begin{equation}\label{eq:rho-beta-minus-rho-betazeroq-bounded}
    \sup_{\vq \in \mathcal{Q}} \sup_{\norm{ \bbeta }_{\infty} \leq R} \abs[\big]{ \riskfunc[\vq]{y}{\eta(\bp(\vx)\trans\bbeta)} - \riskfunc[\vq]{y}{\eta(\bp(\vx)\trans\betazero[\vq])} } \leq \oC{deltaenvelopeboundauxtwo} \prn{ 1 + \riskderenvel{\vx}{y} }.
  \end{equation}
\end{lemma}

\begin{proof}
We prove~\eqref{eq:betazeroq-is-bounded-in-inf-norm} first. By \cref{ass:gl-approximation-error} ($\betazerotrans[\vq] \bp(\vx)$ is close to $\muzero[\vq]{\vx}$) and \cref{ass:dgp-holder}\ref{ass:mu-bounded} ($\muzero[\vq]{\vx}$ is uniformly bounded), $\abs{ \betazerotrans[\vq] \bp(\vx) }$ is bounded uniformly over $\vq \in \mathcal{Q}$ and $\vx \in \mathcal{X}$. By \cref{ass:gl-local-basis}, we can bound the $k$th coordinate of $\betazero[\vq]$
\begin{align*}
  &\abs{ \prn{\betazero[\vq]}_k } \lesssim h^{-d / 2} \prn[\bigg]{\int_{\mathcal{H}_{k}} \prn[\big]{ \betazerotrans[\vq]\bp(\vx) }^2\,\mathrm{d}\vx}^{1 / 2}\\
  &\quad \leq h^{-d / 2} \sup_{\vx \in \mathcal{X}} \abs{ \betazerotrans[\vq]\bp(\vx) } \cdot \prn{\Leb{\mathcal{H}_k}}^{1 / 2} \lesssim \sup_{\vx \in \mathcal{X}} \abs{ \betazerotrans[\vq]\bp(\vx) } \lesssim 1,
\end{align*}
where the constants in $\lesssim$ do not depend on $k$.

Now we prove~\eqref{eq:rho-beta-minus-rho-betazeroq-bounded}. Note that
\begin{align*}
  &\riskfunc[\vq]{y}{\eta\prn{\bp(\vx) \trans \bbeta}} - \riskfunc[\vq]{y}{\eta\prn{\bp(\vx) \trans \betazero[\vq]}}\\
  &\quad = \int_0^{\bp(\vx) \trans (\bbeta - \betazero[\vq])} \prn{ \riskder[\vq]{y}{\eta \prn{ \bp(\vx)\trans \betazero[\vq] + t }} - \riskder[\vq]{y}{\eta \prn{ \muzero[\vq]{\vx} }} }\\
  &\quad\qquad \times \eta^{(1)} \prn{ \bp(\vx)\trans \betazero[\vq] + t }\,\mathrm{d} t\\
  &\qquad + \riskder[\vq]{y}{\eta \prn{ \muzero[\vq]{\vx} }} \int_0^{\bp(\vx) \trans (\bbeta - \betazero[\vq])} \eta^{(1)} \prn{ \bp(\vx)\trans \betazero[\vq] + t }\,\mathrm{d} t
\end{align*}

By \cref{ass:gl-holder}\ref{ass:psi-holder-jumps}, we have a bound
\begin{equation*}
  \abs{ \riskder[\vq]{y}{\eta \prn{ \bp(\vx)\trans \betazero[\vq] + t }} - \riskder[\vq]{y}{\eta \prn{ \muzero[\vq]{\vx} }} } \lesssim \riskderenvel{\vx}{y} + 1.
\end{equation*}
Since both $\bp(\vx)\trans \betazero[\vq]$ and $\bp(\vx) \trans \bbeta$ lie in a fixed compact interval (not depending on $\vx$ or $\vq$), $\eta^{(1)} \prn{ \bp(\vx)\trans \betazero[\vq] + t }$ is uniformly bounded in absolute value. This means that for any $\vx \in \mathcal{X}$, $\vq \in \mathcal{Q}$, $\norm{ \bbeta }_{\infty} \leq R$, we have for some positive constants $\nC{deltaenvelopeboundaux}$ and $\oC{deltaenvelopeboundauxtwo}$ depending only on $R$
\begin{align*}
  &\abs{ \riskfunc[\vq]{y}{\eta\prn{\bp(\vx) \trans \bbeta}} - \riskfunc[\vq]{y}{\eta\prn{\bp(\vx) \trans \betazero[\vq]}} } \leq \oC{deltaenvelopeboundaux} \prn{ 1 + \riskderenvel{\vx}{y} } \cdot \abs{ \bp(\vx) \trans\prn{\bbeta - \betazero[\vq]} }\\
  &\quad \leq \oC{deltaenvelopeboundaux} \prn{ 1 + \riskderenvel{\vx}{y} } \cdot \norm{ \bp(\vx) }_1 \cdot \norm{ \bbeta - \betazero[\vq] }_{\infty}\\
  &\quad \leq \oC{deltaenvelopeboundaux} \prn{ 1 + \riskderenvel{\vx}{y} } \cdot \norm{ \bp(\vx) }_1 \cdot \prn{ \norm{ \bbeta }_{\infty} + \norm{ \betazero[\vq] }_{\infty} }\\
  &\quad \leq \oC{deltaenvelopeboundauxtwo} \prn{ 1 + \riskderenvel{\vx}{y} },
\end{align*}
concluding the proof.
\end{proof}

We are now ready to prove a general consistency result for an estimator under constraints $\norm{ \bbeta }_{\infty} \leq R$ for some large enough constant $R$. This is Lemma~\mainref{lem:consistency, nonconvex} in the paper.

\begin{lemma}[Consistency, nonconvex case]\label{lem:consistency-nonconvex}
  Assume the following conditions.
  \begin{conditions}[wide]
  \item \cref{ass:dgp-holder,ass:gl-quasi-uniform-partition,ass:gl-local-basis,ass:gl-holder,ass:gl-approximation-error} hold.
  \item\label{cond:cons-nonconvex-m-greater-half-d} $m > d / 2$.
  \item\label{cond:cons-nonconvex-rate-condition} The following rate condition holds:
    \begin{align}
      &\frac{\brk[\big]{\log(1 / h)}^{\frac{\nu}{\nu - 1}}}{n h^{\frac{2 \nu}{\nu - 1} d}} = o(1),\text{ or}\label{eq:consistency-nonconvex-moment-case}\\
      &\riskderenvel{\vx_i}{y_i} \text{ is \subgauss\ conditionally on $\vx_i$ and }\frac{\sqrt{\log n} \log(1 / h)}{n h^{2 d}} = o(1).\label{eq:consistency-nonconvex-subgaussian-case}
    \end{align}
  \item $\mathcal{B} = \crl{\vbeta \in \mathbb{R}^K: \norm{\vbeta}_{\infty} \leq R}$ in \eqref{eq:m-estimator}, where $R > 0$ is a fixed number (not depending on $n$) such that $\sup_{\vq \in \mathcal{Q}} \norm{\betazero[\vq]}_{\infty} \leq R / 2$ (existing by \cref{lem:assumptions-of-nonconvex-consistency-are-reasonable}).
  \item\label{cond:lower-bound-of-psi-one} There is a positive constant $c$ such that we have $\inf \Psione[\vq]{\vx}{\zeta} > c$, where the infimum is over $\vx \in \mathcal{X}$, $\vq \in \mathcal{Q}$, $\zeta$ between $\eta \prn{ \bp(\vx)\trans\bbeta }$ and $\eta \prn{ \muzero[\vq]{\vx} }$, and $\bbeta \in \mathcal{B}$.
  \item\label{cond:nonconvex-consistency-additional-complexity-assumption} The class of functions
    \begin{equation*}
      \crl{ (\vx, y) \mapsto \riskfunc[\vq]{y}{\eta(\bp(\vx)\trans\bbeta)} - \riskfunc[\vq]{y}{\eta(\bp(\vx)\trans\betazero[\vq])}: \norm{ \bbeta }_{\infty} \leq R, \vq \in \mathcal{Q} }
    \end{equation*}
    with envelope $\oC{deltaenvelopeboundauxtwo} \prn{ 1 + \riskderenvel{\vx}{y} }$ (by \cref{lem:assumptions-of-nonconvex-consistency-are-reasonable}) satisfies the uniform entropy bound~\eqref{eq:uniform-entropy-bound} with $A \lesssim 1$ and $V \lesssim K$.
  \end{conditions}

  Then \cref{eq:beta-q-close-to-beta0,eq:mu-est-close-to-mu0,eq:l2-rate} hold.
\end{lemma}

\begin{proof}
For $\bbeta$ satisfying the constraint $\norm{\bbeta}_{\infty} \leq R$, define
\renewcommand{\d}[1]{\ensuremath{\operatorname{d}\!{#1}}}
\begin{align*}
  &\delta_{\vq, i}(\bbeta) := \riskfunc[\vq]{y_i}{\eta\prn{\bp(\vx_i) \trans \bbeta}} - \riskfunc[\vq]{y_i}{\eta\prn{\bp(\vx_i) \trans \betazero[\vq]}}\\
  &\quad = \int_0^{\bp(\vx_i) \trans (\bbeta - \betazero[\vq])} \riskder[\vq]{y_i}{\eta \prn{ \bp(\vx_i)\trans \betazero[\vq] + t }} \eta^{(1)} \prn{ \bp(\vx_i)\trans \betazero[\vq] + t }\,\mathrm{d} t.
\end{align*}

Note that
\begin{align*}
  &\Expect*{\delta_{\vq, i}(\bbeta) \given \vx_i} = \int_0^{\bp(\vx_i) \trans (\bbeta - \betazero[\vq])} \expriskder[\vq]{\vx_i}{\eta( \bp(\vx_i)\trans \betazero[\vq] + t)} \eta^{(1)} \prn{ \bp(\vx_i)\trans \betazero[\vq] + t }\,\mathrm{d} t\\
  &\quad \overset{\text{(a)}}{=} \int_0^{\bp(\vx_i) \trans (\bbeta - \betazero[\vq])} \Psione[\vq]{\vx_i}{\xi_{\vq, i, t}}\\
  &\qquad \times \eta^{(1)} \prn{ \zeta_{\vq, i, t} } \eta^{(1)} \prn{ \bp(\vx_i)\trans \betazero[\vq] + t } \crl{ \bp(\vx_i)\trans \betazero[\vq] - \muzero[\vq]{\vx_i} + t } \,\mathrm{d} t\\
  &\quad \overset{\text{(b)}}{\geq} \oC{deltacondauxone} \crl{ \bp(\vx_i) \trans (\bbeta - \betazero[\vq]) }^2 - \nC{} \sup_{\vq, \vx} \abs{ \bp(\vx)\trans \betazero[\vq] - \muzero[\vq]{\vx} } \cdot \abs{ \bp(\vx_i) \trans (\bbeta - \betazero[\vq]) }\\
  &\quad \overset{\text{(c)}}{\geq} \oC{deltacondauxone} \crl{ \bp(\vx_i) \trans (\bbeta - \betazero[\vq]) }^2 - \oC{deltacondauxtwo} h^m \abs{ \bp(\vx_i) \trans (\bbeta - \betazero[\vq]) },
\end{align*}
with some positive constants $\nC{deltacondauxone}$ and $\nC{deltacondauxtwo}$ (depending on $R$), where in (a) we used
\begin{align*}
  &\expriskder[\vq]{\vx_i}{\eta(\bp(\vx_i)\trans \betazero[\vq] + t)} = \expriskder[\vq]{\vx_i}{\eta( \bp(\vx_i)\trans \betazero[\vq] + t)} - \expriskder[\vq]{\vx_i}{\eta( \muzero[\vq]{\vx_i})}\\
  &\quad = \Psione[\vq]{\vx_i}{\xi_{\vq, i, t}} \crl{ \eta \prn{ \bp(\vx_i)\trans \betazero[\vq] + t } - \eta \prn{ \muzero[\vq]{\vx_i} } }\\
  &\quad = \Psione[\vq]{\vx_i}{\xi_{\vq, i, t}} \eta^{(1)} \prn{ \zeta_{\vq, i, t} } \crl{ \bp(\vx_i)\trans \betazero[\vq] - \muzero[\vq]{\vx_i} + t }
\end{align*}
for some $\xi_{\vq, i, t}$ between $\eta \prn{ \bp(\vx_i)\trans \betazero[\vq] + t }$ and $\eta \prn{ \muzero[\vq]{\vx_i} }$, and some $\zeta_{\vq, i, t}$ between $\bp(\vx_i)\trans \betazero[\vq] + t$ and $\muzero[\vq]{\vx_i}$ by the mean-value theorem applied twice; in (b) we used Condition~\ref{cond:lower-bound-of-psi-one}, \cref{ass:gl-holder}, in particular that $\eta(\cdot)$ is strictly monotonic giving the positivity of the product $\eta^{(1)} \prn{ \zeta_{\vq, i, t} } \eta^{(1)} \prn{ \bp(\vx_i)\trans \betazero[\vq] + t }$; in (c) we used \cref{ass:gl-approximation-error}. By Lyapunov's inequality, $\Expectlet_n \brk{\abs{ \bp(\vx_i) \trans (\bbeta - \betazero[\vq]) }} \leq \Expectlet_n \brk{\prn{ \bp(\vx_i) \trans (\bbeta - \betazero[\vq]) }^2}^{1 / 2}$. We conclude
\begin{align*}
  &\Expectlet_n \brk[\big]{\Expect{\delta_{\vq, i}(\bbeta) \given \vx_i}}\\
  &\quad \geq \oC{deltacondauxone} \prn{\bbeta - \betazero[\vq]}\trans \Expectlet_n \brk{\bp(\vx_i) \bp(\vx_i)\trans} \prn{\bbeta - \betazero[\vq]} - \oC{deltacondauxtwo} h^m \Expectlet_n \brk{\prn{ \bp(\vx_i) \trans \prn{\bbeta - \betazero[\vq]} }^2}^{1 / 2}\\
  &\quad \labrel{1}[\geq] \oC{deltacondone} h^d \norm{ \bbeta - \betazero[\vq] }^2 - \oC{deltacondtwo} h^{m + d / 2} \norm{ \bbeta - \betazero[\vq] }
\end{align*}
with probability approaching one for some other positive constants $\nC{deltacondone}$ and $\nC{deltacondtwo}$ (depending on $R$), where \labrel{1} is by \cref{lem:holder-gram}.

Fix $\varepsilon > 0$ smaller than $R / 2$. In this case
\begin{equation*}
  \crl{ \bbeta: \norm{ \bbeta - \betazero[\vq] } \leq \varepsilon } \subset \crl{ \bbeta: \norm{ \bbeta - \betazero[\vq] }_{\infty} \leq \varepsilon } \subset \crl{ \bbeta: \norm{ \bbeta }_{\infty} \leq R}
\end{equation*}
because
\begin{equation*}
  \norm{ \bbeta }_{\infty} \leq \norm{ \bbeta - \betazero[\vq] }_{\infty} + \norm{ \betazero[\vq] }_{\infty} \leq \norm{ \bbeta - \betazero[\vq] }_{\infty} + R / 2.
\end{equation*}

Define the class of functions
\begin{align*}
  &\mathcal{G} := \crl*{(\vx, y) \mapsto \norm{ \bbeta - \betazero[\vq] }^{-1} \prn{\riskfunc[\vq]{y}{\eta\prn{\bp(\vx) \trans \bbeta}} - \riskfunc[\vq]{y}{\eta\prn{\bp(\vx) \trans \betazero[\vq]}}}:\\
  &\qquad \norm{ \bbeta }_{\infty} \leq R, \norm{ \bbeta - \betazero[\vq] } > \varepsilon, \vq \in \mathcal{Q}}.
\end{align*}

It is a product of a subclass of the class
\begin{equation*}
  \crl{ (\vx, y) \mapsto a: 0 < a < 1 / \varepsilon }
\end{equation*}
with envelope $1 / \varepsilon$, obviously satisfying the uniform entropy bound~\eqref{eq:uniform-entropy-bound} with $A \lesssim 1$, $V \lesssim 1$, and a subclass of the class
\begin{equation*}
  \crl{ (\vx, y) \mapsto \riskfunc{y}{\eta(\bp(\vx)\trans\bbeta)} - \riskfunc{y}{\eta(\bp(\vx)\trans\betazero[\vq])}: \norm{ \bbeta }_{\infty} \leq R, \vq \in \mathcal{Q} }
\end{equation*}
with envelope $\oC{deltaenvelopeboundauxtwo} \prn{ 1 + \riskderenvel{\vx}{y} }$, satisfying the uniform entropy bound~\eqref{eq:uniform-entropy-bound} with $A \lesssim 1$ and $V \lesssim K$ by the conditions of \cref{lem:consistency-nonconvex}.

By \cref{lem:permanence-properties}, class $\mathcal{G}$ with envelope $\oC{deltaenvelopeboundauxtwo} / \varepsilon \prn{ 1 + \riskderenvel{\vx}{y} }$ satisfies the uniform entropy bound~\eqref{eq:uniform-entropy-bound} with $A \lesssim 1$ and $V \lesssim K$.

Next, under~\eqref{eq:consistency-nonconvex-moment-case}
\begin{align*}
  &\Expect[\big]{\max_{1 \leq i \leq n} \prn{ \oC{deltaenvelopeboundauxtwo} / \varepsilon }^2 \prn{ 1 + \riskderenvel{\vx}{y} }^2 \given \mX_n}^{1 / 2}\\
  &\quad \leq \Expect[\big]{\max_{1 \leq i \leq n} \prn[\big]{ \oC{deltaenvelopeboundauxtwo} / \varepsilon }^\nu \prn[\big]{ 1 + \riskderenvel{\vx}{y} }^\nu \given \mX_n}^{1 / \nu}\\
  &\quad \leq \Expect[\big]{\sum_{i = 1}^{n} \prn[\big]{ \oC{deltaenvelopeboundauxtwo} / \varepsilon }^\nu \prn[\big]{ 1 + \riskderenvel{\vx}{y} }^\nu \given \mX_n}^{1 / \nu}\\
  &\quad \lesssim \frac{1}{\varepsilon} \prn[\bigg]{ \sum_{i = 1}^n \Expect[\big]{\prn{1 + \riskderenvel{\vx_i}{y_i}}^{\nu} \given \mX_n} }^{1 / \nu}\\
  &\quad \lesssim \frac{1}{\varepsilon} \prn[\bigg]{ \sum_{i = 1}^n 1 }^{1 / \nu} = \frac{n^{1 / \nu}}{\varepsilon}\quad\text{a.\,s.}
\end{align*}
with constants in $\lesssim$ depending on $R$ but not on $n$ or $\varepsilon$, and under~\eqref{eq:consistency-nonconvex-subgaussian-case} by \cref{lem:subgaussian-max-bound}
\begin{align*}
  \Expect*{\max_{1 \leq i \leq n} \prn{ \oC{deltaenvelopeboundauxtwo} / \varepsilon }^2 \prn{ 1 + \riskderenvel{\vx}{y} }^2 \given \mX_n}^{1 / 2} \lesssim \frac{\sqrt{\log n}}{\varepsilon}\quad\text{a.\,s.}
\end{align*}

Moreover,
\begin{align*}
  &\Expectlet_n \brk[\big]{\Expect{g(\vx_i, y_i)^2 \given \vx_i}}\\
  &\quad \leq \norm{ \bbeta - \betazero[\vq] }^{-2} \oC{deltaenvelopeboundaux}^2 \Expectlet_n \brk[\big]{\prn[\big]{ \bp(\vx_i) \trans\prn{\bbeta - \betazero[\vq]} }^2 \Expect{\prn{ 1 + \riskderenvel{\vx_i}{y_i} }^2 \given \vx_i}}\\
  &\quad \lesssim \norm{ \bbeta - \betazero[\vq] }^{-2} \prn{\bbeta - \betazero[\vq]}\trans \Expectlet_n \brk{\bp(\vx_i) \bp(\vx_i)\trans} \prn{\bbeta - \betazero[\vq]}\\
  &\quad \labrel{1}[\lesssim] h^d,
\end{align*}
where \labrel{1} is by \cref{lem:holder-gram}.

By \cref{lem:ueb-maximal-ineq}, we have
\begin{equation*}
  \sup_{g \in \mathcal{G}} \abs[\big]{ \Expectlet\brk[\big]{g(\vx_i, y_i) - \Expect{g(\vx_i, y_i) \given \vx_i}} } \lesssim_{\Problet}
  \begin{cases}
    \sqrt{\frac{\log(1 / h)}{n}} + \frac{\log(1 / h)}{n^{1 - 1 / \nu} h^d} = o (h^d),\\
    \sqrt{\frac{\log(1 / h)}{n}} + \frac{\sqrt{\log n} \log(1 / h)}{n h^d} = o (h^d)
  \end{cases}
\end{equation*}
under~\eqref{eq:consistency-nonconvex-moment-case} and~\eqref{eq:consistency-nonconvex-subgaussian-case} respectively (since $\varepsilon$ is fixed).

Combining, we infer from the previous results that with probability approaching one for all $\norm{ \bbeta }_{\infty} \leq R$, $\norm{ \bbeta - \betazero[\vq] } > \varepsilon$, $\vq \in \mathcal{Q}$
\begin{align*}
  &\Expectlet_n \brk{\delta_{\vq, i}(\bbeta)} \geq \oC{deltacondone} h^d \norm{ \bbeta - \betazero[\vq] }^2 - \oC{deltacondtwo} h^{m + d / 2} \norm{ \bbeta - \betazero[\vq] } + \norm{ \bbeta - \betazero[\vq] } \cdot o (h^d)\\
  &\quad = \norm{ \bbeta - \betazero[\vq] } \cdot \crl{ \oC{deltacondone} h^d \norm{ \bbeta - \betazero[\vq] } - \oC{deltacondtwo} h^{m + d / 2} + o (h^d)}\\
  &\quad > \norm{ \bbeta - \betazero[\vq] } \cdot \crl{ \oC{deltacondone} \varepsilon h^d - \oC{deltacondtwo} h^{m + d / 2} + o (h^d) }\\
  &\quad \labrel{1}[=] \norm{ \bbeta - \betazero[\vq] } \cdot h^d \crl{ \oC{deltacondone} \varepsilon + o \prn{ 1 } } > 0,
\end{align*}
where in \labrel{1} we used $m > d / 2$.

It follows that the constrained minimizer under the constraint $\norm{ \bbeta }_{\infty} \leq R$ has to lie inside the ball $\norm{ \bbeta - \betazero[\vq] } \leq \varepsilon$ for all $\vq \in \mathcal{Q}$ with probability approaching one. Since $\varepsilon$ was arbitrary smaller than $R / 2$, it is equivalent to $\sup_{\vq \in \mathcal{Q}} \norm{ \betahatqconstr - \betazero[\vq] } = o_{\Problet}(1)$. \Cref{eq:mu-est-close-to-mu0,eq:l2-rate} follow (as in \cref{lem:holder-consistency}).
\end{proof}

\subsection{Weaker conditions for special cases}

\Cref{lem:holder-consistency-unconnected,lem:nonconvex-cons-unconnected} consider the special case of unconnected basis functions.

\begin{lemma}[Consistency, convex case, unconnected basis functions]\label{lem:holder-consistency-unconnected}
  Assume the following.

  \begin{conditions}[wide]
  \item \cref{ass:dgp-holder,ass:gl-quasi-uniform-partition,ass:gl-local-basis,ass:gl-holder,ass:complexity,ass:gl-approximation-error} hold.

  \item $\riskfunc[\vq]{y}{\eta(\theta)}$ is convex with respect to $\theta$, and $\riskder[\vq]{y}{\eta(\theta)} \eta^{(1)}(\theta)$ is its left or right derivative, and $\mathcal{B} = \mathbb{R}^K$ in \eqref{eq:m-estimator}.

  \item For all $k \in \crl{ 1, \ldots, K }$ the $k$th basis function $p_k(\cdot)$ is only active on one of the cells of $\Delta$.

  \item The rate of convergence of $h$ to zero is restricted by
    \begin{align*}
      &\frac{\brk[\big]{ \log (1 / h) }^{\frac{\nu}{\nu - 1}}}{n h^{\frac{\nu}{\nu - 1} d}} = o(1),\text{ or}\\
      &\riskderenvel{\vx_i}{y_i} \text{ is \subgauss\ conditionally on $\vx_i$ and }\frac{\sqrt{\log n} \log(1 / h)}{n h^{d}} = o(1).
    \end{align*}
  \end{conditions}

  Then \cref{eq:beta-q-close-to-beta0,eq:mu-est-close-to-mu0,eq:l2-rate} hold.
\end{lemma}

\begin{proof}
As in \cref{lem:holder-consistency}, \cref{eq:mu-est-close-to-mu0,eq:l2-rate} follow from \cref{eq:beta-q-close-to-beta0}.

For $l \in \crl{ 1, \ldots, \kappa }$, the number $M_l$ of basis functions in $\vp(\cdot)$ which are active on the $l$th cell of $\Delta$ is bounded by a constant. Denote the vector of such basis functions $\vp_l := \prn{ p_{l, 1}, \ldots, p_{l, M_l} }\trans$. Define the matrices $\mQ_{0, \vq, l}$ and $\bar{\mQ}_{\vq, l}$ as before with $\vp$ replaced by $\vp_l$ (for different $l$, the dimensions of these square matrices may vary but are bounded from above). By a simple modification of the argument in \cref{lem:holder-gram}, the analogues of~\eqref{eq:q0q-eigenvalues-bound} and \eqref{eq:qbarq-eigenvalues-bound} continue to hold: uniformly over $\vq \in \mathcal{Q}$ and $l$
\begin{equation}\label{eq:gram-matrices-restr-to-one-cell-eigen-bounds}
  \begin{aligned}
    &h^{d} \lesssim \lambda_{\min}(\mQ_{0, \vq, l}) \leq \lambda_{\max}(\mQ_{0, \vq, l}) \lesssim h^{d},\\
    &h^{d} \lesssim \lambda_{\min}(\bar{\mQ}_{\vq, l}) \leq \lambda_{\max}(\bar{\mQ}_{\vq, l}) \lesssim h^{d} \quad\text{w.\,p.\,a. 1}.
  \end{aligned}
\end{equation}

By the assumption of the lemma, we can write $\betazero[\vq] = \prn{ \vbeta_{0, \vq, 1}, \ldots, \vbeta_{0, \vq, \kappa} }\trans$, where $\vbeta_{0, \vq, l}$ is a subvector of dimension $M_l$ corresponding to the elements in $\vp$ active on the $l$th cell.

Fix a sufficiently small $\gamma > 0$. Denote for $l \in \crl{ 1, \ldots, \kappa }$, $i \in \crl{ 1, \ldots, n }$ and $\alphavec_l \in \mathcal{S}^{M_l - 1}$
\begin{equation*}
  \delta_{\vq, i, l}\prn{\alphavec_l} := \alphavec_l\trans \vp_l(\vx_i) \riskder[\vq]{y_i}{\eta\prn{\vp_l(\vx_i)\trans\prn{\vbeta_{0, \vq, l} + \gamma \alphavec_l } }} \eta^{(1)}\prn{\vp_l(\vx_i)\trans\prn{\vbeta_{0, \vq, l} + \gamma \alphavec_l } }.
\end{equation*}
Proceeding in the same way as in \cref{lem:holder-consistency}, we will show
\begin{equation}\label{eq:holder-consistency-strategy-unconnected}
  \inf_{\vq, \alphavec_l, l} \Expectlet_n \brk{\delta_{\vq, i, l}(\alphavec_l)} > 0\quad\text{with probability approaching 1},
\end{equation}
which is again enough to prove the lemma by convexity. It will follow that with probability approaching one the minimizer $\widehat{\bbeta}_{\vq, l}$ of $\Expectlet_n \brk[\big]{\riskfunc[\vq]{y_i}{\eta(\bp_l(\vx_i) \trans \bbeta_l)}}$ with respect to $\bbeta_l$ has to lie inside the ball $\norm{\bbeta_l - \bbeta_{0, \vq, l}} \leq \gamma$, and in particular inside the cube $\norm{\bbeta_l - \bbeta_{0, \vq, l}}_{\infty} \leq \gamma$. But note that $\betahat[\vq] = \prn[\big]{\widehat{\bbeta}_{\vq, 1}\trans, \ldots, \widehat{\bbeta}_{\vq, \bar{k}}\trans}\trans$. So, with probability approaching one, for all $\vq \in \mathcal{Q}$, we have $\norm[\big]{\betahat[\vq] - \betazero[\vq]}_{\infty} \leq \gamma$. Since $\gamma$ was arbitrary (small enough), it is equivalent to $\sup_{\vq \in \mathcal{Q}} \norm[\big]{\betahat[\vq] - \betazero[\vq]}_{\infty} = o_{\Problet}(1)$.

\Cref{eq:holder-consistency-strategy-unconnected} is proven analogously to the corresponding argument in \cref{lem:holder-consistency}. The class of functions
\begin{equation*}
  \mathcal{G}_{1} := \crl[\big]{(\vx_i, y_i) \mapsto \delta_{\vq, i, l}\prn{\alphavec_l}\colon \alphavec_l \in \mathcal{S}^{M_l - 1}, l \in \crl{ 1, \ldots, \kappa }, \vq \in \mathcal{Q}}
\end{equation*}
now satisfies the uniform entropy bound~\eqref{eq:uniform-entropy-bound} with $A \lesssim \kappa \asymp h^{-d}$ (since there are $\kappa$ different values of $l$) and $V \lesssim 1$ (since the vectors $\alphavec_l$ are of bounded dimensions). The bound $\sup_{l, \alphavec_l} \Expectlet_n \brk[\big]{|\alphavec_l\trans \bp_l(\vx_i)|} \lesssim h^d$ with probability approaching one can be proven without assuming $m > d / 2$ by using
\begin{equation*}
\sup_{l, \alphavec_l} \Expectlet_n \brk[\big]{|\alphavec_l\trans \bp_l(\vx_i)|} \leq \sup_{l, \alphavec_l} \norm{\alphavec_l}_{\infty} \Expectlet_n \brk[\big]{\norm{\bp_l(\vx_i)}_1} \lesssim h^d\quad\text{w.\,p.\,a.~1,}
\end{equation*}
since the dimension of $\bp_l(\cdot)$ is uniformly bounded.
\end{proof}

\begin{lemma}[Consistency in the nonconvex case: unconnected basis]\label{lem:nonconvex-cons-unconnected}
  Assume that for all $k \in \crl{ 1, \ldots, K }$ the $k$th basis function $p_k(\cdot)$ is only active on one of the cells of $\Delta$, and define $\bp_l(\cdot)$, $M_l$, $\bbeta_{0, \vq, l}$ as in the proof of \cref{lem:holder-consistency-unconnected}. Further, assume the conditions of \cref{lem:consistency-nonconvex} with \cref{cond:cons-nonconvex-m-greater-half-d} removed, \cref{cond:cons-nonconvex-rate-condition} replaced by
  \begin{align}
    &\frac{\brk[\big]{\log(1 / h)}^{\frac{\nu}{\nu - 1}}}{n h^{\frac{\nu}{\nu - 1} d}} = o(1),\text{ or}\\
    &\riskderenvel{\vx_i}{y_i} \text{ is \subgauss\ conditionally on $\vx_i$ and }\frac{\sqrt{\log n} \log(1 / h)}{n h^{d}} = o(1),
  \end{align}
  and \cref{cond:nonconvex-consistency-additional-complexity-assumption} replaced by the following: the class
  \begin{equation*}
    \crl{ (\vx, y) \mapsto \riskfunc{y}{\eta(\bp_l(\vx)\trans\bbeta_l)} - \riskfunc{y}{\eta(\bp_l(\vx)\trans\bbeta_{0, \vq, l})}: \norm{ \bbeta_l }_{\infty} \leq R, \vq \in \mathcal{Q}, l \in \crl{ 1, \ldots, \kappa } }
  \end{equation*}
  satisfies the uniform entropy bound \eqref{eq:uniform-entropy-bound} with $A \lesssim \kappa \asymp h^{-d}$ and $V \lesssim 1$. Then \cref{eq:beta-q-close-to-beta0,eq:mu-est-close-to-mu0,eq:l2-rate} hold.
\end{lemma}

\begin{proof}
As in \cref{lem:holder-consistency}, \cref{eq:mu-est-close-to-mu0,eq:l2-rate} follow from \cref{eq:beta-q-close-to-beta0}.

Define matrices $\mQ_{0, \vq, l}$ and $\bar{\mQ}_{\vq, l}$ as in the proof of \cref{lem:holder-consistency-unconnected}, and recall that the asymptotic bounds on their eigenvalues are the same as in the general (not restricted to one cell) case, i.\,e. \eqref{eq:gram-matrices-restr-to-one-cell-eigen-bounds} holds.

For any $M_l$-dimensional vector $\bbeta_l$ satisfying the constraint $\norm{\bbeta_l}_{\infty} \leq R$, define
\begin{align*}
  &\delta_{\vq, i, l}(\bbeta_l) := \riskfunc[\vq]{y_i}{\eta(\bp_l(\vx_i) \trans \bbeta_{l})} - \riskfunc[\vq]{y_i}{\eta(\bp_l(\vx_i) \trans \bbeta_{0, \vq, l})}\\
  &\quad = \int_0^{\bp_l(\vx_i) \trans (\bbeta_l - \bbeta_{0, \vq, l})} \riskder[\vq]{y_i}{\eta( \bp_l(\vx_i)\trans \bbeta_{0, \vq, l} + t)} \eta^{(1)}(\bp_l(\vx_i)\trans \bbeta_{0, \vq, l} + t)\,\mathrm{d} t.
\end{align*}
By the same argument as in the proof of \cref{lem:consistency-nonconvex}, we have
\begin{align*}
  &\Expectlet_n \brk[\big]{\Expectlet\brk{\delta_{\vq, i, l}(\bbeta_l) | \vx_i}}\\
  &\quad \geq \oC{celldeltacondone} (\bbeta_l - \bbeta_{0, \vq, l})\trans \Expectlet_n \brk{\bp_l(\vx_i) \bp_l(\vx_i)\trans} (\bbeta_l - \bbeta_{0, \vq, l}) - \oC{celldeltacondauxtwo} h^m \Expectlet_n \brk[\big]{\bigl|\bp_l(\vx_i) \trans (\bbeta_l - \bbeta_{0, \vq, l})\bigr|}\\
  &\quad \labrel{1}[\geq] \oC{celldeltacondone} h^d \norm{\bbeta_l - \bbeta_{0, \vq, l}}^2 - \oC{celldeltacondtwo} h^{m + d} \norm{\bbeta_l - \bbeta_{0, \vq, l}}
\end{align*}
with probability approaching one for some positive constants $\nC{celldeltacondone}$, $\nC{celldeltacondauxtwo}$ and $\nC{celldeltacondtwo}$ (depending on $R$, but not on $\vq$ or $l$), where in \labrel{1} we used $\Expectlet_n \brk[\big]{\big|\bp_l(\vx_i) \trans (\bbeta_l - \bbeta_{0, \vq, l})\big|} \leq \norm{\bbeta_l - \bbeta_{0, \vq, l}}_{\infty} \Expectlet_n \brk[\big]{\norm{\bp_l(\vx_i)}_1} \lesssim \norm{\bbeta_l - \bbeta_{0, \vq, l}}_{\infty} h^d \leq \norm{\bbeta_l - \bbeta_{0, \vq, l}}_{\infty} h^d$ with probability approaching one since the dimension of $\bp_l(\cdot)$ is bounded.

Next, proceeding with the same concentration argument as in \cref{lem:consistency-nonconvex}, we will obtain  that with probability approaching one for all $l \in \crl{ 1, \ldots, \kappa }$, $\vq \in \mathcal{Q}$, $\norm{ \bbeta_l }_{\infty} \leq R$, $\| \bbeta_l - \bbeta_{0, \vq, l} \| > \varepsilon$,
\begin{align*}
  &\Expectlet_n \brk{\delta_{\vq, i, l}(\bbeta_l)} \geq \oC{celldeltacondone} h^d \| \bbeta_l - \bbeta_{0, \vq, l} \|^2 - \oC{celldeltacondtwo} h^{m + d} \| \bbeta_l - \bbeta_{0, \vq, l} \| + \| \bbeta_l - \bbeta_{0, \vq, l} \| \cdot o(h^d)\\
  &\quad = \| \bbeta_l - \bbeta_{0, \vq, l} \| \cdot \bigl\{ \oC{celldeltacondone} h^d \| \bbeta - \betazero[\vq] \| - \oC{celldeltacondtwo} h^{m + d} + o( h^d)\bigr\}\\
  &\quad > \| \bbeta_l - \bbeta_{0, \vq, l} \| \cdot \bigl\{ \oC{celldeltacondone} \varepsilon h^d - \oC{celldeltacondtwo} h^{m + d} + o( h^d)\bigr\}\\
  &\quad \labrel{1}[=]  \| \bbeta_l - \bbeta_{0, \vq, l} \| \cdot h^d \crl{ \oC{celldeltacondone} \varepsilon + o \prn{ 1 } } > 0.
\end{align*}

It follows that the constrained minimizer $\widehat{\bbeta}_{\vq, \mathrm{constr}, l}$ of $\Expectlet_n \brk[\big]{\riskfunc[\vq]{y_i}{\eta(\bp_l(\vx_i) \trans \bbeta_l)}}$ with respect to $\bbeta_l$ under the constraint $\norm{\bbeta_l}_{\infty} \leq R$ has to lie inside the ball $\norm{\bbeta_l - \bbeta_{0, \vq, l}} \leq \varepsilon$, and in particular inside the cube $\norm{\bbeta_l - \bbeta_{0, \vq, l}}_{\infty} \leq \varepsilon$. But this optimization can be solved separately for all $l$, i.\,e. $\betahatqconstr = \prn[\big]{\widehat{\bbeta}_{\vq, \mathrm{constr}, 1}, \ldots, \widehat{\bbeta}_{\vq, \mathrm{constr}, \bar{k}}}\trans$. So, with probability approaching one, for all $\vq \in \mathcal{Q}$, we have $\norm[\big]{\betahatqconstr - \betazero[\vq]}_{\infty} \leq \varepsilon$. Since $\varepsilon$ was arbitrary smaller than $R / 2$, it is equivalent to $\sup_{\vq \in \mathcal{Q}} \norm[\big]{\betahatqconstr - \betazero[\vq]}_{\infty} = o_{\Problet}(1)$.
\end{proof}

The following lemma considers the special case of strongly convex and strongly smooth loss function.

\begin{lemma}[Consistency: strongly convex and strongly smooth loss case]\label{lem:holder-consistency-strong-shape} Assume the following conditions.

  \begin{conditions}[wide]
  \item \cref{ass:dgp-holder,ass:gl-quasi-uniform-partition,ass:gl-local-basis,ass:gl-holder,ass:complexity,ass:gl-approximation-error} hold.

  \item The rate of convergence of $h$ to zero is restricted by
    \begin{align*}
      &\frac{\brk[\big]{\log(1 / h)}^{\frac{\nu}{\nu - 1}}}{n h^{\frac{\nu}{\nu - 1} d}} = o(1),\text{ or}\\
      &\riskderenvel{\vx_i}{y_i} \text{ is \subgauss\ conditionally on $\vx_i$ and }\frac{\log n \log(1 / h)}{n h^d} = o(1).
    \end{align*}

  \item The function $\eta \mapsto \riskder[\vq]{y}{\eta}$ is continuously differentiable on $\mathbb{R}$ (for all $y, \vq$), and there exist fixed (not depending on $n$, $\vq$ or $\theta$) numbers $\lambda, \Lambda$ such that
    \begin{equation*}
      0 < \lambda \leq \frac{\partial}{\partial \theta} \prn[\big]{ \riskder[\vq]{y_i}{\eta(\theta)} \eta^{(1)}(\theta) } \leq \Lambda,
    \end{equation*}
    and $\mathcal{B} = \mathbb{R}^K$.
  \end{conditions}

  Then \cref{eq:beta-q-close-to-beta0,eq:mu-est-close-to-mu0,eq:l2-rate} hold.
\end{lemma}

\begin{proof}
As in \cref{lem:holder-consistency}, \cref{eq:mu-est-close-to-mu0,eq:l2-rate} follow from \cref{eq:beta-q-close-to-beta0}.

Denote for $\vbeta \in \mathbb{R}^K$
\begin{equation*}
  G_n(\vbeta) := \Expectlet_n \brk{\riskder[\vq]{y_i}{\eta\prn{\vp(\vx_i)\trans\vbeta }} \eta^{(1)}\prn{\vp(\vx_i)\trans\vbeta } \vp(\vx_i)},
\end{equation*}
which is the gradient of the convex (by \cref{ass:gl-holder}\ref{ass:psi-holder-jumps}) function $\Expectlet_n \brk{ \riskfunc[\vq]{y_i}{\eta\prn{\vp(\vx_i)\trans\vbeta }}}$ of $\vbeta$. By definition of $\betahat[\vq]$ and differentiability, $G_n\prn[\big]{\betahat[\vq]} = 0$. By the mean value theorem,
\begin{equation}\label{eq:strong-shape-gradient-diff}
  G_n \prn{ \betazero[\vq] } = G_n \prn{ \betazero[\vq] } - G_n \prn[\big]{ \betahat[\vq] } = \Expectlet_n \brk{\mu_i \vp(\vx_i) \vp(\vx_i)\trans} \prn[\big]{ \betazero[\vq] - \betahat[\vq] },
\end{equation}
where
\begin{equation*}
  \mu_i := \frac{\partial}{\partial \theta} \prn[\big]{ \riskder[\vq]{y_i}{\eta(\theta)} \eta^{(1)}(\theta) } \Big|_{\theta = \widetilde{\theta}_i}\quad\text{for some $\widetilde{\theta}_i$ between $\vp(\vx_i)\trans \betazero[\vq]$ and $\vp(\vx_i)\trans \betahat[\vq]$}.
\end{equation*}
By the assumption of the lemma, $0 < \lambda \leq \mu_i \leq \Lambda$. Therefore, for any vector $\boldsymbol{a} \in \mathbb{R}^K$
\begin{align*}
  &\lambda \cdot \boldsymbol{a}\trans \Expectlet_n \brk{\vp(\vx_i) \vp(\vx_i)\trans} \boldsymbol{a}
    \leq \boldsymbol{a}\trans \Expectlet_n \brk{\mu_i \vp(\vx_i) \vp(\vx_i)\trans} \boldsymbol{a}\\
  &\quad = \Expectlet_n \brk{\mu_i \prn{\vp(\vx_i)\trans \boldsymbol{a}}^2} \leq \Lambda \cdot \boldsymbol{a}\trans \Expectlet_n \brk{\vp(\vx_i) \vp(\vx_i)\trans} \boldsymbol{a}.
\end{align*}
Moreover, the matrix
\begin{equation*}
\Expectlet_n \brk{\mu_i \vp(\vx_i) \vp(\vx_i)\trans}
\end{equation*}
has the same multi-banded structure as
\begin{equation*}
\Expectlet_n \brk{\vp(\vx_i) \vp(\vx_i)\trans}.
\end{equation*}
That means that by the same argument as that in \cref{lem:holder-gram} we have
\begin{equation}\label{eq:strong-shape-matrix-inf-norm-bound}
  \norm[\big]{ \Expectlet_n \brk{\mu_i \vp(\vx_i) \vp(\vx_i)\trans}^{-1} }_{\infty} \lesssim_{\Problet} h^{-d}.
\end{equation}

It is shown in the proofs of \cref{lem:holder-uniform-convergence-variance} and \cref{lem:holder-projection-approximation-error} (see also \cref{rem:rates-in-uniform-conv-var,rem:rates-in-proj-approx-error}) that
\begin{equation}\label{eq:strong-shape-grad-inf-norm-bound}
  \norm[\big]{ G_n \prn{ \betazero[\vq] } }_{\infty} = o_{\Problet} (h^d).
\end{equation}

From~\eqref{eq:strong-shape-gradient-diff}
\begin{equation*}
  \norm{ \betazero[\vq] - \betahat[\vq] }_{\infty} \leq \norm[\big]{ \Expectlet_n \brk{\mu_i \vp(\vx_i) \vp(\vx_i)\trans}^{-1} }_{\infty} \cdot \norm{ G_n \prn{ \betazero[\vq] } }_{\infty},
\end{equation*}
which in combination with~\eqref{eq:strong-shape-matrix-inf-norm-bound} and \eqref{eq:strong-shape-grad-inf-norm-bound} gives
\begin{equation*}
  \norm{ \betazero[\vq] - \betahat[\vq] }_{\infty} = o_{\Problet}(1)
\end{equation*}
uniformly over $\vq \in \mathcal{Q}$.
\end{proof}

\section{Bahadur representation}\label{sec:bahadur-representation}

We will now prove our first main result, the novel Bahadur representation which is Theorem~\mainref{theorem: bahadur-repres} in the paper. Recall the notation
\begin{equation}\label{eq:bahterm-def}
  \bahterm[\vq][(\vv)]{\vx} := - \vp^{(\vv)}(\vx)\trans \qzero[\vq][-1] \mathbb{E}_{n}\brk{\vp(\vx_i) \eta^{(1)}\prn{\muzero[\vq]{\vx_i}} \riskder[\vq]{y_i}{\eta\prn{\muzero[\vq]{\vx_i}}}}.
\end{equation}

\begin{theorem}[Bahadur representation]\label{th:bahadur-repres}\leavevmode
  Suppose \cref{ass:dgp-holder,ass:gl-quasi-uniform-partition,ass:gl-local-basis,ass:gl-holder,ass:complexity,ass:gl-approximation-error} hold. Furthermore, assume the following conditions:
  \begin{conditions}[wide]
\item $\sup_{\vq \in \mathcal{Q}}\norm[\big]{\betahatqconstr - \betazero[\vq]}_{\infty} = o_{\Problet}(1)$;
\item there exists a constant $c > 0$ such that $\{\vbeta \in \mathbb{R}^K: \norm{\vbeta - \betazero[\vq]}_{\infty} \leq c, \vq \in \mathcal{Q}\} \subseteq \mathcal{B}$;
\item $\frac{\log^{d + 2} n}{n h^d} = o(1)$;
\item either $\frac{\prn{h^{-1} \log n }^{{\frac{\nu}{\nu - 2} d}}}{n} = o(1)$, or $\riskderenvel{\vx_i}{y_i}$ is \subgauss\ conditionally on $\vx_i$;
\item either $\theta \mapsto \riskfunc[\vq]{y}{\eta(\theta)}$ is convex with left or right derivative $\riskder[\vq]{y}{\eta(\theta)} \eta^{(1)}(\theta)$ and $\mathcal{B} = \mathbb{R}^K$, or the additional complexity \cref{ass:complexity-addition} holds.
\end{conditions}

\begin{assertions}
\item\label{th:bahadur-repres-normal} Then
  \begin{equation}\label{eq:bahadur-representation-conclusion-expr}
    \sup _{\vq \in \mathcal{Q}, \vx \in \mathcal{X}} \abs[\big]{\muhat[\vq][(\vv)]{\vx}
      - \muzero[\vq][(\vv)]{\vx} - \bahterm[\vq][(\vv)]{\vx}} \lesssim_{\mathbb{P}} h^{- \multindnorm{\vv}} \rbah
  \end{equation}
  with
  \begin{equation}
    \label{eq:bahadur-rates-for-statement}
    \rbah := \prn[\bigg]{\frac{\log^d n}{n h^d}}^{1 / 2 + \prn{\sfrac{\alpha}{2} \wedge \sfrac{1}{4}}} \log n + h^{\prn{\alpha \wedge \sfrac{1}{2}} m} \prn[\bigg]{\frac{\log^{d + 1} n}{n h^d}}^{1 / 2} + h^{m}.
  \end{equation}

\item\label{th:bahadur-repres-cont} If, in addition to the previous conditions, \eqref{eq:holder-assumption-without-cases-nonconvex} holds (without any restrictions on $y$), then
  \begin{equation}\label{eq:bahadur-representation-conclusion-expr-cont}
    \sup _{\vq \in \mathcal{Q}, \vx \in \mathcal{X}} \abs[\big]{\muhat[\vq][(\vv)]{\vx} - \muzero[\vq][(\vv)]{\vx} - \bahterm[\vq][(\vv)]{\vx}} \lesssim_{\mathbb{P}} h^{- \multindnorm{\vv}} \rbahcont
  \end{equation}
  with
    \begin{equation}
    \label{eq:bahadur-rates-for-statement-cont}
    \rbahcont := \prn[\bigg]{\frac{\log^d n}{n h^d}}^{(1 + \alpha) / 2} \log n + h^{\alpha m}\prn[\bigg]{\frac{\log^{d + 1} n}{n h^d}}^{1 / 2} + h^{m}.
  \end{equation}
\end{assertions}
\end{theorem}

\begin{remark}\label{rem:qzero-vs-qbar}
As will be clear in the proof, the matrix $\qzero[\vq][-1]$ in \cref{eq:bahadur-representation-conclusion-expr,eq:bahadur-representation-conclusion-expr-cont} can be replaced by $\qbar[\vq][-1]$.
\end{remark}

\subsection{Proof: convex case}
\label{sec:bahadur-representation-convex-proof}

\setcounter{proofstepcounter}{0}
We will now prove \cref{th:bahadur-repres} under the assumption that $\theta \mapsto \riskfunc[\vq]{y}{\eta(\theta)}$ is convex with left or right derivative $\riskder[\vq]{y}{\eta(\theta)} \eta^{(1)}(\theta)$ and $\mathcal{B} = \mathbb{R}^K$. We only show the \ref{th:bahadur-repres-normal} part, since the argument for \ref{th:bahadur-repres-cont} is very similar with minor changes in obvious places.

\paragraph{Notation} In this proof, we will denote
\begin{equation*}
\mathbb{G}_n^i \brk{g(\vx_i, y_i)} := \sqrt{n} \Expectlet_n \brk[\big]{g(\vx_i, y_i) - \Expect{g(\vx_i, y_i) \given \vx_i}}.
\end{equation*}

\subsubsection{Strategy}
By $\sup_{\vx} \abs{\vp^{(\vv)}(\vx)} \lesssim h^{- \multindnorm{\vv}}$ (\cref{ass:gl-local-basis}), $\abs{\vp^{(\vv)}(\vx)\trans \betazero[\vq] - \muzero[\vq][(\vv)]{\vx}} \lesssim h^{m - \multindnorm{\vv}}$ (\cref{ass:gl-approximation-error}),
and \cref{lem:holder-projection-approximation-error}, the inequality
\begin{equation}\label{eq:psuSJW}
  \begin{aligned}
    &\sup _{\vq \in \mathcal{Q}, \vx \in \mathcal{X}} \abs[\big]{\vp^{(\vv)}(\vx)\trans \betahat[\vq]
      - \muzero[\vq][(\vv)]{\vx} \vphantom{\mathbb{E}_{n}\brk{\eta^{(1)}}}\\
    &\quad\quad + \vp^{(\vv)}(\vx)\trans \qbar[\vq][-1] \mathbb{E}_{n}\brk{\vp(\vx_i) \eta^{(1)}\prn{\muzero[\vq]{\vx_i}} \riskder[\vq]{y_i}{\eta\prn{\muzero[\vq]{\vx_i}}}}}\\
    &\quad\lesssim_{\mathbb{P}} h^{- \multindnorm{\vv}} \rbah
  \end{aligned}
\end{equation}
is implied by
\begin{equation}\label{eq:SdSOtN}
  \sup _{\vq}
    \norm[\big]{
      \betahat[\vq]
      - \betazero[\vq]
    - \bahbetabar[\vq]
  }_{\infty} \lesssim_{\mathbb{P}} \prn[\bigg]{\frac{\log^d n}{n h^d}}^{1 / 2 + \prn{\sfrac{\alpha}{2} \wedge \sfrac{1}{4}}} \log n + h^{\prn{\alpha \wedge \sfrac{1}{2}} m} \prn[\bigg]{\frac{\log^{d + 1} n}{n h^d}}^{1 / 2},
\end{equation}
where we put
\begin{align*}
  &\vp_i := \vp(\vx_i),\\
  &\bahbetabarNoLink[\vq] :=
    -\qbar[\vq][-1] \Expectlet_n \brk[\big]{\vp_{i} \eta^{(1)}(\vp_{i}\trans \betazero[\vq]) \riskder[\vq]{y_{i}}{\eta(\vp_{i}\trans\betazero[\vq])}}\numberthis\label{eq:bahbetabar-def}
\end{align*}
to declutter notation, and noted that for $h < 1$
\begin{equation}\label{eq:hgDEXQ}
h^{\prn{\alpha \wedge \sfrac{1}{2}} m} \prn[\bigg]{\frac{\log^{d + 1} n}{n h^d}}^{1 / 2} \lesssim \prn[\bigg]{\frac{\log^d n}{n h^d}}^{1 / 2 + \prn{\sfrac{\alpha}{2} \wedge \sfrac{1}{4}}} \log n + h^{m}.
\end{equation}
Indeed, if $\alpha \leq 1 / 2$, either $\prn[\big]{\frac{\log^{d + 1} n}{n h^d}}^{1 / 2} \leq h^m$,
in which case (for $h < 1$)
\begin{equation*}
h^{\alpha m} \prn[\bigg]{\frac{\log^{d + 1} n}{n h^d}}^{1 / 2} \leq \prn[\bigg]{\frac{\log^{d + 1} n}{n h^d}}^{1 / 2} \leq h^m,
\end{equation*}
or $\prn[\big]{\frac{\log^{d + 1} n}{n h^d}}^{1 / 2} > h^m$, in which case
\begin{equation*}
h^{\alpha m} \prn[\bigg]{\frac{\log^{d + 1} n}{n h^d}}^{1 / 2} \leq \prn[\bigg]{\frac{\log^{d + 1} n}{n h^d}}^{(1 + \alpha) / 2} \leq \prn[\bigg]{\frac{\log^d n}{n h^d}}^{(1 + \alpha) / 2} \log n.
\end{equation*}
If $\alpha > 1 / 2$, \eqref{eq:hgDEXQ} follows from $a b \lesssim a^2 + b^2$.

From \eqref{eq:psuSJW}, \eqref{eq:bahadur-representation-conclusion-expr} follows, because $\sup_{\vx} \norm{\vp^{(\vv)}(\vx)}_1 \lesssim h^{-\multindnorm{\vv}}$ and
\begin{equation}\label{eq:UMbHrv}
\begin{aligned}
  &\sup_{\vq} \norm[\big]{
  \prn{\qbar[\vq][-1] - \qzero[\vq][-1]} \mathbb{E}_{n}\brk{\vp(\vx_i) \eta^{(1)}\prn{\muzero[\vq]{\vx_i}} \riskder[\vq]{y_i}{\eta\prn{\muzero[\vq]{\vx_i}}}}
  }_{\infty}\\
  &\quad \lesssim \sup_{\vq} \norm{\qbar[\vq][-1] - \qzero[\vq][-1]}_{\infty} \sup_{\vq} \norm[\big]{\mathbb{E}_{n}\brk{\vp(\vx_i) \eta^{(1)}\prn{\muzero[\vq]{\vx_i}} \riskder[\vq]{y_i}{\eta\prn{\muzero[\vq]{\vx_i}}}}}_{\infty}\\
  &\quad \labrel{1}[\lesssim_{\Problet}]
    h^{-d} \prn[\bigg]{\frac{\log (1 / h)}{n h^d}}^{1 / 2}
    \prn[\bigg]{\frac{h^d \log (1 / h)}{n}}^{1 / 2}
  = \frac{\log (1 / h)}{n h^d} = o(\rbah),
\end{aligned}
\end{equation}
where \labrel{1} is by \cref{lem:holder-gram} and \cref{eq:uniform-conv-variance-enough}. So, we will be showing \eqref{eq:SdSOtN}.

Note that for any vector $\valpha$
\begin{equation}\label{eq:bahadur-repres-last-step-subgradient-inequality}
  \Expectlet_n \brk[\big]{
    \valpha\trans \vp_i \riskder[\vq]{y_i}{\eta(\vp_{i}\trans (\betahat[\vq] - \valpha))} \eta^{(1)}(\vp_{i}\trans (\betahat[\vq] - \valpha))
  } \leq 0,
\end{equation}
because $g\prn{\betahat[\vq] - \valpha} := \Expectlet_n \brk[\big]{\riskder[\vq]{y_i}{\eta( \vp_{i}\trans(\betahat[\vq] - \valpha))} \eta^{(1)}(\vp_{i}\trans (\betahat[\vq] - \valpha)) \vp_i}$ is a subgradient of the function $f(\vbeta) := \Expectlet_n \brk[\big]{\riskfunc[\vq]{y_i}{\eta\prn{ \vp_{i}\trans \vbeta }}}$ at $\betahat[\vq] - \valpha$, and $\betahat[\vq]$ is the minimizer of this function, giving $\valpha\trans g\prn[\big]{\betahat[\vq] - \valpha} \leq f\prn{\betahat[\vq]} - f\prn{\betahat[\vq] - \valpha} \leq 0$. Apply \eqref{eq:bahadur-repres-last-step-subgradient-inequality} with a particular vector $\valpha_{\vq}$ that will be chosen later, and decompose
\begin{align*}
  0 &\geq \Expectlet_n \brk[\big]{
  \valpha_{\vq}\trans \vp_i \riskder[\vq]{y_i}{\eta(\vp_{i}\trans (\betahat[\vq] - \valpha_{\vq}))} \eta^{(1)}(\vp_{i}\trans (\betahat[\vq] - \valpha_{\vq}))
  }\\
  &= \Expectlet_n \brk[\big]{
    \valpha_{\vq}\trans \vp_i
    \crl{\riskder[\vq]{y_i}{\eta(\vp_{i}\trans (\betahat[\vq] - \valpha_{\vq}))} - \riskder[\vq]{y_i}{\eta(\vp_{i}\trans \betazero[\vq])}}
    \eta^{(1)}(\vp_{i}\trans (\betahat[\vq] - \valpha_{\vq}))
    }\numberthis\label{eq:lqjlhq}
  \\
    &\quad + \termdeltaetaonebound
    + \Expectlet_n \brk[\big]{\valpha_{\vq}\trans \vp_i \riskder[\vq]{y_i}{\eta(\vp_{i}\trans \betazero[\vq])} \eta^{(1)}(\vp_{i}\trans \betazero[\vq])},
\end{align*}
where
\begin{equation}\label{eq:term-delta-eta-one-bound-def}
  \begin{aligned}
    &\termdeltaetaonebound := \Expectlet_n \brk[\big]{\valpha_{\vq}\trans \vp_{i} \riskder[\vq]{y_i}{\eta\prn{\vp_{i}\trans \betazero[\vq]}} \eta^{(1)}\prn{\vp_{i}\trans \prn{ \betahat[\vq] - \valpha_{\vq} }}}\\
    &\qquad - \Expectlet_n \brk[\big]{\valpha_{\vq}\trans \vp_{i} \riskder[\vq]{y_i}{\eta\prn{\vp_{i}\trans \betazero[\vq]}} \eta^{(1)}\prn{\vp_{i}\trans \betazero[\vq]}}
  \end{aligned}
\end{equation}
will be bounded later.
Define for simplicity
\begin{align*}
  \delta_{\vq, i}(\betaone, \betatwo, \valpha) &:= \valpha \trans \vp_{i} \crl[\big]{\riskder[\vq]{y_{i}}{\eta\prn{\vp_{i}\trans \prn{\betazero[\vq] + \betaone + \betatwo - \valpha}}} - \riskder[\vq]{y_i}{\eta\prn{\vp_{i}\trans \betazero[\vq] }}}\\
  &\hphantom{:=\ } \times \eta^{(1)}\prn{\vp_{i}\trans \prn{\betazero[\vq] + \betaone + \betatwo - \valpha}},
\end{align*}
so that
\begin{equation*}
\text{\eqref{eq:lqjlhq}} = \Expectlet_n \brk[\big]{\delta_{\vq, i}(\bahbetabar[\vq], \betahat[\vq] - \betazero[\vq] - \bahbetabar[\vq], \valpha_{\vq})}.
\end{equation*}
Now, add and subtract the conditional mean of this term, continuing
\begin{align*}
0 &\geq n^{-1 / 2} \mathbb{G}_n^i\brk[\big]{
    \delta_{\vq, i}(\bahbetabar[\vq], \betahat[\vq] - \betazero[\vq] - \bahbetabar[\vq], \valpha_{\vq})
    }\numberthis\label{eq:nnqhabbs0}\\
  &\quad + \Expectlet_n \brk[\big]{
    \Expect[\big]{
    \delta_{\vq, i}(\bahbetabar[\vq], \betahat[\vq] - \betazero[\vq] - \bahbetabar[\vq], \valpha_{\vq})
    \given \vx_i}
      }\numberthis\label{eq:nnqhabbs}\\
    &\quad + \termdeltaetaonebound
    + \Expectlet_n \brk[\big]{\valpha_{\vq}\trans \vp_i \riskder[\vq]{y_i}{\eta(\vp_{i}\trans \betazero[\vq])} \eta^{(1)}(\vp_{i}\trans \betazero[\vq])}.
\end{align*}
The difference \eqref{eq:nnqhabbs0} will be shown to be small by the usual concentration argument.
Using Taylor expansion, we will show that the term in \eqref{eq:nnqhabbs} is close to $\valpha_{\vq}\trans \qtilde[\vq] (\betahat[\vq] - \betazero[\vq])$ with the matrix
\begin{equation}\label{eq:qtilde-def}
\qtildeNoLink[\vq] := \Expectlet_n \brk[\big]{\vp_{i}\vp_{i}\trans \Psione[\vq]{\vx_i}{\eta(\vp_{i}\trans \betazero[\vq])} \eta^{(1)}\prn{\vp_{i}\trans\betazero[\vq]}^2},
\end{equation}
that is, we will bound
\begin{equation}\label{eq:term-cond-exp-bound-def}
\termCondExpBound := \Expectlet_n \brk[\big]{\Expect[\big]{\delta_{\vq, i} \prn{ \bahbetabar[\vq], \betahat[\vq] - \betazero[\vq] - \bahbetabar[\vq], \valpha_{\vq} } \given \vx_i}} - \valpha_{\vq}\trans \qtilde[\vq] \prn{ \betahat[\vq] - \betazero[\vq] }.
\end{equation}

To deal with the remaining term $\valpha_{\vq}\trans \qtilde[\vq] \prn{ \betahat[\vq] - \betazero[\vq] }$, start with replacing the approximation $\vp_i\trans \betazero[\vq]$ in the definition of $\qtilde[\vq]$ with the true function $\muzero[\vq]{\vx_i}$ leaving us with $\qbar[\vq]$; the error introduced by this operation
\begin{equation}\label{eq:term-delta-q-def}
\termdeltaq := \valpha_{\vq}\trans \prn{\qtilde[\vq] - \bar{\boldsymbol{Q}}_{\vq}} \prn{ \betahat[\vq] - \betazero[\vq] }
\end{equation}
will be bounded. Next, write
\begin{align*}
  &\valpha_{\vq}\trans \qbar[\vq] \prn{ \betahat[\vq] - \betazero[\vq] } + \Expectlet_n \brk[\big]{\valpha_{\vq}\trans \vp_i \riskder[\vq]{y_i}{\eta(\vp_{i}\trans \betazero[\vq])} \eta^{(1)}(\vp_{i}\trans \betazero[\vq])}\\
  &\quad = \valpha_{\vq}\trans \qbar[\vq] \prn{ \betahat[\vq] - \betazero[\vq] - \bahbetabar[\vq] }.
\end{align*}

We obtained
\begin{equation}\label{eq:qqqhabb}
0 \geq n^{- 1 / 2} \mathbb{G}_n^i\brk[\big]{
    \delta_{\vq, i}(\bahbetabar[\vq], \betahat[\vq] - \betazero[\vq] - \bahbetabar[\vq], \valpha_{\vq})
    } + \termdeltaetaonebound + \termCondExpBound + \termdeltaq + \valpha_{\vq}\trans \qbar[\vq] \prn{ \betahat[\vq] - \betazero[\vq] - \bahbetabar[\vq] }.
\end{equation}

At this point, it would be convenient if we could choose $\valpha_{\vq}$ so that the last term $\valpha_{\vq}\trans \qbar[\vq] \prn{ \betahat[\vq] - \betazero[\vq] - \bahbetabar[\vq] }$ is proportional (with a positive coefficient) to $\norm{\betahat[\vq] - \betazero[\vq] - \bahbetabar[\vq]}_{\infty}$, because then in combination with all the other terms being small in absolute value, we can conclude that $\norm{\betahat[\vq] - \betazero[\vq] - \bahbetabar[\vq]}_{\infty}$ is small in absolute value as well (otherwise it is impossible to obtain a nonnegative quantity). This is essentially what we will do, with one caveat: for the bounds on the other terms to work out, it is helpful if $\valpha_{\vq}$ is very sparse. Following these considerations, we introduce another vector $\bar{\valpha}_{\vq}$ proportional to $\brk[\big]{\qbar[\vq][-1]}_{\vk, k} \sign\prn[\big]{(\betahat[\vq] - \betazero[\vq] - \bahbetabar[\vq])_{\vk, k}}$ (with a positive coefficient), where $(\vk, k)$ is the index of the largest component of $\betahat[\vq] - \betazero[\vq] - \bahbetabar[\vq]$ (in the blockwise index notation as in \cref{lem:holder-gram}), and $\brk[\big]{\qbar[\vq][-1]}_{\vk, k}$ is the $(\vk, k)$th row of $\qbar[\vq][-1]$: this way,
$\bar{\valpha}_{\vq}\trans \qbar[\vq] \prn{ \betahat[\vq] - \betazero[\vq] - \bahbetabar[\vq] }$ is indeed proportional to $\norm{\betahat[\vq] - \betazero[\vq] - \bahbetabar[\vq]}_{\infty}$. Then we choose $\valpha_{\vq}$ to be a sparse vector that is close enough to $\bar{\valpha}_{\vq}$ for an appropriate bound on
\begin{equation}\label{eq:term-delta-v-bound-def}
  \termdeltavbound := \prn{\valpha_{\vq} - \bar{\valpha}_{\vq}}\trans \bar{\boldsymbol{Q}}_{\vq} \prn{ \betahat[\vq] - \betazero[\vq] - \bahbetabar[\vq] },
\end{equation}
which can be done because of the structure of the rows of $\qbar[\vq][-1]$, specifically the exponential decay in \cref{eq:inv-multiband-exp-decay}.

\subsubsection{Main argument}
We will now give some specifics. To show a precise bound in probability in \eqref{eq:bahadur-representation-conclusion-expr} (as opposed to a $o_{\Problet}(\cdot)$ bound), it will be convenient to multiply $\rbah$ by another positive sequence $\gamma_n$ that arbitrarily slowly diverges to infinity:
\begin{equation}\label{eq:bahadur-rtwo-def}
\rtwon := \rbah \gamma_n.
\end{equation}
We will also put
\begin{equation}\label{eq:bahadur-rone-def}
\rone := \brk[\bigg]{\prn[\bigg]{\frac{\log^d n}{n h^d}}^{1 / 2} + h^m} \gamma_n^{\sfrac{1}{(1 + \alpha \wedge (1 / 2))}},
\end{equation}
so that in particular
\begin{equation}\label{eq:two-lemmas-restrict-rone-strong}
  \prn[\bigg]{\frac{\log^d n}{n h^d}}^{1 / 2} + h^m = o(\rone)
\end{equation}
and by \cref{lem:holder-uniform-convergence-variance,lem:holder-projection-approximation-error},
we have on the event $\betabarboundedeventDef$ of $1 - o(1)$ probability
\begin{equation}\label{eq:bahadur-bound-betabar}
  \sup_{\vq \in \mathcal{Q}} \norm{\bahbetabar[\vq]}_{\infty} \leq \rone.
\end{equation}
Let $\gamma_n$ diverge slowly enough so that $\rone + \rtwon = o(1)$; in particular, $\sup_{\vq \in \mathcal{Q}} \norm{\bahbetabar[\vq]}_{\infty}$ is smaller than any positive constant with probability approaching one.

Fix two small enough constants $\nc{space}, \nc{vbar} > 0$ (the restrictions on which will be discussed later).
Since $\oc{space}$ is a constant, with probability approaching one both $\sup_{\vq \in \mathcal{Q}} \norm{\betahat[\vq] - \betazero[\vq]}_{\infty} \leq \oc{space}/2$ (by consistency) and $\sup_{\vq \in \mathcal{Q}} \norm{\bahbetabar[\vq]}_{\infty} \leq \oc{space}/2$ (as just discussed); hence, with probability approaching one $\sup_{\vq \in \mathcal{Q}} \norm{\betahat[\vq] - \betazero[\vq] - \bahbetabar[\vq]}_{\infty} \leq \oc{space}$.
This means that the probability of the event $\bigcap_{\vq \in \mathcal{Q}} \crl{\betahat[\vq] \in \mathcal{O}_{\vq}}$ approaches one, where we use the partitioning
\begin{equation*}
\mathcal{O}_{\vq} := \bigcup_{\ell = -\infty}^{\bar{L}_n} \mathcal{O}_{\vq, \ell},\quad \mathcal{O}_{\vq, \ell} := \crl[\big]{\vbeta \in \mathbb{R}^K \colon 2^{\ell - 1} \rtwon \leq \norm{\vbeta - \betazero[\vq] - \bahbetabar[\vq]}_{\infty} < 2^{\ell} \rtwon}
\end{equation*}
with $\bar{L}_n$ defined as the smallest integer such that $2^{\bar{L}_n} \rtwon \geq \oc{space}$. (The sequence $\bar{L}_n$ diverges to infinity.)

Put
\begin{equation*}
\bar{\valpha}_{\vq} = \oc{vbar} 2^L \rtwon h^{d} \brk[\big]{\qbar[\vq][-1]}_{\vk, k} \sign\prn[\big]{(\betahat[\vq] - \betazero[\vq] - \bahbetabar[\vq])_{\vk, k}},
\end{equation*}
where $(\vk, k)$ is such an index that $\abs{(\betahat[\vq] - \betazero[\vq] - \bahbetabar[\vq])_{\vk, k}} = \norm{\betahat[\vq] - \betazero[\vq] - \bahbetabar[\vq]}_{\infty}$ (as we already discussed above) and $L$ is a constant integer chosen later. The vector $\bar{\valpha}_{\vq}$ is not sparse, but the components decay exponentially with the ``distance'' to the multi-index $(\vk, k)$ according to \cref{eq:inv-multiband-exp-decay}. Therefore, we can zero out all components except a logarithmic neighborhood around $(\vk, k)$ of this vector and control the error introduced by this operation: specifically, take $\valpha_{\vq} \in \mathbb{R}^K$ with components $v_{\vq, \boldsymbol{j}, j}= \bar{v}_{\vq, \boldsymbol{j}, j}$ for $\norm{\vj - \vk}_{\infty} \leq \nc{neighborhood} \log n$ and zero otherwise, where $\oc{neighborhood}$ is some constant. On the event
\begin{equation}\label{eq:eventi-def}
  \eventi := \crl[\Big]{\sup_{\delta \in \Delta} \Expectlet_n \brk{\Indlet(\vx_{i} \in \delta)} \leq \oC{bahadur-var-bound} h^d}
\end{equation}
of $1 - o(1)$ probability (for some large enough constant $\nC{bahadur-var-bound} > 0$), we have
\begin{equation*}
\norm{\qbar[\vq]}_{\infty} \lesssim h^d,
\end{equation*}
so on $\eventi \cap \crl[\big]{\betahat[\vq] \in \mathcal{O}_{\vq, \ell}}$
\begin{equation*}
  \abs{\termdeltavbound} \lesssim (2^L \rtwon) h^{d} (2^{\ell} \rtwon) n^{-\oc{difvvbar}} = o \prn{ h^{d} 2^{L + \ell} \rtwon[2] },
\end{equation*}
where $\nc{difvvbar} > 0$ is a constant (depending on the constant controlling the neighborhood size $\oc{neighborhood}$).

Bounding $\termdeltaetaonebound$, $\termCondExpBound$, $\termdeltaq$ is deferred to \cref{lem:bahadur-delta-eta1-bound,lem:bahadur-conditional-exp-bound,lem:delta-q-bound-claim}.  Specifically, on an event $\eventallDef := \betabarboundedevent \cap \eventi \cap \eventiii$ with probability $1 - o(1)$, using \cref{lem:bahadur-delta-eta1-bound} (where $\eventiii$ is defined) and the restriction
\begin{equation}\label{eq:delta-eta1-bound-restr}
  \rone[2] = o(\rtwon),
\end{equation}
we bound
\begin{equation*}
\abs{\termdeltaetaonebound}
    = O\prn{h^{d} \rone (\rone + 2^{\ell} \rtwon + 2^L \rtwon) 2^L \rtwon}
    = o \prn{ h^{d} 2^{L + \ell} \rtwon[2] }.
  \end{equation*}

Using \cref{lem:bahadur-conditional-exp-bound}, we bound also on $\crl[\big]{\betahat[\vq] \in \mathcal{O}_{\vq, \ell}} \cap \eventall$ with $L \leq \ell \leq \bar{L}_n$:
\begin{equation}\label{eq:JrjuDG}
\abs{\termCondExpBound} \leq \frac{\oc{vbar}}{8} h^d 2^{L + \ell} \rtwon[2].
\end{equation}
Indeed,
\begin{align*}
  &2^L \rtwon h^d \prn[\big]{\rone + 2^{\ell} \rtwon + 2^L \rtwon}^2 \lesssim 2^L \rtwon h^d \rone[2] + 2^L \rtwon h^d 2^{2 \ell} \rtwon[2]\\
  &\quad \labrel{1}[=] 2^L \rtwon h^d 2^{2 \ell} \rtwon[2] + o(h^d 2^L \rtwon[2]),
\end{align*}
where in \labrel{1} we used $\rone[2] = o(\rtwon)$ again. Since $2^{\ell} \rtwon \leq 2 \oc{space}$, we can therefore ensure
\begin{equation*}
\abs[\big]{\termCondExpBound + \valpha_{\vq}\trans \qtilde[\vq] \valpha_{\vq}} \leq \frac{\oc{vbar}}{16} h^d 2^{L + \ell} \rtwon[2]
\end{equation*}
by taking $\oc{space}$ sufficiently small relative to $\oc{vbar}$. Next, taking $\oc{vbar}$ sufficiently small (in particular, making $\oc{vbar}^2$ much smaller than $\oc{vbar}$), we can upper-bound $\lambda_{\max}(\qtilde[\vq]) \norm{\valpha_{\vq}}^2$ by $(\oc{vbar} / 16) h^d 2^{L + \ell} \rtwon[2]$ as well, giving \eqref{eq:JrjuDG}.

Finally, using \cref{lem:delta-q-bound-claim}, we bound on $\eventall$
\begin{equation*}
  \abs{\termdeltaq}
  = O \prn{h^{m + d} (2^L \rtwon) \prn{ \rone + 2^{\ell} \rtwon } }
  \labrel{1}[=] o\prn{h^{d} 2^{L + \ell} \rtwon[2]},
\end{equation*}
where in \labrel{1} we used
\begin{equation}\label{eq:delta-q-bound-rate-restr}
h^m \rone = o(\rtwon).
\end{equation}

Combining, we have
\begin{equation*}
\abs{\termdeltaetaonebound} + \abs{\termCondExpBound} + \abs{\termdeltaq} + \abs{\termdeltavbound} \leq \frac{\oc{vbar}}{4} h^d 2^{L + \ell} \rtwon[2].
\end{equation*}
Since
\begin{equation*}
\bar{\valpha}_{\vq}\trans \bar{\boldsymbol{Q}}_{\vq} \prn{ \betahat[\vq] - \betazero[\vq] - \bahbetabar[\vq] } = \oc{vbar} 2^L \rtwon h^{d} \norm{\betahat[\vq] - \betazero[\vq] - \bahbetabar[\vq]}_{\infty} \geq \oc{vbar} 2^L \rtwon h^{d} 2^{\ell - 1} \rtwon,
\end{equation*}
we conclude from \eqref{eq:qqqhabb} that on $\crl[\big]{\betahat[\vq] \in \mathcal{O}_{\vq, \ell}} \cap \eventall$, with $L \leq \ell \leq \bar{L}_n$,
\begin{align*}
  - n^{- 1 / 2} \mathbb{G}_n^i\brk[\big]{
  \delta_{\vq, i}(\bahbetabar[\vq], \betahat[\vq] - \betazero[\vq] - \bahbetabar[\vq], \valpha_{\vq})
  } &\geq \oc{vbar} 2^L \rtwon h^{d} 2^{\ell - 1} \rtwon - \frac{\oc{vbar}}{4} h^d 2^{L + \ell} \rtwon[2]\\
  &= \frac{\oc{vbar}}{4} h^d 2^{L + \ell} \rtwon[2].
\end{align*}
We will prove that the probability of this event is small enough, by using a concentration argument.
\begin{lemma}[Uniform concentration]\label{lem:bahadur-concentration-of-g}
  Define
\begin{align*}
  &\mathllap{\mathcal{V}} :=
    \begin{aligned}[t]
    \crl[\big]{
    &\valpha \in \mathbb{R}^K: \exists \text{ $d$-dimensional multi-index $\vk$},\\
  &\abs{v_{\boldsymbol{\ell}, l}} \leq \varrho^{\norm{\vk - \vl}_{\infty}} 2^L \rtwon \text{ for $\norm{\vk - \vl}_{\infty} \leq \oc{neighborhood} \log n$ and $v_{\boldsymbol{\ell}, l} = 0$ otherwise}
    },
  \end{aligned}\\
  \mathcal{H}_1 &:= \crl{\vbeta \in \mathbb{R}^K: \norm{\vbeta}_{\infty} \leq r_{1, n}},\\
  \mathcal{H}_{2, \ell} &:= \crl{\vbeta \in \mathbb{R}^K: \norm{\vbeta}_{\infty} \leq 2^{\ell} \rtwon},
\end{align*}
where $\varrho$ is the constant from \cref{lem:holder-gram}. On the event $\eventi$ defined in \cref{eq:eventi-def} we have
  \begin{equation}\label{eq:bahadur-concentration-of-g}
    \begin{aligned}
      &\Expect[\Big]{\sup_{\vq \in \mathcal{Q}, \betaone \in \mathcal{H}_1, \betatwo \in \mathcal{H}_{2, \ell}, \valpha \in \mathcal{V}} \abs[\big]{
        \Expectlet_n \brk{
        \delta_{\vq, i}(\betaone, \betatwo, \valpha) - \Expect{\delta_{\vq, i}(\betaone, \betatwo, \valpha) | \vx_i}
        }
        } \given \crl{\vx_i}_{i = 1}^n}\\
      &\quad\leq \oC{conc} \prn[\bigg]{\frac{h^{\sfrac{d}{2}} 2^L \rtwon \prn{\rone + 2^{\ell} \rtwon + 2^L \rtwon} ^{\alpha \wedge (1 / 2)}}{\sqrt{n}} \log^{\sfrac{(d+1)}{2}} n + \frac{2^L \rtwon \log^{d + 1} n}{n}},
    \end{aligned}
  \end{equation}
  where the constant $\nC{conc}$ does not depend on $L$ or $\ell$.
\end{lemma}

\Cref{sec:bOVPpP} is devoted to the proof of this fact.

On $\crl[\big]{\betahat[\vq] \in \mathcal{O}_{\vq, \ell}} \cap \eventall$, we have (as long as $\oc{vbar}$ is small enough)
\begin{equation*}
  \valpha_{\vq} \in \mathcal{V},\quad \bahbetabar[\vq] \in \mathcal{H}_1,\quad \betahat[\vq] - \betazero[\vq] - \bahbetabar[\vq] \in \mathcal{H}_{2, \ell}.
\end{equation*}
Therefore,
\begin{align*}
  &\Prob[\Big]{\sup_{\vq \in \mathcal{Q}} \norm[\big]{\betahat[\vq] - \betazero[\vq] - \bahbetabar[\vq]}_{\infty} \geq 2^L \rtwon \given \crl{\vx_i}_{i = 1}^n}\\
  &\quad = \Prob[\bigg]{\bigcap_{\vq \in \mathcal{Q}} \bigcup_{\ell = L + 1}^{\infty} \crl[\big]{\betahat[\vq] \in \mathcal{O}_{\vq, \ell}} \given \crl{\vx_i}_{i = 1}^n}\\
  &\quad = \Prob[\bigg]{\bigcap_{\vq \in \mathcal{Q}} \bigcup_{\ell = L + 1}^{\bar{L}_n} \crl[\big]{\betahat[\vq] \in \mathcal{O}_{\vq, \ell}} \cap \eventall \given \crl{\vx_i}_{i = 1}^n} + o_{\Problet}(1)\\
  &\quad \leq \Prob[\bigg]{
    \bigcup_{\ell = L + 1}^{\bar{L}_n}
    \crl[\Big]{
    \sup_{\vq \in \mathcal{Q}, \vbeta_1 \in \mathcal{H}_1, \vbeta_2 \in \mathcal{H}_{2, \ell}, \valpha \in \mathcal{V}} \abs[\big]{n^{- 1 / 2} \mathbb{G}_n^i\brk[\big]{
    \delta_{\vq, i}(\vbeta_1, \vbeta_2, \valpha)
    }} \geq \frac{\oc{vbar}}{4} h^d 2^{L + \ell} \rtwon[2]
    }
    \given \crl{\vx_i}_{i = 1}^n
    } \Indlet(\eventi) + o_{\Problet}(1)\\
  &\quad \leq \sum_{\ell = L + 1}^{\bar{L}_n} \Prob[\bigg]{
    \sup_{\vq \in \mathcal{Q}, \vbeta_1 \in \mathcal{H}_1, \vbeta_2 \in \mathcal{H}_{2, \ell}, \valpha \in \mathcal{V}} \abs[\big]{n^{- 1 / 2} \mathbb{G}_n^i\brk[\big]{
    \delta_{\vq, i}(\vbeta_1, \vbeta_2, \valpha)
    }} \geq \frac{\oc{vbar}}{4} h^d 2^{L + \ell} \rtwon[2]
    \given \crl{\vx_i}_{i = 1}^n
    } \Indlet(\eventi) + o_{\Problet}(1)\\
  &\quad \leq \sum_{\ell = L + 1}^{\bar{L}_n} \prn[\Big]{\frac{\oc{vbar}}{4} h^d 2^{L + \ell} \rtwon[2]}^{-1} \Expect[\Big]{
    \sup_{\vq \in \mathcal{Q}, \vbeta_1 \in \mathcal{H}_1, \vbeta_2 \in \mathcal{H}_{2, \ell}, \valpha \in \mathcal{V}} \abs[\big]{n^{- 1 / 2} \mathbb{G}_n^{i}\brk[\big]{
    \delta_{\vq, i}(\vbeta_1, \vbeta_2, \valpha)
    }}
    \given \crl{\vx_i}_{i = 1}^n
    } \Indlet(\eventi) + o_{\Problet}(1)\\
  &\quad \leq \sum_{\ell = L + 1}^{\bar{L}_n}
    \oC{conc}
    \prn[\Big]{\frac{\oc{vbar}}{4} h^d 2^{L + \ell} \rtwon[2]}^{-1}
    \prn[\bigg]{\frac{h^{\sfrac{d}{2}} 2^L \rtwon \prn{\rone + 2^{\ell} \rtwon + 2^L \rtwon} ^{\alpha \wedge (1 / 2)}}{\sqrt{n}} \log^{\sfrac{(d+1)}{2}} n + \frac{2^L \rtwon \log^{d + 1} n}{n}}\\
  &\quad \lesssim \sum_{\ell = L + 1}^{\bar{L}_n} 2^{- \ell} \rtwon[-1] \prn[\bigg]{
    \frac{\rone[\alpha \wedge (1 / 2)] + \prn[\big]{2^{\ell} \rtwon}^{\alpha \wedge (1 / 2)}}{\sqrt{n h^d}} \log^{\sfrac{(d+1)}{2}} n
    +
    \frac{\log^{d + 1} n}{n h^d}
    },\\
  &\quad \labrel{1}[\lesssim] \sum_{\ell = L + 1}^{\bar{L}_n} \prn[\big]{2^{- \ell} + 2^{\alpha \wedge (1 / 2) - \ell}} \lesssim \sum_{\ell = L + 1}^{\infty} 2^{-\ell} = 2^{-L},
\end{align*}
where in \labrel{1} we used
\begin{align}
  &\frac{\rone[\alpha \wedge (1 / 2)]}{\sqrt{n h^d}} \log^{\sfrac{(d + 1)}{2}} n \lesssim \rtwon,\label{eq:last-restr-first}\\
  &\frac{\log^{\sfrac{(d + 1)}{2}} n}{\sqrt{n h^d}} \lesssim \rtwon[1 - \alpha \wedge (1 / 2)],\label{eq:last-restr-second}\\
  &\frac{\log^{d + 1} n}{n h^d} \lesssim \rtwon.\label{eq:last-restr-third}
\end{align}

We conclude that
\begin{equation*}
\sup_{\vq \in \mathcal{Q}} \norm[\big]{\betahat[\vq] - \betazero[\vq] - \bahbetabar[\vq]}_{\infty} \lesssim_{\Problet} \rtwon.
\end{equation*}
Since $\gamma_n$ was arbitrarily slowly diverging, it follows that in fact
\begin{equation*}
\sup_{\vq \in \mathcal{Q}} \norm[\big]{\betahat[\vq] - \betazero[\vq] - \bahbetabar[\vq]}_{\infty} \lesssim_{\Problet} \rbah.
\end{equation*}

\subsubsection{Concentration argument (proof of \cref{lem:bahadur-concentration-of-g})}\label{sec:bOVPpP}
We will first argue that $\vbeta_1$ and $\vbeta_2$ can be considered effectively low-dimensional (of $\polylog n$ dimension).

Note that $\delta_{\vq, i}(\betaone, \betatwo, \valpha) \neq 0$ only if $\valpha \trans \vp_{i} \neq 0$. For each $\valpha \in \mathcal{V}$, let $\mathcal{J}_{\valpha} := \crl{j: v_j \neq 0}$. By construction, the cardinality of $\mathcal{J}_{\valpha}$ is bounded by $(2 \oc{neighborhood} \log n + 1)^d$. We have $\delta_{\vq, i}(\betaone, \betatwo, \valpha) \neq 0$ only if $p_j(\vx_i) \neq 0$ for some $j \in \mathcal{J}_{\valpha}$, which happens only if $\vx_i \in \mathcal{I}_{\valpha}$ where
\begin{equation*}
  \mathcal{I}_{\valpha} := \bigcup \crl{\delta \in \Delta \colon \delta \cap \supp p_j \neq \varnothing \text{ for some $j \in \mathcal{J}_{\valpha}$}}.
\end{equation*}
The family $\mathcal{I}_{\valpha}$ includes at most $\nc{locationl2} (\oc{neighborhood} \log n)^d$ cells. Moreover, at most $\nc{locationl3} (\oc{neighborhood} \log n)^d$ basis functions in $\vp$ have supports overlapping with $\mathcal{I}_{\valpha}$. Denote the set of indices of such basis functions by $\bar{\mathcal{J}}_{\valpha}$. Based on the above observations, we have $\delta_{\vq, i}(\betaone, \betatwo, \valpha) = \delta_{\vq, i}\prn{\vbeta_{1, \bar{\mathcal{J}}_{\valpha}}, \vbeta_{2, \bar{\mathcal{J}}_{\valpha}}, \valpha}$, where
\begin{align*}
  &\delta_{\vq, i}\prn{\vbeta_{1, \bar{\mathcal{J}}_{\valpha}}, \vbeta_{2, \bar{\mathcal{J}}_{\valpha}}, \valpha}\numberthis \label{eq:delta-q-i-after-dimension-considerations}\\
  &\quad := \sum_{j \in \mathcal{J}_{\valpha}} p_{i, j} v_{j} \brk[\bigg]{
    \riskder[\vq][\bigg]{y_i}{\eta\prn[\bigg]{\sum_{l \in \bar{\mathcal{J}}_{\valpha}} p_{i, l}\prn{\beta_{0, \vq, l}+ \beta_{1, l} + \beta_{2, l}} - \sum_{j \in \mathcal{J}_{\valpha}} p_{i, j} v_{j} }}\\
  &\qquad - \riskder[\vq][\bigg]{y_i}{\eta\prn[\bigg]{\sum_{l \in \bar{\mathcal{J}}_{\valpha}} p_{i, l} \beta_{0, \vq, l}}}
    }\\
  &\qquad \times \eta^{(1)}\prn[\bigg]{\sum_{l \in \bar{\mathcal{J}}_{\valpha}} p_{i, l}\prn{\beta_{0, \vq, l}+\beta_{1, l}+\beta_{2, l}} - \sum_{j \in \mathcal{J}_{\valpha}} p_{i, j} v_{j} } \Indlet\{\vx_i \in \mathcal{I}_{\valpha}\}.
\end{align*}
Accordingly, for $\widetilde{\vbeta}_1 \in \mathbb{R}^{\oc{locationl3} (\oc{neighborhood} \log n)^d}, \widetilde{\vbeta}_2 \in \mathbb{R}^{\oc{locationl3} (\oc{neighborhood} \log n)^d}$, define the function class
\begin{equation*}
  \mathcal{G} := \crl[\big]{
    \prn{\vx_i, y_i} \mapsto \delta_{\vq, i}\prn[\big]{\widetilde{\vbeta}_1, \widetilde{\vbeta}_2, \valpha}:
    \vq \in \mathcal{Q}, \valpha \in \mathcal{V},
    \norm[\big]{\widetilde{\vbeta}_1}_{\infty} \leq \rone,
    \norm[\big]{\widetilde{\vbeta}_2}_{\infty} \leq 2^{\ell} \rtwon
  }.
\end{equation*}

We will now bound $\sup_{g \in \mathcal{G}}\abs{\mathbb{E}_n[g\prn{\vx_i, y_{i}}] - \Expect[\big]{g\prn{\vx_i, y_{i}} \given \vx_i}}$. As usual, the strategy is to check conditions of \cref{lem:ueb-maximal-ineq}.

\begin{lemma}[Bounding variance]\label{lem:bahadur-g-variance-bound}
  On $\eventi$, class $\mathcal{G}$ satisfies the following variance bound:
  \begin{equation*}
    \sup _{g \in \mathcal{G}} \Expectlet_n \brk{\Var{g\prn{\vx_i, y_{i}} \mid \vx_i}}
    \lesssim 2^{2 L} \rtwon[2] h^{d} \prn{\rone + 2^{\ell} \rtwon  + 2^L \rtwon}^{(2 \alpha) \wedge 1}.
  \end{equation*}
\end{lemma}

\begin{proof}
We will now proceed under the assumption that $\measureonY$ is Lebesgue measure, so~\eqref{eq:holder-assumption-with-cases-nonconvex} holds; the argument under~\eqref{eq:holder-assumption-without-cases-nonconvex} is similar (and leads to an even stronger variance bound), so it is omitted.

By the same argument as in the proof of \cref{lem:holder-projection-approximation-error}, using $\abs[\big]{\riskder[\vq]{y}{\eta \prn{p(\vx)\trans \betazero[\vq]}}} \lesssim \riskderenvel{\vx}{y} + 1$, for $y_i$ \textit{outside} the closed segment between the two points
\begin{equation*}
\eta \prn[\bigg]{\sum_{l \in \bar{\mathcal{J}}_{v}} p_{i, l} \beta_{0, \vq, l}},\quad\text{and}\quad \eta\prn[\bigg]{\sum_{l \in \bar{\mathcal{J}}_{v}} p_{i, l} \prn{\beta_{0, \vq, l} + \beta_{1, l} + \beta_{2, l} }- \sum_{j \in \mathcal{J}_{v}} p_{i, j} v_{j}},
\end{equation*}
we have
\begin{align*}
  &\abs[\bigg]{
    \riskder[\vq][\bigg]{y_{i}}{\eta\prn[\bigg]{\sum_{l \in \bar{\mathcal{J}}_{v}} p_{i, l} \prn{\beta_{0, \vq, l} + \beta_{1, l} + \beta_{2, l} }- \sum_{j \in \mathcal{J}_{v}} p_{i, j} v_{j} }} - \riskder[\vq][\bigg]{y_i}{\eta\prn[\bigg]{\sum_{l \in \bar{\mathcal{J}}_{v}} p_{i, l} \beta_{0, \vq, l}}}
    }\\
    &\quad \lesssim \prn[\big]{\riskderenvel{\vx_i}{y_i} + 1} \prn{\rone + 2^{\ell} \rtwon  + 2^L \rtwon}^{\alpha},
\end{align*}
and for $y_i$ \textit{in} this segment we have
\begin{align*}
  &\abs[\bigg]{\riskder[\vq][\bigg]{y_i}{\eta\prn[\bigg]{\sum_{l \in \bar{\mathcal{J}}_{v}} p_{i, l} \prn{\beta_{0, \vq, l} + \beta_{1, l} + \beta_{2, l} }- \sum_{j \in \mathcal{J}_{v}} p_{i, j} v_{j}}} - \riskder[\vq][\bigg]{y_i}{\eta\prn[\bigg]{\sum_{l \in \bar{\mathcal{J}}_{v}} p_{i, l} \beta_{0, \vq, l}}}}\\
  &\quad \lesssim \prn{\riskderenvel{\vx_i}{y_i} + 1} \prn{\rone + 2^{\ell} \rtwon + 2^L \rtwon} + 1
\end{align*}
uniformly over $\vq$.

By construction, for each $\valpha \in \mathcal{V}$, there exists some $\vk_{\valpha}$ such that $|v_{\vl, l}| \leq \varrho^{\|\vl - \vk_{\valpha}\|_{\infty}} 2^L \rtwon$ if $\|\vl - \vk_{\valpha}\|_{\infty} \leq M_{n}$, and otherwise $v_{\vl, l} = 0$. The above facts imply (cf. the proof of \cref{lem:holder-projection-approximation-error}) that for any $\vx_i \in \delta \subset \mathcal{I}_{\valpha}$,
\begin{align*}
  &\Var{\delta_{\vq, i}\prn{\vbeta_{1}, \vbeta_{2}, \valpha} \given \vx_i} \lesssim 2^{2 L} \rtwon[2] \prn{\rone +2^{\ell} \rtwon + 2^L \rtwon}^{(2 \alpha) \wedge 1} \sum_{(\vl, l) \in \mathcal{L}_{\delta}} \varrho^{2 \norm{\vl-\vk_{\valpha}}_{\infty}}\\
  &\text{for}\quad \mathcal{L}_{\delta} := \crl{(\vl, l)\colon \supp p_{\vl, l} \cap \delta \neq \varnothing}.
\end{align*}
In addition, since $\delta_{\vq, i}\prn{\vbeta_{1}, \vbeta_{2}, \valpha} \neq 0$ only if $\vx_i \in \mathcal{I}_{\valpha}$, for all $g \in \mathcal{G}$,
\begin{align*}
  &\mathbb{E}_{n}\brk{\Var{g\prn{\vx_i, y_{i}} \given \vx_i}}
    \lesssim 2^{2 L} \rtwon[2] \prn{\rone + 2^{\ell} \rtwon + 2^L \rtwon}^{(2 \alpha) \wedge 1} \sum_{\delta \subset \mathcal{I}_{\valpha}} \mathbb{E}_{n}\brk{\Indlet(\vx_i \in \delta)} \sum_{(\vl, l) \in \mathcal{L}_{\delta}} \varrho^{2 \norm{\vl - \vk_{\valpha}}_{\infty}}\\
  &\quad = 2^{2 L} \rtwon[2] \prn{\rone+2^{\ell} \rtwon +2^L \rtwon}^{(2 \alpha) \wedge 1} \sum_{\vl, l} \varrho^{2 \norm{\vl - \vk_{\valpha}}_{\infty}} \sum_{\delta \in \mathcal{L}_{\vl, l}^{\star}} \mathbb{E}_{n}\brk{\Indlet(\vx_i \in \delta)}\\
  &\text{for}\quad \mathcal{L}_{\vl, l}^{\star} := \crl{\delta \subset \mathcal{I}_{\valpha} \colon \supp p_{\vl, l} \cap \delta \neq \varnothing}.
\end{align*}

Note that $\mathcal{L}_{\vl, l}^{\star}$ contains a bounded number of elements. Then on $\eventi$,
\begin{align*}
  &\sup _{g \in \mathcal{G}} \Expectlet_n \brk{\Var{g\prn{\vx_i, y_{i}} \given\vx_i}}
    \lesssim 2^{2 L} \rtwon[2] h^{d} \prn{\rone+2^{\ell} \rtwon + 2^L \rtwon}^{(2 \alpha) \wedge 1} \sum_{\vl, l} \varrho^{2 \norm{\vl - \vk_{\valpha}}_{\infty}}\\
  &\quad \lesssim 2^{2 L} \rtwon[2] h^{d} \prn{\rone+2^{\ell} \rtwon + 2^L \rtwon}^{(2 \alpha) \wedge 1} \sum_{\vl} \varrho^{2 \norm{\vl - \vk_{\valpha}}_{\infty}} \tag*{since $l$ is bounded}\\
  &\quad \lesssim 2^{2 L} \rtwon[2] h^{d} \prn{\rone+2^{\ell} \rtwon + 2^L \rtwon}^{(2 \alpha) \wedge 1} \sum_{\vt \in \mathbb{Z}^d} \varrho^{2 \norm{\vt}_{\infty}} \lesssim 2^{2 L} \rtwon[2] h^{d} \prn{\rone + 2^{\ell} \rtwon + 2^L \rtwon}^{(2 \alpha) \wedge 1},
\end{align*}
concluding the proof of \cref{lem:bahadur-g-variance-bound}.
\end{proof}

\begin{lemma}[Complexity of class \texorpdfstring{$\mathcal{G}$}{G}]\label{lem:bahadur-g-complexity}
  Class $\mathcal{G}$ with envelope $2^L \rtwon$ multiplied by a large enough constant satisfies the uniform entropy bound~\eqref{eq:uniform-entropy-bound} with $A \lesssim 1$ and $V \lesssim \log^d n$.
\end{lemma}

\begin{proof}[Proof of \cref{lem:bahadur-g-complexity}]
First, indeed $\sup_{\vx, y} \sup_{g \in \mathcal{G}} \abs{g(\vx, y)} \lesssim 2^L \rtwon$.

Next, the class of functions $\mathcal{W}_1 := \crl{(\vx_i, y_i) \mapsto \valpha\trans \vp(\vx_i) \colon \valpha \in \mathcal{V}}$ is a union of $O\prn{h^{-d}}$ classes $\mathcal{W}_{1, \vk} := \crl{(\vx_i, y_i) \mapsto \valpha\trans \vp(\vx_i) \colon \valpha \in \mathcal{V}_{\vk}}$, where
\begin{equation*}
  \mathcal{V}_{\vk} = \crl[\big]{\valpha \in \mathbb{R}^{K}: \abs{v_{\boldsymbol{\ell}, l}} \leq \varrho^{\|\vk-\boldsymbol{\ell}\|_{\infty}} 2^L \rtwon \text{ for $\|\vk-\boldsymbol{\ell}\|_{\infty} \leq M_n$ and $v_{\boldsymbol{\ell}, l} = 0$ otherwise}}.
\end{equation*}
Since $\mathcal{W}_{1, \vk}$ is a subclass of a vector space of functions of dimension $O\prn{\log^d n}$, by Lemma~2.6.15 in~\parencite{vaart96_weak_conver_empir_proces} it is VC with index $O\prn{\log^d n}$. This implies that $\mathcal{W}_{1, \vk}$ with envelope $O(2^L \rtwon)$ satisfies the uniform entropy bound~\eqref{eq:uniform-entropy-bound} with $A \lesssim 1$ and $V \lesssim \log^d n$. Since there are $O\prn{ h^{-d} }$ such classes and $\log (1 / h) \lesssim \log n$, using the chain
\begin{equation}\label{eq:log-d-entropy-chain}
  O\prn{h^{-d}} \prn[\bigg]{\frac{A}{\varepsilon} }^{O\prn{\log^d n}} = e^{O\prn{\log n}} \prn[\bigg]{\frac{A}{\varepsilon} }^{O\prn{\log^d n}} \leq \prn[\bigg]{\frac{A}{\varepsilon} }^{O\prn{\log n} + O\prn{\log^d n}} = \prn[\bigg]{\frac{A}{\varepsilon} }^{O\prn{\log^d n}}
\end{equation}
(recall that $A \geq e$), we get that $\mathcal{W}_{1}$ also satisfies the uniform entropy bound~\eqref{eq:uniform-entropy-bound} with $A \lesssim 1$ and $V \lesssim \log^d n$.

By \cref{ass:complexity}, the class of functions
\begin{align*}
  \mathcal{W}_2 := \crl[\big]{&(\vx_i, y_i) \mapsto\\
  &\brk{\riskder[\vq]{y_{i}}{\eta\prn{\vp_{i}\trans \prn{\betazero[\vq] + \betaone + \betatwo - \valpha}}} - \riskder[\vq]{y_{i}}{\eta\prn{\vp_{i}\trans \betazero[\vq] }}} \Indlet\{\vx_i \in \mathcal{I}_{\valpha}\} \colon\\
  &\vbeta_{1} \in \mathcal{H}_1, \vbeta_2 \in \mathcal{H}_{2, \ell}, \valpha \in \mathcal{V}, \vq \in \mathcal{Q}}
\end{align*}
with a large enough constant envelope satisfies the uniform entropy bound~\eqref{eq:uniform-entropy-bound} with $A \lesssim 1$ and $V \lesssim \log^d n$.

The class of functions
\begin{align*}
  \mathcal{W}_3 := \crl[\big]{&(\vx_i, y_i) \mapsto \eta^{(1)} \prn{\vp(\vx_i) \trans \prn{\betazero[\vq] + \betaone + \betatwo - \valpha}} \Indlet\{\vx_i \in \mathcal{I}_{\valpha}\} \colon\\
  &\vbeta_1 \in \mathcal{H}_1, \vbeta_2 \in \mathcal{H}_{2, \ell}, \valpha \in \mathcal{V}, \vq \in \mathcal{Q}}
\end{align*}
is a subset of the union over $\delta \in \Delta$ of classes (for some fixed positive constants $c$ and $r$, $n$ large enough)
\begin{equation*}
  \mathcal{W}_{3, \delta} := \crl[\big]{(\vx_i, y_i) \mapsto \eta^{(1)} \prn{\vp(\vx_i) \trans \vbeta} \Indlet\{\vx_i \in \mathcal{N}_{[c \log n]}(\delta)\} \colon \norm{\vbeta - \betazero[\vq]}_{\infty} \leq r, \vq \in \mathcal{Q}}.
\end{equation*}
Note that $\vbeta$ can be assumed to lie in a fixed vector space $\mathcal{B}_{\delta}$ of dimension $\dim \mathcal{B}_{\delta} = O\prn{\log^d n}$. Again applying Lemma~2.6.15 in~\parencite{vaart96_weak_conver_empir_proces} and noting that $\eta^{(1)}$ on a fixed (bounded) interval is Lipschitz, we have that $\mathcal{W}_{3, \delta}$ with a large enough constant envelope satisfies the uniform entropy bound~\eqref{eq:uniform-entropy-bound} with $A \lesssim 1$ and $V \lesssim \log^d n$. Similarly to the argument for $\mathcal{W}_1$, this implies that the same is true for $\mathcal{W}_3$.

Applying \cref{lem:permanence-properties} concludes the proof of \cref{lem:bahadur-g-complexity}.
\end{proof}

\begin{proof}[Proof of \cref{lem:bahadur-concentration-of-g}]
Apply \cref{lem:ueb-maximal-ineq} conditionally on $\crl{\vx_i}_{i = 1}^n$ on $\eventi$.
\end{proof}

\subsubsection{Bounding \texorpdfstring{$\termdeltaetaonebound$, $\termCondExpBound$, $\termdeltaq$}{T1, T2, T3}}

\begin{lemma}[Bounding $\termdeltaetaonebound$]\label{lem:bahadur-delta-eta1-bound}
  There exists an event $\newtarget{def:eventiii}{\eventiiiNoLink}$ whose probability converges to one such that on $\eventiii$
  \begin{align*}
    &\sup_{\vq \in \mathcal{Q}, \betaone \in \mathcal{H}_1, \betatwo \in \mathcal{H}_{2, \ell}, \valpha \in \mathcal{V}} \abs{
      \Expectlet_n \brk{\valpha\trans\vp_{i} \riskder[\vq]{y_i}{\eta\prn{\vp_{i}\trans \betazero[\vq]}} \eta^{(1)}\prn{\vp_{i}\trans\prn{\betazero[\vq] + \betaone + \betatwo - \valpha}}}\\
    &\qquad- \Expectlet_n \brk{\valpha\trans\vp_{i} \riskder[\vq]{y_i}{\eta\prn{\vp_{i}\trans \betazero[\vq]}} \eta^{(1)}\prn{\vp_{i}\trans\betazero[\vq]}}} \lesssim h^{d} \rone (\rone + 2^{\ell} \rtwon  + 2^L \rtwon) 2^L \rtwon.
  \end{align*}
\end{lemma}

\begin{proof}[Proof of \cref{lem:bahadur-delta-eta1-bound}]
Note that
\begin{align*}
  &\Expectlet_n \brk[\big]{\valpha\trans\vp_{i} \riskder[\vq]{y_i}{\eta\prn{\vp_i\trans \betazero[\vq]}} \eta^{(1)}\prn{\vp_{i}\trans\prn{\betazero[\vq] + \betaone + \betatwo - \valpha}}}\\
  &\qquad- \Expectlet_n \brk{\valpha\trans\vp_i \riskder[\vq]{y_i}{\eta\prn{\vp_i\trans \betazero[\vq]}} \eta^{(1)}\prn{\vp_i\trans\betazero[\vq]}}\\
  &\quad= \valpha\trans \Expectlet_n \brk{\riskder[\vq]{y_i}{\eta\prn{\vp_{i}\trans \betazero[\vq]}} \eta^{(2)}(\xi_{\vq, i}) \vp_{i} \vp_{i}\trans} \prn{\betaone + \betatwo - \valpha}\\
  &\quad \lesssim h^{d} \rone (\rone + 2^{\ell} \rtwon  + 2^L \rtwon) 2^L \rtwon,
\end{align*}
where $\xi_{\vq, i}$ is between $\vp_{i}\trans(\betazero[\vq] + \betaone + \betatwo - \valpha)$ and $\vp_{i}\trans\betazero[\vq]$. The bound holds on the event
\begin{align*}
  \eventiii := \crl[\big]{\sup
  \norm[\big]{
  \Expectlet_n \brk{\riskder[\vq]{y_i}{\eta(\vp_{i}\trans \betazero[\vq])} \eta^{(2)}(\xi_{\vq, i}) \vp_{i} \vp_{i}\trans}
  }_{\infty} \leq h^{d} \rone
  },
\end{align*}
where the supremum is over $\betaone \in \mathcal{H}_1$, $\betatwo \in \mathcal{H}_{2, \ell}$, $\valpha \in \mathcal{V}$, $\vq \in \mathcal{Q}$ and $\xi_{\vq, i}$ between $\vp_{i}\trans(\betazero[\vq] + \betaone + \betatwo - \valpha)$ and $\vp_{i}\trans\betazero[\vq]$.
By the same argument as \cref{lem:holder-uniform-convergence-variance,lem:holder-projection-approximation-error}, $\Prob{\eventiii} \to 1$.
\end{proof}

\begin{lemma}[Bounding $\termCondExpBound$]\label{lem:bahadur-conditional-exp-bound}
  On $\eventi$, we have
  \begin{align*}
    &\sup_{\vq \in \mathcal{Q}, \betaone \in \mathcal{H}_1, \betatwo \in \mathcal{H}_{2, \ell}, \valpha \in \mathcal{V}} \abs[\big]{\Expectlet_n \brk{\Expect{\delta_{\vq, i}(\betaone, \betatwo, \valpha) | \vx_i}} - \valpha\trans \qtilde[\vq] \prn{\betaone + \betatwo} + \valpha\trans \qtilde[\vq] \valpha}\\
    &\quad \lesssim 2^L \rtwon h^{d} \prn{\rone + 2^{\ell} \rtwon  + 2^L \rtwon}^2,
  \end{align*}
  and in addition for all $\vq \in \mathcal{Q}$, $\valpha \in \mathcal{V}$
\begin{equation*}
\abs{\valpha\trans \qtilde[\vq] \valpha} \lesssim h^d \norm{\valpha}^2.
\end{equation*}
\end{lemma}

\begin{proof}[Proof of \cref{lem:bahadur-conditional-exp-bound}]
First, on $\eventi$ the largest eigenvalue of $\qtilde[\vq]$ is bounded by $h^d$ up to a constant factor (uniformly in $\vq$):
\begin{align*}
  &\lambda_{\mathrm{max}}\prn[\big]{\qtilde[\vq]} = \sup_{\|\balpha\| = 1} \balpha\trans \qtilde[\vq] \balpha = \sup_{\|\balpha\| = 1} \Expectlet_n \brk[\big]{(\balpha\trans \vp_{i})^2 \Psione[\vq]{\vx_i}{\eta(\vp_{i}\trans \betazero[\vq])} \eta^{(1)}(\vp_{i}\trans\betazero[\vq])^2}\\
  &\quad \leq \sup_{\|\balpha\| = 1} \Expectlet_n \brk[\big]{(\balpha\trans \vp_{i})^2 \abs[\big]{\Psione[\vq]{\vx_i}{\eta(\vp_{i}\trans \betazero[\vq])}} \eta^{(1)}(\vp_{i}\trans\betazero[\vq])^2}\\
  &\quad \lesssim \sup_{\|\balpha\| = 1} \Expectlet_n \brk[\big]{(\balpha\trans \vp_{i})^2}\\
  \shortintertext{(because $\abs[\big]{\Psione[\vq]{\vx_i}{\eta(\vp_{i}\trans \betazero[\vq])}} \eta^{(1)}(\vp_{i}\trans\betazero[\vq])^2 \lesssim 1$ by \cref{ass:gl-holder,ass:gl-approximation-error})}
  &\quad = \sup_{\|\balpha\| = 1} \sum_{l = 1}^{\kappa} \Expectlet_n \brk[\big]{(\balpha\trans \vp_{i})^2 \Indlet\{\vx_i \in \delta_l\}}\\
  &\quad \lesssim \sup_{\|\balpha\| = 1} \sum_{l = 1}^{\kappa} \prn[\bigg]{\sum_{k = 1}^{M_l} \alpha_{l, k}^2} \Expectlet_n \brk[\big]{\Indlet\{\vx_i \in \delta_l\}}\\
  \shortintertext{(because $\sup_{\vx \in \delta_{l}} (\balpha\trans \bp(\vx))^2 \lesssim \sum_{k = 1}^{M_l} \alpha_{l, k}^2$, where $\{\alpha_{l, k}\}_{k = 1}^{M_l}$ are the components of $\balpha$ corresponding to the $M_l$ basis functions supported on $\delta_l$)}
  &\quad \leq \sup_{\delta \in \Delta} \Expectlet_n \brk[\big]{\Indlet\{\vx_i \in \delta\}} \sup_{\|\balpha\| = 1} \sum_{l = 1}^{\kappa} \prn[\bigg]{\sum_{k = 1}^{M_l} \alpha_{l, k}^2}\\
  &\quad \lesssim \sup_{\delta \in \Delta} \Expectlet_n \brk[\big]{\Indlet\{\vx_i \in \delta\}} \sup_{\|\balpha\| = 1} \norm{\balpha}^2 = \sup_{\delta \in \Delta} \Expectlet_n \brk[\big]{\Indlet\{\vx_i \in \delta\}}.
\end{align*}

\sloppy Next, by Taylor expansion,
\begin{align*}
  &\Expect{\delta_{\vq, i}(\betaone, \betatwo, \valpha) | \vx_i}\\
  &\quad= \valpha \trans \vp_{i} \brk{\expriskder[\vq]{\vx_i}{\eta(\vp_{i}\trans (\betazero[\vq] + \betaone + \betatwo - \valpha))} - \expriskder[\vq]{\vx_i}{\eta(\vp_{i}\trans \betazero[\vq])}}\\
  &\quad\quad \times \eta^{(1)}\prn{\vp_{i}\trans \prn{\betazero[\vq] + \betaone + \betatwo - \valpha}}\\
  &\quad= \valpha \trans \vp_{i} \brk[\big]{\Psione[\vq]{\vx_i}{\eta(\vp_{i}\trans \betazero[\vq])}
    \crl[\big]{
    \eta^{(1)}\prn{ \vp_{i}\trans \betazero[\vq] } \vp_{i}\trans \prn{\betaone + \betatwo - \valpha}
    +\\
  &\quad\quad\quad + (1 / 2) \eta^{(2)}(\xi_{\vq, i}) \prn{\vp_{i}\trans \prn{\betaone + \betatwo - \valpha}}^2
    }\\
  &\quad\quad+ (1 / 2) \Psitwo[\vq]{\vx_i}{\tilde{\xi}_{\vq, i}} \crl[\big]{
    \eta\prn{\vp_{i}\trans \prn{\betazero[\vq] + \betaone + \betatwo - \valpha}} - \eta\prn{\vp_{i}\trans \betazero[\vq]}
    }^2}\\
  &\quad\quad\times \eta^{(1)}\prn{\vp_i\trans \prn{\betazero[\vq] + \betaone + \betatwo - \valpha}}
\end{align*}
for some $\xi_{\vq, i}$ between $\vp_{i}\trans \betazero[\vq]$ and $\vp_{i}\trans \prn{\betazero[\vq] + \betaone + \betatwo - \valpha}$, $\tilde{\xi}_{\vq, i}$ is between $\eta\prn{\vp_{i}\trans \betazero[\vq]}$ and $\eta\prn{\vp_{i}\trans \prn{\betazero[\vq] + \betaone + \betatwo - \valpha}}$.
This gives
\begin{equation*}
  \Expectlet_n \brk{\Expect{\delta_{\vq, i}(\betaone, \betatwo, \valpha) | \vx_i}} = \valpha\trans \qtilde[\vq] \prn{\betaone + \betatwo} - \valpha\trans \qtilde[\vq] \valpha + \text{I} + \text{II} + \text{III},
\end{equation*}
where for some $\check{\xi}_{\vq, i}$ between $\vp_{i}\trans \betazero[\vq]$ and $\vp_{i}\trans \prn{ \betazero[\vq] + \betaone + \betatwo - \valpha }$
\begin{align*}
  &\text{I} := \Expectlet_n \brk[\big]{\valpha\trans\vp_{i} \Psione[\vq]{\vx_i}{\eta(\vp_{i}\trans \betazero[\vq])} \eta^{(1)}\prn{\vp_{i}\trans \betazero[\vq]}\eta^{(2)}\prn{\check{\xi}_{\vq, i}} \prn{\vp_{i}\trans \prn{\betaone + \betatwo - \valpha}}^2},\\
  &\text{II} := \frac{1}{2} \Expectlet_n \brk[\big]{
    \valpha\trans\vp_{i} \Psione[\vq]{\vx_i}{\eta(\vp_{i}\trans \betazero[\vq])} \eta^{(2)}\prn{\xi_{\vq, i}}\\
  &\quad\times \eta^{(1)}\prn{\vp_{i}\trans \prn{\betazero[\vq] + \betaone + \betatwo - \valpha}} \prn{\vp_{i}\trans \prn{\betaone + \betatwo - \valpha} }^2},\\
  &\text{III} := \frac{1}{2} \Expectlet_n \brk[\big]{
    \valpha\trans\vp_{i} \Psitwo[\vq]{\vx_i}{\tilde{\xi}_{\vq, i}}
    \crl{\eta\prn{\vp_{i}\trans \prn{\betazero[\vq] + \betaone + \betatwo - \valpha}} - \eta\prn{\vp_{i}\trans \betazero[\vq]}}^2\\
  &\qquad\times \eta^{(1)}\prn{\vp_{i}\trans \prn{\betazero[\vq] + \betaone + \betatwo - \valpha}}}
\end{align*}
and
\begin{align*}
  &\text{I} \lesssim 2^L \rtwon h^{d} \prn{\rone + 2^{\ell} \rtwon  + 2^L \rtwon}^2,\\
  &\text{II} \lesssim 2^L \rtwon h^{d} \prn{\rone + 2^{\ell} \rtwon  + 2^L \rtwon}^2,\\
  &\text{III} \lesssim 2^L \rtwon h^{d} \prn{\rone + 2^{\ell} \rtwon  + 2^L \rtwon}^2
\end{align*}
on the event $\eventi$.
\end{proof}

\begin{lemma}[Bounding $\termdeltaq$]\label{lem:delta-q-bound-claim}
  For the matrix
  $\qtilde[\vq]$ defined in \cref{eq:qtilde-def},
  we have the following bound on the event $\eventi$:
  \begin{equation}\label{eq:delta-q-bound}
    \sup_{\vq \in \mathcal{Q}, \betaone \in \mathcal{H}_{1}, \betatwo \in \mathcal{H}_{2, \ell}, \valpha \in \mathcal{V}} \abs[\big]{\valpha\trans \prn{\qtilde[\vq] - \qbar[\vq]}\prn{\betaone + \betatwo}} \lesssim h^{m + d} 2^L \rtwon\prn{\rone+2^{\ell} \rtwon }.
  \end{equation}
\end{lemma}

\begin{proof}[Proof of \cref{lem:delta-q-bound-claim}]
By the same logic as in \cref{lem:holder-gram}, we have on $\eventi$
\begin{equation*}
  \norm{\qbar[\vq]-\qtilde[\vq]}_{\infty} \vee \norm{\qbar[\vq]-\qtilde[\vq]} \lesssim h^{m} h^{d}
\end{equation*}
uniformly over $\vq$ with probability approaching one. This gives~\eqref{eq:delta-q-bound}, proving \cref{lem:delta-q-bound-claim}.
\end{proof}

We have now proved the deferred lemmas, and the proof of \cref{th:bahadur-repres} is concluded.

\begin{remark}[Rate restrictions]
  The rates in the proof are determined by four restrictions: \cref{eq:two-lemmas-restrict-rone-strong,eq:delta-eta1-bound-restr,eq:last-restr-first,eq:last-restr-second}.
  \Cref{eq:last-restr-third} follows from \cref{eq:last-restr-first,eq:two-lemmas-restrict-rone-strong}. \Cref{eq:delta-q-bound-rate-restr} follows from \cref{eq:two-lemmas-restrict-rone-strong} and \cref{eq:delta-eta1-bound-restr}.
\end{remark}

\subsection{Proof: general case}\label{sec:bahadur-representation-nonconvex-proof}

We will state another version of the Bahadur Representation theorem, where $\theta \mapsto \riskfunc[\vq]{y}{\eta(\theta)}$ is not assumed to be convex. Here, we will use the additional complexity \cref{ass:complexity-addition}.

The argument is almost the same as for \cref{th:bahadur-repres}, so we will only describe the changes that need to be made.

The setup is the same as in the convex case except the definition of $\delta_{\vq, i}(\bbeta_1, \bbeta_2, \valpha)$ is replaced with
\begin{align*}
  &\delta_{\vq, i}(\betaone, \betatwo, \valpha)\\
  &\quad := \riskfunc[\vq]{y_i}{\eta\prn{\bp_i\trans \prn{\betazero[\vq] + \bbeta_1 + \bbeta_2 }}} - \riskfunc[\vq]{y_i}{\eta\prn{\bp_i\trans \prn{\betazero[\vq] + \bbeta_1 + \bbeta_2 - \valpha }}}\\
  &\qquad - \brk[\big]{\eta \prn{ \bp_i\trans \prn{ \betazero[\vq] + \bbeta_1 + \bbeta_2 } } - \eta \prn{ \bp_i\trans \prn{ \betazero[\vq] + \bbeta_1 + \bbeta_2 - \valpha } }} \riskder[\vq]{y_i}{\eta \prn{ \bp_i\trans \betazero[\vq] }}\\
  &\quad = \int_{-\bp_i\trans\valpha}^0 \brk[\big]{\riskder[\vq]{y_i}{\eta\prn{\bp_i\trans \prn{\betazero[\vq] + \bbeta_1 + \bbeta_2 } + t}} - \riskder[\vq]{y_i}{\eta\prn{\bp_i\trans \betazero[\vq]}}}\\
  &\qquad \times \eta^{(1)} \prn{ \bp_i\trans \prn{ \betazero[\vq] + \bbeta_1 + \bbeta_2 } + t }\,\mathrm{d} t,
\end{align*}
and \eqref{eq:delta-q-i-after-dimension-considerations} is changed to fit this definition.

Instead of~\eqref{eq:bahadur-repres-last-step-subgradient-inequality}, we have for any vector $\valpha$
\begin{equation*}
  \Expectlet_n \brk[\big]{\riskfunc[\vq]{y_i}{\eta \prn{ \bp_i\trans \betahatqconstr }} - \riskfunc[\vq]{y_i}{\eta \prn{ \bp_i\trans \prn{ \betahatqconstr - \valpha } }}} \leq 0
\end{equation*}
by the definition of $\betahatqconstr$ as long as $\norm{\valpha}_{\infty}$ is small enough (so that $\betahatqconstr - \valpha$ satisfies the constraints). Accordingly, the definition of $\termdeltaetaonebound$ becomes
\begin{align*}
  &\termdeltaetaonebound := \Expectlet_n \brk[\big]{\prn{ \eta \prn{ \bp_i\trans\betahatqconstr } - \eta \prn{ \bp_i\trans \prn{ \betahatqconstr - \valpha_{\vq} } } } \riskder[\vq]{y_i}{\eta \prn{ \bp_i\trans \betazero[\vq] }}}\\
  &\qquad - \Expectlet_n \brk[\big]{\valpha_{\vq}\trans \bp_i \riskder[\vq]{y_i}{\eta \prn{ \bp_i\trans \betazero[\vq] }} \eta^{(1)} \prn{ \bp_i\trans \betazero[\vq] }}.
\end{align*}

\begin{lemma}[Bounding variance]\label{lem:bahadur-nonconvex-g-variance-bound}
  On $\eventi$ as in \cref{lem:bahadur-g-variance-bound}, class $\mathcal{G}$ satisfies the following variance bound:
  \begin{equation*}
    \sup _{g \in \mathcal{G}} \Expectlet_n \brk{\Var{g\prn{\vx_i, y_{i}} \mid \vx_i}} \lesssim 2^{2 L} \rtwon[2] h^{d} \prn{\rone+2^{\ell} \rtwon +2^L \rtwon}^{(2 \alpha) \wedge 1}.
  \end{equation*}
\end{lemma}

\begin{proof}[Proof of \cref{lem:bahadur-nonconvex-g-variance-bound}]
This is proven by the same argument as in the proof of \cref{th:bahadur-repres}.
\end{proof}

\begin{lemma}[Complexity of class \texorpdfstring{$\mathcal{G}$}{G}]\label{lem:bahadur-nonconvex-g-complexity}
  Class $\mathcal{G}$ with envelope $2^L \rtwon$ multiplied by a large enough constant satisfies the uniform entropy bound~\eqref{eq:uniform-entropy-bound} with $A \lesssim (2^L \rtwon)^{-1}$ and $V \lesssim \log^d n$.
\end{lemma}

Note that $A$ is not constant in this statement but it will not matter since $\log((2^L \rtwon)^{-1}) \lesssim \log n$.

\begin{proof}[Proof of \cref{lem:bahadur-nonconvex-g-complexity}]
This is directly assumed in \cref{ass:complexity-addition}.
\end{proof}

\begin{lemma}[Uniform concentration in \texorpdfstring{$\mathcal{G}$}{G}]\label{lem:bahadur-nonconvex-concentration-of-g}
  On the event $\eventi$, \eqref{eq:bahadur-concentration-of-g} holds.
\end{lemma}

\begin{proof}[Proof of \cref{lem:bahadur-nonconvex-concentration-of-g}]
This is proven by the same argument as in the proof of \cref{th:bahadur-repres}.
\end{proof}

\begin{lemma}\label{lem:nonconvex-delta-q-bound-claim}
  For $\newtarget{def:qtildenonconv}{\qtildeNoLink[\vq]} := \Expectlet_n \brk[\big]{\vp_{i}\vp_{i}\trans \Psione[\vq]{\vx_i}{\eta(\vp_{i}\trans \betazero[\vq])} \eta^{(1)}\prn{\vp_{i}\trans\betazero[\vq]}^2}$, \eqref{eq:delta-q-bound} holds.
\end{lemma}

\begin{proof}[Proof of~\cref{lem:nonconvex-delta-q-bound-claim}]
This is proven by the same argument as in the proof of \cref{th:bahadur-repres}.
\end{proof}

\begin{lemma}\label{lem:nonconvex-bahadur-conditional-exp-bound}
  On $\eventi$, we have
  \begin{align*}
    &\sup_{\vq \in \mathcal{Q}, \betaone \in \mathcal{H}_1, \betatwo \in \mathcal{H}_{2, \ell}, \valpha \in \mathcal{V}} \abs[\big]{\Expectlet_n \brk[\big]{\Expect{\delta_{\vq, i}(\betaone, \betatwo, \valpha) | \vx_i}} - \valpha\trans \qtildenonconv[\vq] \prn{\betaone + \betatwo} + \valpha\trans \qtildenonconv[\vq] \valpha}\\
    &\quad \lesssim 2^L \rtwon h^{d} \prn{\rone + 2^{\ell} \rtwon  + 2^L \rtwon}^2,
  \end{align*}
and in addition for all $\vq \in \mathcal{Q}$, $\valpha \in \mathcal{V}$
\begin{equation*}
\abs{\valpha\trans \qtildenonconv[\vq] \valpha} \lesssim h^d \norm{\valpha}^2.
\end{equation*}
\end{lemma}

\begin{proof}[Proof of \cref{lem:nonconvex-bahadur-conditional-exp-bound}]
First, on $\eventi$ the largest eigenvalue of $\qtildenonconv[\vq]$ is bounded by $h^d$ up to a constant factor (uniformly in $\vq$), which is proven in \cref{lem:bahadur-conditional-exp-bound}.

\sloppy Next, by the Taylor expansion,
\begin{align*}
  &\Expect{\delta_{\vq, i}(\betaone, \betatwo, \valpha) | \vx_i}\\
  &\quad= \int_{-\bp_i\trans\valpha}^0 \brk[\big]{\expriskder[\vq]{\vx_i}{\eta(\vp_{i}\trans (\betazero[\vq] + \betaone + \betatwo) + t)} - \expriskder[\vq]{\vx_i}{\eta(\vp_{i}\trans \betazero[\vq])}}\\
  &\quad\quad \times \eta^{(1)}\prn{\vp_{i}\trans \prn{\betazero[\vq] + \betaone + \betatwo} + t} \,\mathrm{d} t\\
  &\quad= \int_{-\bp_i\trans\valpha}^0 \brk[\big]{
    \Psione[\vq]{\vx_i}{\eta(\vp_{i}\trans \betazero[\vq])}
    \crl{
    \eta^{(1)}\prn{ \vp_{i}\trans \betazero[\vq] } \prn{ \vp_{i}\trans \prn{\betaone + \betatwo} + t}
    \\
  &\quad\quad\quad +\frac{1}{2}\eta^{(2)}(\xi_{\vq, i, t}) \prn{\vp_{i}\trans \prn{\betaone + \betatwo} + t}^2}\\
  &\quad\quad+ \frac{1}{2} \Psitwo[\vq]{\vx_i}{\tilde{\xi}_{\vq, i, t}} \crl{\eta\prn{\vp_{i}\trans \prn{\betazero[\vq] + \betaone + \betatwo} + t} - \eta\prn{\vp_{i}\trans \betazero[\vq]}}^2
    }\\
  &\quad\quad\times \eta^{(1)}\prn{\vp_{i}\trans \prn{\betazero[\vq] + \betaone + \betatwo} + t} \,\mathrm{d} t
\end{align*}
for some $\xi_{\vq, i, t}$ between $\vp_{i}\trans \betazero[\vq]$ and $\vp_{i}\trans \prn{\betazero[\vq] + \betaone + \betatwo} + t$, $\tilde{\xi}_{\vq, i, t}$ between $\eta\prn{\vp_{i}\trans \betazero[\vq]}$ and $\eta\prn{\vp_{i}\trans \prn{\betazero[\vq] + \betaone + \betatwo} + t}$.
This gives
\begin{equation*}
  \Expectlet_n \brk{\Expect{\delta_{\vq, i}(\betaone, \betatwo, \valpha) | \vx_i}} = \valpha\trans \qtildenonconv[\vq] \prn{\betaone + \betatwo} - \valpha\trans \qtildenonconv[\vq] \valpha + \text{I} + \text{II} + \text{III},
\end{equation*}
where for some $\check{\xi}_{\vq, i, t}$ between $\vp_{i}\trans \betazero[\vq]$ and $\vp_{i}\trans \prn{ \betazero[\vq] + \betaone + \betatwo} + t$ again
\begin{align*}
  &\text{I} := \Expectlet_n \brk[\bigg]{\int_{-\bp_i\trans\valpha}^0 \Psione[\vq]{\vx_i}{\eta(\vp_{i}\trans \betazero[\vq])} \eta^{(1)}\prn{\vp_{i}\trans \betazero[\vq]}\eta^{(2)}\prn{\check{\xi}_{\vq, i}} \prn{\vp_{i}\trans \prn{\betaone + \betatwo} + t}^2\,\mathrm{d} t},\\
  &\text{II} := \frac{1}{2} \Expectlet_n \brk[\bigg]{
    \int_{-\bp_i\trans\valpha}^0 \Psione[\vq]{\vx_i}{\eta(\vp_{i}\trans \betazero[\vq])} \eta^{(2)}\prn{\xi_{\vq, i, t}}\\
  &\qquad \times \eta^{(1)}\prn{\vp_{i}\trans \prn{\betazero[\vq] + \betaone + \betatwo } + t} \prn{\vp_{i}\trans \prn{\betaone + \betatwo} + t}^2} \,\mathrm{d} t,\\
  &\text{III} := \frac{1}{2} \Expectlet_n \brk[\bigg]{
    \int_{-\bp_i\trans\valpha}^0 \Psitwo[\vq]{\vx_i}{\tilde{\xi}_{\vq, i, t}}
    \crl{\eta\prn{\vp_{i}\trans \prn{\betazero[\vq] + \betaone + \betatwo} + t} - \eta\prn{\vp_{i}\trans \betazero[\vq]}}^2\\
  &\qquad\times
    \eta^{(1)}\prn{\vp_{i}\trans \prn{\betazero[\vq] + \betaone + \betatwo} + t}\,\mathrm{d} t}
\end{align*}
and
\begin{align*}
  &\text{I} \lesssim 2^L \rtwon h^{d} \prn{\rone + 2^{\ell} \rtwon  + 2^L \rtwon}^2,\\
  &\text{II} \lesssim 2^L \rtwon h^{d} \prn{\rone + 2^{\ell} \rtwon  + 2^L \rtwon}^2,\\
  &\text{III} \lesssim 2^L \rtwon h^{d} \prn{\rone + 2^{\ell} \rtwon  + 2^L \rtwon}^2
\end{align*}
on the event $\eventi$.
\end{proof}

\begin{lemma}\label{lem:nonconvex-bahadur-delta-eta1-bound}
  There exists an event $\newtarget{def:nonconvexeventiii}{\nonconvexeventiiiNoLink}$ whose probability converges to one such that on $\nonconvexeventiii$
  \begin{align*}
    &\sup_{\vq \in \mathcal{Q}, \betaone \in \mathcal{H}_1, \betatwo \in \mathcal{H}_{2, \ell}, \valpha \in \mathcal{V}}
      \abs[\big]{\Expectlet_n \brk{\riskder[\vq]{y_i}{\eta\prn{\vp_{i}\trans \betazero[\vq]}}\\
    &\qquad \times \prn{\eta \prn{ \bp_i\trans \prn{ \betazero[\vq] + \betaone + \betatwo } } - \eta\prn{\vp_{i}\trans\prn{\betazero[\vq] + \betaone + \betatwo - \valpha}}}}\\
    &\qquad - \Expectlet_n \brk{\valpha\trans\vp_{i} \riskder[\vq]{y_i}{\eta\prn{\vp_{i}\trans \betazero[\vq]}} \eta^{(1)}\prn{\vp_{i}\trans\betazero[\vq]}}} \lesssim h^{d} \rone (\rone + 2^{\ell} \rtwon  + 2^L \rtwon) 2^L \rtwon.
  \end{align*}
\end{lemma}

\begin{proof}[Proof of \cref{lem:nonconvex-bahadur-delta-eta1-bound}]
By the Taylor expansion,
\begin{align*}
  &\abs[\big]{\Expectlet_n \brk{\riskder[\vq]{y_i}{\eta\prn{\vp_{i}\trans \betazero[\vq]}} \prn{\eta \prn{ \bp_i\trans \prn{ \betazero[\vq] + \betaone + \betatwo } } - \eta\prn{\vp_{i}\trans\prn{\betazero[\vq] + \betaone + \betatwo - \valpha}}}}\\
  &\qquad - \Expectlet_n \brk{\valpha\trans\vp_{i} \riskder[\vq]{y_i}{\eta\prn{\vp_{i}\trans \betazero[\vq]}} \eta^{(1)}\prn{\vp_{i}\trans\betazero[\vq]}}}\\
  &\quad = \abs[\big]{\Expectlet_n \brk[\big]{\bp_i\trans\valpha \riskder[\vq]{y_i}{\eta\prn{\vp_{i}\trans \betazero[\vq]}} \crl[\big]{\eta^{(2)} \prn{ \xi_{\vq, i} } \bp_i\trans \prn{ \betaone + \betatwo - \valpha } + (1 / 2) \eta^{(2)} \prn{ \tilde{\xi}_{\vq, i} } \bp_i\trans\valpha}}}\\
  &\quad \lesssim h^{d} \rone (\rone + 2^{\ell} \rtwon + 2^L \rtwon) 2^L \rtwon,
\end{align*}
where $\xi_{\vq, i}$ is between $\vp_{i}\trans(\betazero[\vq] + \betaone + \betatwo - \valpha)$ and $\vp_{i}\trans\betazero[\vq]$, $\tilde{\xi}_{\vq, i}$ between $\vp_{i}\trans(\betazero[\vq] + \betaone + \betatwo)$ and $\vp_{i}\trans(\betazero[\vq] + \betaone + \betatwo - \valpha)$. The bound holds on the event $\nonconvexeventiii := \nonconvexeventiiiprime \cap \nonconvexeventiiidoubleprime$, where
\begin{equation*}
  \nonconvexeventiiiprime := \crl[\Big]{
    \sup_{\betaone \in \mathcal{H}_1, \betatwo \in \mathcal{H}_{2, \ell}, \valpha \in \mathcal{V}, \vq \in \mathcal{Q}} \norm[\big]{
      \Expectlet_n \brk{\riskder[\vq]{y_i}{\eta\prn{\vp_{i}\trans \betazero[\vq]}} \eta^{(2)}(\xi_{\vq, i}) \vp_{i} \vp_{i}\trans}
    }_{\infty} \leq h^{d} \rone
  },
\end{equation*}
and $\nonconvexeventiiidoubleprime$ is defined the same way as $\nonconvexeventiiiprime$ with $\xi_{\vq, i}$ replaced by $\tilde{\xi}_{\vq, i}$.

By the same argument as \cref{lem:holder-uniform-convergence-variance}, $\Prob{\nonconvexeventiii} \to 1$.
\end{proof}

\subsection{Rates of convergence}

\begin{corollary}[Uniform rate of convergence]\label{cor:uniform-rate-of-convergence}\leavevmode
  \begin{assertions}[wide]
  \item\label{th:rates-of-conv-i} If the conditions of \cref{th:bahadur-repres}\ref{th:bahadur-repres-normal} hold, then
  \begin{equation}\label{eq:rates-lemma-uniform-rate}
    \sup _{\vq \in \mathcal{Q}, \vx \in \mathcal{X}} \abs[\big]{\muhat[\vq][(\vv)]{\vx} - \muzero[\vq][(\vv)]{\vx}}
    \lesssim_{\Problet} h^{-\multindnorm{\vv}} \brk[\bigg]{\prn[\bigg]{\frac{\log^d n}{n h^d}}^{1 / 2} \log n + h^m}
  \end{equation}

  \item If the conditions of \cref{th:bahadur-repres}\ref{th:bahadur-repres-normal} hold and
  \begin{equation}\label{eq:rates-lemma-additional-rate-restriction}
    \brk{\log n}^{(d + 1) / \prn{ \alpha \wedge \sfrac{1}{2} } + d} = O\prn{n h^d},\quad h^{\prn{ \alpha \wedge \sfrac{1}{2} } m} \log^{d / 2} n = O(1),
  \end{equation}
  then
  \begin{equation}\label{eq:rates-lemma-uniform-rate-modified}
    \sup _{\vq \in \mathcal{Q}, \vx \in \mathcal{X}} \abs[\big]{\muhat[\vq][(\vv)]{\vx} - \muzero[\vq][(\vv)]{\vx}} \lesssim_{\Problet} h^{-\multindnorm{\vv}} \brk[\bigg]{\prn[\bigg]{ \frac{\log n}{n h^d} }^{1 / 2} + h^m}.
  \end{equation}

\item If the conditions of \cref{th:bahadur-repres}\ref{th:bahadur-repres-cont} hold and
    \begin{equation*}
    \brk{\log n}^{(d + 1) / \alpha + d} = O\prn{n h^d},\quad h^{\alpha m} \log^{d / 2} n = O(1),
  \end{equation*}
  then~\eqref{eq:rates-lemma-uniform-rate-modified} is also true.
  \end{assertions}
\end{corollary}

\begin{proof}
By \cref{th:bahadur-repres}, \cref{lem:holder-uniform-convergence-variance} and triangle inequality,
\begin{align*}
  &\sup _{\vq \in \mathcal{Q}, \vx \in \mathcal{X}} \abs[\big]{\vp^{(\vv)}(\vx)\trans \betahat[\vq] - \muzero[\vq][(\vv)]{\vx}}\\
  &\quad \lesssim_{\Problet} h^{-\multindnorm{\vv}} \brk[\bigg]{\prn[\bigg]{ \frac{\log(1 / h)}{n h^d} }^{1 / 2} + \prn[\bigg]{ \frac{\log^d n}{n h^d} }^{1 / 2 + \sfrac{\alpha}{2} \wedge \sfrac{1}{4}} \log n + h^{\prn{\alpha \wedge \sfrac{1}{2}} m}\prn[\bigg]{\frac{\log^{d + 1} n}{n h^d}}^{1 / 2} + h^{m}}.
\end{align*}
Using $\log(1 / h) \lesssim \log n$ and simplifying the right-hand side, we obtain~\eqref{eq:rates-lemma-uniform-rate}.
Additional restrictions~\eqref{eq:rates-lemma-additional-rate-restriction} allow us to get a slightly stronger result~\eqref{eq:rates-lemma-uniform-rate-modified}.
\end{proof}

\begin{corollary}[Mean square rate of convergence]\label{cor:msq-rate-of-convergence}\leavevmode
  \begin{assertions}[wide]
  \item If the conditions of \cref{th:bahadur-repres}\ref{th:bahadur-repres-normal} hold and
    \begin{equation}\label{eq:rates-lemma-additional-rate-restriction-for-l2}
      \brk{\log n}^{(d + 2) / \prn{ \alpha \wedge \sfrac{1}{2} } + d} = o\prn{n h^d},\quad h^{\prn{ \alpha \wedge \sfrac{1}{2} } m} \log^{(d + 1) / 2} n = o(1),
    \end{equation}
    then
    \begin{equation}\label{eq:rates-lemma-l2-rate-modified}
      \sup_{\vq \in \mathcal{Q}} \prn[\Big]{\int_{\mathcal{X}} \abs[\big]{\vp^{(\vv)}(\vx)\trans \betahat[\vq] - \muzero[\vq][(\vv)]{\vx}}^2 \fx{\vx}\,\mathrm{d} \vx }^{1 / 2} \lesssim_{\Problet} h^{-\multindnorm{\vv}} \brk[\bigg]{\frac{1}{\sqrt{n h^d}} + h^m}.
    \end{equation}

  \item If the conditions of \cref{th:bahadur-repres}\ref{th:bahadur-repres-cont} hold and
    \begin{equation*}
      \brk{\log n}^{(d + 2) / \alpha + d} = o\prn{n h^d},\quad h^{\alpha m} \log^{(d + 1) / 2} n = o(1),
    \end{equation*}
    then~\eqref{eq:rates-lemma-l2-rate-modified} is also true.
  \end{assertions}
\end{corollary}

\begin{proof}
To prove \eqref{eq:rates-lemma-l2-rate-modified}, note that
\begin{align*}
  &\sup_{\vq} \int_{\mathcal{X}} \abs[\big]{\vp^{(\vv)}(\vx)\trans \betahat[\vq] - \vp^{(\vv)}(\vx)\trans \betazero[\vq]}^2 \fx{\vx}\,\mathrm{d} \vx\\
  &\quad = \sup_{\vq} \prn[\big]{ \betahat[\vq] - \betazero[\vq] }\trans \Expect{\vp^{(\vv)}(\vx) \vp^{(\vv)}(\vx)\trans} \prn[\big]{ \betahat[\vq] - \betazero[\vq] } \labrel{1}[\lesssim] h^{d - 2 \multindnorm{\vv}} \norm[\big]{\betahat[\vq] - \betazero[\vq]}^2,
\end{align*}
where inequality \labrel{1} is true because the largest eigenvalue of $\Expect{\vp^{(\vv)}(\vx_i) \vp^{(\vv)}(\vx_i)\trans}$ is bounded from above by $h^{d - 2 \multindnorm{\vv}}$ up to a multiplicative coefficient, which is proven by the same argument as for $\vv = 0$ in \cref{lem:holder-gram} in combination with \cref{ass:gl-local-basis}.

It is left to prove
\begin{equation}\label{eq:l2-beta-diff-norm-bound}
  \norm[\big]{\betahat[\vq] - \betazero[\vq]} \lesssim_{\Problet} \frac{1}{h^d \sqrt{n}}.
\end{equation}
By the triangle inequality,
\begin{equation}\label{eq:l2-rates-diff-2norm-bound}
  \begin{split}
    &\norm[\big]{\betahat[\vq] - \betazero[\vq]} \leq \norm[\big]{\betahat[\vq] - \betazero[\vq] + \qzero[\vq][-1] \Expectlet_n \brk{\vp(\vx_i)  \eta^{(1)}\prn{\muzero[\vq]{\vx_i}} \riskder[\vq]{y_i}{\eta\prn{\muzero[\vq]{\vx_i}}}}}\\
    &\quad + \norm[\big]{\qzero[\vq][-1] \Expectlet_n \brk{\vp(\vx_i) \eta^{(1)}\prn{\muzero[\vq]{\vx_i}} \riskder[\vq]{y_i}{\eta\prn{\muzero[\vq]{\vx_i}}}}}.
  \end{split}
\end{equation}
To bound the second term in~\eqref{eq:l2-rates-diff-2norm-bound} on the right-hand side, consider the expectation
\begin{align*}
  &\Expect[\big]{\norm[\big]{\Expectlet_n \brk[\big]{\vp(\vx_i) \eta^{(1)}\prn{\muzero[\vq]{\vx_i}} \riskder[\vq]{y_i}{\eta\prn{\muzero[\vq]{\vx_i}}}}}^2}\\
  &\quad = \frac{1}{n^2} \sum_{i = 1}^n \Expect[\big]{\eta^{(1)}\prn{\muzero[\vq]{\vx_i}}^2 \riskder[\vq]{y_i}{\eta\prn{\muzero[\vq]{\vx_i}}}^2 \norm{\vp(\vx_i)}^2}\\
  &\quad = \frac{1}{n^2} \sum_{i = 1}^n \Expect[\big]{\sigma_{\vq}^2(\vx_i) \eta^{(1)}\prn{\muzero[\vq]{\vx_i}}^2 \norm{\vp(\vx_i)}^2} \overset{\text{(a)}}{\lesssim} \frac{1}{n^2} \sum_{i = 1}^n 1 = \frac{1}{n},
\end{align*}
where in (a) we used uniform boundedness of~$\sigma_{\vq}^2(\vx)$ by \cref{ass:gl-holder}\ref{ass:psi-moments}, uniform boundedness of $\muzero[\vq]{\vx}$ and $\norm{\vp(\vx)}$. By Markov's inequality and \cref{lem:holder-gram}, this immediately implies
\begin{equation*}
  \norm[\big]{\qzero[\vq][-1] \Expectlet_n \brk[\big]{\vp(\vx_i) \eta^{(1)}\prn{\muzero[\vq]{\vx_i}} \riskder[\vq]{y_i}{\eta\prn{\muzero[\vq]{\vx_i}}}}} \lesssim_{\Problet} \norm{\qzero[\vq][-1]} \cdot \frac{1}{\sqrt{n}} \lesssim \frac{1}{h^d \sqrt{n}}.
\end{equation*}

Concerning the first term in~\eqref{eq:l2-rates-diff-2norm-bound}, it is proven in \cref{th:bahadur-repres} (see \cref{eq:UMbHrv}) that
\begin{equation*}
  \norm[\big]{\betahat[\vq] - \betazero[\vq] + \qzero[\vq][-1] \Expectlet_n \brk[\big]{\vp(\vx_i) \eta^{(1)}\prn{\muzero[\vq]{\vx_i}} \riskder[\vq]{y_i}{\eta\prn{\muzero[\vq]{\vx_i}}}}}_{\infty} \lesssim_{\Problet} \rbah
\end{equation*}
for $\rbah$ defined in~\eqref{eq:bahadur-rates-for-statement}, so that
\begin{align*}
  &\norm[\big]{\betahat[\vq] - \betazero[\vq] + \qzero[\vq][-1] \Expectlet_n \brk[\big]{\vp(\vx_i) \eta^{(1)}\prn{\muzero[\vq]{\vx_i}} \riskder[\vq]{y_i}{\eta\prn{\muzero[\vq]{\vx_i}}}}}\\
  &\quad \lesssim_{\Problet} \sqrt{K} \rbah \lesssim \frac{\rbah}{h^{d / 2}} \overset{\text{(a)}}{=} o \prn[\bigg]{ \frac{1}{h^d \sqrt{n}} },
\end{align*}
where in equality (a) we used $\rbah = o \prn[\big]{ 1 / \sqrt{n h^d}}$ under the assumptions. This concludes the proof of~\eqref{eq:l2-beta-diff-norm-bound}.
\end{proof}

\section{Strong approximation}\label{sec:strong-approximation}

\Cref{sec:conditional-yurinskii-coupling} collects results that may be of independent theoretical interest.
Specifically, we use the conditional Strassen's theorem \cref{th:cond-strassen}
to prove conditional Yurinskii coupling in $d$-norm, \cref{th:conditional-yurinski}. We then use it to prove \cref{lem:yurinski-coupling-K-dimensional}, a version of Gaussian approximation for a $K$-dimensional empirical process. \Cref{sec:SA-main-result} proves the main result of this \lcnamecref{sec:strong-approximation}, \cref{th:strong-approximation-binscatter-yurinski}, which is essentially a corollary of \cref{lem:yurinski-coupling-K-dimensional} after some additional technical work done in
\cref{lem:asympt-variance,lem:closeness-of-linear-terms,lem:hats-off}.

\subsection{Yurinskii coupling}\label{sec:conditional-yurinskii-coupling}

The three theorems, and their proofs, in this subsection are self-contained, and hence all variables, functions, and stochastic processes, should be treated as defined within each of the theorems and their proofs, and independently of all other statements elsewhere in the supplemental appendix.

The following \lcnamecref{th:cond-strassen} is due to \parencite{monrad1991nearby}. We use the statement from \parencite{chen2020jackknife} making it explicit that the supremum over Borel sets may not be a random variable (a direct proof may also be found in that work). Let $(S, d)$ be a Polish space (where $d$ is its metric), and $\mathcal{B}(S)$ its Borel sigma-algebra.

\begin{theorem}[Conditional Strassen's theorem]\label{th:cond-strassen}
  Let $X$ be a random variable defined on some probability space $(\Omega, \mathcal{F}, \Problet)$ and with values in some Polish $(S, d)$. Let $\mathcal{J}$ be a countably generated sub-sigma algebra of $\mathcal{F}$ and assume that this probability space is rich enough: there exists a random variable $U$ that is independent of the sigma-algebra $\mathcal{J} \vee \sigma(X)$. Let $\mathcal{B}(S) \times \Omega \ni (A, \omega) \mapsto G(A| \mathcal{J})(\omega)$ be a regular conditional distribution on $\mathcal{B}(S)$, i.\,e., for each $A \in \mathcal{B}(S)$, $G(A| \mathcal{J})$ is $\mathcal{J}$-measurable, and for each $\omega \in \Omega$, $G(\cdot|\mathcal{J})(\omega)$ is a probability measure on $\mathcal{B}(S)$. Suppose that for some nonnegative numbers $\alpha$ and $\beta$
  \begin{equation*}
    \Expectlet^{*} \brk[\Big]{\sup_{A \in \mathcal{B}} \crl{\Prob{X \in A \given \mathcal{J}} - G(A^{\alpha} | \mathcal{J})}} \leq \beta,
  \end{equation*}
  where $\Expectlet^{*}$ denotes outer expectation. Then on this probability space there exists an $S$-valued random element $Y$ such that $G(\cdot | \mathcal{J})$ is its regular conditional distribution given $\mathcal{J}$ and $\Prob{d(X, Y) > \alpha} \leq \beta$.
\end{theorem}

The following \lcnamecref{th:conditional-yurinski} is a conditional version of Lemma~39 in \parencite{Belloni-Chernozhukov-Chetverikov-FernandezVal_2019_JoE}; its proof carefully leverages \cref{th:cond-strassen}. See also \parencite{cattaneo2024yurinskii} for a related, but different, conditional Yurinskii's coupling result.

\begin{theorem}[Conditional Yurinskii coupling, \texorpdfstring{$d$}{d}-norm]\label{th:conditional-yurinski}
  Let a random vector $\bC$ with values in $\mathbb{R}^m$, sequence of random vectors $\crl{\bxi_i}_{i = 1}^n$ with values in $\mathbb{R}^k$ and sequence of random vectors $\crl{\bg_i}_{i = 1}^n$ with values in $\mathbb{R}^k$ be defined on the same probability space and be such that all the $2 n$ vectors $\crl{\bxi_i}_{i = 1}^n \cup \crl{\bg_i}_{i = 1}^n$ are independent conditionally on $\bC$, are mean zero conditionally on $\bC$ and for each $i \in \crl{1, \ldots, n}$ the distribution of $\bg_i$ conditionally on $\bC$ is $\mathcal{N}\prn{0, \Var{\bxi_i \given \bC}}$. Assume that this probability space is rich enough: there exists a random variable $U$ that is independent of $\crl{\bxi_i}_{i = 1}^n \cup \crl{\bg_i}_{i = 1}^n$. Denote
  \begin{equation*}
    \bS := \bxi_1 + \ldots + \bxi_n,\quad \bT := \bg_1 + \ldots + \bg_n,
  \end{equation*}
  and let
  \begin{equation*}
    \beta := \sum_{i = 1}^n \Expect[\big]{\norm{\bxi_i}^2 \norm{\bxi_i}_d} + \sum_{i = 1}^n \Expect[\big]{\norm{\bg_i}^2 \norm{\bg_i}_d}
  \end{equation*}
  be finite. Then for each $\delta > 0$, on this probability space there exists a random vector $\bT'$ such that $\Problet_{\bT | \sigma(\bC)}(\cdot)$ is its regular conditional distribution given $\sigma(\bC)$, and
  \begin{equation*}
    \Prob[\big]{\norm{\bS - \bT'}_d > 3 \delta} \leq \min_{t \geq 0} \prn[\bigg]{ 2 \Prob{\norm{\bZ}_d > t} + \frac{\beta}{\delta^3} t^2 },
  \end{equation*}
  where $\bZ \sim \mathcal{N}(0, I_k)$.
\end{theorem}

\begin{proof}
By the conditional Strassen's theorem, it is enough to show
\begin{equation*}
  \Expectlet^{*} \brk[\Big]{\sup_{A \in \mathcal{B}\prn{\mathbb{R}^k}} \prn[\big]{\Prob{\bS \in A \given \sigma(\bC)} - \Problet_{\bT | \sigma(\bC)}\prn{A^{3 \delta}}}} \leq \min_{t \geq 0} \prn[\bigg]{ 2 \Prob{\norm{\bZ}_d > t} + \frac{\beta}{\delta^3} t^2},
\end{equation*}
or, equivalently, for any $t > 0$
\begin{equation}\label{eq:cond-yur-strassen-enough}
  \Expectlet^{*} \brk[\Big]{\sup_{A \in \mathcal{B}\prn{\mathbb{R}^k}} \prn[\big]{\Prob{\bS \in A \given \bC} - \Prob{\bT \in A^{3 \delta} \given \bC}}} \leq 2 \Prob{\norm{bZ}_d > t} + \frac{\beta}{\delta^3} t^2.
\end{equation}

Fix $t > 0$ and $A \in \mathcal{B}\prn{\mathbb{R}^k}$. Let $f\colon \mathbb{R}^k \to \mathbb{R}$ be the same as in the proof of Lemma~39 in~\parencite{Belloni-Chernozhukov-Chetverikov-FernandezVal_2019_JoE}, namely, it is such that for all $\vx, \vy \in \mathbb{R}^k$,
\begin{align*}
  & \abs[\big]{f(\vx + \vy) - f(\vx) - \vy\trans \nabla f(\vx) - (1 / 2) \vy\trans \nabla^2 f(\vx) \vy} \leq \frac{\|\vy\|^2\|\vy\|_d}{\sigma^2 \delta},\\
  & (1-\epsilon) 1\{\vx \in A\} \leq f(\vx) \leq \epsilon+(1-\epsilon) 1\crl{\vx \in A^{3 \delta}},
\end{align*}
with $\sigma := \delta / t$ and $\epsilon := \Prob*{\norm{\bZ}_d > t}$.

Then note that
\begin{align}
  &\Prob{\bS \in A \given \bC} = \Expect[\big]{\Indlet\{\bS \in A\} - f(\bS) \given \bC} + \Expect[\big]{f(\bS) - f(\bT) \given \bC} + \Expect*{f(\bT) \given \bC}\nonumber\\
  &\quad \leq \epsilon \Expect[\big]{\Indlet\{\bS \in A\} \given \bC} + \Expect[\big]{f(\bS) - f(\bT) \given \bC} + \epsilon + (1 - \epsilon) \Expect[\big]{\Indlet\crl{\bT \in A^{3 \delta}} \given \bC}\nonumber\\
  &\quad \leq 2 \epsilon + \Expect[\big]{\Indlet\crl{\bT \in A^{3 \delta}} \given \bC} + \Expect[\big]{f(\bS) - f(\bT) \given \bC}.\label{eq:cond-yur-rhs-with-f-dif}
\end{align}

Now we bound $\Expect{f(\bS) - f(\bT) \given \bC}$:
\begin{align*}
  &\Expect{f(\bS) - f(\bT) \given \bC} = \sum_{i = 1}^n \Expect{f\prn{\bX_i + \bY_i} - f\prn{\bX_i + \bW_i} \given \bC}\\
  &\quad\leq \sum_{i = 1}^n \Expect[\bigg]{f \prn{ \bX_i } + \bY_i\trans \nabla f \prn{ \bX_i } + \frac{1}{2} \bY_i\trans \nabla^2 f \prn{ \bX_i } \bY_i + \frac{\norm{ \bY_i }^2 \norm{ \bY_i }_d}{\sigma^2 \delta} \given \bC}\\
  &\qquad - \sum_{i = 1}^n \Expect[\bigg]{f \prn{ \bX_i } + \bW_i\trans \nabla f \prn{ \bX_i } + \frac{1}{2} \bW_i\trans \nabla^2 f \prn{ \bX_i } \bW_i - \frac{\norm{ \bW_i }^2 \norm{ \bW_i }_d}{\sigma^2 \delta} \given \bC}\\
  &\quad \labrel{1}[=] \sum_{i = 1}^n \Expect{\frac{\norm{ \bY_i }^2 \norm{ \bY_i }_d + \norm{ \bW_i }^2 \norm{ \bW_i }_d}{\sigma^2 \delta} \given \bC}\text{ a.\,s.}
\end{align*}
for $\bX_i := \bxi_1 + \ldots + \bxi_{i - 1} + \bg_{i + 1} + \ldots + \bg_n$, $\bY_i := \bxi_i$, $\bW_i := \bg_i$, in \labrel{1} we used the conditional independence of the family $\crl{ \bxi_i }_{i = 1}^n \cup \crl{ \bg_i }_{i = 1}^n$, that they are conditionally mean zero and the equality of the corresponding conditional second moments.

We conclude that almost surely
\begin{equation*}
  \sup_{A \in \mathcal{B} \prn{ \mathbb{R}^k }}\prn[\big]{\Prob{\bS \in A \given \bC} - \Prob{\bT \in A^{3 \delta} \given \bC}} \leq 2 \epsilon + \sum_{i = 1}^n \Expect[\bigg]{\frac{\norm{ \bxi_i }^2 \norm{ \bxi_i}_d + \norm{ \bg_i }^2 \norm{ \bg_i }_d}{\sigma^2 \delta} \given \bC}.
\end{equation*}
By the definition of outer expectation, this implies
\begin{equation*}
  \Expectlet^{*} \brk[\Big]{\sup_{A \in \mathcal{B} \prn{ \mathbb{R}^k }}\prn[\big]{\Prob{\bS \in A \given \bC} - \Prob{\bT \in A^{3 \delta} \given \bC}}} \leq 2 \epsilon + \frac{\beta}{\sigma^2 \delta},
\end{equation*}
which is~\eqref{eq:cond-yur-strassen-enough}.
\end{proof}

The following \lcnamecref{lem:yurinski-coupling-K-dimensional} generalizes Lemma 36 in \parencite{Belloni-Chernozhukov-Chetverikov-FernandezVal_2019_JoE}, and also builds on the argument for Lemma SA27 in \parencite{Cattaneo-Feng-Underwood_2024_JASA}.

\begin{lemma}[Yurinskii coupling: \texorpdfstring{$K$}{K}-dimensional process]\label{lem:yurinski-coupling-K-dimensional}
  Let $\crl{ \vx_i, y_i }_{i = 1}^n$ be a random sample, where $\vx_i$ has compact support $\mathcal{X} \subset \mathbb{R}^d$, $y_i \in \mathcal{Y} \subset \mathbb{R}$ is a scalar. Also let $\mathcal{Q} \subseteq \mathbb{R}^{d_{\mathcal{Q}}}$ be a fixed compact set.

  Let $A_n\colon \mathcal{Q} \times \mathcal{X} \times \mathcal{Y} \to \mathbb{R}$ be a Borel measurable function satisfying $\sup_{\vq \in \mathcal{Q}}\abs{A_n(\vq, \vx_i, y_i)} \leq \bar{A}_n(\vx_i, y_i)$, where $\bar{A}_n(\vx_i, y_i)$ is a Borel measurable envelope, $\Expect{A_n(\vq, \vx_i, y_i) \given \vx_i} = 0$ for all $\vq \in \mathcal{Q}$, $\sup_{\vx \in \mathcal{X}}\Expect{|\bar{A}(\vx_i, y_i)|^{\nu} \given \vx_i = \vx} \leq \mu_n < \infty$ for some $\nu \geq 3$ with $\mu_n \gtrsim 1$ and $\log \mu_n \lesssim \log n$, which satisfies the Lipschitz condition
  \begin{equation*}
    \sup_{\vx \in \mathcal{X}} \Expect[\big]{\abs{
        A_n\prn{\vq, \vx_i, y_i} - A_n\prn{\tilde{\vq}, \vx_i, y_i}
      }^2 \given \vx_i = \vx} \lesssim \| \vq - \tilde{\vq} \|
  \end{equation*}
  for all $\vq, \tilde{\vq} \in \mathcal{Q}$. Also, the (regular) conditional variance $\Expect{A_n\prn{\vq, \vx_i, y_i}^2 \given \vx_i = \vx}$ is continuous in $\vx \in \mathcal{X}$. Moreover, assume that the class of functions $\crl{ (\vx, y) \mapsto A_n(\vq, \vx, y)\colon \vq \in \mathcal{Q} }$ is VC-subgraph with an index bounded from above by a constant not depending on $n$.

  Let $\bb(\cdot)$ be a Borel measurable function $\mathcal{X} \to \mathbb{R}^K$ (where $K = K_n$ is some sequence of positive integers tending to infinity and satisfying $\log K \lesssim \log n$) such that $\sup_{\vx \in \mathcal{X}} \norm{ \bb(\vx) } \leq \zeta_K$ and the probability of the event $\mathcal{A} := \crl{ \sup_{\norm{ \balpha } = 1} \Expectlet_n \brk[\big]{\prn{ \balpha\trans \bb(\vx_i) }^2} \leq C_{\mathrm{Gr}} }$ approaches one, where $C_{\mathrm{Gr}}$ is some positive constant. Assume $\zeta_K$ satisfies $1 / \zeta_K \lesssim 1$, $\abs{\log \zeta_K} \lesssim \log n$.

  Let $r_{n, \mathrm{yur}} = \newtarget{def:ryur-in-yurinski-lemma}{\ryurNoLink} \to 0$ be a sequence of positive numbers satisfying
  \begin{equation}\label{eq:yurinski-lemma-rate-restriction}
    \prn[\bigg]{ \frac{\zeta^3_K \mu_n^{3 / \nu}}{\sqrt{n}} }^{\frac{1}{3 + 2 d_{\mathcal{Q}}}} \sqrt{\log n} + \frac{\zeta_K \mu_n^{1 / \nu}}{n^{1 / 2 - 1 / \nu}} \log n = o \prn{ \ryur }.
  \end{equation}

  Assume also that the probability space is rich enough, and denote
\begin{equation*}
\bG_n(\vq) := \sqrt{n} \Expectlet_n\brk{A_n(\vq, \vx_i, y_i) \bb(\vx_i)}.
\end{equation*}
Then, on the same probability space, there exists a $K$-dimensional, conditionally on $\mX_n$ mean-zero Gaussian, process $\bZ_n(\vq)$ on $\mathcal{Q}$ with a.\,s. continuous trajectories, satisfying
\begin{gather*}
  \Expect[\big]{\bG_n(\vq) \bG_n(\tilde{\vq})\trans \given \mX_n} = \Expect[\big]{\bZ_n(\vq) \bZ_n(\tilde{\vq})\trans \given \mX_n},\quad \vq, \tilde{\vq} \in \mathcal{Q}; \\
  \sup_{\vq \in \mathcal{Q}} \norm{ \bG_n(\vq) - \bZ_n(\vq) }_{\infty} = o_{\Problet}(\ryur).
\end{gather*}

  Moreover, if $\bar{A}_n(\vq, \vx_i, y_i)$ is \subgauss, then \eqref{eq:yurinski-lemma-rate-restriction} can be replaced with
  \begin{equation*}
    \prn[\bigg]{ \frac{\zeta^3_K}{\sqrt{n}} }^{\frac{1}{3 + 2 d_{\mathcal{Q}}}} \sqrt{\log n} + \frac{\zeta_K}{\sqrt{n}} \log^{3 / 2} n = o(\ryur).
  \end{equation*}
\end{lemma}

\begin{proof}

\newcommand{\MDimYurinskiApproximationDelta}{a_n^{-2} \left( \log n \right)^{-1 - \gamma}}

Let $\mathcal{Q}_n^{\delta} := \crl{\vq_1, \ldots, \vq_{|\mathcal{Q}_n^{\delta}|}}$ be an internal $\delta_n$-covering of $\mathcal{Q}$ with respect to the 2-norm $\|\cdot\|$ of cardinality $|\mathcal{Q}_n^{\delta}| \lesssim 1 / \delta_n^{d_{\mathcal{Q}}}$, where $\delta_n$ is chosen later. Denote $\pi_n^{\delta} \colon \mathcal{Q} \to \mathcal{Q}$ a sequence of projections associated with this covering: it maps each point in $\mathcal{Q}$ to the center of the ball containing this point (if such a ball is not unique, choose one by an arbitrary rule).

\paragraph{Strategy}
The plan of attack is to
\begin{enumerate}
\item show that $\bG_n(\vq)$ does not deviate too much in $\sup$-norm from its projected version, i.\,e. bound the tails of $\sup_{\vq \in \mathcal{Q}} \norm{ \bG_n(\vq) - \bG_n\circ\pi_n^{\delta}(\vq) }_{\infty}$,\sloppy
\item apply Yurinskii coupling to the finite-dimensional vector $\prn[\big]{\bG_n\circ\pi_n^{\delta}(\vq)}_{\vq \in \mathcal{Q}_n^{\delta}}$ and obtain a conditionally on $\mX_n$ Gaussian vector $\bZ_n^{\delta}$ with the right structure that is close enough, i.\,e. with a bound on the tails of $\norm{ \bG_n\circ\pi_n^{\delta} - \bZ_n^{\delta} }_{\infty}$,
\item extend this conditionally Gaussian vector to a $K$-dimensional conditionally Gaussian process $\bZ_n$,
\item and finally show that $\bZ_n(\vq)$ does not deviate too much from its projected version, i.\,e. bound the tails of $\sup_{\vq \in \mathcal{Q}} \norm{ \bZ_n(\vq) - \bZ_n \circ \pi_n^{\delta}(\vq) }_{\infty}$.
\end{enumerate}

If we complete these steps, it will prove the theorem by the triangle inequality.

\paragraph{Discretization of $G_n$} Consider the class of functions
\begin{equation*}
  \mathcal{G}_n' := \crl{ \mathcal{X} \times \mathcal{Y} \ni (\vx, y) \mapsto A_n(\vq, \vx, y) b_l(\vx)\colon 1 \leq l \leq K, \vq \in \mathcal{Q}}
\end{equation*}
with envelope $\zeta_K \bar{A}_n(X)$. Since $\crl{A_n(\vq, \vx, y)}$ is a VC class with $O(1)$ index and envelope $\bar{A}_n(\vx, y)$, $\mathcal{G}_n'$ satisfies the uniform entropy bound \eqref{eq:uniform-entropy-bound} with $A \lesssim K$ and $V \lesssim 1$.

Next, consider the class of functions
\begin{equation*}
  \mathcal{G}_n^{\delta} := \crl{ \mathcal{X} \times \mathcal{Y} \ni (\vx, y) \mapsto \prn[\big]{A_n(\vq, \vx, y) - A_n(\tilde{\vq}, \vx, y)} b_l(\vx) \colon 1 \leq l \leq K, \vq, \tilde{\vq} \in \mathcal{Q}, \| \vq - \tilde{\vq} \| \lesssim \delta_n }
\end{equation*}
with envelope $(\vx, y) \mapsto 2 \zeta_K \bar{A}_n(\vx, y)$. Using \cref{lem:permanence-properties}, we get that this class satisfies the uniform entropy bound~\eqref{eq:uniform-entropy-bound} with $A \lesssim K$ and $V \lesssim 1$.

Now we apply \cref{lem:ueb-maximal-ineq} conditionally on $\mX_n$ on $\mathcal{A}$ with $\|F\|_{\bar{\Problet}, 2} \leq 2 \zeta_K \mu_n^{1 / \nu}$ since
\begin{equation*}
  \|\bar{A}_n(\vx, y)\|_{\bar{\Problet}, 2}^2 = \frac{1}{n} \sum_{i = 1}^n \Expect{\bar{A}_n(\vx_i, y_i)^2 \given \vx_i} \leq \frac{1}{n} \sum_{i = 1}^n \mu_n^{2 / \nu} = \mu_n^{2 / \nu},
\end{equation*}
$\|M\|_{\Problet, 2} \leq 2 \zeta_K \prn{\mu_n n}^{1 / \nu}$ since
\begin{align*}
  &\Expect[\Big]{\prn[\Big]{\max_{1 \leq i \leq n} \bar{A}_n(\vx_i, y_i)}^2 \given \mX_n} \leq \Expect[\Big]{\prn[\Big]{\max_{1 \leq i \leq n} \bar{A}_n(\vx_i, y_i)}^\nu \given \mX_n}^{2 / \nu}\\
  &\quad \leq \Expect[\bigg]{\sum_{i = 1}^n \bar{A}_n(\vx_i, y_i)^{\nu} \given \mX_n}^{2 / \nu} \leq (\mu_n n)^{2 / \nu},
\end{align*}
and $\sigma \lesssim \sqrt{\delta_n}$ since on $\mathcal{A}$
\begin{align*}
  &\frac{1}{n} \sum_{i = 1}^n \Expect[\big]{\prn[\big]{A_n(\vq, \vx_i, y_i) - A_n(\tilde{\vq}, \vx_i, y_i)}^2  b_{l}(\vx_i)^{2} \given \mX_n}\\
  &\quad = \frac{1}{n} \sum_{i = 1}^n b_{l}(\vx_i)^{2} \Expect[\big]{\prn[\big]{A_n(\vq, \vx_i, y_i) - A_n(\tilde{\vq}, \vx_i, y_i)}^2 \given \mX_n}\\
  &\quad \lesssim \|\vq - \tilde{\vq}\| \frac{1}{n} \sum_{i = 1}^n b_{l}(\vx_i)^2 \lesssim \|\vq - \tilde{\vq}\|.
\end{align*}
This gives that on $\mathcal{A}$
\begin{align*}
  &\Expect[\Big]{\sup_{\|\vq - \tilde{\vq}\| \leq \delta_n} \norm{ \bG_n(\vq) - \bG_n(\tilde{\vq}) }_{\infty} \given \mX_n}\\
  &\quad \lesssim \sqrt{ \delta_n \log \frac{K \zeta_K \mu_n^{1 / \nu}}{\sqrt{\delta_n}}} + \frac{\zeta_K \prn{\mu_n n}^{1 / \nu}}{\sqrt{n}} \log \frac{K \zeta_K \mu_n^{1 / \nu}}{\sqrt{\delta_n}}.
\end{align*}

By Markov's inequality and since $\Prob*{\mathcal{A}} \to 1$, for any sequence $t_n > 0$
\begin{equation}\label{eq:discretization-gn-final-bound}
  \begin{aligned}
    &\Prob*{\sup_{\|\vq - \tilde{\vq}\| \leq \delta_n} \norm{ \bG_n(\vq) - \bG_n(\tilde{\vq}) }_{\infty} > t_n}\\
    &\quad \lesssim \frac{1}{t_n} \sqrt{ \delta_n \log \frac{K \zeta_K \mu_n^{1 / \nu}}{\sqrt{\delta_n}}} + \frac{\zeta_K \prn{\mu_n n}^{1 / \nu}}{t_n \sqrt{n}} \log \frac{K \zeta_K \mu_n^{1 / \nu}}{\sqrt{\delta_n}} + o(1),
    \end{aligned}
\end{equation}
where the constant in $\lesssim$ does not depend on $n$.

\paragraph{Coupling}
Define a $K |\mathcal{Q}_n^{\delta}|$-dimensional vector $\bxi_i := \prn[\big]{A_n(\vq, \vx_i, y_i) b_l(\vx_i) / \sqrt{n}}_{1 \leq l \leq K, \vq \in \mathcal{Q}_n^{\delta}}$, so that we have $\bG_n \circ \pi_n^{\delta} = \sum_{i = 1}^n \bxi_i$. We make some preparations before applying~\cref{th:conditional-yurinski}.

Firstly, we bound $\Expect[\big]{\norm{ \bxi_i }^2 \norm{ \bxi_i }_{\infty} \given \mX_n}$.
\begin{align*}
  &\sum_{i = 1}^n \Expect[\big]{\norm{ \bxi_i }^2 \norm{ \bxi_i }_{\infty} \given \mX_n}\\
  &\quad = \frac{1}{n^{3/2}} \sum_{i = 1}^n \|\bb(\vx_i)\|^2 \|\bb(\vx_i)\|_{\infty} \Expect[\bigg]{
    \sum_{\vq \in \mathcal{Q}_n^{\delta}} A_n(\vq, \vx_i, y_i)^2 \cdot \max_{\vq \in \mathcal{Q}_n^{\delta}}\abs{ A_n(\vq, \vx_i, y_i) }
    \given \vx_i}\\
  &\quad \leq \frac{1}{n^{3/2}} |\mathcal{Q}_n^{\delta}| \sum_{i = 1}^n \|\bb(\vx_i)\|^2 \|\bb(\vx_i)\|_{\infty} \Expect[\big]{\bar{A}_n(\vx_i, y_i)^3 \given \vx_i}\\
  &\quad \leq \frac{|\mathcal{Q}_n^{\delta}| \mu_n^{3 / \nu}}{n^{3/2}} \sum_{i = 1}^n \|\bb(\vx_i)\|^2 \|\bb(\vx_i)\|_{\infty} \lesssim \frac{|\mathcal{Q}_n^{\delta}| \mu_n^{3 / \nu}}{n^{1 / 2}} \zeta_K^3.
\end{align*}

Secondly, for $i \in \{1, \ldots, n\}$ let $\bg_i \sim \mathcal{N}(0, \Sigma_i)$ be independent vectors, where $\Sigma_i = \Var{\bxi_i \given \vx_i}$. Since there is an independent random variable $U_1$ distributed uniformly on $[0, 1]$, we can construct the family $\crl{\bg_i}$ on the same probability space. Then by Jensen's inequality for any $\lambda > 0$ we have
\begin{align*}
  &\Expect[\big]{\norm{ \bg_i }_{\infty}^2 \given \vx_i} \leq \frac{1}{\lambda} \log \Expect[\big]{e^{\lambda \norm{\bg_i}_{\infty}^2} \given \vx_i} \leq \frac{1}{\lambda} \log \sum_{t = 1}^{K \abs{\mathcal{Q}_n^{\delta}}} \Expect[\big]{e^{\lambda \prn{g_{i t}}^2} \given \vx_i}\\
  &\quad \leq \frac{-\frac{1}{2} \log\prn[\bigg]{1 - \frac{2 \lambda}{n} \zeta_K^2 \mu_n^{2 / \nu}} + \log K + \log \abs{\mathcal{Q}_n^{\delta}}}{\lambda} \lesssim \frac{\zeta_K^2 \mu_n^{2 / \nu}}{n} \prn[\big]{ \log K + \log |\mathcal{Q}_n^{\delta}| },
\end{align*}
where we used the moment-generating function of $\chi_1^2$: $\Expect{\exp\crl{\alpha \chi_1^2}} = \prn{ 1 - 2 \alpha }^{- 1 / 2}$ for $\alpha < 1 / 2$, the bound $\Var{\xi_{i t} \given \vx_i} \leq \left.\zeta_K^2 \mu_n^{2 / \nu} \right/ n$, and put $\lambda := \prn[\big]{4 \zeta_K^2 \mu_n^{2 / \nu} / n}^{-1}$. Also,
\begin{align*}
  &\Expect[\big]{\norm{ \bg_i }^4 \given \vx_i}^{1 / 2} = \Expect[\bigg]{\prn[\bigg]{\sum_{t = 1}^{K |\mathcal{Q}_n^{\delta}|} g_{i t}^2}^2 \given \vx_i}^{1 / 2} \leq \sum_{t = 1}^{K |\mathcal{Q}_n^{\delta}|} \Expect{g_{i t}^4 \given \vx_i}^{1 / 2}\\
  &\quad = \sum_{t = 1}^{K |\mathcal{Q}_n^{\delta}|} \sqrt{3} \Expect{g_{i t}^2 \given \vx_i} \lesssim \Expect[\big]{\norm{ \bg_i }^2 \given \vx_i} = \Expect{\norm{ \bxi_i }^2 \given \vx_i},
\end{align*}
which gives
\begin{align*}
  \sum_{i = 1}^n \Expect[\big]{\norm{ \bg_i }^4 \given \vx_i}^{1 / 2} \lesssim \sum_{i = 1}^n \Expect[\big]{\norm{ \bxi_i }^2 \given \vx_i} \leq \mu_n^{2 / \nu} |\mathcal{Q}_n^{\delta}| \Expectlet_n \brk[\big]{ \|\bb(\vx_i)\|^2 } \lesssim \zeta_K^2 \mu_n^{2/\nu} |\mathcal{Q}_n^{\delta}|.
\end{align*}
Therefore, by the Cauchy-Schwarz inequality,
\begin{align*}
  &\sum_{i = 1}^n \Expect[\big]{\norm{ \bg_i }^2 \norm{ \bg_i }_{\infty} \given \vx_i} \leq \sum_{i = 1}^n \Expect[\big]{\norm{ g_i }^4 \given \vx_i}^{1 / 2} \Expect[\big]{\norm{ g_i }^2_{\infty} \given \vx_i}^{1 / 2}\\
  &\quad \lesssim \frac{\zeta^3_K \mu_n^{3 / \nu}}{\sqrt{n}} |\mathcal{Q}_n^{\delta}| \sqrt{\log \prn[\big]{K |\mathcal{Q}_n^{\delta}|}}.
\end{align*}

Now, since there exists a random variable $U_2$ independent of $\crl{ \bxi_i }_{i = 1}^n \cup \crl{ \bg_i }_{i = 1}^n$, applying~\cref{th:conditional-yurinski} with
\begin{align*}
  &\beta := \sum_{i = 1}^n \Expect[\big]{\norm{ \bxi_i }^2 \norm{ \bxi_i }_{\infty}} + \sum_{i = 1}^n \Expect[\big]{\norm{ \bg_i }^2 \norm{ \bg_i }_{\infty}}\\
  &\quad \lesssim \frac{\zeta^3_K \mu_n^{3 / \nu}}{\sqrt{n}} |\mathcal{Q}_n^{\delta}| \sqrt{\log \prn[\big]{K |\mathcal{Q}_n^{\delta}|}}\\
  &\quad \lesssim \frac{\zeta^3_K \mu_n^{3 / \nu}}{n^{1 / 2} \delta_n^{d_{\mathcal{Q}}}} \sqrt{\log \frac{K}{\delta_n^{d_{\mathcal{Q}}}}}
\end{align*}
gives that for any $t_n > 0$, on the same probability space there exists a vector $\bZ_n^{\delta} \sim \mathcal{N}\prn[\big]{0, \Var{\bG_n \circ \pi_n^{\delta} \given \mX_n}}$, generally different for different $t_n$, such that
\begin{align*}
  \Prob[\big]{\norm[\big]{ \bG_n\circ\pi_n^{\delta} - \bZ_n^{\delta} }_{\infty} > 3 t_n} \leq \min_{s \geq 0} \crl[\big]{ 2 \Prob*{\|\bN\|_{\infty} > s} + \frac{\beta}{t_n^3} s^2 },
\end{align*}
where $\bN$ is a $K |\mathcal{Q}_n^{\delta}|$-dimensional standard Gaussian vector. By the union bound,
\begin{equation*}
  \Prob{\|\bN\|_{\infty} > s} \leq 2 K |\mathcal{Q}_n^{\delta}| e^{-s^2 /2},
\end{equation*}
so by taking $s := C \sqrt{\log \prn[\big]{K |\mathcal{Q}_n^{\delta}|}}$ for a positive constant $C$ not depending on $n$ (chosen later), we have
\begin{equation*}
  \begin{split}
    &\Prob[\big]{\norm[\big]{ \bG_n\circ\pi_n^{\delta} - \bZ_n^{\delta} }_{\infty} > 3 t_n} \lesssim \prn[\big]{K |\mathcal{Q}_n^{\delta}|}^{1 - C^2 / 2} + \frac{\beta}{t_n^3} C^2 \log \prn[\big]{K |\mathcal{Q}_n^{\delta}|}\\
    &\quad \lesssim \prn[\bigg]{\frac{K}{\delta_n^{d_{\mathcal{Q}}}}}^{1 - C^2 / 2} + \frac{\beta}{t_n^3} C^2 \log \frac{K}{\delta_n^{d_{\mathcal{Q}}}}\\
    &\quad \lesssim \prn[\bigg]{\frac{K}{\delta_n^{d_{\mathcal{Q}}}}}^{1 - C^2 / 2} + C^2 \frac{\zeta^3_K \mu_n^{3 / \nu}}{t_n^3 \sqrt{n} \delta_n^{d_{\mathcal{Q}}}} \log^{3 / 2} \frac{K}{\delta_n^{d_{\mathcal{Q}}}},
  \end{split}
\end{equation*}
where the constant in $\lesssim$ does not depend on $n$.

\paragraph{Embedding a conditionally Gaussian vector into a conditionally Gaussian process}
For a fixed vector in $\bX \in \mathcal{X}^n$, by standard existence results for Gaussian processes, there exists a mean-zero $K$-dimensional Gaussian process whose covariance structure is the same as that of $\bG_n(\vq)$ given $\mX_n = \bX$. It follows from Kolmogorov's continuity criterion that this process can be defined on $\cont{\mathcal{Q}, \mathbb{R}^K}$. The laws of such processes define a family of Gaussian probability measures on the Borel $\sigma$-algebra $\mathcal{B} \prn{ \cont{\mathcal{Q}, \mathbb{R}^K} }$ of the space $\cont{\mathcal{Q}, \mathbb{R}^K}$, denoted $\crl{\Problet_{\bX}}_{\bX \in \mathcal{X}^n}$. In order to construct one process that is conditionally on $\mX_n'$, a mean zero Gaussian process with the same conditional covariance structure as that of $\bG_n(\vq)$, where $\mX_n' \overset{d}{=} \mX_n$, we need to show that this family of measures is a probability kernel as a function $\mathcal{X}^n \times \mathcal{B} \prn{ \cont{\mathcal{Q}, \mathbb{R}^K} } \to [0, 1]$. This follows by a standard argument: we can take a $\pi$-system of sets, generating $\mathcal{B} \prn{ \cont{\mathcal{Q}, \mathbb{R}^K} }$, of the form $B = \crl{\bof \in \cont{\mathcal{Q}, \mathbb{R}^K}: \bof(\vq_1) \in B_1, \ldots, \bof(\vq_m) \in B_m}$, where $m \in \{1, 2, \ldots\}$, $\vq_j \in \mathcal{Q}$, and each $B_j$ is a parallelepiped in $\mathbb{R}^K$ with edges parallel to the coordinate axes, and notice that for such sets $\bX \mapsto \Problet_{\bX}(B)$ is a Borel function (since a mean-zero Gaussian vector is a linear transformation of a standard Gaussian vector). The sets $A \in \mathcal{B} \prn{ \cont{\mathcal{Q}, \mathbb{R}^K} }$ such that $\bX \mapsto \Problet_{\bX}(A)$ is a Borel function form a $\lambda$-system. It is left to apply the monotone class theorem.

We have shown that there exists a law on $\mathcal{X}^n \times \cont{\mathcal{Q}, \mathbb{R}^K}$ which is the joint law of $\prn[\big]{\mX_n', \crl{\bZ_n'(\vq)}_{\vq \in \mathcal{Q}}}$, where $\crl{\bZ_n'(\vq)}_{\vq \in \mathcal{Q}}$ is a conditionally on $\mX_n'$ mean zero Gaussian process with the same conditional covariance structure as that of $\bG_n'(\vq)$ (where $\bG_n'(\vq)$ is the same function of $\mX_n$ as $\bG_n(\vq)$ except $\mX_n$ is replaced by $\mX_n'$). Projecting this $\cont{\mathcal{Q}, \mathbb{R}^K}$-process onto $\mathcal{Q}_n^{\delta}$, we obtain a vector $\bZ_n^{\delta\prime}$ such that
\begin{equation*}
\prn[\big]{\mX_n', \bZ_n^{\delta\prime}} \overset{d}{=} \prn[\big]{\mX_n, \bZ_n^{\delta}}.
\end{equation*}
Since $\cont{\mathcal{Q}, \mathbb{R}^K}$ is Polish and there exists a uniformly distributed random variable $U_3$ independent of $\crl{ (\vx_i, y_i) }_{i = 1}^n \cup \crl{\bZ_n^{\delta}}$, we can apply Theorem 8.17 (transfer) in~\parencite{Kallenberg2021}, and obtain that there exists (on the same probability space) a random element $\crl{\bZ_n(\vq)}_{\vq \in \mathcal{Q}} \in \cont{\mathcal{Q}, \mathbb{R}^K}$ such that
\begin{equation*}
  \prn[\big]{\mX_n, \bZ_n^{\delta}, \crl{\bZ_n(\vq)}_{\vq \in \mathcal{Q}}} = \prn[\big]{\mX_n', \bZ_n^{\delta\prime}, \crl{\bZ_n'(\vq)}_{\vq \in \mathcal{Q}}}.
\end{equation*}
In particular, $\crl{\bZ_n(\vq)}_{\vq \in \mathcal{Q}}$ is the conditionally Gaussian process whose projection on $\mathcal{Q}_n^{\delta}$ is the vector $\bZ_n^{\delta}$ a.\,s.

\paragraph{Discretization of $Z_n$} Consider the stochastic process $X_{n}$ defined for $t = (l, \vq, \tilde{\vq}) \in T$ with
\begin{equation*}
  T := \crl{ (l, \vq, \tilde{\vq})\colon l \in \{1, 2, \ldots, K\}, \vq, \tilde{\vq} \in \mathcal{Q}, \|\vq - \tilde{\vq}\| \leq \delta_n }
\end{equation*}
as $X_{n, t} := Z_{n, l}(\vq) - Z_{n, l}(\tilde{\vq})$. It is a separable (because each $Z_{n, l}(\cdot)$ has a.\,s. continuous trajectories) mean-zero Gaussian conditionally on $\mX_n$ process with the index set $T$ considered a metric space: $\dist \prn[\big]{(l, \vq, \tilde{\vq}), (l', \vq', \tilde{\vq}')} = \|\vq - \vq'\| + \|\tilde{\vq} - \tilde{\vq}'\| + \Indlet\{l \neq l'\}$.

We will apply \cref{lem:maximal-gaussian} to this process. Note that on the event $\mathcal{A}$
\begin{align*}
  &\sigma(X_n)^2 := \sup_{t \in T} \Expect[\big]{X_{n, t}^2 \given \mX_n} = \sup_{\|\vq - \tilde{\vq}\| \leq \delta_n} \max_l \Expect[\big]{ \prn[\big]{Z_{n, l}(\vq) - Z_{n, l}(\tilde{\vq})}^2 \given \mX_n}\\
  &\quad = \sup_{\|\vq - \tilde{\vq}\| \leq \delta_n} \max_l \Expect[\big]{ \prn[\big]{G_{n, l}(\vq) - G_{n, l}(\tilde{\vq})}^2 \given \mX_n}\\
  &\quad = \sup_{\|\vq - \tilde{\vq}\| \leq \delta_n} \max_l \frac{1}{n} \sum_{i = 1}^n b_l(\vx_i)^2 \Expect[\big]{\prn[\big]{ A_n(\vq, \vx_i, y_i) - A_n(\tilde{\vq}, \vx_i, y_i) }^2 \given \vx_i}\\
  &\quad \lesssim \sup_{\|\vq - \tilde{\vq}\| \leq \delta_n} \max_l \frac{1}{n} \sum_{i = 1} b_l(\vx_i)^2 \|\vq - \tilde{\vq}\| \leq \delta_n \max_l \Expectlet_n \brk{b_l(\vx_i)^2}\\
  &\quad \leq \delta_n \underbrace{\sup_{\norm{ \balpha } = 1} \Expectlet_n \brk[\big]{\prn[\big]{\balpha \trans \bb(\vx_i)}^2}}_{\leq C_{\mathrm{Gr}}} \lesssim \delta_n.
\end{align*}

Next, we define and bound the semimetric $\rho(t, t')$:
\begin{align*}
  &\rho(t, t')^2 := \Expect[\big]{\prn{ X_{n, t} - X_{n, t'} }^2 \given \mX_n}\\
  &\quad = \Expect[\big]{\prn[\big]{\prn[\big]{ Z_{n, l}(\vq) - Z_{n, l}(\tilde{\vq}) } - \prn[\big]{ Z_{n, l'}(\vq') - Z_{n, l'}(\tilde{\vq}') }}^2 \given \mX_n}\\
  &\quad \lesssim \Expect[\big]{\prn[\big]{ Z_{n, l}(\vq) - Z_{n, l'}(\vq') }^2 \given \mX_n} + \Expect[\big]{\prn[\big]{ Z_{n, l}(\tilde{\vq}) - Z_{n, l'}(\tilde{\vq}') }^2 \given \mX_n}.
\end{align*}
The first term on the right is bounded the following way: if $l \neq l'$,
\begin{align*}
  &\Expect[\big]{\prn[\big]{ Z_{n, l}(\vq) - Z_{n, l'}(\vq') }^2 \given \mX_n} = \Expect[\big]{\prn[\big]{ G_{n, l}(\vq) - G_{n, l'}(\vq') }^2 \given \mX_n}\\
  &\quad = \frac{1}{n} \sum_{i = 1}^n \Expect[\big]{\prn[\big]{ A_n(\vq, \vx_i, y_i) b_l(\vx_i) - A_n(\vq', \vx_i, y_i) b_{l'}(\vx_i) }^2 \given \vx_i}\\
  &\quad \lesssim \frac{1}{n} \sum_{i = 1}^n b_l(\vx_i)^2 \Expect[\big]{\prn[\big]{ A_n(\vq, \vx_i, y_i) - A_n(\vq', \vx_i, y_i) }^2 \given \vx_i}\\
  &\qquad + \frac{1}{n} \sum_{i = 1}^n \prn[\big]{b_l(\vx_i)^2 + b_{l'}(\vx_i)^2} \Expect[\big]{A_n(\vq', \vx_i, y_i)^2 \given \vx_i}\\
  &\quad \lesssim \|\vq - \vq'\| \Expectlet_n \brk[\big]{b_l(\vx_i)^2} + \mu_n^{2 / \nu} \prn[\big]{\Expectlet_n \brk{b_l(\vx_i)^2} + \Expectlet_n \brk{b_{l'}(\vx_i)^2}} \lesssim \|\vq - \vq'\| + \mu_n^{2 / \nu}.
\end{align*}
Similarly, if $l = l'$, $\Expect[\big]{\prn[\big]{ Z_{n, l}(\vq) - Z_{n, l'}(\vq') }^2 \given \mX_n} \lesssim \|\vq - \vq'\|$.

The term $\Expect[\big]{\prn[\big]{ Z_{n, l}(\tilde{\vq}) - Z_{n, l'}(\tilde{\vq}') }^2 \given \mX_n}$ is bounded the same way, and we conclude
\begin{align*}
  \rho(t, t')^2 \lesssim \mu_n^{2 / \nu} \prn[\big]{ \|\vq - \vq'\| + \|\tilde{\vq} - \tilde{\vq}'\| + \Indlet\{l \neq l'\}} = \mu_n^{2 / \nu} \mathrm{dist}\prn[\big]{ (l, \vq, \tilde{\vq}), (l', \vq', \tilde{\vq}') }.
\end{align*}

In other words, we have proven that for some positive constant $\nC{covering}$ we have on $\mathcal{A}$
\begin{equation*}
  \rho(t, t')^2 \leq \oC{covering} \mu_n^{2 / \nu} \mathrm{dist}\prn[\big]{ (l, \vq, \tilde{\vq}), (l', \vq', \tilde{\vq}') }.
\end{equation*}

This means that an $\prn[\big]{\varepsilon / \mu_n^{1 / \nu} \sqrt{\oC{covering}}}^2$-covering of $T$ with respect to $\dist(\cdot)$ induces an $\varepsilon$-covering of $T$ with respect to $\rho$, and hence
\begin{equation}\label{eq:prelim-n-bound}
  N(T, \rho, \varepsilon) \leq N\prn[\bigg]{T, \dist(\cdot), \prn[\bigg]{\frac{\varepsilon}{\mu_n^{1 / \nu} \sqrt{\oC{covering}}}}^2}.
\end{equation}
Therefore,
\begin{equation*}
  \log N(T, \rho, \varepsilon) \lesssim \log\prn{K \mu_n^{1 / \nu} / \varepsilon}.
\end{equation*}

Applying \cref{lem:maximal-gaussian} gives that on the event $\mathcal{A}$
\begin{align*}
  &\Expect[\Big]{\sup_{t \in T} \abs{ X_{n, t} } \given \mX_n} \lesssim \sigma(X_n) + \int_0^{2 \sigma(X_n)} \sqrt{\log N(T, \rho, \varepsilon)}\,\mathrm{d}\varepsilon\\
  &\quad \lesssim \sigma(X_n) + \sigma(X_n) \sqrt{\log\prn[\big]{K \mu_n^{1 / \nu} / \sigma(X_n)}}  \labrel{1}[\lesssim] \prn[\bigg]{ \delta_n \log \prn[\bigg]{K \frac{\mu_n^{1 / \nu}}{\delta_n}} }^{1 / 2},
\end{align*}
where in \labrel{1} we used our bound $\sigma(X_n) \lesssim \sqrt{\delta_n}$ above and that $x \mapsto x \log \frac{1}{x}$ is increasing for sufficiently small $x$. Rewriting and applying Markov's inequality, we obtain that on $\mathcal{A}$
\begin{equation*}
  \Prob[\bigg]{\sup_{\|\vq - \tilde{\vq}\| \leq \delta_n} \norm{ \bZ_n(\vq) - \bZ_n(\tilde{\vq}) }_{\infty} > t_n \given \mX_n} \lesssim \frac{1}{t_n} \sqrt{\delta_n \log \frac{K \mu_n^{1 / \nu}}{\delta_n}},
\end{equation*}
where the constant in $\lesssim$ does not depend on $n$. Since $\Prob*{\mathcal{A}} \to 1$, this implies
\begin{equation}\label{eq:discretization-zn-final-bound}
  \Prob[\bigg]{\sup_{\|\vq - \tilde{\vq}\| \leq \delta_n} \norm{ \bZ_n(\vq) - \bZ_n(\tilde{\vq}) }_{\infty} > t_n} \lesssim \frac{1}{t_n} \sqrt{\delta_n \log \frac{K \mu_n^{1 / \nu}}{\delta_n}} + o(1).
\end{equation}

\paragraph{Choosing $\delta_n$ and conclusion} Combining the bounds obtained above, for any given positive sequence $t_n$ and any constant $C > 0$ of our choice
\begin{align*}
  &\Prob*{\sup_{\vq \in \mathcal{Q}} \norm{ \bG_n(\vq) - \bZ_n(\vq) }_{\infty} > t_n} \leq \Prob*{\sup_{\|\vq - \tilde{\vq}\| \leq \delta_n} \norm{ \bG_n(\vq) - \bG_n(\tilde{\vq}) }_{\infty} > t_n / 3}\\
  &\qquad + \Prob*{\norm[\big]{ \bG_n\circ\pi_n^{\delta} - \bZ_n^{\delta} }_{\infty} > t_n / 3} + \Prob*{\sup_{\|\vq - \tilde{\vq}\| \leq \delta_n} \norm{ Z_n(\vq) - Z_n(\tilde{\vq}) }_{\infty} > t_n / 3}\\
  &\quad \lesssim \frac{1}{t_n} \sqrt{ \delta_n \log \frac{K \zeta_K \mu_n^{1 / \nu}}{\sqrt{\delta_n}}} + \frac{\zeta_K \prn{\mu_n n}^{1 / \nu}}{t_n \sqrt{n}} \log \frac{K \zeta_K \mu_n^{1 / \nu}}{\sqrt{\delta_n}}\\
  &\qquad + \prn[\bigg]{\frac{K}{\delta_n^{d_{\mathcal{Q}}}}}^{1 - C^2 / 2} + C^2 \frac{\zeta^3_K \mu_n^{3 / \nu}}{t_n^3 \sqrt{n} \delta_n^{d_{\mathcal{Q}}}} \log^{3 / 2} \frac{K}{\delta_n^{d_{\mathcal{Q}}}}\\
  &\qquad + \frac{1}{t_n} \sqrt{\delta_n \log \frac{K \mu_n^{1 / \nu}}{\delta_n}} + o(1),
\end{align*}
where the constant in $\lesssim$ does not depend on $n$.

Take, for example, $C = 2$ (so that $1 - C^2 / 2$ is negative). Now we approximately (assuming that each $\log(\cdot)$ on the right is $O(\log n)$ and ignoring constant coefficients) optimize this over $\delta_n$. This gives
\begin{equation*}
  \delta_n := \prn[\bigg]{ \frac{\zeta^3_K \mu_n^{3 / \nu}}{t_n^2 \sqrt{n}} \log n }^{\frac{2}{1 + 2 d_{\mathcal{Q}}}}.
\end{equation*}

Let $\ell_n \to 0$ be a positive sequence satisfying
\begin{equation*}
  \prn[\bigg]{ \frac{\zeta^3_K \mu_n^{3 / \nu}}{\sqrt{n}} }^{\frac{1}{3 + 2 d_{\mathcal{Q}}}} \sqrt{\log n} + \frac{\zeta_K \mu_n^{1 / \nu}}{n^{1 / 2 - 1 / \nu}} \log n = o \prn{ \ell_n \ryur }.
\end{equation*}
Clearly, whenever \eqref{eq:yurinski-lemma-rate-restriction} holds, such a sequence exists. Putting $t_n := \ell_n r_{\mathrm{yur}}$, we get
\begin{equation*}
  \Prob[\bigg]{\sup_{\vq \in \mathcal{Q}} \norm{ \bG_n(\vq) - \bZ_n(\vq) }_{\infty} > \ell_n r_{\mathrm{yur}}} = o(1).
\end{equation*}
Fix $\varepsilon > 0$. For $n$ large enough, $\ell_n < \varepsilon$. Then for these $n$ we have
\begin{equation*}
  \Prob[\bigg]{\sup_{\vq \in \mathcal{Q}} \norm{ \bG_n(\vq) - \bZ_n(\vq) }_{\infty} > \varepsilon r_{\mathrm{yur}}} \leq \Prob[\bigg]{\sup_{\vq \in \mathcal{Q}} \norm{ \bG_n(\vq) - \bZ_n(\vq) }_{\infty} > \ell_n r_{\mathrm{yur}}} = o(1).
\end{equation*}

\Cref{lem:yurinski-coupling-K-dimensional} is proven.
\end{proof}

\subsection{Main Result}\label{sec:SA-main-result}

We begin by presenting our main strong approximation result, which  is a special case of more technical and lengthy \cref{lem:yurinski-coupling-K-dimensional}. To simplify exposition, the following notation will be helpful from this point onwards:
\begin{align}
  \ellhatqx[\vq] &:= - h^{d / 2} \frac{\qhat[\vq][-1] \vp^{(\vv)}(\vx)}{\sqrt{\bigomegahat[\vq]}},\label{eq:ellhatqx-def}\\
  \elltwobarqx[\vq] &:= - h^{d / 2} \frac{\qbar[\vq][-1] \vp^{(\vv)}(\vx)}{\sqrt{\bigomegabar[\vq]}},\label{eq:elltwobarqx-def}\\
  \ellqx[\vq] &:= - h^{d / 2} \frac{\qzero[\vq][-1] \vp^{(\vv)}(\vx)}{\sqrt{\bigomega[\vq]}},\label{eq:ellqx-def}\\
  \bart[\vq] &:= h^{-d / 2} \elltwobarqx[\vq]\trans \sqrt{n} \Expectlet_n \brk[\big]{ \riskder[\vq]{y_i}{\eta \prn{ \muzero[\vq]{\vx_i} }} \eta^{(1)} \prn{ \muzero[\vq]{\vx_i} } \vp(\vx_i) },\label{eq:bartqx-def}\\
  \hatTqx[\vq] &:= \frac{\vp^{(\vv)}(\vx)\trans \betahat[\vq] - \muzero[\vq][(\vv)]{\vx}}{\sqrt{\bigomegahat[\vq] / n}},\label{eq:hatTqx-def}
\end{align}
where $\newtarget{def:bigomegahat}{\bigomegahat[\vq]}$ is some feasible estimator of variance.

\begin{theorem}[Strong approximation]\label{th:strong-approximation-binscatter-yurinski}\leavevmode
  \begin{assertions}[wide]
  \item\label{th:yurinski-normal} Suppose \cref{ass:dgp-holder,ass:gl-quasi-uniform-partition,ass:gl-local-basis,ass:gl-holder,ass:complexity,ass:gl-approximation-error} hold with $\nu \geq 3$. Furthermore, assume the following conditions hold:
    \begin{gather}
      \sup_{\vq \in \mathcal{Q}} \sup_{\vx \in \mathcal{X}} \abs{\muhat[\vq][(\vv)]{\vx} - \muzero[\vq][(\vv)]{\vx}} \lesssim_{\Problet} h^{- \multindnorm{\vv}} \ruc;\label{eq:gRcohQ}\\
      \sup_{\vq \in \mathcal{Q}} \sup_{\vx \in \mathcal{X}} \abs{\muhat[\vq][(\vv)]{\vx} - \muzero[\vq][(\vv)]{\vx} - \bahterm[\vq][(\vv)]{\vx}} \lesssim_{\Problet} h^{- \multindnorm{\vv}} \rbah\text{ with } \frac{\log n}{n h^d} \lesssim \rbah;\label{eq:yXRXnK}\\
      \sup_{\vq \in \mathcal{Q}} \sup_{\vx \in \mathcal{X}} \abs[\big]{\bigomegahat[\vq] - \bigomegabar[\vq]} \lesssim_{\Problet} h^{-2 \multindnorm{\vv} - d} \romega \text{ with } \romega = o(1);\label{eq:variance-closeness}\\
      \begin{multlined}
      \Expect{
        (\riskder[\vq]{y_i}{\eta\prn{\muzero[\vq]{\vx_i}}} \eta^{(1)}\prn{\muzero[\vq]{\vx_i}} - \riskder[\tilde{\vq}]{y_i}{\eta\prn{\muzero[\tilde{\vq}]{\vx_i}}} \eta^{(1)}\prn{\muzero[\tilde{\vq}]{\vx_i}}
        )^2 \given \vx_i
      }\\
      \lesssim \| \vq - \tilde{\vq} \| \text{ for all } \vq, \tilde{\vq} \in \mathcal{Q}.
      \end{multlined}\label{eq:l2-norm-holder-type-restriction}
    \end{gather}

    Let $\newtarget{def:yurinski-corollary-rate}{\rstr}$ be any positive sequence of numbers converging to zero\footnote{In particular, the left-hand side of \eqref{eq:MFPKjN} is automatically assumed to be $o(1)$.} such that
\begin{equation}\label{eq:MFPKjN}
  \prn[\bigg]{\frac{1}{n h^{3 d}}}^{\frac{1}{6 + 4 d_{\mathcal{Q}}}} \sqrt{\log n} + \frac{1}{h^{\sfrac{d}{2}} n^{1 / 2 - 1 / \nu}} \log n = o(\rstr).
\end{equation}

    Denote
\begin{equation*}
\newtarget{def:gn}{\gnNoLink}(\vq) := h^{- d / 2} \sqrt{n} \Expectlet_n \brk[\big]{\riskder[\vq]{y_i}{\eta\prn{ \muzero[\vq]{\vx_i}}} \eta^{(1)} \prn{\muzero[\vq]{\vx_i}} \vp(\vx_i)}.
\end{equation*}

    Then (provided the probability space is rich enough) there exists a $K$-dimensional,
    conditionally on $\mX_n$ mean-zero Gaussian, process $\newtarget{def:zn}{\znNoLink}(\vq)$ on $\mathcal{Q}$ with a.\,s. continuous trajectories, satisfying
    \begin{gather}
      \Expect{\gn(\vq) \gn(\tilde{\vq})\trans \given \mX_n} = \Expect{\zn(\vq) \zn(\tilde{\vq})\trans \given \mX_n};\\
      \sup_{\vq} \norm{ \gn(\vq) - \zn(\vq) }_{\infty} = o_{\Problet}\prn{\rstr};\numberthis\label{eq:XcWykx}\\
      \sup_{\vx \in \mathcal{X}} \sup_{\vq \in \mathcal{Q}} \abs[\big]{ \elltwobarqx[\vq]\trans \gn(\vq) - \elltwobarqx[\vq]\trans \zn(\vq) } = o_{\Problet}(\rstr),\numberthis\label{eq:yurinski-corollary-statement}\\
      \sup_{\vq, \vx} \abs[\big]{ \hatTqx[\vq] - \elltwobarqx[\vq]\trans \zn(\vq) } \lesssim_{\Problet} \sqrt{n h^d} (\ruc \romega + \rbah) + o(\rstr).\label{eq:strong-approximation-t-process}
    \end{gather}

  \item\label{th:yurinski-subgauss} If, in addition to the previous conditions, $\riskderenvel{\vx_i}{y_i}$ is \subgauss\ conditionally on $\vx_i$, and $\newtarget{def:yurinski-corollary-rate-subgauss}{\rstrsubgaussNoLink} \to 0$ is a positive sequence of numbers satisfying
  \begin{equation*}
    \prn[\bigg]{\frac{1}{n h^{3 d}}}^{\frac{1}{6 + 4 d_{\mathcal{Q}}}} \sqrt{\log n} + \frac{\log^{3 / 2} n}{\sqrt{n h^d}} = o(\rstrsubgauss),
  \end{equation*}
  then
  \begin{equation}\label{eq:yurinski-corollary-statement-subgauss}
    \sup_{\vx \in \mathcal{X}} \sup_{\vq \in \mathcal{Q}} \abs[\big]{ \elltwobarqx[\vq]\trans \gn(\vq) - \elltwobarqx[\vq]\trans \zn(\vq) } = o_{\Problet}(\rstrsubgauss).
  \end{equation}
  \end{assertions}
\end{theorem}

\begin{proof}
Note that \eqref{eq:XcWykx} implies \eqref{eq:yurinski-corollary-statement} because $\sup_{\vq, \vx} \norm[\big]{\elltwobarqx[\vq]}_1 \lesssim_{\Problet} 1$ by \cref{lem:closeness-of-linear-terms} below. To prove \eqref{eq:XcWykx}, apply \cref{lem:yurinski-coupling-K-dimensional} with $\gn(\cdot)$ as in the statement,
\begin{align*}
  &A_{n}(\vq, y) = A_n(\vq, \vx, y) := \riskder[\vq]{y}{\eta\prn{\muzero[\vq]{\vx}}} \eta^{(1)}\prn{\muzero[\vq]{\vx}},\\
  &\bb(\cdot) := h^{- d / 2} \vp(\cdot),\\
  &\zeta_K \lesssim K^{1 / 2} \lesssim h^{- d / 2},\\
  &\mu_n \lesssim 1.
\end{align*}
Finally, \eqref{eq:strong-approximation-t-process} follows from combining \eqref{eq:yurinski-corollary-statement} with \cref{lem:hats-off}.
\end{proof}

\begin{remark}\label{rem:strong-approx-correspondence}
  This theorem contains Theorem~\mainref{thm: strong approximation} in the main paper; the notation there is simplified for better readability. Specifically, $T(\vx,q)$ in the main paper corresponds to $\hatTqx[\vq]$ here and the approximating process $Z(\vx,q)$ corresponds to $\elltwobarqx[\vq]\trans \zn(\vq)$ here. The reason for differing notation is that we require more precision here: given the form of \cref{lem:yurinski-coupling-K-dimensional}, it is more natural to use $K$-dimensional processes such as $\zn(\vq)$, also making for simpler presentation of precise results in \cref{sec:uniform-inference} later.
\end{remark}

The rest of this subsection will be devoted to proving the lemmas used in the proof of \cref{th:strong-approximation-binscatter-yurinski}. First, we prove bounds on the asymptotic variance and its approximation. Similar bounds were proven in \parencite{Cattaneo-Crump-Farrell-Feng_2024_AER,Cattaneo-Farrell-Feng_2020_AOS}.

\begin{lemma}[Asymptotic variance]\label{lem:asympt-variance}
  Suppose Assumptions \cref{ass:dgp-holder,ass:gl-quasi-uniform-partition,ass:gl-local-basis,ass:gl-holder,ass:complexity,ass:gl-approximation-error} hold, and $\frac{\log (1 / h)}{n h^{d}} = o(1)$. Then
  \begin{align}
    &h^{-d - 2 \multindnorm{\vv}} \lesssim \inf_{\vq \in \mathcal{Q}} \inf_{\vx \in \mathcal{X}} \bigomega[\vq] \leq \sup_{\vq \in \mathcal{Q}} \sup_{\vx \in \mathcal{X}} \bigomega[\vq] \lesssim h^{-d - 2 \multindnorm{\vv}},\label{eq:bigomega-bounds}\\
    &\sup_{\vq \in \mathcal{Q}} \sup_{\vx \in \mathcal{X}} \abs[\big]{\bigomegabar[\vq] - \bigomega[\vq]} \lesssim_{\Problet} h^{- d - 2 \multindnorm{\vv}} \sqrt{\frac{\log(1 / h)}{n h^d}},\label{eq:bigomega-and-bigomegabar-closeness}\\
    &h^{-d - 2 \multindnorm{\vv}} \lesssim_{\Problet} \inf_{\vq \in \mathcal{Q}} \inf_{\vx \in \mathcal{X}} \bigomegabar[\vq] \leq \sup_{\vq \in \mathcal{Q}} \sup_{\vx \in \mathcal{X}} \bigomegabar[\vq] \lesssim_{\Problet} h^{-d - 2 \multindnorm{\vv}}.\label{eq:bigomegabar-bounds}
  \end{align}
\end{lemma}

\begin{proof}

For the lower bound in~\eqref{eq:bigomega-bounds}, we have
\begin{equation}\label{eq:bigomega-lower-bound}
  \begin{split}
    &\bigomega[\vq] \geq \lambda_{\min}(\sigmazero[\vq]) \cdot \norm{ \qzero[\vq][-1]  \vp^{(\vv)}(\vx) }^2 \geq \lambda_{\min}(\sigmazero[\vq]) \cdot \norm{ \qzero[\vq]}^{-2} \norm{ \vp^{(\vv)}(\vx) }^2\\
    &\quad \gtrsim \lambda_{\min}(\sigmazero[\vq]) \cdot h^{- 2 d} \cdot h^{- 2 \multindnorm{\vv}} \gtrsim h^{d} \cdot h^{- 2 d} \cdot h^{- 2 \multindnorm{\vv}} = h^{- d - 2 \multindnorm{\vv}}
  \end{split}
\end{equation}
by Assumptions~\ref{ass:gl-local-basis}, \ref{ass:gl-holder} ($\sigma^2_{\vq}(\vx)$ and $\eta^{(1)} \prn[\big]{ \muzero[\vq]{\vx} }^2$ are uniformly over $\vx$ bounded away from zero) and \cref{lem:holder-gram}.

For the upper bound in~\eqref{eq:bigomega-bounds}, we have
\begin{align*}
  &\bigomega[\vq] \leq \lambda_{\max}(\sigmazero[\vq]) \cdot \big\| \qzero[\vq][-1] \vp^{(\vv)}(\vx) \big\|^2 \leq \lambda_{\max}(\sigmazero[\vq]) \cdot \big[\lambda_{\max}(\qzero[\vq][-1])\big]^{2} \big\| \bp^{(\vv)}(\vx) \big\|^2\\
  &\quad \overset{\text{(a)}}{\lesssim} \lambda_{\max}(\sigmazero[\vq]) \cdot h^{-2 d} \cdot h^{- 2 \multindnorm{\vv}} \overset{\text{(b)}}{\lesssim} h^d \cdot h^{-2 d} \cdot h^{- 2 \multindnorm{\vv}} = h^{- d - 2 \multindnorm{\vv}},
\end{align*}
where (a) is by \cref{ass:gl-local-basis} and \cref{lem:holder-gram}, (b) is by \cref{ass:gl-holder} ($\sigma^2_{\vq}(\vx)$ and $\eta^{(1)} \prn[\big]{ \muzero[\vq]{\vx} }^2$ are uniformly over $\vx$ bounded) and \cref{lem:holder-gram}.

We will now prove~\eqref{eq:bigomega-and-bigomegabar-closeness}. We start by noticing
\begin{align*}
  &\sup_{\vq, \vx} \abs[\big]{ \vp^{(\vv)}(\vx)\trans \qzero[\vq][-1] \prn{ \sigmabar[\vq] - \sigmazero[\vq] } \qzero[\vq][-1] \vp^{(\vv)}(\vx) }\\
  &\quad \lesssim_{\Problet} \sup_{\vq, \vx} \norm{ \sigmabar[\vq] - \sigmazero[\vq] } \cdot \norm{ \qzero[\vq][-1] }^2 \cdot \norm{ \vp^{(\vv)}(\vx) }^2\\
  &\quad \overset{\text{(a)}}{\lesssim_{\Problet}} h^{-2 d} \cdot h^{-2 \multindnorm{\vv}} \sup_{\vq} \norm{ \sigmabar[\vq] - \sigmazero[\vq] } \overset{\text{(b)}}{\lesssim_{\Problet}} h^{-2 d - 2 \multindnorm{\vv}} h^d \sqrt{\frac{\log(1 / h)}{n h^d}} = h^{- d - 2 \multindnorm{\vv}} \sqrt{\frac{\log(1 / h)}{n h^d}},
\end{align*}
where in (a) we used $\sup_{\vq} \norm{ \qzero[\vq][-1] } \lesssim_{\Problet} h^{-d}$ by \cref{lem:holder-gram} and $\sup_{\vx} \norm{ \vp^{(\vv)}(\vx) } \lesssim h^{-\multindnorm{\vv}}$ by \cref{ass:gl-local-basis}, in (b) we used $\sup_{\vq} \norm{ \sigmabar[\vq] - \sigmazero[\vq] } \lesssim_{\Problet} h^d \sqrt{\frac{\log(1 / h)}{n h^d}}$ which is proven by the same argument as $\sup_{\vq} \norm{ \qbar[\vq] - \qzero[\vq] } \lesssim_{\Problet} h^d \sqrt{\frac{\log(1 / h)}{n h^d}}$ in \cref{lem:holder-gram}. Similarly,
\begin{align*}
  &\sup_{\vq, \vx} \abs[\big]{ \vp^{(\vv)}(\vx)\trans (\qbar[\vq][-1] - \qzero[\vq][-1]) \sigmabar[\vq] \qbar[\vq][-1] \vp^{(\vv)}(\vx) }\\
  &\quad \leq \sup_{\vq, \vx} \norm{\qbar[\vq][-1] - \qzero[\vq][-1]} \cdot \norm{\qbar[\vq][-1]} \cdot \norm{\sigmabar[\vq]} \cdot \norm{\vp^{(\vv)}(\vx)}^2\\
  &\quad \lesssim_{\Problet} h^{-d} \sqrt{\frac{\log(1 / h)}{n h^d}} \cdot h^{-d} \cdot h^d \cdot h^{-2 \multindnorm{\vv}} = h^{- d - 2 \multindnorm{\vv}} \sqrt{\frac{\log(1 / h)}{n h^d}}.
\end{align*}
Finally,
\begin{equation*}
  \sup_{\vq, \vx} \abs[\big]{ \vp^{(\vv)}(\vx)\trans \qzero[\vq][-1] \sigmabar[\vq] (\qzero[\vq][-1] - \qbar[\vq][-1]) \vp^{(\vv)}(\vx) } \lesssim_{\Problet} h^{- d - 2 \multindnorm{\vv}} \sqrt{\frac{\log(1 / h)}{n h^d}}.
\end{equation*}
Combining the bounds above gives $\sup_{\vq, \vx} \abs[\big]{\bigomegabar[\vq] - \bigomega[\vq]} \lesssim_{\Problet} h^{- d - 2 \multindnorm{\vv}} \sqrt{\frac{\log(1 / h)}{n h^d}}$.

\eqref{eq:bigomegabar-bounds} is an immediate consequence of~\eqref{eq:bigomega-bounds} and~\eqref{eq:bigomega-and-bigomegabar-closeness}, since $\frac{\log (1 / h)}{n h^{d}} = o(1)$.
\end{proof}

\begin{lemma}[Closeness of linear terms]\label{lem:closeness-of-linear-terms}
  Suppose \cref{ass:dgp-holder,ass:gl-quasi-uniform-partition,ass:gl-local-basis,ass:gl-holder,ass:complexity,ass:gl-approximation-error}, \eqref{eq:variance-closeness} hold, and $\frac{\log(1 / h)}{n h^d} = o(1)$. Then
  \begin{align}
    &\sup_{\vq, \vx} \abs[\big]{ \bigomegahat[\vq]^{-1 / 2} - \bigomegabar[\vq]^{-1 / 2} } \lesssim_{\Problet} h^{d / 2 + \multindnorm{\vv}} \romega,\label{eq:omega-inverse-root-diff-bound}\\
    &\sup_{\vq, \vx} \abs[\big]{ \bigomegabar[\vq]^{-1 / 2} - \bigomega[\vq]^{-1 / 2} } \lesssim_{\Problet} h^{d / 2 + \multindnorm{\vv}} \sqrt{\frac{\log(1 / h)}{n h^d}},\label{eq:omega-bar-inverse-root-diff-bound}\\
    &\sup_{\vq, \vx} \norm{ \ellqx[\vq] }_1 \lesssim 1,\label{eq:ellqx-bound}\\
    &\sup_{\vq, \vx} \norm{ \elltwobarqx[\vq] - \ellqx[\vq] }_1 \lesssim_{\Problet} \sqrt{\frac{\log(1 / h)}{n h^d}}, \label{eq:elltwobardiff-bound}\\
    &\sup_{\vq, \vx} \norm{ \elltwobarqx[\vq] }_1 \lesssim 1\quad\text{w.\,p.\,a. 1}.\label{eq:elltwobar-bound}
  \end{align}

  If, in addition, \cref{ass:gram-estimate} holds, then
\begin{align}
  &\sup_{\vq, \vx} \norm{ \ellhatqx[\vq] - \elltwobarqx[\vq] }_1 \lesssim_{\Problet} \rqhat + \romega,\label{eq:ellhatdiff-bound}\\
  &\sup_{\vq, \vx} \norm{ \ellhatqx[\vq] }_1 \lesssim 1\quad\text{w.\,p.\,a. 1}.\label{eq:ellhat-bound}
\end{align}
\end{lemma}

\begin{proof}
\Cref{eq:omega-inverse-root-diff-bound} follows from the following chain of inequalities:
\begin{align*}
  &\sup_{\vq, \vx} \abs[\big]{ \bigomegahat[\vq]^{-1 / 2} - \bigomegabar[\vq]^{-1 / 2}}\\
  &\quad \leq \sup_{\vq, \vx} \brk[\big]{\bigomegahat[\vq] \bigomegabar[\vq]}^{- 1 / 2} \brk[\big]{ \bigomegahat[\vq]^{1 / 2} + \bigomegabar[\vq]^{1 / 2} }^{-1} \sup_{\vq, \vx} \abs[\big]{ \bigomegahat[\vq]  - \bigomegabar[\vq] }\\
  &\quad \overset{\text{(a)}}{\lesssim_{\Problet}} \prn{h^{d + 2 \multindnorm{\vv}}}^{\sfrac{3}{2}} \sup_{\vq, \vx} \abs[\big]{ \bigomegahat[\vq]  - \bigomegabar[\vq] } \overset{\text{(b)}}{\lesssim_{\Problet}} h^{d / 2 + \multindnorm{\vv}} \romega,
\end{align*}
where in (a) we used that $\bigomegahat[\vq]$ is close to $\bigomegabar[\vq]$ by \eqref{eq:variance-closeness} and the lower bound on $\bigomegabar[\vq]$ from \cref{lem:asympt-variance}, and in (b) we used \eqref{eq:variance-closeness}.

To prove \eqref{eq:omega-bar-inverse-root-diff-bound},
recall that
\begin{align*}
  &\sup_{\vq, \vx} \abs[\big]{\bigomegabar[\vq] - \bigomega[\vq]} \lesssim_{\Problet} h^{- d - 2 \multindnorm{\vv}} \sqrt{\frac{\log(1 / h)}{n h^d}},\quad\text{and}\\
&\inf_{\vq, \vx} \bigomegabar[\vq] \gtrsim_{\Problet} h^{- d - 2 \multindnorm{\vv}}
\end{align*}
by \cref{lem:asympt-variance}. Therefore,
  \begin{align*}
    &\sup_{\vq, \vx} \abs[\big]{ \bigomegabar[\vq]^{-1 / 2} - \bigomega[\vq]^{-1 / 2}}\\
    &\quad \leq \sup_{\vq, \vx} \brk[\big]{\bigomegabar[\vq] \bigomega[\vq]}^{- 1 / 2} \brk[\big]{\bigomegabar[\vq]^{1 / 2} + \bigomega[\vq]^{1 / 2}}^{-1} \sup_{\vq, \vx} \abs[\big]{ \bigomegabar[\vq]  - \bigomega[\vq] }\\
    &\quad \lesssim_{\Problet} \prn{h^{d + 2 \multindnorm{\vv}}}^{\sfrac{3}{2}} \sup_{\vq, \vx} \abs[\big]{ \bigomegabar[\vq]  - \bigomega[\vq] } \lesssim_{\Problet} h^{3 d / 2 + 3 \multindnorm{\vv}} h^{-d - 2 \multindnorm{\vv}} \sqrt{\frac{\log(1 / h)}{n h^d}}\\
    &\quad = h^{d / 2 + \multindnorm{\vv}} \sqrt{\frac{\log(1 / h)}{n h^d}}.
  \end{align*}

\Cref{eq:ellqx-bound}
follows from
\begin{equation*}
  h^{d / 2} \norm{ \qbar[\vq][-1] \vp^{(\vv)}(\vx) }_1 \leq h^{d / 2} \norm{ \qbar[\vq][-1] }_1 \cdot \norm{ \vp^{(\vv)}(\vx) }_1 \lesssim_{\Problet} h^{- d / 2 - \multindnorm{\vv}}
\end{equation*}
and \cref{lem:asympt-variance}.

We will now prove~\eqref{eq:elltwobardiff-bound}.
We have that
\begin{align*}
  &\sup_{\vq} \norm{ \qbar[\vq][-1] }_1 \lesssim_{\Problet} h^{-d}\\
  \intertext{by \cref{lem:holder-gram}, and}
  &\sup_{\vx} \norm{ \vp^{(\vv)}(\vx) }_1 \lesssim h^{-\multindnorm{\vv}}
\end{align*}
by \cref{ass:gl-local-basis}. Combining this with~\eqref{eq:omega-bar-inverse-root-diff-bound} gives
  \begin{equation*}
    \norm*{\frac{h^{d / 2} \qbar[\vq][-1] \vp^{(\vv)}(\vx)}{\sqrt{\bigomegabar[\vq]}} - \frac{h^{d / 2} \qbar[\vq][-1] \vp^{(\vv)}(\vx)}{\sqrt{\vphantom{\bigomegabar[\vq]} \bigomega[\vq]}}}_1 \lesssim_{\Problet} \sqrt{\frac{\log(1 / h)}{n h^d}}.
  \end{equation*}
It is left to bound
  \begin{align*}
    &\sup_{\vq, \vx} \frac{h^{d / 2}}{\sqrt{\bigomega[\vq]}} \norm[\big]{ \prn{ \qbar[\vq][-1] - \qzero[\vq][-1] } \vp^{(\vv)}(\vx) }_1\\
    &\quad \overset{\text{(a)}}{\lesssim} h^{d + \multindnorm{\vv}} \sup_{\vq, \vx} \norm{ \qbar[\vq][-1] - \qzero[\vq][-1] }_1 \norm{ \vp^{(\vv)}(\vx) }_1\\
    &\quad \overset{\text{(b)}}{\lesssim}_{\Problet} h^{d + \multindnorm{\vv}} \cdot h^{-d} \sqrt{\frac{\log(1 / h)}{n h^d}} \cdot h^{-\multindnorm{\vv}} = \sqrt{\frac{\log(1 / h)}{n h^d}},
  \end{align*}
  where in (a) we used $\inf_{\vq, \vx} \bigomega[\vq] \gtrsim h^{-d - 2 \multindnorm{\vv}}$ by \cref{lem:asympt-variance}, in (b) we used $\norm{ \qbar[\vq][-1] - \qzero[\vq][-1] }_1 \lesssim_{\Problet} h^{-d} \sqrt{\log(1 / h) / \prn{n h^d}}$ by \cref{lem:holder-gram} and $\norm{ \vp^{(\vv)}(\vx) }_1 \lesssim h^{-\multindnorm{\vv}}$ by \cref{ass:gl-local-basis}.

\Cref{eq:elltwobar-bound} follows from~\eqref{eq:elltwobardiff-bound}, \eqref{eq:ellqx-bound} and $\log(1 / h) / \prn{n h^d} = o(1)$.

We will now prove~\eqref{eq:ellhatdiff-bound}. By the triangle inequality,
\begin{align}
  &\sup_{\vq, \vx} \norm[\bigg]{ \frac{\qhat[\vq][-1] \vp^{(\vv)}(\vx)}{\bigomegahat[\vq]^{1 / 2}} - \frac{\qbar[\vq][-1] \vp^{(\vv)}(\vx)}{\bigomegabar[\vq]^{1 / 2}}}_1\\
  &\quad \lesssim \sup_{\vq, \vx} \norm[\bigg]{ \frac{\prn[\big]{ \qhat[\vq][-1] - \qbar[\vq][-1] } \vp^{(\vv)}(\vx)}{\bigomegahat[\vq]^{1 / 2}} }_1 + \sup_{\vq, \vx} \norm[\big]{ \qbar[\vq][-1] \vp^{(\vv)}(\vx)
    \prn[\big]{
    \bigomegahat[\vq]^{-1 / 2}
    -
    \bigomegabar[\vq]^{-1 / 2}
    } }_1. \label{eq:plug-in-two-terms}
\end{align}

To bound the first term in~\eqref{eq:plug-in-two-terms}, recall that $\bigomegahat[\vq] \gtrsim_{\Problet} h^{- d - 2 \multindnorm{\vv}}$ by \cref{lem:asympt-variance} and \eqref{eq:variance-closeness}. Then
\begin{align*}
  &\sup_{\vq, \vx} \norm[\bigg]{ \frac{\prn{ \qhat[\vq][-1] - \qbar[\vq][-1] } \vp^{(\vv)}(\vx)}{\sqrt{\bigomegahat[\vq]}} }_1 \lesssim_{\Problet} h^{d / 2 + \multindnorm{\vv}} \sup_{\vq, \vx} \norm[\big]{\prn{ \qhat[\vq][-1] - \qbar[\vq][-1] } \vp^{(\vv)}(\vx) }_1\\
  &\quad \lesssim h^{d / 2 + \multindnorm{\vv}} \sup_{\vq, \vx} \norm{\qhat[\vq][-1] - \qbar[\vq][-1]}_1 \cdot \norm{ \vp^{(\vv)}(\vx) }_{1}\\
  &\quad \lesssim h^{d / 2} \norm{\qhat[\vq][-1] - \qbar[\vq][-1]}_\infty \lesssim_\Problet h^{d / 2} \cdot h^{-d} \rqhat = h^{-d / 2} \rqhat,
\end{align*}
where in the last inequality we used $\norm{\qhat[\vq] - \qbar[\vq]}_\infty \lesssim_{\Problet} h^d \rqhat$ by assumption, and
\begin{equation*}
  \norm{\qhat[\vq][-1] - \qbar[\vq][-1]}_\infty \leq \norm{\qhat[\vq][-1]}_\infty \cdot \norm{\qhat[\vq] - \qbar[\vq]}_\infty \cdot \norm{\qbar[\vq][-1]}_\infty \lesssim_{\Problet} h^{-d} \prn{ h^d \rqhat } h^{-d} = h^{-d} \rqhat
\end{equation*}
also by assumption and \cref{lem:holder-gram}.

It is left to bound the second term in~\eqref{eq:plug-in-two-terms}:
\begin{align*}
  &\sup_{\vq, \vx} \norm[\big]{ \qbar[\vq][-1] \vp^{(\vv)}(\vx) \prn[\big]{\bigomegahat[\vq]^{-1 / 2} - \bigomegabar[\vq]^{-1 / 2}} }_1\\
  &\quad \leq \sup_{\vq} \norm{ \qbar[\vq][-1] }_1 \sup_{\vx} \norm{ \vp^{(\vv)}(\vx) }_1 \sup_{\vq, \vx} \abs[\big]{ \bigomegahat[\vq]^{-1 / 2} - \bigomegabar[\vq]^{-1 / 2} }\\
  &\quad \overset{\text{(a)}}{\lesssim_{\Problet}} h^{-d} \cdot h^{-\multindnorm{\vv}} \cdot \sup_{\vq, \vx} \abs[\big]{ \bigomegahat[\vq]^{-1 / 2} - \bigomegabar[\vq]^{-1 / 2} }\\
  &\quad \overset{\text{(b)}}{\lesssim_{\Problet}} h^{-d - \multindnorm{\vv}} h^{d / 2 + \multindnorm{\vv}} \romega = h^{-d / 2} \romega,
\end{align*}
where (a) is by \cref{lem:holder-gram} and \cref{ass:gl-local-basis}, (b) is by~\eqref{eq:omega-inverse-root-diff-bound}. \Cref{eq:ellhatdiff-bound} is proven.

\Cref{eq:ellhat-bound} follows from~\eqref{eq:ellhatdiff-bound}, \eqref{eq:elltwobar-bound} and $\rqhat + \romega = o(1)$.
\end{proof}

\begin{lemma}[Hats off]\label{lem:hats-off}
Suppose \cref{ass:dgp-holder,ass:gl-quasi-uniform-partition,ass:gl-local-basis,ass:gl-holder,ass:complexity,ass:gl-approximation-error}, \eqref{eq:gRcohQ}, \eqref{eq:yXRXnK} and \eqref{eq:variance-closeness} hold, $\nu \geq 3$, and $\frac{\log^3 n}{n h^{3 d}} = o(1)$. Define (and fix for all further arguments) $\rhatsoff$ as an arbitrary positive sequence satisfying
  \begin{equation}\label{eq:rhatsoff-def}
    \sqrt{n h^d} \prn{\ruc \romega + \rbah} = o(\rhatsoffNoLink).
  \end{equation}
  Then
  \begin{equation}\label{eq:hats-off}
    \sup_{\vx \in \mathcal{X}} \sup_{\vq \in \mathcal{Q}} \abs{ \hatTqx[\vq] - \bart[\vq] } \lesssim_{\Problet} \sqrt{n h^d} \prn{\ruc \romega + \rbah} = o(\rhatsoffNoLink).
  \end{equation}
\end{lemma}

\begin{proof}
First, note that by \cref{lem:asympt-variance}, \eqref{eq:gRcohQ} and \eqref{eq:variance-closeness}
\begin{align*}
  &\sup_{\vq, \vx} \abs[\bigg]{ \frac{\muhat[\vq][(\vv)]{\vx} - \muzero[\vq][(\vv)]{\vx}}{\sqrt{\bigomegahat[\vq] / n}} - \frac{\muhat[\vq][(\vv)]{\vx} - \muzero[\vq][(\vv)]{\vx}}{\sqrt{\vphantom{\widehat{\boldsymbol{\Omega}}_{\vq}}\bigomegabar[\vq] / n}} }\\
  &\quad \leq \sqrt{n} \sup_{\vq, \vx} \abs[\big]{ \muhat[\vq][(\vv)]{\vx} - \muzero[\vq][(\vv)]{\vx} } \sup_{\vq, \vx} \frac{1}{\sqrt{\bigomegahat[\vq] \bigomegabar[\vq]} \prn[\Big]{ \sqrt{\bigomegahat[\vq]} + \sqrt{\vphantom{\widehat{\boldsymbol{\Omega}}_{\vq}}\bigomegabar[\vq]} }}\\
  &\qquad \cdot \sup_{\vq, \vx} \abs[\big]{ \bigomegahat[\vq]  - \bigomegabar[\vq] }\\
  &\quad \lesssim_{\Problet} \sqrt{n} \cdot h^{-\multindnorm{\vv}} \ruc \cdot \prn{h^{d + 2 \multindnorm{\vv}}}^{3 / 2} \cdot h^{- 2 \multindnorm{\vv} - d} \romega\\
  &\quad = \sqrt{n} h^{d / 2} \ruc \romega.
\end{align*}

By \eqref{eq:yXRXnK}, \cref{ass:gl-local-basis,lem:holder-gram}, \cref{eq:uniform-conv-variance-enough} in \cref{lem:holder-uniform-convergence-variance},
\begin{align*}
  &\sup_{\vq, \vx} \abs[\big]{
    \muhat[\vq][(\vv)]{\vx}
    - \muzero[\vq][(\vv)]{\vx}
    + \vp^{(\vv)}(\vx)\trans \qbar[\vq][-1] \mathbb{E}_{n}\brk{\vp(\vx_i) \eta^{(1)}\prn{\muzero[\vq]{\vx_i}} \riskder[\vq]{y_i}{\eta\prn{\muzero[\vq]{\vx_i}}}}
    }\\
  &\quad \leq
    \sup_{\vq, \vx} \abs{\muhat[\vq][(\vv)]{\vx} - \muzero[\vq][(\vv)]{\vx} - \bahterm[\vq][(\vv)]{\vx}}\\
  &\qquad + \sup_{\vq, \vx} \norm{\vp^{(\vv)}(\vx)}_1 \norm[\big]{\qbar[\vq][-1] - \qzero[\vq][-1]}_{\infty} \norm[\big]{\mathbb{E}_{n}\brk{\vp(\vx_i) \eta^{(1)}\prn{\muzero[\vq]{\vx_i}} \riskder[\vq]{y_i}{\eta\prn{\muzero[\vq]{\vx_i}}}}}_{\infty}\\
  &\quad \lesssim h^{-\multindnorm{\vv}} \rbah + h^{-d - \multindnorm{\vv}} \sqrt{\frac{\log (1 / h)}{n h^d}} \sqrt{\frac{h^{d} \log (1 / h)}{n}}\\
  &\quad = h^{-\multindnorm{\vv}} \rbah + h^{-\multindnorm{\vv}} \frac{\log (1 / h)}{n h^d}
    \labrel{1}[\lesssim] h^{-\multindnorm{\vv}} \rbah,\numberthis\label{eq:TvfoLC}
\end{align*}
where in \labrel{1} we used $\frac{\log (1 / h)}{n h^d} \lesssim \frac{\log n}{n h^d} \lesssim \rbah$.
Combining this with $h^{-2 \multindnorm{\vv} - d} \lesssim_{\Problet} \inf_{\vq, \vx} |\bigomegabar[\vq]|$ (by \cref{lem:asympt-variance})
\begin{equation}
  \sup_{\vq, \vx} \abs[\bigg]{ \frac{\vp^{(\vv)}(\vx)\trans \betahat[\vq] - \muzero[\vq][(\vv)]{\vx}}{\sqrt{\bigomegabar[\vq] / n}} - \bart[\vq] } \lesssim_{\Problet} \sqrt{n} \prn[\big]{ h^{d + 2 \multindnorm{\vv}} }^{1 / 2} h^{- \multindnorm{\vv}} \rbah = \sqrt{n h^d} \rbah.
\end{equation}

It is left to apply the triangle inequality.
\end{proof}

\section{Uniform inference}\label{sec:uniform-inference}

\subsection{Plug-in approximation}\label{sec:plug-in}

Recall that the conditional covariance structure of $\zn(\cdot)$ is
\begin{equation*}
\Expect[\big]{\zn(\vq) \zn(\tilde{\vq})\trans \given \mX_n} = h^{-d} \Sigmacov{\vq}{\tilde{\vq}}
\end{equation*}
with
\begin{equation}\label{eq:jvgoXA}
\begin{aligned}
  \Sigmacov{\vq}{\tilde{\vq}} &:= \Expectlet_n \brk[\big]{\Scov{\vq}{\tilde{\vq}}{\vx_i} \eta^{(1)}(\muzero[\vq]{\vx_i}) \eta^{(1)}(\muzero[\tilde{\vq}]{\vx_i}) \vp(\vx_i) \vp(\vx_i)\trans}
    \in \cont{\mathcal{Q} \times \mathcal{Q}, \mathbb{R}^{K \times K}},\\
  \Scov{\vq}{\tilde{\vq}}{\vx} &:= \Expect[\big]{\riskder[\vq]{y_i}{\eta(\muzero[\vq]{\vx_i})} \riskder[\tilde{\vq}]{y_i}{\eta(\muzero[\tilde{\vq}]{\vx_i})} \given \vx_i = \vx}.
\end{aligned}
\end{equation}

\begin{theorem}[Plug-in approximation]\label{th:plug-in}
  Suppose that all of the following is true:
  \begin{conditions}[wide]
  \item\label{th:plug-in-cond-1} All the conditions of \cref{th:strong-approximation-binscatter-yurinski}\ref{th:yurinski-normal} and \cref{ass:gram-estimate} hold.
  \item\label{th:plug-in-cond-3} The function $\vx \mapsto \Scov{\vq}{\tilde{\vq}}{\vx}$ defined in \eqref{eq:jvgoXA}
    is continuous for all $\vq, \tilde{\vq}$, and
    \begin{equation}\label{eq:tJdMYD}
      \sup_{\vx, \vq, \vq_1 \neq \vq_2} \frac{\abs{\Scov{\vq}{\vq_1}{\vx} - \Scov{\vq}{\vq_2}{\vx}}}{\norm{\vq_1 - \vq_2}} \lesssim 1.
    \end{equation}
  \item\label{th:plug-in-cond-4} There is an estimator
    $\ScovEstDef{\vq}{\tilde{\vq}}{\vx}$ which is a known (not depending on $\Problet_{\mD_n}$) measurable function of $(\mD_n, \vq, \tilde{\vq}, \vx)$, satisfying the bound
    \begin{equation}\label{eq:vJCqsc}
      \sup_{\vq, \tilde{\vq}, \vx} \abs[\big]{\ScovEst{\vq}{\tilde{\vq}}{\vx} - \Scov{\vq}{\tilde{\vq}}{\vx}} \lesssim_{\Problet} \rs\quad\text{with $\rs = o(1)$.}
    \end{equation}
  \end{conditions}

  Then, on the same probability space, there is a mean-zero Gaussian conditionally on $\mD_n$ process $\newtarget{def:znhatst}{\znhatstNoLink}(\vq)$ in $\cont{\mathcal{Q}, \mathbb{R}^K}$, whose distribution is known (does not depend on $\Problet_{\mD_n}$), such that
  \begin{equation}\label{eq:feasible-vector-process}
    \Expect[\Big]{\sup_{\vq \in \mathcal{Q}} \norm[\big]{\znst(\vq) - \znhatst(\vq)}_{\infty} \given \mD_n}
  \lesssim_{\Problet} (\ruc + \rs)^{1 / (2 d_{\mathcal{Q}} + 2)} \sqrt{\log n},
  \end{equation}
  where $\newtarget{def:znst}{\znstNoLink}(\vq)$ is a process in $\cont{\mathcal{Q}, \mathbb{R}^K}$ satisfying
  \begin{equation*}
    (\mX_n, \zn(\cdot)) \overset{d}{=} (\mX_n, \znst(\cdot)),\quad \znst(\cdot) \underset{\mX_n}{\independent} \vy_n.
  \end{equation*}
  Further,
  \begin{align}
    &h^{d / 2} \sup_{\vq \in \mathcal{Q}} \sup_{\vx \in \mathcal{X}} \abs[\bigg]{\frac{\vp^{(\vv)}(\vx)\trans \qhat[\vq][-1]}{\sqrt{\bigomegahat[\vq]}} \znhatst(\vq) - \frac{\vp^{(\vv)}(\vx)\trans \qbar[\vq][-1]}{\sqrt{\bigomegabar[\vq]}} \znst(\vq)}
    \lesssim_{\Problet} \rplugin,\label{eq:agGxAc}\\
    &\text{where}\quad \rplugin := \brk[\big]{(\ruc + \rs)^{1 / (2 d_{\mathcal{Q}} + 2)} + \rqhat + \romega} \sqrt{\log n}.\label{eq:uBkMSP}
  \end{align}
  Equivalently (see \cref{lem:conv-to-zero-in-cond-p-equiv-to-conv-to-zero-in-p}), fixing an arbitrary positive sequence $\RpluginDef$ such that $\rplugin = o(\Rplugin)$, we have
  \begin{equation*}
    \Prob[\bigg]{h^{d / 2} \sup_{\vq \in \mathcal{Q}} \sup_{\vx \in \mathcal{X}} \abs[\bigg]{ \frac{\vp^{(\vv)}(\vx)\trans \qhat[\vq][-1]}{\sqrt{\bigomegahat[\vq]}} \znhatst(\vq) - \frac{\vp^{(\vv)}(\vx)\trans \qbar[\vq][-1]}{\sqrt{\bigomegabar[\vq]}} \znst(\vq) } > \Rplugin \given \data} = o_{\Problet}(1).
  \end{equation*}

\end{theorem}

\begin{remark}\label{rem:plug-in-correspondence}
  One-dimensional processes $Z^{\star}(\vx, \vq)$ and $\widehat{Z}(\vx, \vq)$ in the main paper are defined as follows:
\begin{align*}
  &Z^{\star}(\vx, \vq) := \elltwobarqx[\vq]\trans \znst(\vq) = - h^{d / 2} \frac{\vp^{(\vv)}(\vx)\trans \qbar[\vq][-1]}{\sqrt{\bigomegabar[\vq]}} \znst(\vq),\\
  &\widehat{Z}(\vx, \vq) := \ellhatqx[\vq]\trans \znhatst(\vq) = - h^{d / 2} \frac{\vp^{(\vv)}(\vx)\trans \qhat[\vq][-1]}{\sqrt{\bigomegahat[\vq]}} \znhatst(\vq).
\end{align*}
\end{remark}

The rest of \cref{sec:plug-in} is devoted to proving \cref{th:plug-in}.

\subsubsection{Strategy}
First, $\znst(\cdot)$ exists by Theorem 8.17 (transfer) in \cite{Kallenberg2021} with $(\xi, \eta, \zeta) = (\mX_n, \zn(\cdot), \vy_n)$, $\tilde{\xi} = \xi$.

The rest of the proof will be broken up in two steps.

First, we will prove the existence of $\znhatst(\vq)$. A natural approach is to approximate the covariance structure $\Sigmacov{\vq}{\tilde{\vq}}$ using plug-in, giving an estimated covariance structure $\SigmacovEst{\vq}{\tilde{\vq}}$,
discretize $\mathcal{Q}$ by a suitable $\delta_n$-covering $\mathcal{Q}^{\delta_n}$ and argue that since the discrete versions $\mSigmabar^{\delta_n}$ and $\mSigmahat^{\delta_n}$ of $\Sigmacov{\vq}{\tilde{\vq}}$ and $\SigmacovEst{\vq}{\tilde{\vq}}$ are close, the Gaussian vector $\zn(\vq|_{\mathcal{Q}^{\delta_n}}) \in \mathbb{R}^{K |\mathcal{Q}^{\delta_n}|}$ is close to a ``feasible'' Gaussian vector with covariance $\mSigmahat^{\delta_n}$. Then we can embed this vector into a conditionally Gaussian process $\znhatst(\vq)$. One technical caveat is that a simple plug-in method can lead to $\mSigmahat^{\delta_n}$ not being positive-semidefinite, or to $\SigmacovEst{\vq}{\tilde{\vq}}$ not having a suitable Lipschitz property in $\vq$ and $\tilde{\vq}$ (necessary for controlling the discretization error between a continuous process and a discrete vector). For this reason, we just project $\SigmacovEst{\vq}{\tilde{\vq}}$ onto the space of structures satisfying the desired properties
(program \eqref{eq:bMrBeR}).

The second step is proving \eqref{eq:feasible-vector-process} (given $\znhatst(\vq)$ and $\znst(\vq)$ that are close) by bounding $\sup_{\vq \in \mathcal{Q}} \znst(\vq)$ using the maximal inequality in \cref{lem:maximal-gaussian}, and applying \cref{lem:closeness-of-linear-terms}.

\subsubsection{Existence of the feasible process}\label{sec:existence-of-feasible-process}

As the first step, we will prove that under \cref{th:plug-in-cond-1,th:plug-in-cond-3,th:plug-in-cond-4}, the process $\znhatst(\vq)$ as described in \cref{th:plug-in} exists.

\paragraph{Estimating the covariance structure}
By \cref{lem:holder-gram}, for a large enough constant $\nC{lqkjeh}$ the probability of the event
\begin{equation}\label{eq:large-event-lmax}
\mathcal{A} := \crl{\lambda_{\max}(\Expectlet_n \brk{\vp(\vx_i) \vp(\vx_i)\trans}) \leq \oC{lqkjeh} h^d}
\end{equation}
converges to one. Note that $\Scov{\vq}{\tilde{\vq}}{\vx}$ is bounded uniformly in $\vx, \vq, \tilde{\vq}$ because it is continuous in these arguments and $\mathcal{Q} \times \mathcal{Q} \times \mathcal{X}$ is compact.
Combining with \cref{eq:tJdMYD}, we see that there is a large enough $\nC{ahhQYwqt}$ such that on $\mathcal{A}$
\begin{equation*}
\sup_{\vq, \tilde{\vq}} \norm{\Sigmacov{\vq}{\tilde{\vq}}} \leq \oC{ahhQYwqt} h^d
\end{equation*}
and
\begin{equation*}
  \sup_{\vq, \vq_1 \neq \vq_2} \frac{\abs{\Sigmacov{\vq}{\vq_1} - \Sigmacov{\vq}{\vq_2}}}{\norm{\vq_1 - \vq_2}}
  \leq \oC{ahhQYwqt} h^d.
\end{equation*}

We first approximate $\Sigmacov{\vq}{\tilde{\vq}}$
by
\begin{equation}\label{eq:uylQEQ}
\SigmacovEst{\vq}{\tilde{\vq}} := \Expectlet_n \brk[\big]{\ScovEst{\vq}{\tilde{\vq}}{\vx_i} \eta^{(1)}(\muhat[\vq]{\vx_i}) \eta^{(1)}(\muhat[\tilde{\vq}]{\vx_i}) \vp(\vx_i) \vp(\vx_i)\trans}.
\end{equation}
The latter is known (a function of $\mD_n$ not depending on $\Problet_{\mD_n}$), but may not be positive semidefinite or Lipschitz in the arguments. In order to recover these properties, we solve the following program:
\begin{align*}
  \text{minimize}\quad &\sup_{\vq, \tilde{\vq}} \norm{\mM_{\vq, \tilde{\vq}} - \SigmacovEst{\vq}{\tilde{\vq}}}\numberthis\label{eq:bMrBeR}\\
  \text{over}\quad &\mM_{\vq, \tilde{\vq}} \in \cont{\mathcal{Q} \times \mathcal{Q}, \mathbb{R}^{K \times K}},\\
  &\sup_{\vq, \tilde{\vq}} \norm{\mM_{\vq, \tilde{\vq}}} \leq \oC{ahhQYwqt} h^d,\\
  &\mM_{\vq, \tilde{\vq}} = \mM_{\tilde{\vq}, \vq} = \mM_{\vq, \tilde{\vq}}\trans = \mM_{\tilde{\vq}, \vq}\trans,\quad \vq, \tilde{\vq} \in \mathcal{Q},\\
                   &\prn[\big]{\mM_{\vq_k, \vq_m}}_{k, m = 1}^M \in \mathbb{R}^{M K \times M K} \text{ is symmetric positive semidefinite for any $\prn{\vq_k}_{k = 1}^M$},\\
  &\norm{\mM_{\vq, \vq_1} - \mM_{\vq, \vq_2}} \leq \oC{ahhQYwqt} h^d \norm{\vq_1 - \vq_2}, \quad \vq, \vq_1, \vq_2 \in \mathcal{Q}.
\end{align*}

\begin{lemma}[Solution to the program]\label{lem:solution-to-program}
There exists an exact minimizer $\SigmacovEstPlusDef{\vq}{\tilde{\vq}}$ of \eqref{eq:bMrBeR} which is a (known) measurable function of $\mD_n$ with values in $\cont{\mathcal{Q} \times \mathcal{Q}, \mathbb{R}^{K \times K}}$.
\end{lemma}

\begin{proof}
The space $\cont{\mathcal{Q} \times \mathcal{Q}, \mathbb{R}^{K \times K}}$ with the sup-norm
\begin{equation*}
\sup_{\vq, \tilde{\vq}} \norm{\mM_{\vq, \tilde{\vq}}}
\end{equation*}
is Polish, and the feasible set is non-empty and compact by the multi-dimensional Arzel\`a--Ascoli theorem. Indeed, the uniform bound and Lipschitz constraints give equicontinuity, and the symmetry, positive-semidefiniteness and Lipschitz constraints are closed under uniform convergence.

The function
\begin{equation*}
\prn[\big]{\mathcal{X}^n \times \mathcal{Y}, \cont{\mathcal{Q} \times \mathcal{Q}, \mathbb{R}^{K \times K}}} \ni \prn[\big]{\mD_n, \mM_{\vq, \tilde{\vq}}} \mapsto \sup_{\vq, \tilde{\vq}} \norm{\mM_{\vq, \tilde{\vq}} - \SigmacovEst{\vq}{\tilde{\vq}}} \in \mathbb{R}
\end{equation*}
is a Carath\'eodory function, that is, measurable in $\mD_n$ and continuous in $\mM_{\vq, \tilde{\vq}}$.
By the Measurable Maximum Theorem 18.19 in \cite{aliprantis2006hitchhiker} (here minimum), the result follows.
\end{proof}

Now we prove that $\SigmacovEstPlus{\vq}{\tilde{\vq}}$ is sufficiently close to $\Sigmacov{\vq}{\tilde{\vq}}$.

\begin{lemma}[Consistency of $\SigmacovEstPlus{\vq}{\tilde{\vq}}$]\label{lem:ANyFZe}
The function $\SigmacovEstPlus{\vq}{\tilde{\vq}}$ defined in \cref{lem:solution-to-program} satisfies
\begin{equation*}
\sup_{\vq, \tilde{\vq}} \norm[\big]{\SigmacovEstPlus{\vq}{\tilde{\vq}} - \Sigmacov{\vq}{\tilde{\vq}}} \lesssim_{\Problet} h^d \prn{\ruc + \rs}.
\end{equation*}
\end{lemma}

\begin{proof}
From \cref{eq:vJCqsc,eq:gRcohQ,lem:holder-gram} we have
\begin{equation}\label{eq:VHEuxh}
\sup_{\vq, \tilde{\vq}} \norm[\big]{\Sigmacov{\vq}{\tilde{\vq}} - \SigmacovEst{\vq}{\tilde{\vq}}} \lesssim_{\Problet} h^d \prn{\ruc + \rs}.
\end{equation}
But $\Sigmacov{\vq}{\tilde{\vq}}$ lies in the feasibility set of \eqref{eq:bMrBeR} on $\mathcal{A}$, which implies that
\begin{equation*}
\sup_{\vq, \tilde{\vq}} \norm[\big]{\SigmacovEstPlus{\vq}{\tilde{\vq}} - \SigmacovEst{\vq}{\tilde{\vq}}} \leq \sup_{\vq, \tilde{\vq}} \norm[\big]{\Sigmacov{\vq}{\tilde{\vq}} - \SigmacovEst{\vq}{\tilde{\vq}}}\quad\text{on $\mathcal{A}$}.
\end{equation*}
Since $\Problet(\mathcal{A}) = 1 - o(1)$, we have
\begin{equation*}
\sup_{\vq, \tilde{\vq}} \norm[\big]{\SigmacovEstPlus{\vq}{\tilde{\vq}} - \SigmacovEst{\vq}{\tilde{\vq}}} \lesssim_{\Problet} h^d \prn{\ruc + \rs},
\end{equation*}
and by the triangle inequality the result follows.
\end{proof}

\paragraph{Discretization and controlling the deviations of the infeasible process}
Having estimated the covariance structure, we proceed with the same steps as in \cref{lem:yurinski-coupling-K-dimensional}. First, we discretize $\mathcal{Q}$. Let $\mathcal{Q}^{\delta_n} := \crl{\vq_1, \ldots, \vq_{|\mathcal{Q}^{\delta_n}|}}$ be an internal $\delta_n$-covering\footnote{This means $\mathcal{Q}_n^{\delta} \subset \mathcal{Q}$ and balls with radius $\delta_n$ centered at these points cover $\mathcal{Q}$.} of $\mathcal{Q}$ with respect to the 2-norm $\|\cdot\|$ of cardinality $|\mathcal{Q}^{\delta_n}| \lesssim 1 / \delta_n^{d_{\mathcal{Q}}}$, where $\delta_n$ is chosen later. As already proven in \cref{lem:yurinski-coupling-K-dimensional}, there is a bound on the deviations of $\znst(\cdot)$:
\begin{equation*}
  \Expect[\Big]{\sup_{\norm{\vq - \tilde{\vq}} \leq \delta_n} \norm{\znst(\vq) - \znst(\tilde{\vq})}_{\infty} \given \mX_n}
  = \Expect[\Big]{\sup_{\norm{\vq - \tilde{\vq}} \leq \delta_n} \norm{\znst(\vq) - \znst(\tilde{\vq})}_{\infty} \given \mD_n}
  \lesssim \sqrt{\delta_n \log \frac{K}{\delta_n}}\quad\text{on $\mathcal{A}$,}
\end{equation*}
where in the equality we used that $\znst(\cdot)$ is independent of $\vy_n$ conditionally on $\mX_n$.

\paragraph{Closeness of the discrete vectors}
Consider now the conditionally Gaussian vector $\znst(\vq|_{\mathcal{Q}^{\delta_n}}) \in \mathbb{R}^{K |\mathcal{Q}^{\delta_n}|}$ and denote
\begin{equation*}
\mSigmabar^{\delta_n} := h^d \Var{\znst(\vq|_{\mathcal{Q}^{\delta_n}}) \given \mX_n} = \prn{\mSigmabar_{\vq_k, \vq_m}}_{k, m = 1}^{|\mathcal{Q}^{\delta_n}|} \in \mathbb{R}^{K |\mathcal{Q}^{\delta_n}| \times K |\mathcal{Q}^{\delta_n}|},
\end{equation*}
so that
\begin{equation*}
\znst(\vq|_{\mathcal{Q}^{\delta_n}}) = h^{-d / 2} \prn{\mSigmabar^{\delta_n}}^{1 / 2} \vxi_n,\quad\text{where }\Problet_{\vxi_n | \mX_n} = \mathcal{N}\prn[\big]{0, \mI_{\mathbb{R}^{K |\mathcal{Q}^{\delta_n}|}}},
\end{equation*}
and $\vxi_n$ is independent of $\vy_n$ given $\mX_n$. (Since $\Problet_{\vxi_n | \mX_n}$ does not depend on $\mX_n$, this vector is in fact independent of $\mD_n$.) Discretizing $\SigmacovEstPlus{\vq}{\tilde{\vq}}$ in the same way, we can put
\begin{equation*}
\mSigmahat^{\delta_n} := \prn[\big]{\SigmacovEstPlus{\vq_k}{\vq_m}}_{k, m = 1}^{|\mathcal{Q}^{\delta_n}|} \in \mathbb{R}^{K |\mathcal{Q}^{\delta_n}| \times K |\mathcal{Q}^{\delta_n}|}.
\end{equation*}
Since the functions $\mSigmabar_{\vq, \tilde{\vq}}$ and $\SigmacovEstPlus{\vq}{\tilde{\vq}}$ are close, the two matrices are close as well, which we will make precise in the following \lcnamecref{lem:phBAhW}.

\begin{lemma}[The matrices $\mSigmabar^{\delta_n}$ and $\mSigmahat^{\delta_n}$ are close]\label{lem:phBAhW}
  We have
  \begin{equation*}
    \norm[\big]{\mSigmabar^{\delta_n} - \mSigmahat^{\delta_n}} \lesssim_{\Problet} |\mathcal{Q}^{\delta_n}| h^d \prn{\ruc + \rs}.
  \end{equation*}
\end{lemma}

\begin{proof}
We can bound the Frobenius norm using the bound on each element of the matrix:
\begin{equation*}
\norm[\big]{\mSigmabar^{\delta_n} - \mSigmahat^{\delta_n}} \leq \norm[\big]{\mSigmabar^{\delta_n} - \mSigmahat^{\delta_n}}_F \leq |\mathcal{Q}^{\delta_n}| \sup_{\vq, \tilde{\vq}} \norm[\big]{\SigmacovEstPlus{\vq}{\tilde{\vq}} - \Sigmacov{\vq}{\tilde{\vq}}}.
\end{equation*}
It is left to combine with \cref{lem:ANyFZe}.
\end{proof}

Applying \cref{lem:maximal-ineq-gaussian-vectors}, we then have
\begin{equation}
  \begin{aligned}
    h^{-d / 2} \Expect[\big]{
    \norm[\big]{\prn{\mSigmabar^{\delta_n}}^{1 / 2} \vxi_n - \prn{\mSigmahat^{\delta_n}}^{1 / 2} \vxi_n}_{\infty}
    \given
    \mD_n
    } &\leq 2 h^{-d / 2} \sqrt{\log(K |\mathcal{Q}^{\delta_n}|)} \norm{\mSigmabar^{\delta_n} - \mSigmahat^{\delta_n}}^{1 / 2}\\
    &\lesssim_{\Problet} |\mathcal{Q}^{\delta_n}|^{1 / 2} \prn{\ruc + \rs}^{1 / 2} \sqrt{\log(K |\mathcal{Q}^{\delta_n}|)}.
  \end{aligned}
\end{equation}
\paragraph{Embedding the conditionally Gaussian vector $\prn{\mSigmahat^{\delta_n}}^{1 / 2} \vxi_n$ into a conditionally Gaussian process}
By the same argument as in \cref{lem:yurinski-coupling-K-dimensional}, there exists a law
on $(\mathcal{X} \times \mathcal{Y})^n \times \cont{\mathcal{Q}, \mathbb{R}^K}$ which is the joint law of $\mD_n' \overset{d}{=} \mD_n$ and a conditionally on $\mD_n'$ mean-zero Gaussian process
$\znhatstprDef(\vq)$
with the conditional covariance structure
\begin{equation*}
\Expect{\znhatstpr(\vq) \znhatstpr(\tilde{\vq})\trans \given \mD_n'} = h^{-d} \SigmacovEstPlusPr{\vq}{\tilde{\vq}}
\end{equation*}
(where $\SigmacovEstPlusPr{\vq}{\tilde{\vq}}$ is the same function of the data as $\SigmacovEstPlus{\vq}{\tilde{\vq}}$ with $\mD_n$ replaced by $\mD_n'$).
Continuity follows by the Kolmogorov-Chentsov theorem:
\begin{equation}\label{eq:YHasPT}
\Expect[\big]{\norm[\big]{\znhatstpr(\vq) - \znhatstpr(\tilde{\vq})}^a \given \mD_n'} \leq C_n \norm{\vq - \tilde{\vq}}^{d_{\mathcal{Q}} + b},
\end{equation}
for $a, b > 0$ chosen as follows. The vector $\znhatstpr(\vq) - \znhatstpr(\tilde{\vq})$ is conditionally on $\mD_n'$ Gaussian with covariance
\begin{equation*}
h^{-d} \prn[\big]{\SigmacovEstPlusPr{\vq}{\vq} - 2 \SigmacovEstPlusPr{\vq}{\tilde{\vq}} + \SigmacovEstPlusPr{\tilde{\vq}}{\tilde{\vq}}}.
\end{equation*}
Hence, we have for any $m > 0$
\begin{align*}
  &\Expect[\big]{\norm[\big]{\znhatstpr(\vq) - \znhatstpr(\tilde{\vq})}^{2 m} \given \mD_n'}\\
  &\quad \leq h^{-d m} \norm[\big]{\SigmacovEstPlusPr{\vq}{\vq} - 2 \SigmacovEstPlusPr{\vq}{\tilde{\vq}} + \SigmacovEstPlusPr{\tilde{\vq}}{\tilde{\vq}}}^m \Expect{\norm{\vxi_K}^{2 m}}\\
&\quad \leq \prn{2 \oC{ahhQYwqt}}^{m} \Expect{\norm{\vxi_K}^{2 m}} \norm{\vq - \tilde{\vq}}^m,
\end{align*}
where $\vxi_K \sim \mathcal{N}(0, \mI_K)$.
So we can take, for example, $m = d_{\mathcal{Q}} + 1$, $a = 2 m$, $b = 1$ in \eqref{eq:YHasPT}.

Projecting this $\cont{\mathcal{Q}, \mathbb{R}^K}$-process onto $\mathcal{Q}^{\delta_n}$, we obtain a vector $\znhatstpr\prn[\big]{\vq|_{\mathcal{Q}^{\delta_n}}} \in \mathbb{R}^{K |\mathcal{Q}^{\delta_n}|}$
such that
\begin{equation*}
\prn[\big]{\mD_n', \znhatstpr\prn[\big]{\vq|_{\mathcal{Q}^{\delta_n}}}} \overset{d}{=} \prn[\big]{\mD_n, h^{-d / 2} \prn{\mSigmahat^{\delta_n}}^{1 / 2} \vxi_n}.
\end{equation*}
Since $\cont{\mathcal{Q}, \mathbb{R}^K}$ is Polish, by Theorem~8.17 (transfer) in \cite{Kallenberg2021} on our probability space there exists a random element $\znhatst(\vq) \in \cont{\mathcal{Q}, \mathbb{R}^K}$ such that
\begin{equation*}
\prn[\big]{\mD_n', \znhatstpr\prn[\big]{\vq|_{\mathcal{Q}^{\delta_n}}}, \znhatstpr(\vq)} \overset{d}{=} \prn[\big]{\mD_n, h^{-d / 2} \prn{\mSigmahat^{\delta_n}}^{1 / 2} \vxi_n, \znhatst(\vq)}.
\end{equation*}
In particular, almost surely $h^{- d / 2} \prn{\mSigmahat^{\delta_n}}^{1 / 2} \vxi_n$ is the projection of $\znhatst$ on $\mathcal{Q}^{\delta_n}$.

\paragraph{Controlling the deviations of $\znhatstNoLink(\vq)$}
Consider the stochastic process $X_n$ with index set
\begin{equation*}
  T := \crl{ (l, \vq, \tilde{\vq})\colon l \in \{1, 2, \ldots, K\}, \vq, \tilde{\vq} \in \mathcal{Q}, \|\vq - \tilde{\vq}\| \leq \delta_n }
\end{equation*}
whose value at $t = (l, \vq, \tilde{\vq})$ is $X_{n, t} := \znhatstComp{l}(\vq) - \znhatstComp{l}(\tilde{\vq})$. It is a separable (because each $\znhatstComp{l}(\cdot)$ has continuous trajectories) mean-zero Gaussian conditionally on $\mD_n$ process with the index set $T$ considered a metric space: $\dist \prn[\big]{(l, \vq, \tilde{\vq}), (l', \vq', \tilde{\vq}')} = \|\vq - \vq'\| + \|\tilde{\vq} - \tilde{\vq}'\| + \Indlet\{l \neq l'\}$.

The quantity $\sigma(X_n)^2$ is defined and bounded as follows:
\begin{equation*}
  \sigma(X_n)^2 := \sup_{t \in T} \Expect{X_{n, t}^2 \given \mD_n} = \sup_{\norm{\vq - \tilde{\vq}} \leq \delta_n} \max_l \Expect[\big]{ \prn[\big]{\znhatstComp{l}(\vq) - \znhatstComp{l}(\tilde{\vq})}^2 \given \mD_n}
  \lesssim \delta_n,
\end{equation*}
because $\znhatstComp{l}(\vq) - \znhatstComp{l}(\tilde{\vq})$ is conditionally on $\mD_n$ mean-zero normal variable with variance
\begin{align*}
  &h^{-d} \prn[\big]{\prn[\big]{\SigmacovEstPlusSc{\vq}{\vq}}_{l l} - 2 \prn[\big]{\SigmacovEstPlusSc{\vq}{\tilde{\vq}}}_{l l} + \prn[\big]{\SigmacovEstPlusSc{\tilde{\vq}}{\tilde{\vq}}}_{l l}}\\
  &\quad \leq h^{-d} \abs[\big]{\prn[\big]{\SigmacovEstPlusSc{\vq}{\vq}}_{l l} - \prn[\big]{\SigmacovEstPlusSc{\vq}{\tilde{\vq}}}_{l l}} + h^{-d} \abs[\big]{\prn[\big]{\SigmacovEstPlusSc{\vq}{\tilde{\vq}}}_{l l} - \prn[\big]{\SigmacovEstPlusSc{\tilde{\vq}}{\tilde{\vq}}}_{l l}}\\
  &\quad \leq h^{-d} \norm[\big]{\SigmacovEstPlus{\vq}{\vq} - \SigmacovEstPlus{\vq}{\tilde{\vq}}}
    + h^{-d} \norm[\big]{\SigmacovEstPlus{\vq}{\tilde{\vq}} - \SigmacovEstPlus{\tilde{\vq}}{\tilde{\vq}}}\\
  &\quad \leq 2 \oC{ahhQYwqt} \norm{\vq - \tilde{\vq}}.
\end{align*}

Next, we define and bound the semimetric $\rho(t, t')$:
\begin{align*}
  \rho(t, t')^2 &:= \Expect{(X_{n, t} - X_{n, t'})^2 \given \mD_n}\\
             &= \Expect[\big]{\prn[\big]{\prn[\big]{ \znhatstComp{l}(\vq) - \znhatstComp{l}(\tilde{\vq}) } - \prn[\big]{ \znhatstComp{l'}(\vq') - \znhatstComp{l'}(\tilde{\vq}') }}^2 \given \mD_n}\\
             &\lesssim \Expect[\big]{\prn[\big]{ \znhatstComp{l}(\vq) - \znhatstComp{l'}(\vq') }^2 \given \mD_n}
               -
               \Expect[\big]{\prn[\big]{ \znhatstComp{l}(\tilde{\vq}) - \znhatstComp{l'}(\tilde{\vq}') }^2 \given \mD_n}.
\end{align*}
The first term on the right is bounded the following way: if $l \neq l'$,
\begin{align*}
  &\Expect[\big]{\prn[\big]{ \znhatstComp{l}(\vq) - \znhatstComp{l'}(\vq') }^2 \given \mD_n}\\
  &\quad \lesssim \Expect[\big]{\prn[\big]{ \znhatstComp{l}(\vq) - \znhatstComp{l}(\vq') }^2 \given \mD_n}
    + \Expect[\big]{\znhatstComp{l}(\vq')^2 \given \mD_n} + \Expect[\big]{\znhatstComp{l'}(\vq')^2 \given \mD_n}\\
  &\quad \lesssim \norm{\vq - \vq'} + 1,
\end{align*}
and if $l = l'$,
\begin{equation*}
\Expect[\big]{\prn[\big]{ \znhatstComp{l}(\vq) - \znhatstComp{l'}(\vq') }^2 \given \mD_n} \lesssim \norm{\vq - \vq'}.
\end{equation*}
The term $\Expect[\big]{\prn[\big]{ \znhatstComp{l}(\tilde{\vq}) - \znhatstComp{l'}(\tilde{\vq}') }^2 \given \mD_n}$ is bounded the same way, and we conclude that
\begin{equation*}
\rho(t, t')^2 \leq \oC{lkjwh} \prn[\big]{\norm{\vq - \vq'} + \norm{\tilde{\vq} - \tilde{\vq}'} + \Indlet\crl{l \neq l'}} = \oC{lkjwh} \mathrm{dist}\prn[\big]{(l, \vq, \tilde{\vq}), (l', \vq', \tilde{\vq}')}.
\end{equation*}
for some positive constant $\nC{lkjwh}$. This means that, for any $\varepsilon > 0$, an $\oC{lkjwh}^{-1} \varepsilon^2$-covering of $T$ with respect to $\mathrm{dist}(\cdot)$ induces an $\varepsilon$-covering of $T$ with respect to $\rho$, and hence
\begin{equation*}
N(T, \rho, \varepsilon) \leq N\prn[\big]{T, \mathrm{dist}(\cdot), \oC{lkjwh}^{-1} \varepsilon^2}.
\end{equation*}
But the right-hand side is clearly $O(K \varepsilon^{-2})$ because $T$ is a subset of $\crl{1, \ldots, K} \times \mathcal{Q} \times \mathcal{Q}$, so for $\varepsilon \in (0, 1)$
\begin{equation*}
\log N(T, \rho, \varepsilon) \lesssim \log \frac{K}{\varepsilon}.
\end{equation*}

Applying \cref{lem:maximal-gaussian} gives
\begin{align*}
  &\Expect[\Big]{\sup_{t \in T} \abs{X_{n, t}} \given \mD_n} = \Expect[\Big]{\sup_{\norm{\vq - \tilde{\vq}} \leq \delta_n} \norm[\big]{\znhatst(\vq) - \znhatst(\tilde{\vq})}_{\infty} \given \mD_n} \lesssim \sigma(X_n) + \int_0^{2 \sigma(X_n)} \sqrt{\log N(T, \rho, \varepsilon)}\,\mathrm{d}\varepsilon\\
  &\quad \lesssim \sigma(X_n) + \int_0^{2 \sigma(X_n)} \sqrt{\log (K / \varepsilon)}\,\mathrm{d}\varepsilon \lesssim \sigma(X_n) + \int_0^{2 \sigma(X_n)} \prn[\big]{\sqrt{\log K} + \sqrt{\log (\varepsilon^{-1})}}\,\mathrm{d}\varepsilon\\
  &\quad \lesssim \sigma(X_n) \sqrt{\log K} + \frac{1}{\sqrt{\log (\sigma(X_n)^{-1})}} \int_0^{2 \sigma(X_n)} \log (\varepsilon^{-1})\,\mathrm{d}\varepsilon\\
  &\quad \lesssim \sigma(X_n) \sqrt{\log K} + \frac{1 + \log (\sigma(X_n)^{-1})}{\sqrt{\log (\sigma(X_n)^{-1})}} \sigma(X_n) \lesssim \sigma(X_n) \sqrt{\log (K \sigma(X_n)^{-1})} \lesssim \sqrt{\delta_n \log \frac{K}{\delta_n}},
\end{align*}
where in the last inequality we used the bound $\sigma(X_n) \lesssim \sqrt{\delta_n}$ from above and that $x \mapsto x \log (1 / x)$ is increasing for sufficiently small $x$.

\paragraph{Choosing $\delta_n$ and conclusion}
Combining the bounds obtained above, we can write
\begin{align*}
  &\Expect[\Big]{\sup_{\vq \in \mathcal{Q}} \norm[\big]{\znst(\vq) - \znhatst(\vq)}_{\infty} \given \mD_n}\\
  &\quad \leq \Expect[\Big]{\sup_{\norm{\vq - \tilde{\vq}} \leq \delta_n} \norm[\big]{\znst(\vq) - \znst(\tilde{\vq})}_{\infty} \given \mD_n}
    +
    h^{- d / 2}  \Expect[\big]{\norm[\big]{\prn{\mSigmabar^{\delta_n}}^{1 / 2} \vxi_n - \prn{\mSigmahat^{\delta_n}}^{1 / 2} \vxi_n}_{\infty} \given \mD_n}\\
  &\qquad + \Expect[\Big]{\sup_{\norm{\vq - \tilde{\vq}} \leq \delta_n} \norm[\big]{\znhatst(\vq) - \znhatst(\tilde{\vq})}_{\infty} \given \mD_n}\\
  &\quad \lesssim_{\Problet} \prn[\bigg]{\sqrt{\delta_n} + \frac{(\ruc + \rs)^{1 / 2}}{\delta_n^{d_{\mathcal{Q}} / 2}}} \sqrt{\log \frac{K}{\delta_n}}.
\end{align*}
Choose the approximately optimal
\begin{equation*}
\delta_n := (\ruc + \rs)^{1 / (d_{\mathcal{Q}} + 1)},
\end{equation*}
giving
\begin{equation*}
  \Expect[\Big]{\sup_{\vq \in \mathcal{Q}} \norm[\big]{\znst(\vq) - \znhatst(\vq)}_{\infty} \given \mD_n}
  \lesssim_{\Problet} (\ruc + \rs)^{1 / (2 d_{\mathcal{Q}} + 2)} \sqrt{\log n}.
\end{equation*}

\subsubsection{Proof of \eqref{eq:agGxAc}}

Having proven the existence of $\znhatst(\cdot)$,
we will now prove that under \cref{th:plug-in-cond-1,th:plug-in-cond-3,th:plug-in-cond-4}, \cref{eq:agGxAc} holds. We build on the proof of Theorem~6.3 in \parencite{Cattaneo-Farrell-Feng_2020_AOS}.

\begin{lemma}[Bounding the supremum of the Gaussian process]\label{lem:bounding-the-sup-of-gaussian-process}
  In the setting of \cref{th:plug-in}, we have $\sup_{\vq \in \mathcal{Q}} \norm{ \zn(\vq) }_{\infty} \lesssim_{\Problet} \sqrt{\log K}$.
\end{lemma}

\begin{proof}
Consider the stochastic process $X_n$ defined for $t = (l, \vq) \in T$ with
\begin{equation*}
  T := \crl{ (l, \vq): l \in \crl{ 1, 2, \ldots, K }, \vq \in \mathcal{Q} }
\end{equation*}
as $X_{n, t} := Z_{n, l}(\vq)$. It is a separable mean-zero Gaussian process with the index set $T$ considered a metric space: $\dist \prn[\big]{ (l, \vq), (l', \vq') } = \|\vq - \vq'\| + \Indlet\{l \neq l'\}$. Note that
\begin{align*}
  &\sigma(X_n)^2 := \sup_{t \in T} \Expect*{X_{n, t}^2 \given \mX_n} = \sup_{\vq \in \mathcal{Q}} \max_l \Expect*{Z_{n, l}^2(\vq) \given \mX_n}\\
  &\quad = h^{-d} \sup_{\vq \in \mathcal{Q}} \max_l \frac{1}{n} \sum_{i = 1}^n \Expect[\big]{\riskder[\vq]{y_i}{\eta \prn{ \muzero[\vq]{\vx_i}}}^2 \given \vx_i} \eta^{(1)}\prn{ \muzero[\vq]{\vx_i} }^2 p_l(\vx_i)^2\\
  &\quad \lesssim h^{-d} \Expectlet_n \brk{p_l(\vx_i)^2} \leq h^{-d} \sup_{\alphavec \in \mathcal{S}^{K - 1}} \Expectlet_n \brk[\big]{\alphavec\trans \vp(\vx_i) \vp(\vx_i)\trans \alphavec} \lesssim 1\quad\text{w.\,p.\,a. 1}.
\end{align*}
Next, we will bound
\begin{equation*}
  \rho(t, t')^{2} := \Expect[\big]{\prn{ X_{n, t} - X_{n, t'} }^2 \given \mX_n} = \Expect[\big]{\prn{ Z_{n, l}(\vq) - Z_{n, l'}(\vq') }^2 \given \mX_n}.
\end{equation*}
If $l \neq l'$,
\begin{align*}
  &\Expect[\big]{\prn{ Z_{n, l}(\vq) - Z_{n, l'}(\vq') }^2 \given \mX_n}\\
  &\quad = \frac{1}{n} \sum_{i = 1}^n h^{-d} \Expectlet\left[ \left( \riskder[\vq]{y_i}{\eta \left( \muzero[\vq]{\vx_i} \right)} \eta^{(1)}\left( \muzero[\vq]{\vx_i} \right) p_l(\vx_i) \right. \right.\\
  &\qquad \left. \left. - \riskder[\vq']{y_i}{\eta \left( \muzero[\vq']{\vx_i} \right)} \eta^{(1)}\left( \muzero[\vq']{\vx_i} \right) p_{l'}(\vx_i)\right)^2 \given \mX_n \right]\\
  &\quad = h^{-d} \Expectlet_n \brk[\big]{\Expect[\big]{\prn{A_n(\vq, y_i, \vx_i) - A_n(\vq', y_i, \vx_i)}^2 \given \mX_n} p_l(\vx_i)^2}\\
  &\qquad + h^{-d} \Expectlet_n \brk[\big]{\Expect[\big]{A_n(\vq', y_i, \vx_i)^2 \given \mX_n} \prn[\big]{ p_{l}(\vx_i)^2 + p_{l'}(\vx_i)^2}}\\
  &\quad \lesssim h^{-d} \norm{ \vq - \vq' } \Expectlet_n \brk[\big]{p_l(\vx_i)^2} + h^{-d} \Expectlet_n \brk[\big]{p_{l}(\vx_i)^2 + p_{l'}(\vx_i)^2}\\
  &\quad \lesssim \norm{ \vq - \vq' } + 1\quad\text{w.\,p.\,a. 1,}
\end{align*}
where we denoted $A_n(\vq, y_i, \vx_i)  := \riskder[\vq]{y_i}{\eta \prn{ \muzero[\vq]{\vx_i} }} \eta^{(1)}\prn{ \muzero[\vq]{\vx_i} }$ to simplify notation. Similarly, if $l = l'$,
\begin{equation*}
  \Expect[\big]{\prn{ Z_{n, l}(\vq) - Z_{n, l'}(\vq') }^2 \given \mX_n} \lesssim \norm{ \vq - \vq' }\quad\text{w.\,p.\,a. 1.}
\end{equation*}
We conclude
\begin{equation*}
  \rho(t, t')^{2} \lesssim \dist \prn[\big]{ (l, \vq), (l', \vq') }.
\end{equation*}
This means that an $\varepsilon^2$-covering of $T$ with respect to $\dist(\cdot)$ induces an $\varepsilon$-covering of $T$ with respect to $\rho$, and hence
\begin{equation}\label{eq:prelim-n-bound-plug-in}
  N(T, \rho, \varepsilon) \leq N\prn{T, \dist(\cdot), \varepsilon^2}.
\end{equation}
On the other hand, since $\mathcal{Q}$ does not depend on $n$, clearly for (sufficiently small) $\tilde{\varepsilon} > 0$, $N(T, \dist(\cdot), \tilde{\varepsilon}) \lesssim K \prn{ C_1 / \tilde{\varepsilon} }^{C_2}$, where $C_1$ and $C_2$ are both constants. Combining this with \eqref{eq:prelim-n-bound-plug-in} we get
\begin{equation*}
  \log N(T, \rho, \varepsilon) \lesssim \log\prn{K / \varepsilon}.
\end{equation*}

Now we apply \cref{lem:maximal-gaussian}:
\begin{align*}
  &\Expect[\Big]{\sup_{t \in T} \abs{ X_{n, t} } \given \mX_n} \lesssim \sigma(X_n) + \int_0^{2 \sigma(X_n)} \sqrt{\log N(T, \rho, \varepsilon)}\,\mathrm{d}\varepsilon\\
  &\quad \lesssim \sigma(X_n) + \sigma(X_n) \sqrt{\log\prn{K / \sigma(X_n)}}  \lesssim_1 \sqrt{\log K}\quad \text{w.\,p.\,a. 1},
\end{align*}
where in $\lesssim_1$ we used our bound $\sigma(X_n) \lesssim 1$ above. Rewriting, we obtain
\begin{equation}\label{eq:Zq-Gaussian-maximal-bound}
  \Expect[\Big]{\sup_{\vq \in \mathcal{Q}} \norm{ \zn(\vq) }_{\infty} \given \mX_n} \lesssim \sqrt{\log K}\quad\text{w.\,p.\,a. 1}.
\end{equation}
By Markov's inequality this gives $\sup_{\vq \in \mathcal{Q}} \norm{ \zn(\vq) }_{\infty} \lesssim_{\Problet} \sqrt{\log K}$.
\end{proof}

By the triangle inequality, \cref{eq:agGxAc} will follow from the two bounds
\begin{align}
  &h^{d / 2} \sup_{\vq, \vx} \abs[\Bigg]{ \frac{\vp^{(\vv)}(\vx)\trans \qhat[\vq][-1]}{\sqrt{\bigomegahat[\vq]}} \znhatst(\vq) - \frac{\vp^{(\vv)}(\vx)\trans \qhat[\vq][-1]}{\sqrt{\bigomegahat[\vq]}} \znst(\vq) } \lesssim_{\Problet} (\ruc + \rs)^{1 / (2 d_{\mathcal{Q}} + 2)} \sqrt{\log n},\label{eq:plug-in-first-term-bound}\\
  &h^{d / 2} \sup_{\vq, \vx} \abs[\Bigg]{ \frac{\vp^{(\vv)}(\vx)\trans \qhat[\vq][-1]}{\sqrt{\bigomegahat[\vq]}} \znst(\vq) - \frac{\vp^{(\vv)}(\vx)\trans \qbar[\vq][-1]}{\sqrt{\bigomegabar[\vq]}} \znst(\vq) } \lesssim_{\Problet} (\rqhat + \romega) \sqrt{\log n}.\label{eq:plug-in-second-term-bound}
\end{align}

To prove~\eqref{eq:plug-in-first-term-bound},
combine $\sup_{\vq, \vx} \norm[\big]{ \ellhatqx[\vq] }_1 \lesssim_{\Problet} 1$ by \cref{lem:closeness-of-linear-terms} with~\eqref{eq:feasible-vector-process}.

To prove~\eqref{eq:plug-in-second-term-bound}, combine \cref{lem:bounding-the-sup-of-gaussian-process} with $\sup_{\vq, \vx} \norm[\big]{ \ellhatqx[\vq] - \elltwobarqx[\vq] }_1 \lesssim_{\Problet} \rqhat + \romega$ by \cref{lem:closeness-of-linear-terms}.

This completes the proof of \cref{th:plug-in}.

\subsection{Confidence bands}

As a corollary of \cref{th:plug-in}, we can obtain the following result, whose proof is based on \cite{Belloni-Chernozhukov-Chetverikov-FernandezVal_2019_JoE}.

\begin{theorem}[Confidence bands]\label{th:confidence-bands}
  Suppose all conditions of \cref{th:plug-in} are true, and $\prn{\Rplugin + \rstr + \rhatsoff} \sqrt{\log n} = o(1)$. Then $\Prob[\big]{\sup_{\vq, \vx} | \hatTqx[\vq] | < \kst(1 - \alpha)} = 1 - \alpha + o(1)$, where $\newtarget{def:k-star}{\kstNoLink}(\eta)$ is the conditional $\eta$-quantile of $\sup_{\vq, \vx} \abs[\big]{\ellhatqx[\vq]\trans \znhatst(\vq)}$ given the data. Equivalently, $\Prob[\big]{\muzero[\vq][(\vv)]{\vx} \in \cb[\vv]{1 - \alpha}{\vx}{\vq}} = 1 - \alpha + o(1)$ with
\begin{equation}
\label{eq:cb-def}
\cb[\vv]{1 - \alpha}{\vx}{\vq} := \prn[\Big]{\muhat[\vq][(\vv)]{\vx} - \kst(1 - \alpha) \sqrt{\bigomegahat[\vq][\vv] / n}, \muhat[\vq][(\vv)]{\vx} + \kst(1 - \alpha) \sqrt{\bigomegahat[\vq][\vv] / n}}.
\end{equation}
\end{theorem}

To prove this result we first introduce some auxiliary lemmas. To simplify the exposition, let
\begin{align}
  &\bigvNoLink := \sup_{\vq, \vx} \abs{ \hatTqx[\vq] },\label{eq:v-def}\\
  &\bigvstNoLink := \sup_{\vq, \vx} \abs[\bigg]{ h^{d / 2} \frac{\vp^{(\vv)}(\vx)\trans \qhat[\vq][-1]}{\sqrt{\bigomegahat[\vq]}} \znhatst(\vq) },\label{eq:v-star-def}\\
  &\bigvtildeNoLink := \sup_{\vq, \vx} \abs[\bigg]{ h^{d / 2} \frac{\vp^{(\vv)}(\vx)\trans \qbar[\vq][-1]}{\sqrt{\bigomegabar[\vq]}} \zn(\vq) },\label{eq:v-tilde-def}\\
  &\bigvtildestNoLink := \sup_{\vq, \vx} \abs[\bigg]{ h^{d / 2} \frac{\vp^{(\vv)}(\vx)\trans \qbar[\vq][-1]}{\sqrt{\bigomegabar[\vq]}} \znst(\vq) }.\label{eq:v-tilde-star-def}
\end{align}

Let $\newtarget{def:k-star-tilde}{\ksttildeNoLink}(\eta)$ denote the conditional $\eta$-quantile of $\bigvtildest$ given $\mX_n$.

\begin{lemma}[Closeness rates]\label{lem:v-rates} Random variables $\bigv$, $\bigvtilde$, $\bigvst$, $\bigvtildest$ satisfy the following:
  \begin{assertions}[wide]
  \item\label{lem:v-rates-v-v-tilde} $\abs{ \bigv - \bigvtilde } = o_{\Problet} \prn{ \rstr + \rhatsoff }$;
  \item\label{lem:v-rates-v-star-v-tilde} $\abs{ \bigvst - \bigvtildest } = o_{\Problet}(\Rplugin)$, where $\Rplugin$ is defined in \cref{th:plug-in};
  \item\label{lem:v-rates-v-bar-star-depends} $\bigvtildest \underset{\mX_n}{\independent} \vy_n$.
  \end{assertions}
\end{lemma}

\begin{proof}
By \cref{th:strong-approximation-binscatter-yurinski}, we have
\begin{equation*}
  \abs[\bigg]{ \bigv - \sup_{\vq, \vx} \abs[\big]{\elltwobarqx[\vq]\trans \zn(\vq)} } = o_{\Problet}\prn{\rstr + \rhatsoff},
\end{equation*}
which is \cref{lem:v-rates-v-v-tilde}.

By \cref{th:plug-in}, we have
\begin{equation*}
  \abs[\bigg]{ \bigvst - \sup_{\vq, \vx} \left| \elltwobarqx[\vq]\trans \znst(\vq) \right| } =  o_{\Problet}(\Rplugin),
\end{equation*}
which is \cref{lem:v-rates-v-star-v-tilde}.

\Cref{lem:v-rates-v-bar-star-depends} follows from the definition of the process $\znst(\cdot)$ and the fact that $\elltwobarqx[\vq]$ only depends on the data via $\mX_n$.
\end{proof}




\begin{lemma}[First sequence]\label{lem:conf-bands-quantiles-are-close}
  There exists a sequence of positive numbers $\nu_{n, 1} \to 0$ such that w.\,p.\,a.~1
  \begin{equation*}
    \kst(1 - \alpha) \leq \ksttilde(1 - \alpha + \nu_{n, 1}) + \Rplugin\quad\text{and}\quad \kst(1 - \alpha) \geq \ksttilde(1 - \alpha - \nu_{n, 1}) - \Rplugin.
  \end{equation*}
\end{lemma}

\begin{proof}
This follows from $| \bigvst - \bigvtildest | = o_{\Problet}\prn{\Rplugin}$ by \cref{lem:v-rates}, directly applying \cref{lem:closeness-in-prob-implies-closeness-of-cond-quantiles}.
\end{proof}

\begin{lemma}[Second sequence]\label{lem:conf-bands-anti-concentration-of-v-tilde} There exists a constant $C_{\bigvtildest} > 0$ such that with $1 - o(1)$ probability
  \begin{equation*}
    \sup_{u \in \mathbb{R}} \Prob[\big]{\abs[\big]{ \bigvtildest - u } \leq \Rplugin + \rstr + \rhatsoff \given \data} \leq C_{\bigvtildest} \prn{\Rplugin + \rstr + \rhatsoff} \sqrt{\log(1 / h)} =: \nu_{n, 2},
  \end{equation*}
  Moreover, for the sequence $\nu_{n, 2} \to 0$ just defined, the following holds w.\,p.\,a.~1:
  \begin{align}
    &\ksttilde\prn{ 1 - \alpha - \nu_{n, 1} } - \ksttilde\prn{1 - \alpha - \nu_{n, 1} - \nu_{n, 2} } \geq \Rplugin + \rstr + \rhatsoff,\label{eq:quantile-diff-eq-one}\\
    &\ksttilde\prn{ 1 - \alpha + \nu_{n, 1} + \nu_{n, 2} } - \ksttilde\prn{1 - \alpha + \nu_{n, 1}} \geq \Rplugin + \rstr + \rhatsoff.\label{eq:quantile-diff-eq-two}
  \end{align}
\end{lemma}

\begin{proof} By \cref{lem:anti-concentration-for-suprema-of-separable-gaussian-processes}, using that $\elltwobarqx[\vq]\trans \znst(\vq)$ is a separable mean-zero Gaussian conditionally on $\mX_n$ process on $\mathcal{Q} \times \mathcal{X}$ with $\Expect[\big]{\prn[\big]{\elltwobarqx[\vq]\trans \znst(\vq)}^2 \given \mX_n} = 1$, we have with $1 - o(1)$ probability
\begin{align*}
  &\sup_{u \in \mathbb{R}} \Prob[\big]{\abs[\big]{ \bigvtildest - u } \leq  \Rplugin + \rstr + \rhatsoff \given \data}\\
  &\quad = \sup_{u \in \mathbb{R}} \Prob[\big]{\abs[\big]{ \bigvtildest - u } \leq  \Rplugin + \rstr + \rhatsoff \given \mX_n}\\
  &\quad \leq 4 \prn{ \Rplugin + \rstr + \rhatsoff } \prn[\bigg]{ \Expect[\Big]{\sup_{\vq, \vx} \abs[\big]{ \elltwobarqx[\vq]\trans \znst(\vq)} \given \mX_n} + 1}\\
  &\quad \leq 4 \prn{ \Rplugin + \rstr + \rhatsoff } \prn[\bigg]{ \sup_{\vq, \vx} \norm[\big]{ \elltwobarqx[\vq] }_1 \Expect[\Big]{\sup_{\vq} \norm[\big]{ \znst(\vq) }_{\infty} \given \mX_n} + 1}\\
  &\quad \overset{\text{(a)}}{\lesssim} \prn{ \Rplugin + \rstr + \rhatsoff } \prn[\Big]{ \Expect[\Big]{\sup_{\vq} \norm[\big]{ \znst(\vq) }_{\infty} \given \mX_n} + 1}\\
  &\quad \overset{\text{(b)}}{\lesssim} \prn{\Rplugin + \rstr + \rhatsoff} \sqrt{\log(1 / h)},
\end{align*}
where in (a) we used $\norm[\big]{ \elltwobarqx[\vq] }_1 \lesssim 1$ w.\,p.\,a.~1 by \cref{lem:closeness-of-linear-terms}; (b) is by \cref{lem:bounding-the-sup-of-gaussian-process}.

We will now prove \eqref{eq:quantile-diff-eq-one}. Note that $\sup_{u \in \mathbb{R}} \Prob[\big]{| \bigvtildest - u | \leq \Rplugin + \rstr + \rhatsoff \given \data} \leq \nu_{n, 2}$ w.\,p.\,a.~1 implies that w.\,p.\,a.~1
\begin{equation*}
  \sup_{u \in \mathbb{R}} \Prob*{u < \bigvtildest \leq u + \Rplugin + \rstr + \rhatsoff \given \data} \leq \nu_{n, 2},
\end{equation*}
that is,
\begin{equation*}
  \sup_{u \in \mathbb{R}} \crl[\big]{ \Prob[\big]{\bigvtildest \leq u + \Rplugin + \rstr + \rhatsoff \given \data} - \Prob[\big]{\bigvtildest \leq u \given \data} } \leq \nu_{n, 2}.
\end{equation*}
Since this is true for any $u$, we can in particular replace $u$ with a random variable $\ksttilde(1 - \alpha - \nu_{n, 1} - \nu_{n, 2})$. Using $\Prob[\big]{\bigvtildest \leq \ksttilde\prn{1 - \alpha - \nu_{n, 1} - \nu_{n, 2}} \given \data} \geq 1 - \alpha - \nu_{n, 1} - \nu_{n, 2}$, this gives w.\,p.\,a.~1
\begin{equation*}
  \Prob[\big]{\bigvtildest \leq \ksttilde\prn{1 - \alpha - \nu_{n, 1} - \nu_{n, 2}} + \Rplugin + \rstr + \rhatsoff \given \data} - \prn{1 - \alpha - \nu_{n, 1} - \nu_{n, 2}} \leq \nu_{n, 2}
\end{equation*}
or
\begin{equation*}
  \Prob[\big]{\bigvtildest \leq \ksttilde\prn{1 - \alpha - \nu_{n, 1} - \nu_{n, 2}} + \Rplugin + \rstr + \rhatsoff \given \data} \leq 1 - \alpha - \nu_{n, 1}.
\end{equation*}
By monotonicity of a (conditional) distribution function, this means that w.\,p.\,a.~1
\begin{equation*}
  \ksttilde\prn{1 - \alpha - \nu_{n, 1} - \nu_{n, 2}} + \Rplugin + \rstr + \rhatsoff \leq \ksttilde\prn{ 1 - \alpha - \nu_{n, 1} }.
\end{equation*}
This proves the inequality~\eqref{eq:quantile-diff-eq-one}. The inequality~\eqref{eq:quantile-diff-eq-two} is proven similarly.
\end{proof}

\begin{proof}[Concluding the proof of \cref{th:confidence-bands}]
Note that
\begin{align*}
  &\Prob*{\bigv > \kst(1 - \alpha)} \overset{\text{(a)}}{\leq} \Prob*{\bigv > \ksttilde(1 - \alpha - \nu_{n, 1}) - \Rplugin}\\
  &\quad \overset{\text{(b)}}{\leq} \Prob*{\bigvtilde > \ksttilde(1 - \alpha - \nu_{n, 1}) - (\Rplugin + \rstr + \rhatsoff)} + o(1)\\
  &\quad \overset{\text{(c)}}{\leq} \Prob*{\bigvtilde > \ksttilde(1 - \alpha - \nu_{n, 1} - \nu_{n, 2})} + o(1)\\
  &\quad \overset{\text{(d)}}{=} \Expect*{\Prob*{\bigvtilde > \ksttilde(1 - \alpha - \nu_{n, 1} - \nu_{n, 2}) \given \mX_n}} + o(1)\\
  &\quad \overset{\text{(e)}}{\leq} \alpha + \nu_{n, 1} + \nu_{n, 2} + o(1) = \alpha + o(1),
\end{align*}
where (a) is by \cref{lem:conf-bands-quantiles-are-close}, (b) is by Assertion~\ref{lem:v-rates-v-v-tilde} in \cref{lem:v-rates}, (c) is by \cref{lem:conf-bands-anti-concentration-of-v-tilde}, (d) is by the law of iterated expectations, (e) is by the definition of a conditional quantile and using that $\bigvtilde$ has the same conditional distribution as $\bigvtildest$.

Similarly,
\begin{align*}
  &\Prob[\big]{\bigv > \kst(1 - \alpha)} \overset{\text{(a)}}{\geq} \Prob[\big]{\bigv > \ksttilde(1 - \alpha + \nu_{n, 1}) + \Rplugin} + o(1)\\
  &\quad \overset{\text{(b)}}{\geq} \Prob[\big]{\bigvtilde > \ksttilde(1 - \alpha + \nu_{n, 1}) + (\Rplugin + \rstr + \rhatsoff)} + o(1)\\
  &\quad \overset{\text{(c)}}{\geq} \Prob[\big]{\bigvtilde > \ksttilde(1 - \alpha + \nu_{n, 1} + \nu_{n, 2})} + o(1)\\
  &\quad \overset{\text{(d)}}{=} \Expect[\big]{\Prob[\big]{\bigvtilde > \ksttilde(1 - \alpha + \nu_{n, 1} + \nu_{n, 2}) \given \mX_n}} + o(1)\\
  &\quad \overset{\text{(e)}}{=} \alpha - \nu_{n, 1} - \nu_{n, 2} + o(1) = \alpha + o(1),
\end{align*}
where (a) is by \cref{lem:conf-bands-quantiles-are-close}, (b) is by Assertion~\ref{lem:v-rates-v-v-tilde} in \cref{lem:v-rates}, (c) is by \cref{lem:conf-bands-anti-concentration-of-v-tilde}, (d) is by the law of iterated expectations. In (e) we used that the distribution function of $\bigvtilde$ conditional on $\mX_n$ is continuous w.\,p.\,a.~1, because by the same anti-concentration argument as in the proof of \cref{lem:conf-bands-anti-concentration-of-v-tilde} there is a positive constant $C$ such that on an event $\mathcal{A}_n$ satisfying $\Prob{\mathcal{A}_n} \to 1$ we have for any $\varepsilon > 0$
\begin{equation*}
  \sup_{u \in \mathbb{R}} \Prob[\big]{|\bigvtilde - u| \leq \varepsilon \given \mX_n} \leq C \varepsilon \sqrt{\log(1 / h)}.
\end{equation*}
In particular, on $\mathcal{A}_n$ all jumps of the distribution function of $\bigvtilde$ conditional on $\mX_n$ are bounded by $C \varepsilon \sqrt{\log(1 / h)}$, which implies that the distribution function is continuous on $\mathcal{A}_n$, since $\varepsilon$ is arbitrary. \cref{th:confidence-bands} is proven.
\end{proof}

The following \lcnamecref{th:k-s-distance} uses a Gaussian anti-concentration result from \parencite{Chernozhukov-Chetverikov-Kato_2014b_AoS} to convert $o_{\Problet}(\cdot)$ bounds obtained above to the bound on the Kolmogorov-Smirnov distance ($\sup$-norm distance of distribution functions), similarly to \parencite{Cattaneo-Farrell-Feng_2020_AOS}.

\subsection{Kolmogorov-Smirnov distance bound}

\begin{theorem}[Kolmogorov-Smirnov Distance]\label{th:k-s-distance}
  Suppose all conditions of \cref{th:plug-in} are true, and $\prn{\Rplugin + \rstr + \rhatsoff} \sqrt{\log n} = o(1)$. Then
  \begin{equation*}
    \sup_{u \in \mathbb{R}} \abs[\Bigg]{ \Prob[\Big]{\sup_{\vq \in \mathcal{Q}} \sup_{\vx \in \mathcal{X}} \left|\hatTqx[\vq]\right| \leq u} - \Prob[\Bigg]{ \sup_{\vq \in \mathcal{Q}}\sup_{\vx \in \mathcal{X}} \abs[\Bigg]{\frac{h^{d / 2}\vp^{(\vv)}(\vx)\trans \qhat[\vq][-1]}{\sqrt{\bigomegahat[\vq]}} \znhatst(\vq)} \leq u \given \data}} = o_{\Problet}(1).
  \end{equation*}
\end{theorem}

\begin{proof}
We will rely on the following \lcnamecref{lem:OWTSYr} which is proven in \cref{sec:lem-OWTSYr-proof} by the same discretization argument as in \cref{th:plug-in}.

\begin{lemma}\label{lem:OWTSYr}
On the same probability space, there exists a mean-zero unconditionally Gaussian process $\tildeZ(\cdot)$ in $\cont{\mathcal{Q}, \mathbb{R}^K}$ such that
\begin{gather}
  \tildeZ(\cdot) \independent \mX_n,\\
  \Expect[\big]{\prn[\big]{\ellqx[\vq]\trans \tildeZ(\vq)}^2} = 1,\label{eq:tildeZq-mean}\\
  \Expect[\Big]{\sup_{\vq \in \mathcal{Q}} \norm[\big]{ \zn(\vq) - \tildeZ(\vq) }_{\infty} \given \mX_n} \lesssim_{\Problet} \prn[\bigg]{\frac{\log n}{n h^d}}^{1 / (4 d_{\mathcal{Q}} + 4)} \sqrt{\log n}.\label{eq:tildeZq-closeness}
\end{gather}
Also, there exists a process $\tildeZst(\cdot)$ in $\cont{\mathcal{Q}, \mathbb{R}^K}$, such that
\begin{gather}
  \tildeZst(\cdot) \independent \mD_n,\\
  \prn[\big]{\mX_n, \zn(\cdot), \tildeZ(\cdot)} \overset{d}{=} \prn[\big]{\mX_n, \znst(\cdot), \tildeZst(\cdot)},
\end{gather}
and in particular
\begin{equation}\label{eq:tildeZqst-closeness}
  \Expect[\Big]{\sup_{\vq \in \mathcal{Q}} \norm[\big]{ \znst(\vq) - \tildeZst(\vq) }_{\infty} \given \mX_n} \lesssim_{\Problet} \prn[\bigg]{\frac{\log n}{n h^d}}^{1 / (4 d_{\mathcal{Q}} + 4)} \sqrt{\log n}.
\end{equation}
\end{lemma}

Since $\tildeZst(\cdot) \overset{d}{=} \tildeZ(\cdot)$,
\begin{equation*}
  \Prob[\Big]{\sup_{\vq, \vx}\abs[\big]{\ellqx[\vq]\trans \tildeZ(\vq)} \leq u} = \Prob[\Big]{\sup_{\vq, \vx}\abs[\big]{\ellqx[\vq]\trans \tildeZst(\vq) } \leq u}.
\end{equation*}

This means that, by the triangle inequality, it is sufficient to prove the following bounds:
\begin{align}
  &\sup_{u \in \mathbb{R}} \abs[\Big]{ \Prob[\Big]{\sup_{\vq, \vx} \abs[\big]{\hatTqx[\vq]} \leq u} - \Prob[\Big]{\sup_{\vq, \vx}\abs[\big]{\ellqx[\vq]\trans \tildeZ(\vq)} \leq u}} = o(1),\label{eq:confidence-bands-t-z-bound}\\
  &\sup_{u \in \mathbb{R}} \abs[\Big]{ \Prob[\Big]{\sup_{\vq, \vx} \abs[\big]{\ellhatqx[\vq]\trans \znhatst(\vq)} \leq u \given \data} - \Prob[\Big]{\sup_{\vq, \vx}\abs[\big]{\ellqx[\vq]\trans \tildeZst(\vq)} \leq u}} = o_{\Problet}(1).\label{eq:confidence-bands-z-z-hat-bound}
\end{align}

We will now prove \eqref{eq:confidence-bands-t-z-bound}.
Note that for any random variables $\xi$ and $\eta$ and any $s > 0$,
\begin{equation}\label{eq:confidence-bands-xi-eta-ineq}
  \sup_{u \in \mathbb{R}} \abs[\big]{ \Prob{\xi \leq u} - \Prob{\eta \leq u} } \leq \sup_{u \in \mathbb{R}} \Prob{\abs{ \eta - u } \leq s} + \Prob{\abs{ \xi - \eta } > s},
\end{equation}
which follows from the two bounds
\begin{align*}
  &\Prob{\xi \leq u} \leq \Prob{\xi \leq u\text{ and } \abs{ \xi - \eta } \leq s} + \Prob{\abs{ \xi - \eta } > s}\\
  &\quad \leq \Prob{\eta \leq u + s} + \Prob{\abs{ \xi - \eta } > s}\\
  &\quad \leq \Prob{\eta \leq u} + \Prob{u < \eta \leq u + s} + \Prob{\abs{ \xi - \eta } > s};\\
  &\Prob{\xi > u} \leq \Prob{\xi > u\text{ and }\abs{ \xi - \eta } \leq s} + \Prob{\abs{ \xi - \eta } > s}\\
  &\quad \leq \Prob{\eta > u - s} + \Prob{\abs{ \xi - \eta } > s}\\
  &\quad \leq \Prob{\eta > u} + \Prob{u - s < \eta \leq u} + \Prob{\abs{ \xi - \eta } > s}.
\end{align*}

\paragraph{Strategy}
To obtain~\eqref{eq:confidence-bands-t-z-bound},
we will apply~\eqref{eq:confidence-bands-xi-eta-ineq} with $\xi = \sup_{\vq, \vx} \abs{ \hatTqx[\vq] }$, $\eta = \sup_{\vq, \vx} \abs[\big]{ \ellqx[\vq]\trans \tildeZ(\vq) }$ and $s = \rstr + \rhatsoff$, use that the term $\Prob{\abs{ \xi - \eta } > s}$ is $o(1)$ because $\abs[\big]{\hatTqx[\vq]}$ and $\abs[\big]{\ellqx[\vq]\trans \tildeZ(\vq)}$ are close, and apply Gaussian anti-concentration
\cref{lem:anti-concentration-for-suprema-of-separable-gaussian-processes} to show that $\sup_{u \in \mathbb{R}} \Prob{\abs{ \eta - u } \leq s}$ is also $o(1)$. The argument for \eqref{eq:confidence-bands-z-z-hat-bound} will be similar.\sloppy

\paragraph{Proof of \eqref{eq:confidence-bands-t-z-bound}}
By $\sup_{\vq, \vx} \norm{ \ellqx[\vq] }_1 \lesssim 1$ (from \cref{lem:closeness-of-linear-terms}) and~\eqref{eq:tildeZq-closeness}, a simple application of \cref{lem:tech-stoch-orders-lemma} gives
\begin{equation}\label{eq:ellqx-z-is-close-to-ellqx-tildez}
  h^{d / 2} \sup_{\vq, \vx} \abs[\Bigg]{ \frac{\vp^{(\vv)}(\vx)\trans \qzero[\vq][-1]}{\sqrt{\bigomega[\vq]}} \zn(\vq) - \frac{\vp^{(\vv)}(\vx)\trans \qzero[\vq][-1]}{\sqrt{\bigomega[\vq]}} \tildeZ(\vq) } = o_{\Problet}\prn{\rstr}.
\end{equation}

By~\eqref{eq:elltwobardiff-bound} in \cref{lem:closeness-of-linear-terms}, we have $\sup_{\vq, \vx} \norm{\elltwobarqx[\vq] - \ellqx[\vq]}_1 \lesssim_{\Problet} \sqrt{\frac{\log(1 / h)}{n h^d}}$, which by \cref{lem:bounding-the-sup-of-gaussian-process} gives $\sup_{\vq, \vx} \abs{\elltwobarqx[\vq]\trans \zn(\vq) - \ellqx[\vq]\trans \zn(\vq)} \lesssim_{\Problet} \frac{\log(1 / h)}{\sqrt{n h^d}} = o(\rstr)$. Combining it with \cref{th:strong-approximation-binscatter-yurinski}, \cref{lem:hats-off} and~\eqref{eq:ellqx-z-is-close-to-ellqx-tildez},
\begin{equation*}
  \Prob[\Big]{\abs[\Big]{ \sup_{\vq, \vx} \abs[\big]{\hatTqx[\vq]} - \sup_{\vq, \vx} \abs[\big]{\ellqx[\vq]\trans \tildeZ(\vq)} } > \rstr + \rhatsoff} = o(1).
\end{equation*}

Then we can apply~\eqref{eq:confidence-bands-xi-eta-ineq} with $\xi = \sup_{\vq, \vx} \abs{ \hatTqx[\vq] }$, $\eta = \sup_{\vq, \vx} \abs[\big]{ \ellqx[\vq]\trans \tildeZ(\vq) }$ and $s = \rstr + \rhatsoff$, and get
\begin{align*}
  &\sup_{u \in \mathbb{R}} \abs[\Big]{ \Prob[\Big]{\sup_{\vq, \vx} \abs{\hatTqx[\vq]} \leq u} - \Prob[\Big]{\sup_{\vq, \vx} \abs[\big]{\ellqx[\vq]\trans \tildeZ(\vq)} \leq u} }\\
  &\quad \leq \sup_{u \in \mathbb{R}} \Prob[\Big]{\abs[\Big]{ \sup_{\vq, \vx} \left| \ellqx[\vq]\trans \tildeZ(\vq) \right| - u } \leq \rstr + \rhatsoff} + o(1).
\end{align*}

For~\eqref{eq:confidence-bands-t-z-bound}, it is left to show that
\begin{equation*}
  \sup_{u \in \mathbb{R}} \Prob[\Big]{\abs[\Big]{ \sup_{\vq, \vx} \abs[\big]{\ellqx[\vq]\trans \tildeZ(\vq)} - u } \leq \rstr + \rhatsoff} = o(1).
\end{equation*}
We will show a stronger version
\begin{equation}\label{eq:confidence-bands-xi-eta-deviation-from-const-bound}
  \sup_{u \in \mathbb{R}} \Prob[\Big]{\abs[\Big]{ \sup_{\vq, \vx} \abs[\big]{\ellqx[\vq]\trans \tildeZ(\vq)} - u } \leq \Rplugin + \rstr + \rhatsoff} = o(1).
\end{equation}
We apply \cref{lem:anti-concentration-for-suprema-of-separable-gaussian-processes}:
\begin{align*}
    &\sup_{u \in \mathbb{R}} \Prob[\Big]{\abs[\Big]{ \sup_{\vq, \vx} \abs[\big]{\ellqx[\vq]\trans \tildeZ(\vq)} - u} \leq \Rplugin + \rstr + \rhatsoff}\\
    &\quad \leq 4 (\Rplugin + \rstr + \rhatsoff) \prn[\Big]{ \Expect[\Big]{\sup_{\vq, \vx} \abs[\big]{ \ellqx[\vq]\trans \tildeZ(\vq) }} + 1 }\\
    &\quad \leq 4 (\Rplugin + \rstr + \rhatsoff) \prn[\Big]{
      \Expect[\Big]{\sup_{\vq, \vx} \norm{ \ellqx[\vq] }_1 \norm[\big]{ \tildeZ(\vq) }_{\infty}} + 1 }\\
    &\quad \overset{\text{(a)}}{\lesssim} (\Rplugin + \rstr + \rhatsoff) \prn[\Big]{
      \Expect[\Big]{\sup_{\vq} \norm[\big]{\tildeZ(\vq)}_{\infty}} + 1 }\\
    &\quad \overset{\text{(b)}}{\lesssim} (\Rplugin + \rstr + \rhatsoff) \sqrt{\log(1 / h)},\numberthis\label{eq:conf-bands-gaussian-concentration-chain}
\end{align*}
where in (a) we used $\norm{ \ellqx[\vq] }_1 \lesssim 1$, see~\eqref{eq:ellqx-bound}; (b) is by
\begin{align*}
  &\Expect[\Big]{\sup_{\vq} \norm[\big]{\tildeZ(\vq)}_{\infty}} = \Expect[\Big]{\sup_{\vq} \norm[\big]{\tildeZ(\vq)}_{\infty} \given \mX_n}\\
  &\quad\leq \Expect[\Big]{\sup_{\vq} \norm[\big]{\tildeZ(\vq) - \zn(\vq)}_{\infty} \given \mX_n} + \Expect[\Big]{\sup_{\vq} \norm[\big]{\zn(\vq)}_{\infty} \given \mX_n}\\
  &\quad \overset{\text{(c)}}{=} o_{\Problet}\prn{\rstr} + O_{\Problet}\prn[\big]{\sqrt{\log (1 / h)}} \overset{\text{(d)}}{=} O_{\Problet}\prn[\big]{\sqrt{\log (1 / h)}}
\end{align*}
using~\eqref{eq:tildeZq-closeness} and \cref{lem:bounding-the-sup-of-gaussian-process} for (c), $\rstr \lesssim 1 \lesssim \sqrt{\log (1 / h)}$ for (d), and noting that
\begin{align*}
  \Expect[\Big]{\sup_{\vq} \norm[\big]{\tildeZ(\vq)}_{\infty}} \lesssim_{\Problet} \sqrt{\log (1 / h)}\\
  \intertext{is equivalent to}
  \Expect[\Big]{\sup_{\vq} \norm[\big]{\tildeZ(\vq)}_{\infty}} \lesssim \sqrt{\log (1 / h)}.
\end{align*}
Since $(\Rplugin + \rstr + \rhatsoff) \sqrt{\log(1 / h)} = o(1)$, the right-hand side in~\eqref{eq:conf-bands-gaussian-concentration-chain} is $o(1)$, proving~\eqref{eq:confidence-bands-xi-eta-deviation-from-const-bound}, which was sufficient for \eqref{eq:confidence-bands-t-z-bound}.

\paragraph{Proof of~\eqref{eq:confidence-bands-z-z-hat-bound}}
By $\sup_{\vq, \vx} \norm{ \ellqx[\vq] }_1 \lesssim 1$ (from \cref{lem:closeness-of-linear-terms}) and~\eqref{eq:tildeZqst-closeness}, a simple application of \cref{lem:tech-stoch-orders-lemma} gives
\begin{equation}\label{eq:ellqx-zst-is-close-to-ellqx-tildezst}
  h^{d / 2} \sup_{\vq, \vx} \abs[\Bigg]{ \frac{\vp^{(\vv)}(\vx)\trans \qzero[\vq][-1]}{\sqrt{\bigomega[\vq]}} \znst(\vq) - \frac{\vp^{(\vv)}(\vx)\trans \qzero[\vq][-1]}{\sqrt{\bigomega[\vq]}} \tildeZst(\vq) } = o_{\Problet}\prn{\rstr}.
\end{equation}

By \cref{th:plug-in}, we have $\sup_{\vq, \vx} \abs[\big]{\ellhatqx[\vq]\trans \znhatst(\vq) - \elltwobarqx[\vq] \znst(\vq)} = o_{\Problet}(\Rplugin)$. By~\eqref{eq:elltwobardiff-bound} in \cref{lem:closeness-of-linear-terms}, we have $\sup_{\vq, \vx} \norm{\elltwobarqx[\vq] - \ellqx[\vq]}_1 \lesssim_{\Problet} \sqrt{\frac{\log(1 / h)}{n h^d}}$, which by \cref{lem:bounding-the-sup-of-gaussian-process} gives $\sup_{\vq, \vx} \abs{\elltwobarqx[\vq]\trans \znst(\vq) - \ellqx[\vq]\trans \znst(\vq)} \lesssim_{\Problet} \frac{\log(1 / h)}{\sqrt{n h^d}} = o(\rstr)$. By the triangle inequality,
\begin{equation*}\sup_{\vq, \vx} \abs{\ellhatqx[\vq]\trans \znhatst(\vq) - \ellqx[\vq]\trans \znst(\vq)} = o_{\Problet}(\Rplugin + \rstr).
\end{equation*}
Combining with \eqref{eq:ellqx-zst-is-close-to-ellqx-tildezst} and applying the triangle inequality again, we obtain
\begin{equation*}
  \Prob[\Big]{\abs[\Big]{\sup_{\vq, \vx} \abs[\big]{\ellhatqx[\vq]\trans \znhatst(\vq)} - \sup_{\vq, \vx}\abs[\big]{\ellqx[\vq]\trans \tildeZst(\vq) }} > \Rplugin + \rstr} = o(1),
\end{equation*}
implying by Markov's inequality that for any $\varepsilon > 0$,
\begin{equation*}
  \Prob[\Big]{\Prob[\Big]{\abs[\Big]{\sup_{\vq, \vx} \abs[\big]{\ellhatqx[\vq]\trans \znhatst(\vq)} - \sup_{\vq, \vx}\abs[\big]{\ellqx[\vq]\trans \tildeZst(\vq) }} > \Rplugin + \rstr \given \data} > \varepsilon} = o(1),
\end{equation*}
that is,
\begin{equation*}
  \Prob[\Big]{\abs[\Big]{\sup_{\vq, \vx} \abs[\big]{\ellhatqx[\vq]\trans \znhatst(\vq)} - \sup_{\vq, \vx}\abs[\big]{\ellqx[\vq]\trans \tildeZst(\vq) }} > \Rplugin + \rstr \given \data} = o_{\Problet}(1).
\end{equation*}

Then we can apply~\eqref{eq:confidence-bands-xi-eta-ineq} with $\Prob*{\cdot \given \data}$ instead of $\Prob{\cdot}$, $\xi = \sup_{\vq, \vx} \abs[\big]{\ellhatqx[\vq]\trans \znhatst(\vq)}$, $\eta = \sup_{\vq, \vx}\abs[\big]{\ellqx[\vq]\trans \tildeZst(\vq) }$ and $s = \Rplugin + \rstr$, and get
\begin{align*}
  &\sup_{u \in \mathbb{R}} \abs[\Big]{ \Prob[\Big]{\sup_{\vq, \vx} \abs[\big]{\ellhatqx[\vq]\trans \znhatst(\vq)} \leq u \given \data} - \Prob[\Big]{\sup_{\vq, \vx} \abs[\big]{\ellqx[\vq]\trans \tildeZst(\vq)} \leq u} }\\
  &\quad \leq \sup_{u \in \mathbb{R}} \Prob[\Big]{\abs[\Big]{ \sup_{\vq, \vx} \abs[\big]{\ellqx[\vq]\trans \tildeZst(\vq)} - u } \leq \Rplugin + \rstr} + o_{\Problet}(1)\\
  &\quad = \sup_{u \in \mathbb{R}} \Prob[\Big]{\abs[\Big]{ \sup_{\vq, \vx} \left|\ellqx[\vq]\trans \tildeZ(\vq) \right| - u } \leq \Rplugin + \rstr} + o_{\Problet}(1),
\end{align*}
where we used that $\tildeZst(\vq)$ is independent of the data allowing us to remove the conditioning on $\data$, and again that $\tildeZ(\cdot)$ and $\tildeZst(\cdot)$ have the same laws.

It is left to use that
\begin{align*}
  &\sup_{u \in \mathbb{R}} \Prob[\Big]{\abs[\Big]{ \sup_{\vq, \vx} \abs[\big]{\ellqx[\vq]\trans \tildeZ(\vq)} - u } \leq \Rplugin + \rstr}\\
  &\quad \leq \sup_{u \in \mathbb{R}} \Prob[\Big]{\abs[\Big]{ \sup_{\vq, \vx} \abs[\big]{\ellqx[\vq]\trans \tildeZ(\vq)} - u } \leq \Rplugin + \rstr + \rhatsoff} \labrel{1}[=] o(1),
\end{align*}
where \labrel{1} is by~\eqref{eq:confidence-bands-xi-eta-deviation-from-const-bound}.

\cref{th:k-s-distance} is proven.
\end{proof}

\subsubsection{Proof of \cref{lem:OWTSYr}: construction of \texorpdfstring{$\tildeZNoLink(\cdot)$}{tildeZ} and \texorpdfstring{$\tildeZstNoLink(\cdot)$}{tildeZ*}}\label{sec:lem-OWTSYr-proof}

Recall that the conditional covariance structure of $\zn(\cdot)$ is
\begin{equation*}
\Expect[\big]{\zn(\vq) \zn(\tilde{\vq})\trans \given \mX_n} = h^{-d} \Sigmacov{\vq}{\tilde{\vq}}
\end{equation*}
with $\Sigmacov{\vq}{\tilde{\vq}}$ defined in \eqref{eq:jvgoXA}. The matrix $\Sigmacov{\vq}{\tilde{\vq}}$ is approximated by the non-random matrix $\mSigma_{\vq, \tilde{\vq}}$ defined by
\begin{equation*}
\mSigma_{\vq, \tilde{\vq}} := \Expect[\big]{\Scov{\vq}{\tilde{\vq}}{\vx_i} \eta^{(1)}(\muzero[\vq]{\vx_i}) \eta^{(1)}(\muzero[\tilde{\vq}]{\vx_i}) \vp(\vx_i) \vp(\vx_i)\trans}.
\end{equation*}
By the same argument as in \cref{lem:holder-gram}, using the Lipschitz continuity of $\Scov{\vq}{\tilde{\vq}}{\vx}$, we have
\begin{equation}\label{eq:mjWpbO}
  \sup_{\vq, \tilde{\vq}} \norm{\Sigmacov{\vq}{\tilde{\vq}} - \mSigma_{\vq, \tilde{\vq}}}
  \lesssim_{\Problet} h^{d} \prn[\bigg]{\frac{\log (1 / h)}{n h^{d}}}^{1 / 2} =: h^d \rsigma.
\end{equation}

\paragraph{Discretization and controlling the deviations of $\znNoLink(\cdot)$}
Then we proceed with the same discretization as in \cref{th:plug-in}. Let $\mathcal{Q}^{\delta_n} := \crl{\vq_1, \ldots, \vq_{|\mathcal{Q}^{\delta_n}|}}$ be an internal $\delta_n$-covering of $\mathcal{Q}$ with respect to the 2-norm $\|\cdot\|$ of cardinality $|\mathcal{Q}^{\delta_n}| \lesssim 1 / \delta_n^{d_{\mathcal{Q}}}$, where $\delta_n$ is chosen later. As already proven in \cref{lem:yurinski-coupling-K-dimensional}, there is a bound on the deviations of $\zn(\cdot)$:
\begin{equation*}
  \Expect[\Big]{\sup_{\norm{\vq - \tilde{\vq}} \leq \delta_n} \norm{\zn(\vq) - \zn(\tilde{\vq})}_{\infty} \given \mX_n}
  \lesssim \sqrt{\delta_n \log \frac{K}{\delta_n}}\quad\text{on $\mathcal{A}$,}
\end{equation*}
where $\mathcal{A}$ is the event from \eqref{eq:large-event-lmax}.

\paragraph{Closeness of the discrete vectors}
Consider now the conditionally Gaussian vector $\zn(\vq|_{\mathcal{Q}^{\delta_n}}) \in \mathbb{R}^{K |\mathcal{Q}^{\delta_n}|}$ and denote
\begin{equation*}
\mSigmabar^{\delta_n} := h^d \Var{\zn(\vq|_{\mathcal{Q}^{\delta_n}}) \given \mX_n} = \prn{\mSigmabar_{\vq_k, \vq_m}}_{k, m = 1}^{|\mathcal{Q}^{\delta_n}|} \in \mathbb{R}^{K |\mathcal{Q}^{\delta_n}| \times K |\mathcal{Q}^{\delta_n}|},
\end{equation*}
so that
\begin{equation*}
\zn(\vq|_{\mathcal{Q}^{\delta_n}}) = h^{-d / 2} \prn{\mSigmabar^{\delta_n}}^{1 / 2} \vxi_n,\quad\text{where }\Problet_{\vxi_n | \mX_n} = \mathcal{N}\prn[\big]{0, \mI_{\mathbb{R}^{K |\mathcal{Q}^{\delta_n}|}}}.
\end{equation*}
(Since $\Problet_{\vxi_n | \mX_n}$ does not depend on $\mX_n$, the vector $\vxi_n$ is independent of $\mX_n$.) Discretizing $\mSigma_{\vq, \tilde{\vq}}$ in the same way, we can put
\begin{equation*}
\mSigma^{\delta_n} := \prn[\big]{\mSigma_{\vq_k, \vq_m}}_{k, m = 1}^{|\mathcal{Q}^{\delta_n}|} \in \mathbb{R}^{K |\mathcal{Q}^{\delta_n}| \times K |\mathcal{Q}^{\delta_n}|}.
\end{equation*}
Since the functions $\mSigmabar_{\vq, \tilde{\vq}}$ and $\mSigma_{\vq, \tilde{\vq}}$ are close, the two matrices are close as well, which we will make precise in the following \lcnamecref{lem:zptGdd}.

\begin{lemma}[The matrices $\mSigmabar^{\delta_n}$ and $\mSigma^{\delta_n}$ are close]\label{lem:zptGdd}
  We have
  \begin{equation*}
    \norm[\big]{\mSigmabar^{\delta_n} - \mSigma^{\delta_n}} \lesssim_{\Problet} |\mathcal{Q}^{\delta_n}| h^d \rsigma.
  \end{equation*}
\end{lemma}

\begin{proof}
We can bound the Frobenius norm using the bound on each element of the matrix:
\begin{equation*}
\norm[\big]{\mSigmabar^{\delta_n} - \mSigma^{\delta_n}} \leq \norm[\big]{\mSigmabar^{\delta_n} - \mSigma^{\delta_n}}_F \leq |\mathcal{Q}^{\delta_n}| \sup_{\vq, \tilde{\vq}} \norm[\big]{\Sigmacov{\vq}{\tilde{\vq}} - \mSigma_{\vq, \tilde{\vq}}}.
\end{equation*}
It is left to combine with \cref{eq:mjWpbO}.
\end{proof}

Applying \cref{lem:maximal-ineq-gaussian-vectors}, we then have
\begin{equation}
  \begin{aligned}
    h^{-d / 2} \Expect[\big]{
    \norm[\big]{\prn{\mSigmabar^{\delta_n}}^{1 / 2} \vxi_n - \prn{\mSigma^{\delta_n}}^{1 / 2} \vxi_n}_{\infty}
    \given
    \mX_n
    } &\leq 2 h^{-d / 2} \sqrt{\log(K |\mathcal{Q}^{\delta_n}|)} \norm{\mSigmabar^{\delta_n} - \mSigma^{\delta_n}}^{1 / 2}\\
    &\lesssim_{\Problet} |\mathcal{Q}^{\delta_n}|^{1 / 2} \rsigma^{1 / 2} \sqrt{\log(K |\mathcal{Q}^{\delta_n}|)}.
  \end{aligned}
\end{equation}

\paragraph{Embedding the conditionally Gaussian vector into a conditionally Gaussian process}
By the same argument as in \cref{lem:yurinski-coupling-K-dimensional}, there exists a law
on $\mathcal{X}^n \times \cont{\mathcal{Q}, \mathbb{R}^K}$ which is the joint law of $\mX_n' \overset{d}{=} \mX_n$ and a conditionally on $\mX_n'$ mean-zero Gaussian process
$\tildeZprDef(\vq)$
with the conditional covariance structure
\begin{equation*}
\Expect[\big]{\tildeZpr(\vq) \tildeZpr(\tilde{\vq})\trans \given \mX_n'} = h^{-d} \mSigma_{\vq, \tilde{\vq}}.
\end{equation*}
Continuity follows by the Kolmogorov-Chentsov theorem:
\begin{equation}\label{eq:kynzOk}
\Expect[\big]{\norm[\big]{\tildeZpr(\vq) - \tildeZpr(\tilde{\vq})}^a \given \mX_n'} \leq C_n \norm{\vq - \tilde{\vq}}^{d_{\mathcal{Q}} + b},
\end{equation}
for $a, b > 0$ chosen as follows. The vector $\tildeZpr(\vq) - \tildeZpr(\tilde{\vq})$ is conditionally on $\mX_n'$ Gaussian with covariance
\begin{equation*}
h^{-d} \prn[\big]{\mSigma_{\vq, \vq} - 2 \mSigma_{\vq, \tilde{\vq}} + \mSigma_{\tilde{\vq}, \tilde{\vq}}}.
\end{equation*}
Hence, we have for any $m > 0$
\begin{align*}
  &\Expect[\big]{\norm[\big]{\tildeZpr(\vq) - \tildeZpr(\tilde{\vq})}^{2 m} \given \mX_n'}\\
  &\quad \leq h^{-d m} \norm[\big]{\mSigma_{\vq, \vq} - 2 \mSigma_{\vq, \tilde{\vq}} + \mSigma_{\tilde{\vq}, \tilde{\vq}}}^m \Expect{\norm{\vxi_K}^{2 m}}\\
  &\quad \leq C_n \norm{\vq - \tilde{\vq}}^m,
\end{align*}
where $\vxi_K \sim \mathcal{N}(0, \mI_K)$, and we used that $\mSigma_{\vq, \tilde{\vq}}$ is Lipschitz in $\tilde{\vq}$ a.\,s. (with the constant allowed to depend on $n$).
So we can take, for example, $m = d_{\mathcal{Q}} + 1$, $a = 2 m$, $b = 1$ in \eqref{eq:kynzOk}.

Projecting this $\cont{\mathcal{Q}, \mathbb{R}^K}$-process onto $\mathcal{Q}^{\delta_n}$, we obtain a vector $\tildeZpr\prn[\big]{\vq|_{\mathcal{Q}^{\delta_n}}} \in \mathbb{R}^{K |\mathcal{Q}^{\delta_n}|}$
such that
\begin{equation*}
\prn[\big]{\mX_n', \tildeZpr\prn[\big]{\vq|_{\mathcal{Q}^{\delta_n}}}} \overset{d}{=} \prn[\big]{\mX_n, h^{-d / 2} \prn{\mSigma^{\delta_n}}^{1 / 2} \vxi_n}.
\end{equation*}
Since $\cont{\mathcal{Q}, \mathbb{R}^K}$ is Polish, by Theorem~8.17 (transfer) in \cite{Kallenberg2021} on our probability space there exists a random element $\tildeZ(\vq) \in \cont{\mathcal{Q}, \mathbb{R}^K}$ such that
\begin{equation*}
\prn[\big]{\mX_n', \tildeZpr\prn[\big]{\vq|_{\mathcal{Q}^{\delta_n}}}, \tildeZpr(\vq)} \overset{d}{=} \prn[\big]{\mX_n, h^{-d / 2} \prn{\mSigma^{\delta_n}}^{1 / 2} \vxi_n, \tildeZ(\vq)},
\end{equation*}
In particular, almost surely $h^{-d / 2} \prn{\mSigma^{\delta_n}}^{1 / 2} \vxi_n$ is the projection of $\tildeZ(\cdot)$ on $\mathcal{Q}^{\delta_n}$. Note that $\Problet_{\tildeZ(\vq) | \mX_n}$ does not depend on $\mX_n$, which means $\tildeZ(\vq)$ is independent of $\mX_n$. In addition, note that \eqref{eq:tildeZq-mean} holds.

\paragraph{Controlling the deviations of $\tildeZNoLink(\vq)$}

It is proven by the same argument as in \cref{lem:yurinski-coupling-K-dimensional,th:plug-in} that
\begin{equation*}
\Expect[\Big]{\sup_{\norm{\vq - \tilde{\vq}} \leq \delta_n} \norm{\tildeZ(\vq) - \tildeZ(\tilde{\vq})}_{\infty} \given \mX_n} = \Expect[\Big]{\sup_{\norm{\vq - \tilde{\vq}} \leq \delta_n} \norm{\tildeZ(\vq) - \tildeZ(\tilde{\vq})}_{\infty}} \lesssim \sqrt{\delta_n \log \frac{K}{\delta_n}}.
\end{equation*}

\paragraph{Choosing $\delta_n$ and conclusion}
Combining the bounds obtained above, we can write
\begin{align*}
  &\Expect[\Big]{\sup_{\vq \in \mathcal{Q}} \norm[\big]{\zn(\vq) - \tildeZ(\vq)}_{\infty} \given \mX_n}\\
  &\quad \leq \Expect[\Big]{\sup_{\norm{\vq - \tilde{\vq}} \leq \delta_n} \norm[\big]{\zn(\vq) - \zn(\tilde{\vq})}_{\infty} \given \mX_n}
    +
    h^{-d / 2} \Expect[\big]{\norm[\big]{\prn{\mSigmabar^{\delta_n}}^{1 / 2} \vxi_n - \prn{\mSigma^{\delta_n}}^{1 / 2} \vxi_n}_{\infty} \given \mX_n}\\
  &\qquad + \Expect[\Big]{\sup_{\norm{\vq - \tilde{\vq}} \leq \delta_n} \norm[\big]{\tildeZ(\vq) - \tildeZ(\tilde{\vq})}_{\infty} \given \mX_n}\\
  &\quad \lesssim_{\Problet} \prn[\bigg]{\sqrt{\delta_n} + \frac{\rsigma^{1 / 2}}{\delta_n^{d_{\mathcal{Q}} / 2}}} \sqrt{\log \frac{K}{\delta_n}}.
\end{align*}
Choose the approximately optimal
\begin{equation*}
\delta_n := (\rsigma)^{1 / (d_{\mathcal{Q}} + 1)},
\end{equation*}
giving
\begin{equation*}
  \Expect[\Big]{\sup_{\vq \in \mathcal{Q}} \norm[\big]{\zn(\vq) - \tildeZ(\vq)}_{\infty} \given \mX_n}
  \lesssim_{\Problet} (\rsigma)^{1 / (2 d_{\mathcal{Q}} + 2)} \sqrt{\log n}.
\end{equation*}

\paragraph{Constructing \texorpdfstring{$\tildeZstNoLink(\cdot)$}{tildeZ*}}

By Theorem~8.17 (transfer) in \cite{Kallenberg2021}, on our probability space there is a random element $\tildeZst(\cdot)$ with values in $\cont{\mathcal{Q}, \mathbb{R}^K}$ such that
\begin{equation}\label{eq:gSXmWZ}
\prn[\big]{\mX_n, \zn(\cdot), \tildeZ(\cdot)} \overset{d}{=} \prn[\big]{\mX_n, \znst(\cdot), \tildeZst(\cdot)},
\end{equation}
and $\tildeZst(\cdot)$ is conditionally on $(\mX_n, \znst(\cdot))$ independent of $\vy_n$. Since also $\znst(\cdot)$ conditionally on $\mX_n$ independent of $\vy_n$, we have by Theorem~8.12 (chain rule) in \cite{Kallenberg2021} that $(\znst(\cdot), \tildeZst(\cdot))$ is conditionally on $\mX_n$ independent of $\vy_n$; in particular, $\tildeZst(\cdot)$ is conditionally on $\mX_n$ independent of $\vy_n$. But by~\eqref{eq:gSXmWZ}, $\tildeZst(\cdot)$ is independent of $\mX_n$. Again by the chain rule, $\tildeZst(\cdot)$ is independent of $(\mX_n, \vy_n) = \mD_n$.

\section{Examples}\label{sec:examples}

This section discusses in detail the four motivating examples introduced in the paper.

\subsection{Quantile regression}\label{sec:example-quantile-regression}

The next result verifies that the quantile regression case under the same conditions on $f_{Y|X}$ as Condition S.2 in~\parencite{Belloni-Chernozhukov-Chetverikov-FernandezVal_2019_JoE} is a special case of our setting.

\begin{proposition}[Verification of \cref{ass:gl-holder} for quantile regression]\label{prop:gl-holder-quantile-regression}
  Suppose \cref{ass:dgp-holder,ass:gl-quasi-uniform-partition,ass:gl-local-basis} hold with $\mathcal{Q} = [\varepsilon_0, 1 - \varepsilon_0]$ for some $\varepsilon_0 \in (0, 1 / 2)$, the loss is given by $\riskfunc[q]{y}{\eta} = \prn{ q - \Indlet\{y < \eta\} } (y - \eta)$, and $\Expect{|y_1|} < \infty$. Assume further that $\eta(\cdot)\colon \mathbb{R} \to \rangeeta$ is strictly monotonic and twice continuously differentiable with $\rangeeta$ an open connected subset of $\mathbb{R}$ containing the conditional $q$-quantile of $y_1 | \vx_1 = \vx$, given by $\eta\prn[\big]{\muzero[q]{\vx}}$ for all $(\vx, q)$, and that $|\eta^{(1)}(\zeta)|$ is bounded away from zero uniformly over $\vx \in \mathcal{X}$, $q \in \mathcal{Q}$, and $\zeta$ satisfying $|\zeta - \muzero[q]{\vx}| \leq r$; $y \mapsto F_{Y|X}(y | \vx)$ is twice continuously differentiable with first derivative $f_{Y|X}(y | \vx)$ (in particular, $\measureonY$ is the Lebesgue measure); $f_{Y|X}\prn[\big]{\eta\prn[\big]{\muzero[q]{\vx}} | \vx}$ is bounded away from zero uniformly over $q \in \mathcal{Q}$, $\vx \in \mathcal{X}$, and the derivative of $y \mapsto f_{Y|X}(y | \vx)$ is continuous and bounded in absolute value from above uniformly over $y \in \mathcal{Y}_{\vx}, \vx \in \mathcal{X}$.
Then \cref{ass:gl-holder,ass:complexity,ass:complexity-addition,eq:l2-norm-holder-type-restriction} are also true.
\end{proposition}

\begin{remark*}
Taking $\eta\prn[\big]{\muzero[q]{\vx}}$ to be the conditional $q$-quantile does not violate~\eqref{eq:holder-rho-problem} by \cref{lem:mu-borel-measurable}.
\end{remark*}

\begin{remark*}
  In the setting of \cref{prop:gl-holder-quantile-regression}, it is not necessary to assume that $\muzero[q]{\vx}$ is Lipschitz in parameter (as we do in \cref{ass:dgp-holder}\ref{ass:mu-bounded}). Since
  \begin{equation*}
    \frac{\partial}{\partial q} \muzero[q]{\vx} = \frac{1}{f_{Y|X} \prn{\eta\prn{\muzero[q]{\vx}} | \vx} \eta^{(1)}\prn{\muzero[q]{\vx}}},
  \end{equation*}
  the Lipschitz property follows from $f_{Y|X} \prn{\eta\prn{\muzero[q]{\vx}} | \vx}$ and $|\eta^{(1)}\prn{\muzero[q]{\vx}}|$ being bounded from below uniformly over $\vx \in \mathcal{X}$, $q \in \mathcal{Q}$.
\end{remark*}

\begin{proof}
We will verify the assumptions one by one.

\paragraph{Verifying \cref{ass:gl-holder}\ref{ass:abs-cont}} It is easy to see that the a.\,e. derivative $\eta \mapsto \riskder[q]{y}{\eta} \equiv \Indlet \crl{ y - \eta < 0 } - q$ of $\eta \mapsto \riskfunc[q]{y}{\eta}$ is Lebesgue integrable and satisfies
\begin{equation*}
\int_a^b \riskder[q]{y}{\eta}\,\mathrm{d}\eta = \riskfunc[q]{y}{b} - \riskfunc[q]{y}{a}
\end{equation*}
for any $[a, b] \subset \rangeeta$.

\paragraph{Verifying \cref{ass:gl-holder}\ref{ass:psi-moments}}
Since $\eta(\muzero[q]{\cdot})$ is the conditional $q$-quantile, we have
\begin{equation*}
  \Expect{\riskder[q]{y_i}{\eta\prn{\muzero[q]{\vx_i}}} \given \vx_i} = \Expect{\Indlet \crl{ y_i < \eta\prn{\muzero[q]{\vx_i}} } - q \given \vx_i} = q - q = 0
\end{equation*}
and
\begin{align*}
  &\sigma_q^2(\vx) = \Expect{\riskder[q]{y_i}{\eta\prn{\muzero[q]{\vx_i}}}^2 \given \vx_i = \vx}\\
  &\quad = \Expect{\prn{ \Indlet \crl{ y_i < \eta\prn{\muzero[q]{\vx_i}} } - q }^2 \given \vx_i = \vx} = q - 2 q^2 + q^2 = q (1 - q)
\end{align*}
is constant in $\vx$ (in particular continuous in $\vx$) and bounded away from zero since both $q$ and $1 - q$ are bounded away from zero. Since $q (1 - q)$ is smooth, $\sigma^2_{q}(\vx)$ is Lipschitz in $q$. The family $\crl{ \riskder[q] {y_i}{\eta \prn{ \muzero[q]{\vx_i} }} \colon q \in \mathcal{Q} }$ has a positive measurable envelope 1 which has uniformly bounded conditional moments of any order.

\paragraph{Verifying \cref{ass:gl-holder}\ref{ass:psi-holder-jumps}}
Clearly, $\riskfunc[q]{y}{\eta}$ is convex in $\eta$ and the a.\,e. derivative in $\eta$ is $\riskder[q]{y}{\eta} \equiv \Indlet \crl{ y - \eta < 0 } - q$, which is piecewise constant with only one jump (therefore piecewise H\"older with $\alpha = 1$). The link function $\eta(\cdot)$ is strictly monotonic and twice continuously differentiable by assumption.

\paragraph{Verifying \cref{ass:gl-holder}\ref{ass:psi-eta}}
The conditional expectation
\begin{align*}
  \expriskder[q]{\vx}{\eta} = \Expect{\Indlet \crl{ y < \eta } - q \given \vx_i = \vx} = \int_{-\infty}^{\eta} f_{Y|X}(y | \vx)\,\mathrm{d} y - q
\end{align*}
is twice continuously differentiable in $\eta$ (the integral on the right is a Riemann integral, possibly improper) and its second derivative $f'_{Y|X}(\eta | \vx)$ is continuous and bounded in absolute value. By the mean value theorem, this means that $f_{Y|X}\prn{\eta\prn{\muzero[q]{\vx}} | \vx}$ being bounded away from zero implies $f_{Y|X}\prn{\eta(\zeta) | \vx}$ is bounded away from zero for $\zeta$ sufficiently close to $\muzero[q]{\vx}$. Together with the assumed lower bound on $|\eta^{(1)}(\zeta)|$, this gives the required uniform lower bound on $\Psione[q]{\vx}{\eta(\zeta)} \eta^{(1)}(\zeta)^2$. The bound on $|\Psione[q]{\vx}{\eta(\zeta)}|$ from above in such a neighborhood (and in fact everywhere) is automatic since $f_{Y|X}(y | \vx)$ is bounded from above (uniformly over $y \in \mathcal{Y}_{\vx}$, $\vx \in \mathcal{X}$).

\paragraph{Verifying \cref{ass:complexity}}
This verification proceeds similarly to Lemmas~25--28 in~\parencite{Belloni-Chernozhukov-Chetverikov-FernandezVal_2019_JoE}.

The class of functions
\begin{equation*}
  \mathcal{W}_1 := \crl[\big]{ (\vx, y) \mapsto \Indlet \crl{ y < \eta(\vp(\vx)\trans \vbeta) } \colon \vbeta \in \mathbb{R}^K }
\end{equation*}
is VC-subgraph with index $O(K)$ by Lemmas 2.6.15 and 2.6.18 in~\parencite{vaart96_weak_conver_empir_proces} (since $\eta(\cdot)$ is monotone).
The class of functions
\begin{equation*}
  \mathcal{G}_2 = \crl[\big]{ (\vx, y) \mapsto \Indlet \crl{ y < \eta\prn{\muzero[q]{\vx}} } \colon q \in \mathcal{Q} }
\end{equation*}
is VC-subgraph with index $2$ since $\eta(\muzero[q]{\vx})$ is increasing in $q$ for any $\vx \in \mathcal{X}$, giving that the class of sets $\crl{ (\vx, y): y < \eta\prn{\muzero[q]{\vx}}}$ with $q \in \mathcal{Q}$ is linearly ordered by inclusion. The VC property of $\mathcal{W}_1$ with envelope 1 implies that it satisfies the uniform entropy bound~\eqref{eq:uniform-entropy-bound} with $A \lesssim 1$ and $V \lesssim K$. The VC property of $\mathcal{G}_2$ with envelope 1 implies that it satisfies the uniform entropy bound~\eqref{eq:uniform-entropy-bound} with $A \lesssim 1$ and $V \lesssim 1$. By \cref{lem:permanence-properties}, for any fixed $r > 0$ the class
\begin{align*}
  &\mathcal{G}_1 = \crl[\big]{ (\vx, y) \mapsto \Indlet \crl{ y < \eta(\vp(\vx)\trans \vbeta) } - \Indlet \crl{ y < \eta\prn{\muzero[q]{\vx}} }: \norm{ \vbeta - \betazero[q] }_{\infty} \leq r, q \in \mathcal{Q}}
\end{align*}
with envelope 2 satisfies the uniform entropy bound~\eqref{eq:uniform-entropy-bound} with $A \lesssim 1$ and $V \lesssim K$ because it is a subclass of $\mathcal{W}_1 - \mathcal{G}_2$.

For a fixed vector space $\mathcal{B}$ of dimension $\dim \mathcal{B}$,
\begin{equation*}
  \mathcal{W}_{\mathcal{B}} := \crl{ (\vx, y) \mapsto \Indlet \crl{ y < \eta(\vp(\vx)\trans \vbeta) }: \vbeta \in \mathcal{B} }
\end{equation*}
is VC with index $O\prn{\dim \mathcal{B}}$ by Lemmas 2.6.15 and 2.6.18 in~\parencite{vaart96_weak_conver_empir_proces} (since $\eta(\cdot)$ is monotone). Therefore, for any fixed $c > 0$, $\delta \in \Delta$, the class
\begin{equation*}
  \mathcal{W}_{2, \delta} := \crl[\big]{ (\vx, y) \mapsto \Indlet \crl{ y < \eta(\vp(\vx)\trans \vbeta) } \Indlet \crl{ \vx \in \mathcal{N}_{[c \log n]}(\delta) }: \vbeta \in \mathbb{R}^K }
\end{equation*}
is VC-subgraph with index $O\prn{ \log^d n }$ by Lemma 2.6.18 in~\parencite{vaart96_weak_conver_empir_proces} because it is contained in the product of $\mathcal{W}_{\mathcal{B}_{2, \delta}}$ for some vector space $\mathcal{B}_{2, \delta}$ of dimension $\dim \mathcal{B}_{2, \delta} \lesssim \log^d n$ and a fixed function $\Indlet \crl{ \vx \in \mathcal{N}_{[c \log n]}(\delta) }$. This means $\mathcal{W}_{2, \delta}$ with envelope 1 satisfies the uniform entropy bound~\eqref{eq:uniform-entropy-bound} with $A \lesssim 1$ and $V \lesssim \log^d n$. Then the union of $O\prn{h^{-d}}$ such classes
\begin{equation*}
  \mathcal{W}_2 := \crl[\big]{ (\vx, y) \mapsto \Indlet \crl{ y < \eta(\vp(\vx)\trans \vbeta) } \Indlet \crl{ \vx \in \mathcal{N}_{[c \log n]}(\delta) }: \vbeta \in \mathbb{R}^K, \delta \in \Delta }
\end{equation*}
with envelope 1 satisfies the uniform entropy bound~\eqref{eq:uniform-entropy-bound} with $A \lesssim 1$ and $V \lesssim \log^d n$, see~\eqref{eq:log-d-entropy-chain}. By \cref{lem:permanence-properties}, the same is true of
\begin{align*}
  \mathcal{G}_3 = \crl[\big]{
  &(\vx, y) \mapsto \\
  &\brk[\big]{ \Indlet \crl{y < \eta\prn{\vp(\vx)\trans \vbeta}} - \Indlet \crl{y < \eta\prn{\vp(\vx)\trans \betazero[q]}} } \Indlet \crl{\vx \in \mathcal{N}_{[c \log n]}(\delta)}:\\
  &\vbeta \in \mathbb{R}^K, \delta \in \Delta, q \in \mathcal{Q}}
\end{align*}
with envelope 2 because it is a subclass of $\mathcal{W}_2 - \mathcal{W}_2$.

The class of functions
\begin{equation*}
  \mathcal{G}_4 = \crl{ \vx \mapsto f_{Y|X} \prn{ \eta\prn{\muzero[q]{\vx}} | \vx }: q \in \mathcal{Q} }
\end{equation*}
has a bounded envelope by the assumptions of the lemma. Moreover, $\mathcal{G}_4$ has the following property: for any $q_1, q_2 \in \mathcal{Q}$ we have for some $\xi_{\vx, q_1, q_2}$ between $\eta\prn{\muzero[q_1]{\vx}}$ and $\eta\prn{\muzero[q_2]{\vx}}$
\begin{align*}
  &\abs[\big]{ f_{Y|X} \prn{ \eta\prn{\muzero[q_1]{\vx}} | \vx } - f_{Y|X} \prn{ \eta\prn{\muzero[q_2]{\vx}} | \vx } }\\
  &\quad = \abs[\big]{ f'_{Y|X}\prn{\xi_{\vx, q_1, q_2}} } \cdot \abs{ \eta\prn{\muzero[q_1]{\vx}} - \eta\prn{\muzero[q_2]{\vx}} }\\
  &\quad \lesssim \abs[\big]{ \eta\prn{\muzero[q_1]{\vx}} - \eta\prn{\muzero[q_2]{\vx}} } \tag*{since $f'_{Y|X}$ is uniformly bounded}\\
  &\quad \lesssim \abs[\big]{ \muzero[q_1]{\vx} - \muzero[q_2]{\vx} } \tag*{since $\eta(\cdot)$ is Lipschitz}\\
  &\quad \lesssim \abs[\big]{ q_1 - q_2 } \tag*{since $\muzero[q]{\vx}$ is Lipschitz in $q$.}
\end{align*}
with constants in $\lesssim$ not depending on $\vx, q_1, q_2$ or $n$ (this is also proven in Lemma~20 in~\cite{Belloni-Chernozhukov-Chetverikov-FernandezVal_2019_JoE}). Since $\mathcal{Q}$ is a fixed one-dimensional segment, this implies that $\mathcal{G}_4$ satisfies the uniform entropy bound~\eqref{eq:uniform-entropy-bound} with $A, V \lesssim 1$.

For a fixed $l$, the class of functions
\begin{equation*}
  \crl{ (\vx, y) \mapsto p_l(\vx) \Indlet \crl{ y < \eta(\vp(\vx)\trans \vbeta) }: \vbeta \in \mathbb{R}^K }
\end{equation*}
is VC-subgraph with index $O(1)$ because it is contained in the product of $\mathcal{W}_{\mathcal{B}_{5, l}}$ for some vector space $\mathcal{B}_{5, l}$ of dimension $\dim \mathcal{B}_{5, l} \lesssim 1$, and a fixed function $p_l(\vx)$. Then this class with envelope $O(1)$ satisfies the uniform entropy bound~\eqref{eq:uniform-entropy-bound} with $A, V \lesssim 1$. Since, as we have shown above, the same is true of $\mathcal{G}_2$, by \cref{lem:permanence-properties}, it is also true of
\begin{align*}
  \mathcal{G}_5 = \crl[\big]{
  &\prn{\vx, y} \mapsto\\
  &p_l(\vx) \brk[\big]{\Indlet \crl{ y < \eta\prn{\muzero[q]{\vx}} } - \Indlet \crl{ y < \eta\prn{\vp(\vx)\trans \betazero[q]} } }: q \in \mathcal{Q}}.
\end{align*}

\paragraph{Verifying \cref{ass:complexity-addition}}
The functions in the class have the form
\begin{align*}
  &\riskfunc[q]{y}{\eta\prn{\bp(\vx) \trans \prn{\betazero[q] + \bbeta }}} \Indlet\prn{\vx \in \mathcal{N}_{[c \log n]}(\delta)}\\
  &\qquad - \riskfunc[q]{y}{\eta\prn{\bp(\vx)\trans \prn{\betazero[q] + \bbeta - \vv }}}\Indlet\prn{\vx \in \mathcal{N}_{[c \log n]}(\delta)} \\
  &\qquad - \brk[\big]{ \eta \prn[\big]{ \bp(\vx)\trans \prn{ \betazero[q] + \bbeta } } - \eta \prn[\big]{ \bp(\vx)\trans \prn{ \betazero[q] + \bbeta - \vv } } } \riskder[q]{y}{\eta \prn{ \bp(\vx)\trans \betazero[q] }}\\
  &\qquad \times \Indlet\prn{\vx \in \mathcal{N}_{[c \log n]}(\delta)}\\
  &\quad =: T_1 + T_2 + T_3 + T_4,
\end{align*}
where
\begin{align*}
  &T_1 := y \brk[\big]{ \Indlet\crl[\big]{y < \eta\prn[\big]{\bp(\vx) \trans \prn{\betazero[q] + \bbeta -\vv }}} - \Indlet\crl[\big]{y < \eta\prn[\big]{\bp(\vx) \trans \prn{\betazero[q] + \bbeta }}} } \\
  &\qquad \times \Indlet\crl{\vx \in \mathcal{N}_{[c \log n]}(\delta)},\\
  &T_2 := \eta\prn{\bp(\vx) \trans \prn{\betazero[q] + \bbeta }} \Indlet\crl[\big]{y < \eta\prn[\big]{\bp(\vx) \trans \prn{\betazero[q] + \bbeta }}} \Indlet\crl{\vx \in \mathcal{N}_{[c \log n]}(\delta)},\\
  &T_3 := - \eta \prn[\big]{ \bp(\vx)\trans \prn{ \betazero[q] + \bbeta - \vv } } \Indlet\crl[\big]{y < \eta \prn[\big]{ \bp(\vx)\trans \prn{ \betazero[q] + \bbeta - \vv } }} \Indlet\crl{\vx \in \mathcal{N}_{[c \log n]}(\delta)},\\
  &T_4 := - \Indlet\crl[\big]{y < \eta \prn[\big]{ \bp(\vx)\trans \betazero[q] }} \brk[\big]{ \eta\prn[\big]{\bp(\vx) \trans \prn{\betazero[q] + \bbeta }} - \eta \prn[\big]{ \bp(\vx)\trans \prn{ \betazero[q] + \bbeta - \vv } }  }\\
  &\qquad \times \Indlet\crl{\vx \in \mathcal{N}_{[c \log n]}(\delta)}.
\end{align*}

Note that for $T_1$ to be nonzero, $y$ has to lie in a fixed interval (not depending on $n$), say $\brk{-\tilde{R}, \tilde{R}}$. The class of functions
\begin{equation*}
  \crl{ (\vx, y) \mapsto y \Indlet\crl{|y| \leq \tilde{R}} \Indlet\crl{y < \eta(\bp(\vx)\trans\bbeta)}: \bbeta \in \mathcal{B} },
\end{equation*}
where $\mathcal{B}$ is any linear subspace of $\mathbb{R}^K$, is VC-subgraph with index $O(\dim \mathcal{B})$ by Lemmas 2.6.15 and 2.6.18 in~\parencite{vaart96_weak_conver_empir_proces} (since $\eta(\cdot)$ is monotone and $y \mapsto y \Indlet\crl{|y| \leq \tilde{R}}$ is one fixed function). The class $\crl{ (\vx, y) \mapsto T_1, \norm{ \vbeta - \betazero[q] }_{\infty} \leq r, \norm{ \vv }_{\infty} \leq \epsninassump, q \in \mathcal{Q} }$ with $\delta$ fixed is a subclass of the difference of two such classes, so by \cref{lem:permanence-properties} this class with a large enough constant envelope satisfies the uniform entropy bound~\eqref{eq:uniform-entropy-bound} with $A \lesssim 1$ and $V \lesssim \log^d n$. Therefore, the union of $O\prn{h^{-d}}$ such classes
\begin{equation*}
  \crl{ (\vx, y) \mapsto T_1, \norm{ \vbeta - \betazero[q] }_{\infty} \leq r, \norm{ \vv }_{\infty} \leq \epsninassump, \delta \in \Delta, q \in \mathcal{Q} }
\end{equation*}
with a large enough constant envelope satisfies the uniform entropy bound~\eqref{eq:uniform-entropy-bound} with $A \lesssim 1$ and $V \lesssim \log^d n$, see~\eqref{eq:log-d-entropy-chain}.

Similarly, the classes
\begin{align*}
  &\crl{ (\vx, y) \mapsto T_2, \norm{ \vbeta - \betazero[q] }_{\infty} \leq r, \norm{ \vv }_{\infty} \leq \epsninassump, \delta \in \Delta, q \in \mathcal{Q} },\\
  &\crl{ (\vx, y) \mapsto T_3, \norm{ \vbeta - \betazero[q] }_{\infty} \leq r, \norm{ \vv }_{\infty} \leq \epsninassump, \delta \in \Delta, q \in \mathcal{Q} },\\
  &\crl{ (\vx, y) \mapsto T_4, \norm{ \vbeta - \betazero[q] }_{\infty} \leq r, \norm{ \vv }_{\infty} \leq \epsninassump, \delta \in \Delta, q \in \mathcal{Q} }\\
\end{align*}
with large enough constant envelopes also satisfy the uniform entropy bound~\eqref{eq:uniform-entropy-bound} with $A \lesssim 1$ and $V \lesssim \log^d n$. We used that $\epsninassump \to 0$ giving that $\vv$ is bounded in $\infty$-norm (like $\vbeta$). Applying \cref{lem:permanence-properties} one more time, we have that there exist some constants $\nC{quant-compl-verification-a} \geq e$, $\nC{quant-compl-verification-v} \geq 1$ and $\nC{quant-compl-verification-envelope} > 0$ such that
\begin{equation}\label{eq:quant-reg-compl-verification-uniform-entropy-bound}
  \sup_{\mathbb{Q}} N\prn{\mathcal{G}, \norm{ \cdot }_{\mathbb{Q}, 2}, \varepsilon \oC{quant-compl-verification-envelope}} \leq \prn[\bigg]{\frac{\oC{quant-compl-verification-a}}{\varepsilon}}^{\oC{quant-compl-verification-v} \log^d n}
\end{equation}
for all $0 < \varepsilon \leq 1$, where the supremum is taken over all finite discrete probability measures $\mathbb{Q}$ and $\mathcal{G}$ is the class defined in \cref{ass:complexity-addition}. Note that the integral representation of $\mathcal{G}$ makes it clear that this class not only has a large enough constant envelope, but is also bounded by $\oC{quant-compl-verification-const-before-eps} \epsninassump$, where $\nC{quant-compl-verification-const-before-eps}$ is a large enough constant.

For large enough $n$ we can replace $\varepsilon$ with $\oC{quant-compl-verification-const-before-eps} \varepsilon \epsninassump / \oC{quant-compl-verification-envelope}$ in~\eqref{eq:quant-reg-compl-verification-uniform-entropy-bound}, giving
\begin{equation*}
  \sup_{\mathbb{Q}} N\prn{\mathcal{G}, \norm{ \cdot }_{\mathbb{Q}, 2}, \oC{quant-compl-verification-const-before-eps} \varepsilon \epsninassump } \leq \prn[\bigg]{\frac{\oC{quant-compl-verification-a} \oC{quant-compl-verification-envelope}}{\oC{quant-compl-verification-const-before-eps} \varepsilon \epsninassump}}^{\oC{quant-compl-verification-v} \log^d n}
\end{equation*}
for all $0 < \varepsilon \leq 1$. For large enough $n$, $\oC{quant-compl-verification-a} \oC{quant-compl-verification-envelope} / (\oC{quant-compl-verification-const-before-eps} \epsninassump) \geq e$. The verification is complete.

\paragraph{Verifying \eqref{eq:l2-norm-holder-type-restriction}}
In this case, $\riskder[q]{y}{\eta} = \Indlet \crl{y < \eta } - q$ and $\eta(\muzero[q]{\vx_i})$ is the $q$-quantile of $y_i$ conditional on $\vx_i$. Without loss of generality, we will assume that $\eta(\cdot)$ is strictly increasing and $q \leq \tilde{q}$ (the other cases are symmetric).
\begin{align*}
  &\Expect[\big]{\abs[\big]{\riskder[q]{y_i}{\eta(\muzero[q]{\vx_i})} \eta^{(1)}(\muzero[q]{\vx_i}) - \riskder[\tilde{q}]{y_i}{\eta(\muzero[\tilde{q}]{\vx_i})} \eta^{(1)}(\muzero[\tilde{q}]{\vx_i}) }^2 \given \vx_i}\\
  &\quad = \Expect[\big]{\riskder[q]{y_i}{\eta(\muzero[q]{\vx_i})}^2 \given \vx_i} \eta^{(1)}(\muzero[q]{\vx_i})^2\\
  &\qquad + \Expect[\big]{\riskder[\tilde{q}]{y_i}{\eta(\muzero[\tilde{q}]{\vx_i})}^2 \given \vx_i} \eta^{(1)}(\muzero[\tilde{q}]{\vx_i})^2\\
  &\qquad - 2 \Expect[\big]{\riskder[q]{y_i}{\eta(\muzero[q]{\vx_i})} \riskder[\tilde{q}]{y_i}{\eta(\muzero[\tilde{q}]{\vx_i})} \given \vx_i} \eta^{(1)}(\muzero[q]{\vx_i}) \eta^{(1)}(\muzero[\tilde{q}]{\vx_i})\\
  &\quad = (q - q^2) \eta^{(1)}(\muzero[q]{\vx_i})^2 + (\tilde{q} - \tilde{q}^2) \eta^{(1)}(\muzero[\tilde{q}]{\vx_i})^2\\
  &\qquad - 2 (q - q \tilde{q}) \eta^{(1)}(\muzero[q]{\vx_i}) \eta^{(1)}(\muzero[\tilde{q}]{\vx_i})\\
  &\quad = q \bigl| \eta^{(1)}(\muzero[q]{\vx_i}) - \eta^{(1)}(\muzero[\tilde{q}]{\vx_i}) \bigr|^2 + (\tilde{q} - q) \eta^{(1)}(\muzero[\tilde{q}]{\vx_i})^2\\
  &\qquad - \bigl| q \eta^{(1)}(\muzero[q]{\vx_i}) - \tilde{q} \eta^{(1)}(\muzero[\tilde{q}]{\vx_i}) \bigr|^2\\
  &\quad \leq q \bigl| \eta^{(1)}(\muzero[q]{\vx_i}) - \eta^{(1)}(\muzero[\tilde{q}]{\vx_i}) \bigr|^2 + (\tilde{q} - q) \eta^{(1)}(\muzero[\tilde{q}]{\vx_i})^2\\
  &\quad \labrel{1}[\lesssim] q (\tilde{q} - q)^2 + \tilde{q} - q \lesssim \tilde{q} - q,
\end{align*}
where in \labrel{1} we used that $\eta^{(1)}(\cdot)$ on a fixed compact is Lipschitz and $\muzero[q]{\vx}$ is Lipschitz in $q$ uniformly over $\vx$, as well as boundedness of $\muzero[q]{\vx}$ uniformly over $q$ and $\vx$.

This concludes the proof of \cref{prop:gl-holder-quantile-regression}.
\end{proof}

\begin{proposition}[Verification of the conditions of \cref{lem:consistency-nonconvex}]
  Suppose all conditions of \cref{prop:gl-holder-quantile-regression} hold. In addition, suppose there is a positive constant $\nC{quantile-bounded-away-from-zero-const}$ such that we have $\inf f_{Y|X}(y|\vx) > \oC{quantile-bounded-away-from-zero-const}$, where the infimum is over $\vx \in \mathcal{X}$, $q \in \mathcal{Q}$, $\norm{ \bbeta }_{\infty} \leq R$ for $R$ described in \cref{lem:consistency-nonconvex}, $y$ between $\eta \prn{ \bp(\vx)\trans\bbeta }$ and $\eta \prn{ \muzero[q]{\vx} }$. Then conditions in \cref{cond:lower-bound-of-psi-one,cond:nonconvex-consistency-additional-complexity-assumption} of \cref{lem:consistency-nonconvex} also hold.
\end{proposition}

\begin{proof}
We only need to verify \cref{lem:consistency-nonconvex}\ref{cond:nonconvex-consistency-additional-complexity-assumption} since \cref{lem:consistency-nonconvex}\ref{cond:lower-bound-of-psi-one} is directly assumed in this lemma ($\Psione[q]{\vx}{\eta} = f_{Y|X}(\eta|\vx)$ in this case).

In this verification, we will use $\theta_1 := \bp(\vx)\trans \bbeta$, $\theta_2 := \bp(\vx)\trans \betazero[q]$ to simplify notations. Rewrite
\begin{align*}
  &\riskfunc[q]{y}{\eta(\theta_1)} - \riskfunc[q]{y}{\eta(\theta_2)}\\
  &\quad = y [\Indlet\crl{y < \eta(\theta_2)} - \Indlet\crl{y < \eta(\theta_1)}] + q [\eta(\theta_1) - \eta(\theta_2)]\\
  &\qquad + \eta(\theta_1) \Indlet\crl{y < \eta(\theta_1)} - \eta(\theta_2) \Indlet\crl{y < \eta(\theta_2)}.
\end{align*}

By the same argument as in the proof of \cref{prop:gl-holder-quantile-regression}, the class
\begin{equation*}
  \crl{(\vx, y) \mapsto y [\Indlet\crl{y < \eta(\theta_2)} - \Indlet\crl{y < \eta(\theta_1)}]: \norm{\bbeta}_{\infty} \leq R, q \in \mathcal{Q}}
\end{equation*}
with a large enough constant envelope satisfies the uniform entropy bound~\eqref{eq:uniform-entropy-bound} with $A \lesssim 1$ and $V \lesssim K$.

The class $\crl{(\vx, y) \mapsto q: q \in \mathcal{Q}}$ is of course VC with a constant index (as a subclass of the class of constant functions), and the class $\crl{(\vx, y) \mapsto \eta(\bp(\vx)\trans \bbeta): \bbeta \in \mathbb{R}^K}$ is VC with index $O(K)$ because the space of functions $\bp(\vx)\trans \bbeta$ is a linear space with $O(K)$ dimension, and $\eta(\cdot)$ is monotone. Applying \cref{lem:permanence-properties}, we see that the class $\crl{(\vx, y) \mapsto q [\eta(\theta_1) - \eta(\theta_2)]: \norm{\bbeta}_{\infty} \leq R, q \in \mathcal{Q}}$ with a large enough constant envelope satisfies the uniform entropy bound~\eqref{eq:uniform-entropy-bound} with $A \lesssim 1$ and $V \lesssim K$.

The class $\crl{(\vx, y) \mapsto \Indlet\crl{y < \eta(\bp(\vx)\trans \bbeta)}: \bbeta \in \mathbb{R}^K}$ is VC with index $O(K)$ because the space of functions $\bp(\vx)\trans \bbeta$ is a linear space with $O(K)$ dimension, and $\eta(\cdot)$ is monotone. Again applying \cref{lem:permanence-properties}, we see that the class $\crl{(\vx, y) \mapsto \eta(\theta_1) \Indlet\crl{y < \eta(\theta_1)}: \norm{\bbeta}_{\infty} \leq R}$ with a large enough constant envelope satisfies the uniform entropy bound~\eqref{eq:uniform-entropy-bound} with $A \lesssim 1$ and $V \lesssim K$. The same is true of its subclass $\crl{(\vx, y) \mapsto \eta(\theta_2) \Indlet\crl{y < \eta(\theta_2)}: q \in \mathcal{Q}}$.

It is left to apply \cref{lem:permanence-properties} again, concluding that the class described in \cref{lem:consistency-nonconvex}\ref{cond:nonconvex-consistency-additional-complexity-assumption} with a large enough constant envelope satisfies the uniform entropy bound~\eqref{eq:uniform-entropy-bound} with $A \lesssim 1$ and $V \lesssim K$.
\end{proof}

\subsection{Distribution regression}\label{sec:example-distribution-regression}

\begin{proposition}[Verification of assumptions for distribution regression]\label{prop:gl-holder-distribution-regression}
  Let $\mathcal{Q} = [-A, A]$ for some $A > 0$. Suppose \cref{ass:dgp-holder} holds with the loss function
  $\riskfunc[q]{y}{\eta} = \prn{ \Indlet\crl{y \leq q} - \eta }^2$, the link function $\eta(\cdot)\colon \mathbb{R} \to (0, 1)$ is strictly monotonic and twice continuously differentiable,
  $q \mapsto F_{Y|X}(q|\vx)$ is continuously differentiable with derivative $f_{Y|X}(q|\vx)$ (in particular, $\measureonY$ is Lebesgue measure), $|\eta^{(1)}(\zeta)|$ is bounded away from zero uniformly over $\vx \in \mathcal{X}$, $q \in \mathcal{Q}$, and $\zeta$ satisfying $|\zeta - \muzero[q]{\vx}| \leq r$,
  $\vx \mapsto F_{Y|X}(q|\vx)$ is a continuous function, and $F_{Y|X}(q|\vx) = \eta\prn{\muzero[q]{\vx}}$ lies in a compact subset of $(0, 1)$ for all $q \in \mathcal{Q}, \vx \in \mathcal{X}$ (this subset does not depend on q or $\vx$).
    Then \cref{ass:gl-holder,ass:complexity,ass:complexity-addition}, \cref{eq:l2-norm-holder-type-restriction}, and \cref{cond:lower-bound-of-psi-one,cond:nonconvex-consistency-additional-complexity-assumption} of \cref{lem:consistency-nonconvex} are also true.
  \end{proposition}

\begin{remark*}
Taking $\eta\prn{\muzero[q]{\vx}}$ to be the conditional distribution function does not violate \cref{eq:holder-rho-problem} by the following standard argument. For any Borel function $\mu(\cdot)\colon \mathcal{X} \to \mathbb{R}$ one has
\begin{align*}
&\Expect[\big]{\prn[\big]{\Indlet\crl{y_1 \leq q} - \eta(\mu(\vx_1))}^2}\\
&\quad = \Expect[\big]{\prn[\big]{\Indlet\crl{y_1 \leq q} - F_{Y|X}(q | \vx_1) + F_{Y|X}(q | \vx_1) - \eta(\mu(\vx_1))}^2}\\
&\quad = \Expect[\big]{\prn[\big]{\Indlet\crl{y_1 \leq q} - F_{Y|X}(q | \vx_1)}^2}\\
&\qquad + 2 \Expect[\big]{\prn[\big]{\Indlet\crl{y_1 \leq q} - F_{Y|X}(q | \vx_1)} \prn[\big]{F_{Y|X}(q | \vx_1) - \eta(\mu(\vx_1))}}\\
&\qquad + \Expect[\big]{\prn[\big]{F_{Y|X}(q|\vx_1) - \eta(\mu(\vx_1))}^2}.
\end{align*}
Since the cross term is zero (proven by conditioning on $\vx_1$), this means
\begin{align*}
&\Expect[\big]{\prn[\big]{\Indlet\crl{y_1 \leq q} - \eta(\mu(\vx_1))}^2}\\
&\quad = \Expect[\big]{\prn[\big]{\Indlet\crl{y_1 \leq q} - F_{Y|X}(q | \vx_1)}^2} + \Expect[\big]{\prn[\big]{F_{Y|X}(q|\vx_1) - \eta(\mu(\vx_1))}^2}\\
&\quad \geq \Expect[\big]{\prn[\big]{\Indlet\crl{y_1 \leq q} - F_{Y|X}(q | \vx_1)}^2}.
\end{align*}
Pointwise in $q \in \mathcal{Q}$, equality holds if and only if $\eta(\mu(\vx_1)) = F_{Y|X}(q | \vx_1)$ almost surely.
\end{remark*}

\begin{remark}\label{rem:feasibility-distr-reg}
  In this case,
  \begin{equation*}
    \Expect[\big]{\zn(q) \zn(\tilde{q})\trans \given \mX_n} = h^{-d} \Expectlet_n \brk[\big]{
      \Scov{q}{\tilde{q}}{\vx_i}
      \eta^{(1)}(\muzero[q]{\vx_i}) \eta^{(1)}(\muzero[\tilde{q}]{\vx_i})
      \bp(\vx_i) \bp(\vx_i)\trans}
  \end{equation*}
  with
  \begin{equation*}
    \Scov{q}{\tilde{q}}{\vx} = 4 F_{Y|X}(q \wedge \tilde{q} | \vx) \prn[\big]{1 - F_{Y|X}(q \vee \tilde{q} | \vx)}.
  \end{equation*}
  This covariance structure is not known, but it can be estimated by
  \begin{equation}\label{eq:dist-reg-covariance-structure}
    h^{-d} \Expectlet_n \brk[\big]{
      \ScovEst{q}{\tilde{q}}{\vx_i}
      \eta^{(1)}(\muhat[q]{\vx_i}) \eta^{(1)}(\muhat[\tilde{q}]{\vx_i})
      \bp(\vx_i) \bp(\vx_i)\trans}.
  \end{equation}
  with
\begin{equation*}
\ScovEst{q}{\tilde{q}}{\vx} = 4 \eta(\muhat[q \wedge \tilde{q}]{\vx}) \prn[\big]{1 - \eta(\muhat[q \vee \tilde{q}]{\vx})}.
\end{equation*}
\end{remark}

\begin{proof} We will verify the assumptions one by one.

\paragraph{Verifying \cref{ass:gl-holder}\ref{ass:abs-cont}} It is easy to see that the a.\,e. derivative $\eta \mapsto \riskder[q]{y}{\eta} \equiv 2 \prn{ \eta - \Indlet\crl{y \leq q} }$ of $\eta \mapsto \riskfunc[q]{y}{\eta}$ is Lebesgue integrable and satisfies
\begin{equation*}
\int_a^b \riskder[q]{y}{\eta}\,\mathrm{d}\eta = \riskfunc[q]{y}{b} - \riskfunc[q]{y}{a}
\end{equation*}
for any $[a, b] \subset (0, 1)$.

\paragraph{Verifying \cref{ass:gl-holder}\ref{ass:psi-moments}}
The first-order optimality condition
\begin{equation*}
  \Expect[\big]{2 \prn{ \eta \prn{ \muzero[q]{\vx_i} } - \Indlet\crl{y_i \leq q} } \given \vx_i} = 2 \prn[\big]{ F_{Y|X}(q | \vx_i) - \Expect{\Indlet\crl{y_i \leq q} \given \vx_i} } = 0
\end{equation*}
indeed holds. The conditional variance
\begin{equation*}
  \sigma^2_q(\vx) := 4 \Expect{\prn{ \eta \prn{ \muzero[q]{\vx_i} } - \Indlet\crl{y_i \leq q} }^2 \given \vx_i = \vx} = 4 F_{Y|X}(q|\vx) \prn{ 1 - F_{Y|X}(q|\vx) }
\end{equation*}
is continuous and bounded away from zero by the assumptions ($F_{Y|X}(q|\vx)$ cannot achieve 0 or 1).

The family $\crl{ 2 \prn{ F_{Y|X}(q|\vx) - \Indlet\crl{y \leq q} }: q \in \mathcal{Q}  }$ is bounded in absolute value by $\riskderenvel{\vx}{y} \equiv 2$.

\paragraph{Verifying \cref{ass:gl-holder}\ref{ass:psi-holder-jumps}}\label{sec:distribution-regression-verifying-psi-holder-jumps}
The loss function $\riskfunc[q]{y}{\eta}$ is infinitely smooth with respect to $\eta$. Its first derivative is
$\riskder[q]{y}{\eta} = 2 \prn{ \eta - \Indlet\crl{y \leq q} }$.
Since the derivative of $\eta(\cdot)$ is bounded on a compact interval, the function $\riskder[q]{y}{\eta \prn{ \theta }}$ is Lipschitz in $\theta$ on a compact interval.

\paragraph{Verifying \cref{ass:gl-holder}\ref{ass:psi-eta}}
The conditional expectation
\begin{equation*}
  \expriskder[q]{\vx}{\eta} := \Expect{2 \prn{ \eta - \Indlet\crl{y_i \leq q} } \given \vx_i = \vx} = 2 \eta - 2 F_{Y|X}(q|\vx)
\end{equation*}
is linear, and in particular infinitely smooth, in $\eta$. Its first partial derivative
\begin{equation*}
  \Psione[q]{\vx}{\eta} = \frac{\partial}{\partial \eta} \expriskder[q]{\vx}{\eta} = 2
\end{equation*}
is a nonzero constant, so it is bounded and bounded away from zero everywhere. The assumed lower bound on $|\eta^{(1)}(\zeta)|$ gives the required uniform lower bound on $\Psione[q]{\vx}{\eta(\zeta)} \eta^{(1)}(\zeta)^2$. The second partial derivative is zero, and so it is also bounded.

\paragraph{Verifying \cref{ass:complexity}}
The class of functions
\begin{equation*}
  \mathcal{G}_{11} := \crl{ (\vx, y) \mapsto F_{Y|X}(q|\vx): q \in \mathcal{Q} }
\end{equation*}
with a constant envelope satisfies the uniform entropy bound~\eqref{eq:uniform-entropy-bound} with $A \lesssim 1$ and $V \lesssim 1$ because it is VC-subgraph with index 1 since the subgraphs are linearly ordered by inclusion (by monotonicity in $q$).

The class
\begin{equation*}
  \mathcal{G}_{1 2} := \crl{ (\vx, y) \mapsto \vp(\vx)\trans \vbeta: \vbeta \in \mathbb{R}^K }
\end{equation*}
is VC with index not exceeding $K + 2$ by Lemma 2.6.15 in~\parencite{vaart96_weak_conver_empir_proces}. Since in a fixed bounded interval $\eta(\cdot)$ is Lipschitz, the class
\begin{equation*}
  \mathcal{G}_{1 3} := \crl{ (\vx, y) \mapsto \eta \prn{\vp(\vx)\trans \vbeta}: \norm{ \bbeta - \betazero[q] }_{\infty} \leq r}
\end{equation*}
with a large enough constant envelope satisfies the uniform entropy bound~\eqref{eq:uniform-entropy-bound} with $A \lesssim 1$ and $V \lesssim K$. By \cref{lem:permanence-properties}, the class
\begin{equation*}
  \mathcal{G}_1 := \crl[\big]{ (\vx, y) \mapsto 2 \prn[\big]{ \eta \prn{\bp(\vx)\trans \bbeta} - F_{Y|X}(q|\vx) }: \norm{ \bbeta - \betazero[q] }_{\infty} \leq r, q \in \mathcal{Q} },
\end{equation*}
with a large enough constant envelope satisfies the uniform entropy bound~\eqref{eq:uniform-entropy-bound} with $A \lesssim 1$ and $V \lesssim K$.

The class of sets $\crl{ (\vx, y): y \leq q }$ with $q \in \mathcal{Q}$ is linearly ordered by inclusion, so it is VC with a constant index and so is the class of functions
\begin{equation*}
  \mathcal{G}_{2 1} := \crl{ (\vx, y) \mapsto \Indlet\crl{y \leq q}: q \in \mathcal{Q} },
\end{equation*}
giving by \cref{lem:permanence-properties} that
\begin{equation*}
  \mathcal{G}_2 := \crl{ (\vx, y) \mapsto 2 \prn{ F_{Y|X}(q|\vx) - \Indlet\crl{y \leq q} }: q \in \mathcal{Q} },
\end{equation*}
which is a subclass of $2 (\mathcal{G}_{11} - \mathcal{G}_{2 1})$, satisfies the uniform entropy bound~\eqref{eq:uniform-entropy-bound} with $A \lesssim 1$ and $V \lesssim 1$.

For a fixed vector space $\mathcal{B}$ of dimension $\dim \mathcal{B}$,
\begin{equation*}
  \mathcal{W}_{\mathcal{B}} := \crl{ (\vx, y) \mapsto \bp(\vx)\trans\bbeta: \bbeta \in \mathcal{B} }
\end{equation*}
is VC-subgraph with index $O(\dim \mathcal{B})$ by Lemmas 2.6.15 and 2.6.18 in~\parencite{vaart96_weak_conver_empir_proces}. Therefore, again using that $\eta(\cdot)$ is Lipschitz in a compact interval, by \cref{lem:permanence-properties} we have that for any fixed $c > 0$, $\delta \in \Delta$, the class
\begin{equation*}
  \mathcal{W}_{3, \delta} := \crl{ (\vx, y) \mapsto \eta \prn{ \bp(\vx)\trans\bbeta } \Indlet\crl{\vx \in \mathcal{N}_{[c \log n]}(\delta)}: \norm{ \bbeta - \betazero[q] }_{\infty} \leq r, q \in \mathcal{Q}}
\end{equation*}
with a large enough constant envelope satisfies the uniform entropy bound~\eqref{eq:uniform-entropy-bound} with $A \lesssim 1$ and $V \lesssim \log^d n$. Then the union of $O(h^{-d})$ such classes
\begin{equation*}
  \mathcal{W}_3 := \crl{ (\vx, y) \mapsto \eta \prn{ \bp(\vx)\trans\bbeta } \Indlet\crl{\vx \in \mathcal{N}_{[c \log n]}(\delta)}: \norm{ \bbeta - \betazero[q] }_{\infty} \leq r, q \in \mathcal{Q}, \delta \in \Delta}
\end{equation*}
satisfies the uniform entropy bound~\eqref{eq:uniform-entropy-bound} with $A \lesssim 1$ and $V \lesssim \log^d n$ (see~\eqref{eq:log-d-entropy-chain}). The same is true of
\begin{align*}
  \mathcal{G}_3 := \crl[\big]{&(\vx, y) \mapsto 2 \prn[\big]{\eta \prn{ \bp(\vx)\trans\bbeta } - \eta \prn{ \bp(\vx)\trans\betazero[q] }} \Indlet\crl{\vx \in \mathcal{N}_{[c \log n]}(\delta)}:\\
  &\norm{ \bbeta - \betazero[q] }_{\infty} \leq r, q \in \mathcal{Q}, \delta \in \Delta}
\end{align*}
because it is a subclass of $2 \mathcal{W}_3 - 2 \mathcal{W}_3$.

The class
\begin{equation*}
  \mathcal{G}_4 := \crl{ \vx \mapsto 2 }
\end{equation*}
consists of just one bounded function, so clearly it satisfies the uniform entropy bound~\eqref{eq:uniform-entropy-bound} with envelope 2, $A \lesssim 1$, $V \lesssim 1$.

Finally, for any fixed $l \in \crl{ 1, \ldots, K }$ the class
\begin{equation*}
  \mathcal{G}_{51} := \crl{ (\vx, y) \mapsto p_l(\vx) \eta \prn{ \bp(\vx)\trans \betazero[q] }: q \in \mathcal{Q} }
\end{equation*}
satisfies the uniform entropy bound~\eqref{eq:uniform-entropy-bound} with a large enough constant envelope, $A \lesssim 1$ and $V \lesssim 1$ because $\eta(\cdot)$ is Lipschitz and $\mathcal{G}_{5 1}$ is contained in a fixed function multiplied by $\eta(\mathcal{W}_{\mathcal{B}})$ for a linear space $\mathcal{B}$ of a constant dimension, and by \cref{lem:permanence-properties}
\begin{equation*}
  \mathcal{G}_{5} := \crl{ (\vx, y) \mapsto 2 p_l(\vx) \prn{ F_{Y|X}(q|\vx) - \eta \prn{ \bp(\vx)\trans \betazero[q] } }: q \in \mathcal{Q} }
\end{equation*}
with a constant envelope satisfies the uniform entropy bound~\eqref{eq:uniform-entropy-bound} with $A \lesssim 1$ and $V \lesssim 1$.

\paragraph{Verifying \cref{ass:complexity-addition}}
\Cref{ass:complexity-addition} holds because the class of functions
\begin{align*}
  \crl[\big]{&(\vx, y) \mapsto \brk[\big]{\eta \prn[\big]{ \bp(\vx)\trans \prn{ \betazero[q] + \bbeta } } - \eta \prn{ \bp(\vx)\trans \prn{ \betazero[q] + \bbeta - \vv } } }\\
  &\times \brk[\big]{ \eta \prn{ \bp(\vx)\trans \prn{ \betazero[q] + \bbeta } } + \eta \prn{ \bp(\vx)\trans \prn{ \betazero[q] + \bbeta - \vv } } - 2 \eta \prn{ \bp(\vx)\trans \betazero[q] } }\\
  &\times \Indlet\prn{\vx \in \mathcal{N}_{[c \log n]}(\delta)}:\\
  &\norm{ \bbeta - \betazero[q] }_{\infty} \leq r, \norm{ \vv }_{\infty} \leq \epsninassump, \delta \in \Delta, q \in \mathcal{Q}}
\end{align*}
is contained in the product of two classes
\begin{align*}
  \mathcal{V}_1 := \crl[\big]{
  &(\vx, y) \mapsto \brk[\big]{\eta \prn{ \bp(\vx)\trans \bbeta } - \eta \prn{ \bp(\vx)\trans \prn{ \bbeta - \vv } } } \Indlet\prn{\vx \in \mathcal{N}_{[c \log n]}(\delta)}:\\
  &\norm{ \bbeta }_{\infty} \leq \tilde{r}, \norm{ \vv }_{\infty} \leq \epsninassump, \delta \in \Delta, q \in \mathcal{Q}},
\end{align*}
and
\begin{align*}
  \mathcal{V}_2 := \crl[\big]{
  &(\vx, y) \mapsto \brk[\big]{ \eta \prn[\big]{ \bp(\vx)\trans \prn{ \betazero[q] + \bbeta } } + \eta \prn[\big]{ \bp(\vx)\trans \prn{ \betazero[q] + \bbeta - \vv } } - 2 \eta \prn{ \bp(\vx)\trans \betazero[q] } }\\
  &\times \Indlet\prn{\vx \in \mathcal{N}_{[c \log n]}(\delta)}: \norm{ \bbeta - \betazero[q] }_{\infty} \leq r, \norm{ \vv }_{\infty} \leq \epsninassump, \delta \in \Delta, q \in \mathcal{Q}}
\end{align*}
for some fixed $\tilde{r} > 0$. Class $\mathcal{V}_1$ with envelope $\epsninassump$ multiplied by a large enough constant (since $\eta$ is Lipschitz) satisfies the uniform entropy bound~\eqref{eq:uniform-entropy-bound} with $A \lesssim 1 / \epsninassump$ and $V \lesssim \log^d n$ (this can be shown by further breaking down $\mathcal{V}_1$ into classes $\crl{ \eta \prn{ \bp(\vx)\trans \bbeta } }$ and $\crl{ \eta \prn{ \bp(\vx)\trans \prn{ \bbeta - \vv } } }$ with constant envelopes, using \cref{lem:permanence-properties} and then replacing $\varepsilon$ in the uniform entropy bound by $\varepsilon \cdot \epsninassump$). Class $\mathcal{V}_2$ with a large enough constant envelope satisfies the uniform entropy bound~\eqref{eq:uniform-entropy-bound} with $A \lesssim 1$ and $V \lesssim \log^d n$ because it is true for each of the three additive terms it can be broken down into. We omit the details since they are the same as in the verification of \cref{ass:complexity}.

\paragraph{Verifying Condition~\ref{cond:lower-bound-of-psi-one} in \cref{lem:consistency-nonconvex}}
This condition holds trivially because $\Psione[q]{\cdot}{\cdot}$ is a positive constant.

\paragraph{Verifying Condition~\ref{cond:nonconvex-consistency-additional-complexity-assumption} in \cref{lem:consistency-nonconvex}}
The class of functions described in this condition is
\begin{align*}
  \crl[\big]{
  &(\vx, y) \mapsto \prn[\big]{ \eta \prn{ \bp(\vx)\trans\betazero[q] } - \eta \prn{ \bp(\vx)\trans\bbeta } } \cdot \prn[\big]{ 2 \Indlet\crl{y \leq q} - \eta \prn{ \bp(\vx)\trans\bbeta } - \eta \prn{ \bp(\vx)\trans\betazero[q] } }:\\
  &\norm{ \bbeta }_{\infty} \leq R, q \in \mathcal{Q}}.
\end{align*}
The assertion follows by \cref{lem:permanence-properties} since
\begin{equation*}
  \crl{ (\vx, y) \mapsto \Indlet\crl{y \leq q}: q \in \mathcal{Q} }
\end{equation*}
is a VC-subgraph class with a constant index and
\begin{equation*}
  \crl{ (\vx, y) \mapsto \eta \prn{ \bp(\vx)\trans\bbeta }: \norm{ \bbeta }_{\infty} \leq R }
\end{equation*}
with a large enough constant envelope satisfies the uniform entropy bound~\eqref{eq:uniform-entropy-bound} with $A \lesssim 1$ and $V \lesssim K$ because $\eta(\cdot)$ is Lipschitz.

\paragraph{Verifying \eqref{eq:l2-norm-holder-type-restriction}}
Without loss of generality, assume $q \leq \tilde{q}$ (the other case is symmetric).
\begin{align*}
  &\Expect[\big]{\abs{\riskder[q]{y_i}{\eta(\muzero[q]{\vx_i})} \eta^{(1)}(\muzero[q]{\vx_i}) - \riskder[\tilde{q}]{y_i}{\eta(\muzero[\tilde{q}]{\vx_i})} \eta^{(1)}(\muzero[\tilde{q}]{\vx_i}) }^2 \given \vx_i}\\
  &\quad = \Expect[\big]{\riskder[q]{y_i}{\eta(\muzero[q]{\vx_i})}^2 \given \vx_i} \eta^{(1)}(\muzero[q]{\vx_i})^2\\
  &\qquad + \Expect[\big]{\riskder[\tilde{q}]{y_i}{\eta(\muzero[\tilde{q}]{\vx_i})}^2 \given \vx_i} \eta^{(1)}(\muzero[\tilde{q}]{\vx_i})^2\\
  &\qquad - 2 \Expect[\big]{\riskder[q]{y_i}{\eta(\muzero[q]{\vx_i})} \riskder[\tilde{q}]{y_i}{\eta(\muzero[\tilde{q}]{\vx_i})} \given \vx_i} \eta^{(1)}(\muzero[q]{\vx_i}) \eta^{(1)}(\muzero[\tilde{q}]{\vx_i})\\
  &\quad = 4 \brk[\big]{F_{Y|X}(q|\vx_i) - F_{Y|X}(q|\vx_i)^2} \eta^{(1)}(\muzero[q]{\vx_i})^2\\
  &\qquad + 4 \brk[\big]{F_{Y|X}(\tilde{q}|\vx_i) - F_{Y|X}(\tilde{q}|\vx_i)^2} \eta^{(1)}(\muzero[\tilde{q}]{\vx_i})^2\\
  &\qquad - 8 \brk[\big]{F_{Y|X}(q|\vx_i) - F_{Y|X}(q|\vx_i) F_{Y|X}(\tilde{q}|\vx_i)} \eta^{(1)}(\muzero[q]{\vx_i}) \eta^{(1)}(\muzero[\tilde{q}]{\vx_i})\\
  &\quad = 4 F_{Y|X}(q|\vx_i) \abs[\big]{\eta^{(1)}(\muzero[q]{\vx_i}) - \eta^{(1)}(\muzero[\tilde{q}]{\vx_i})}^2\\
  &\qquad + 4 \brk[\big]{F_{Y|X}(\tilde{q}|\vx_i) - F_{Y|X}(q|\vx_i)} \eta^{(1)}(\muzero[\tilde{q}]{\vx_i})^2\\
  &\qquad - 4 \abs[\big]{F_{Y|X}(q|\vx_i) \eta^{(1)}(\muzero[q]{\vx_i}) - F_{Y|X}(\tilde{q}|\vx_i) \eta^{(1)}(\muzero[\tilde{q}]{\vx_i})}^2\\
  &\quad \leq 4 F_{Y|X}(q|\vx_i) \abs[\big]{\eta^{(1)}(\muzero[q]{\vx_i}) - \eta^{(1)}(\muzero[\tilde{q}]{\vx_i})}^2\\
  &\qquad + 4 \brk[\big]{F_{Y|X}(\tilde{q}|\vx_i) - F_{Y|X}(q|\vx_i)} \eta^{(1)}(\muzero[\tilde{q}]{\vx_i})^2\\
  &\quad \labrel{1}[\lesssim] (\tilde{q} - q)^2 + (\tilde{q} - q) \lesssim \tilde{q} - q,
\end{align*}
where in \labrel{1} we used that $\eta^{(1)}(\cdot)$ and $\eta(\cdot)$ on a fixed compact are Lipschitz and $\muzero[q]{\vx}$ is Lipschitz in $q$ uniformly over $\vx$ (therefore, $\eta(\muzero[q]{\vx}) = F_{Y|X}(q|\vx)$ is also Lipschitz in $q$ uniformly over $\vx$), as well as boundedness of $\muzero[q]{\vx}$ uniformly over $q$ and $\vx$.

\cref{prop:gl-holder-distribution-regression} is proven.
\end{proof}

\subsection{\texorpdfstring{$L_p$}{Lp} regression}\label{sec:example-lp-regression}

\begin{proposition}[Verification of \cref{ass:gl-holder} for \texorpdfstring{$L_p$}{Lp} regression]\label{prop:lp-regression}
  Suppose \cref{ass:dgp-holder,ass:gl-quasi-uniform-partition,ass:gl-local-basis} hold with $\mathcal{Q}$ a singleton, loss function $\riskfunc{y}{\eta} = \abs{ y - \eta }^p$, $p \in (1, 2]$, $\muzero{\cdot}$ as defined in \eqref{eq:mu-for-lp-reg-definition}, $\measureonY$ the Lebesgue measure. Assume the real inverse link function $\eta(\cdot)\colon \mathbb{R} \to \rangeeta$ is strictly monotonic and twice continuously differentiable with $\rangeeta$ an open connected subset of $\mathbb{R}$, and that $|\eta^{(1)}(\zeta)|$ is bounded away from zero uniformly over $\vx \in \mathcal{X}$ and $\zeta \in B(\vx)$ with $B(\vx)$ defined in~\eqref{eq:ballqx-def}.
Denoting by $a_l$ and $a_r$ the left and the right ends of $\rangeeta$ respectively (possibly $\pm \infty$), assume that $\int_{\mathbb{R}} \riskder{y}{a_l} f_{Y|X}(y|\vx)\,\mathrm{d}y < 0$ if $a_l$ is finite, and $\int_{\mathbb{R}} \riskder{y}{a_r} f_{Y|X}(y|\vx)\,\mathrm{d}y > 0$ if $a_r$ is finite.
Also assume that $\Expect{|y_1|^{\nu (p - 1)}} < \infty$ for some $\nu > 2$, and that $\vx \mapsto f_{Y|X}(y | \vx)$ is continuous for any $y \in \mathcal{Y}$.
In addition, assume that for each $\vx$ the map
\begin{equation*}
u \mapsto \int_{\mathbb{R}} \abs{u - y}^{p - 1} \sign(u - y) f_{Y|X}(y|\vx)\,\mathrm{d}y
\end{equation*}
is twice continuously differentiable with derivatives
\begin{equation*}
\frac{\mathrm{d}^j}{\mathrm{d} u^j} \int_{\mathbb{R}} \abs{u - y}^{p - 1} \sign(u - y) f_{Y|X}(y|\vx)\,\mathrm{d}y = \int_{\mathbb{R}} \abs{u - y}^{p - 1} \sign(u - y) \frac{\partial^j}{\partial y^j} f_{Y|X}(y|\vx)\,\mathrm{d}y,\quad j \in \crl{1, 2}.
\end{equation*}
Moreover, the function $\int_{\mathbb{R}} \abs{ \eta(\zeta) - y }^{p - 1} \sign(\eta(\zeta) - y) \frac{\partial}{\partial y} f_{Y|X}(y|\vx)\,\mathrm{d}y$ is bounded and bounded away from zero uniformly over $\vx \in \mathcal{X}$, $\zeta \in B(\vx)$, and the function $\int_{\mathbb{R}} \abs{ \eta(\zeta) - y }^{p - 1} \sign(\eta(\zeta) - y) \frac{\partial^2}{\partial y^2} f_{Y|X}(y|\vx)\,\mathrm{d}y$ is bounded in absolute value uniformly over $\vx \in \mathcal{X}$, $\zeta \in B(\vx)$. Then \cref{ass:gl-holder,ass:complexity,ass:complexity-addition}, \cref{eq:l2-norm-holder-type-restriction} and \cref{cond:nonconvex-consistency-additional-complexity-assumption} of \cref{lem:consistency-nonconvex} are also true.
\end{proposition}

\begin{proof}
Since $\mathcal{Q}$ is a singleton, we will omit the index $q$ in notations.

\paragraph{Verifying the assumptions of \cref{lem:mu-borel-measurable}}
The fact that
\begin{equation*}
\zeta \mapsto \int_{\mathbb{R}} \riskder{y}{\zeta} f_{Y|X}(y|\vx)\,\mathrm{d}y
\end{equation*}
is continuous is proven below in the verification of \cref{ass:gl-holder}\ref{ass:psi-holder-jumps}. To ensure that it crosses zero if $a_l$ or $a_r$ is not finite, we show
\begin{align}
  &\int_{\mathbb{R}} \riskder{y}{\zeta} f_{Y|X}(y|\vx)\,\mathrm{d}y \to - \infty\text{ as $\zeta \to -\infty$},\label{eq:psi-int-goes-to-minus-inf}\\
  &\int_{\mathbb{R}} \riskder{y}{\zeta} f_{Y|X}(y|\vx)\,\mathrm{d}y \to + \infty\text{ as $\zeta \to +\infty$}.\label{eq:psi-int-goes-to-plus-inf}
\end{align}
To prove~\eqref{eq:psi-int-goes-to-minus-inf}, recall that $\riskder{y}{\zeta} = p \abs{ y - \zeta }^{p - 1} \sign(\zeta- y)$ and therefore
\begin{align*}
  &\int_{\mathbb{R}} \riskder{y}{\zeta} f_{Y|X}(y|\vx)\,\mathrm{d}y = - p \int_\zeta^{+\infty} (y - \zeta)^{p - 1} f_{Y|X}(y|\vx)\,\mathrm{d}y + p \int^{\zeta}_{-\infty} (\zeta - y)^{p - 1} f_{Y|X}(y|\vx)\,\mathrm{d}y\\
  &\quad = -p (- \zeta)^{p - 1} \underbrace{\int_{\mathbb{R}} \prn[\bigg]{ 1 + \frac{y}{-\zeta} }^{p - 1} \Indlet\crl{y > \zeta} f_{Y|X}(y|\vx)\,\mathrm{d}y}_{\to 1}\\
  &\qquad + p \underbrace{\int_{\mathbb{R}} (\zeta - y)^{p - 1} \Indlet\crl{y \leq \zeta} f_{Y|X}(y|\vx)\,\mathrm{d}y}_{\to 0} \to - \infty,
\end{align*}
where we used dominated convergence because for $-\zeta > 1$ we have $1 + y / (-\zeta) \leq 1 + |y|$ in the first integral and $\zeta - y \leq -y = |y|$ in the second integral. \Cref{eq:psi-int-goes-to-plus-inf} is proven similarly.

\paragraph{Verifying \cref{ass:gl-holder}\ref{ass:abs-cont}}
The function $\riskfunc{y}{\eta}$ is continuously differentiable with respect to $\eta \in \mathbb{R}$, and its first derivative is the continuous function $\riskder{y}{\eta} = p \abs{ y - \eta }^{p - 1} \sign(\eta - y)$, therefore $\riskfunc{y}{\eta}$ for any fixed $y$ is absolutely continuous with respect to $\eta$ on bounded intervals.

\paragraph{Verifying \cref{ass:gl-holder}\ref{ass:psi-moments}}
The first-order optimality condition is true because $\muzero{\cdot}$ is defined this way in~\eqref{eq:mu-for-lp-reg-definition}.

The function
\begin{equation*}
  \sigma^2(\vx) := \Expect[\big]{\riskder{y_i}{\eta\prn{\muzero{\vx}}}^2 \given \vx_i = \vx} = p^2 \int_{\mathbb{R}} \abs{ y - \eta\prn{\muzero{\vx}} }^{2 p - 2} f_{Y|X}(y|\vx)\,\mathrm{d}y
\end{equation*}
is continuous because $\vx \mapsto |y - \eta\prn{\muzero{\vx}}|^{2 p - 2} f_{Y|X}(y|\vx)$ is continuous and dominated by $\prn{|y|^{2 p - 2} + C} C'$ for large enough constants $C$ and $C'$. As a continuous function on a compact set, $\sigma^2(\vx)$ is bounded away from zero because it is non-zero since $y_1$ has a conditional density.

The family of functions $\crl{ p \abs{ y - \eta \prn{ \mu_0(\vx) } }^{p - 1} \sign\prn{\eta \prn{ \mu_0(\vx) } - y} }$ only contains one element. Note that $\vx \mapsto \abs{ y - \eta \prn{ \mu_0(\vx) } }^{\nu(p - 1)} f_{Y|X}(y|\vx)$ is continuous and dominated by $(| y |^{\nu(p - 1)} + C) C'$ for large enough constants $C$ and $C'$. Therefore, the function
\begin{equation*}
  \vx \mapsto \int_{\mathbb{R}} \abs{ y - \eta \prn{ \mu_0(\vx) } }^{\nu(p - 1)} f_{Y|X}(y|\vx)\,\mathrm{d}y
\end{equation*}
is also continuous. As a continuous function on a compact set, it is bounded.

\paragraph{Verifying \cref{ass:gl-holder}\ref{ass:psi-holder-jumps}}
The function $x \mapsto |x|^{\alpha} \sign(x)$ for $\alpha \in (0, 1]$ is $\alpha$-H\"older for $x \in \mathbb{R}$ (with constant 2). Therefore, putting $\alpha := p - 1$, for any pair of reals $\zeta_1$ and $\zeta_2$ in a fixed bounded interval we have
\begin{align*}
  &\sup_{\vx} \sup_{\lambda \in [0, 1]} \sup_y \abs{ \riskder{y}{\eta\prn{\zeta_1 + \lambda(\zeta_2 - \zeta_1)}} - \riskder{y}{\eta\prn{\zeta_2}} }\\
  &\quad \leq 2 p \abs{ \eta\prn{\zeta_1 + \lambda(\zeta_2 - \zeta_1)} - \eta\prn{\zeta_2} }^{p - 1} \labrel{1}[\lesssim] \abs{ \zeta_1 - \zeta_2 }^{p - 1},
\end{align*}
where in \labrel{1} we used that the link function $\eta(\cdot)$ in a fixed bounded interval is Lipschitz.

\paragraph{Verifying \cref{ass:gl-holder}\ref{ass:psi-eta}}
The differentiability and boundedness conditions on $\Psione{\vx}{\eta}$ and $\Psitwo{\vx}{\eta}$ in~\cref{ass:gl-holder}\ref{ass:psi-eta} are directly assumed in the statement of \cref{prop:lp-regression}. The lower bound on $\Psione{\vx}{\eta(\zeta)} \eta^{(1)}(\zeta)^2$ follows from the assumed lower bounds on the first derivative expression and on $|\eta^{(1)}(\zeta)|$.

\paragraph{Verifying \cref{ass:complexity}}
Since $\mathcal{G}_2$, $\mathcal{G}_4$, $\mathcal{G}_5$ are just singletons (and the existence of corresponding envelopes holds trivially), it is enough to consider $\mathcal{G}_1$ and $\mathcal{G}_3$.

Assume that $\bbeta$ and $\widetilde{\bbeta}$ are such that $\norm{ \bbeta - \betazero }_{\infty} \leq r$ and $\norm{ \widetilde{\bbeta} - \betazero }_{\infty} \leq r$. Note that
\begin{align*}
  &\abs[\big]{\brk[\big]{\riskder{y}{\eta \prn{ \bp(\vx)\trans \bbeta }} - \riskder{y}{\eta \prn{ \mu_0(\vx) }}} - \brk[\big]{\riskder{y}{\eta\prn{\bp(\vx)\trans \widetilde{\bbeta}}} - \riskder{y}{\eta\prn{\mu_0(\vx)}}}}\\
  &\quad \leq 2 p \abs[\big]{ \eta\prn{\bp(\vx)\trans \bbeta} - \eta\prn{\bp(\vx)\trans \widetilde{\bbeta}} }^{p - 1} \lesssim \norm{ \bbeta - \widetilde{\bbeta} }_{\infty}^{p - 1}.
\end{align*}
The result for $\mathcal{G}_1$ follows.

Similarly,
\begin{equation}\label{eq:lp-reg-ggg-bound}
  \begin{aligned}
    &\abs[\big]{\brk[\big]{\riskder{y}{\eta \prn{ \bp(\vx)\trans \bbeta }} - \riskder{y}{\eta \prn{ \bp(\vx)\trans \betazero }}} - \brk[\big]{\riskder{y}{\eta\prn{\bp(\vx)\trans \widetilde{\bbeta}}} - \riskder{y}{\eta\prn{\bp(\vx)\trans \betazero}}}}\\
    &\quad \leq 2 p \abs[\big]{ \eta\prn{\bp(\vx)\trans \bbeta} - \eta\prn{\bp(\vx)\trans \widetilde{\bbeta}} }^{p - 1} \lesssim \norm{ \bbeta - \widetilde{\bbeta} }_{\infty}^{p - 1}.
  \end{aligned}
\end{equation}
For a fixed cell $\delta \in \Delta$ the class of functions of the form
\begin{equation*}
  \brk[\big]{ \riskder[q]{y}{\eta\prn{\vp(\vx)\trans \vbeta}} - \riskder[q]{y}{\eta\prn{\vp(\vx)\trans \betazero }} } \Indlet\crl{\vx \in \mathcal{N}_{[c \log n]}(\delta)}
\end{equation*}
can be parametrized by $\bbeta$ lying in a fixed vector space $\mathcal{B}_{\delta}$ of dimension $O(\log^d n)$. The result now follows from the bound~\eqref{eq:lp-reg-ggg-bound} (by using \eqref{eq:log-d-entropy-chain} similarly to the proof of \cref{prop:gl-holder-logistic-regression}).

\paragraph{Verifying the addition to \cref{ass:complexity-addition}}
Fix $\delta \in \Delta$. Let $\bbeta$ and $\widetilde{\bbeta}$ be such that $\norm{ \bbeta - \betazero }_{\infty} \leq r$ and $\norm{ \widetilde{\bbeta} - \betazero }_{\infty} \leq r$; let $\vv$ and $\widetilde{\vv}$ be such that $\norm{ \vv }_{\infty} \leq \epsninassump$ and $\norm{ \widetilde{\vv} }_{\infty} \leq \epsninassump$. To declutter notation, put
\begin{align*}
  &g(t) := \brk[\big]{ \riskder{y}{\eta\prn{\bp(\vx)\trans \prn{\betazero + \bbeta } + t}} - \riskder{y}{\eta\prn{\bp(\vx)\trans \betazero}} }\\
  &\qquad \times \eta^{(1)} \prn{ \bp(\vx)\trans \prn{ \betazero + \bbeta } + t },\\
  &\tilde{g}(t) := \brk[\big]{ \riskder{y}{\eta\prn{\bp(\vx)\trans \prn{\betazero + \widetilde{\bbeta} } + t}} - \riskder{y}{\eta\prn{\bp(\vx)\trans \betazero}} }\\
  &\qquad \times \eta^{(1)} \prn{ \bp(\vx)\trans \prn{ \betazero + \widetilde{\bbeta} } + t }.
\end{align*}

Note that
\begin{align*}
  \int_{-\bp(\vx)\trans\vv}^0 g(t)\,\mathrm{d} t - \int_{-\bp(\vx)\trans\widetilde{\vv}}^0 \tilde{g}(t)\,\mathrm{d} t = \int_{-\bp(\vx)\trans\vv}^0 \prn[\big]{g(t) - \tilde{g}(t)} \,\mathrm{d} t + \int_{-\bp(\vx)\trans\vv}^{-\bp(\vx)\trans\widetilde{\vv}} \tilde{g}(t)\,\mathrm{d} t.
\end{align*}
Since $\riskder{y}{\cdot}$ is $(p-1)$-H\"older continuous in the second argument and functions $\eta(\cdot)$, $\eta^{(1)}(\cdot)$ in a fixed bounded interval are Lipschitz, we get that uniformly over $t$ and $\vx$ in these integrals
\begin{equation*}
  \abs[\big]{ g(t) - \tilde{g}(t) } \lesssim \norm[\big]{ \bbeta - \widetilde{\bbeta} }_{\infty}^{p - 1} + \norm[\big]{ \bbeta - \widetilde{\bbeta} }_{\infty},
\end{equation*}
and $\abs[\big]{ \tilde{g}(t) }$ is bounded. This gives
\begin{equation*}
  \abs[\bigg]{ \int_{-\bp(\vx)\trans\vv}^0 g(t)\,\mathrm{d} t - \int_{-\bp(\vx)\trans\widetilde{\vv}}^0 \tilde{g}(t)\,\mathrm{d} t } \lesssim \epsninassump \prn[\big]{ \norm[\big]{ \bbeta - \widetilde{\bbeta} }_{\infty}^{p - 1} + \norm[\big]{ \bbeta - \widetilde{\bbeta} }_{\infty} } + \norm[\big]{ \vv - \widetilde{\vv} }_{\infty}
\end{equation*}
It means that taking an $\varepsilon$-net (for $\varepsilon$ smaller than 1) in the space of $\bbeta$ and an $\epsninassump \varepsilon$-net in the space of $\vv$ induces an $\oC{lp-reg-last-eps} \epsninassump \varepsilon$-net in the space of functions
\begin{equation*}
  \crl[\bigg]{ (\vx, y) \mapsto \int_{-\bp(\vx)\trans\vv}^0 g(t)\,\mathrm{d} t \cdot \Indlet\crl{\vx \in \mathcal{N}_{[c \log n]}(\delta)}: \norm{ \vbeta - \betazero[q] }_{\infty} \leq r, \norm{ \vv }_{\infty} \leq \epsninassump }
\end{equation*}
in terms of the sup-norm, where $\nC{lp-reg-last-eps}$ is some constant. Possibly increasing $\oC{lp-reg-last-eps}$, we can conclude that this class with envelope $\oC{lp-reg-last-eps}$ satisfies the uniform entropy bound \eqref{eq:uniform-entropy-bound} with $A \lesssim 1$ and $V \lesssim \log^d n$ (where we used that $\bbeta$ and $\vv$ can be assumed to lie in a vector space of dimension $O(\log^d n)$). By \eqref{eq:log-d-entropy-chain}, the same can be said about the union of $O \prn{ h^{-d} }$ such classes (corresponding to different $\delta$). The verification is concluded.

\paragraph{Verifying \eqref{eq:l2-norm-holder-type-restriction}} This is obvious because $\mathcal{Q}$ is a singleton.

\paragraph{Verifying \cref{lem:consistency-nonconvex}\ref{cond:nonconvex-consistency-additional-complexity-assumption}}
It follows from the proof of \cref{lem:assumptions-of-nonconvex-consistency-are-reasonable} that
\begin{align*}
  &\abs[\big]{\riskfunc{y}{\eta(\bp(\vx)\trans \bbeta)} - \riskfunc{y}{\eta(\bp(\vx)\trans \widetilde{\bbeta})}} \lesssim \prn[\big]{1 + \riskderenvel{\vx}{y}} \abs[\big]{\bp(\vx)\trans (\bbeta - \widetilde{\bbeta})}\\
  &\quad \lesssim \prn[\big]{1 + \riskderenvel{\vx}{y}} \norm[\big]{\bbeta - \widetilde{\bbeta}}_{\infty}.
\end{align*}
The required uniform entropy bound follows immediately from this.

\cref{prop:lp-regression} is proven.
\end{proof}

\subsection{Logistic regression}\label{sec:example-logistic-regression}

\begin{proposition}[Verification of \cref{ass:gl-holder} and others for logistic regression]\label{prop:gl-holder-logistic-regression}
  Suppose \cref{ass:dgp-holder,ass:gl-quasi-uniform-partition,ass:gl-local-basis} hold with $\mathcal{Q}$ a singleton, $\mathcal{Y} = \crl{ 0, 1 }$, $\eta(\theta) = 1 / \prn{1 + e^{-\theta}}$, $\measureonY$ is the counting measure on $\crl{0, 1}$, and $\riskfunc{y}{\eta} = - y \log(\eta) - (1 - y) \log(1 - \eta)$. Assume also $\newtarget{def:logRegPi}{\logregpiNoLink}(\vx) := \Prob*{y_1 = 1 \given \vx_1 = \vx}$ is continuous and $\logregpi(\vx)$ lies in the interval $(0, 1)$ for any $\vx \in \mathcal{X}$. Then \cref{ass:gl-holder,ass:complexity,ass:complexity-addition}, \eqref{eq:l2-norm-holder-type-restriction},  \cref{cond:lower-bound-of-psi-one,cond:nonconvex-consistency-additional-complexity-assumption} of \cref{lem:consistency-nonconvex} are true.
\end{proposition}

We will prove \cref{prop:gl-holder-logistic-regression} now. Since $\mathcal{Q}$ is a singleton, we will omit the index $q$ in notations.

\paragraph{Verifying \cref{ass:gl-holder}\ref{ass:abs-cont}}
The function $\riskfunc{y}{\eta}$ is infinitely smooth with respect to $\eta \in (0, 1)$, and its first derivative is $\riskder{y}{\eta} = (1 - y) / (1 - \eta) - y / \eta$, therefore $\riskfunc{y}{\eta}$ for any fixed $y$ is absolutely continuous with respect to $\eta$ on bounded intervals.

\paragraph{Verifying \cref{ass:gl-holder}\ref{ass:psi-moments}}
We have $\Expect{\riskder{y_i}{\eta\prn{\muzero{\vx_i}}} \given \vx_i} = 0$ since
\begin{equation*}
\eta(\mu_0(\vx_i)) = \Expect{y_i \given \vx_i} = \Prob{y_i = 1 \given \vx_i} = \logregpi(\vx_i).
\end{equation*}

Next,
\begin{align*}
  &\sigma^2(\vx) = \Expect[\big]{\riskder{y_i}{\eta\prn{\muzero{\vx_i}}}^2 \given \vx_i = \vx} = \frac{\Expect[\big]{\prn{y_i - \logregpi(\vx_i)}^2 \given \vx_i = \vx}}{\eta^{(1)} \prn{ \mu_0(\vx) }^2}\\
  &\quad = \frac{\logregpi(\vx) \prn{ 1 - \logregpi(\vx) }}{\logregpi(\vx)^2 \prn{ 1 - \logregpi(\vx)}^2} = \frac{1}{\logregpi(\vx) \prn{ 1 - \logregpi(\vx)}}
\end{align*}
is continuous and bounded away from zero (is not less than $4$).

Since $\vx_i$ lies in a compact set, $\riskder{y}{\eta\prn{\muzero{\vx_i}}}$ is bounded, so it has moments of any order.

\paragraph{Verifying \cref{ass:gl-holder}\ref{ass:psi-holder-jumps}}
The function $\riskfunc{y}{\eta}$ is infinitely smooth with respect to $\eta \in (0, 1)$, and its first derivative is $\riskder{y}{\eta} = (1 - y) / (1 - \eta) - y / \eta$. Using the famous expression for the derivative of the logistic function $\eta^{(1)}(\theta) = \eta(\theta) \prn{ 1 - \eta(\theta) }$, we get
\begin{align*}
  &\frac{\partial}{\partial \theta} \riskfunc{y}{\eta(\theta)} = \riskder{y}{\eta(\theta)} \eta^{(1)}(\theta) = (1 - y) \eta(\theta) - y \prn{ 1 - \eta(\theta) } = \eta(\theta) - y,\\
  &\frac{\partial^2}{\partial \theta^2} \riskfunc{y}{\eta(\theta)} = \eta^{(1)}(\theta) = \eta(\theta) \prn{ 1 - \eta(\theta) }.
\end{align*}
Since the logistic link maps to $(0, 1)$, the second derivative is positive (and does not depend on $y$). Therefore, $\riskfunc{y}{\eta(\theta)}$ is convex with respect to $\theta$ for any $y$.

Uniformly over $\zeta_1$ and $\zeta_2$ in a fixed bounded interval, we have
\begin{equation*}
  \sup_y \abs[\big]{ \riskder{y}{\eta(\zeta_1)} - \riskder{y}{\eta(\zeta_2)} } \lesssim \abs[\big]{ \eta(\zeta_1) - \eta(\zeta_2) } \overset{\text{(a)}}{\leq} \abs{ \zeta_1 - \zeta_2 },
\end{equation*}
where in (a) we used that the derivative of $\eta(\cdot)$ does not exceed $1$. So in this case the H\"older parameter $\alpha = 1$.

The logistic function $\eta(\cdot)$ is strictly monotonic and infinitely smooth on $\mathbb{R}$.

\paragraph{Verifying \cref{ass:gl-holder}\ref{ass:psi-eta}}
In this case
\begin{equation*}
  \Psi \prn{ \vx; \eta } = \frac{\eta - \Expect{y_i \given \vx_i = \vx}}{\eta (1 - \eta)} = \frac{\eta - \logregpi(\vx)}{\eta (1 - \eta)}.
\end{equation*}

This function is infinitely smooth in $\eta$ for $\eta \in (0, 1)$. Its first derivative
\begin{equation*}
  \Psione{\vx}{\eta} = \frac{\partial}{\partial \eta} \Psi \prn{ \vx; \eta } = \frac{\eta^2 - 2 \eta \logregpi(\vx) + \logregpi(\vx)}{\eta^2 (1 - \eta)^2}.
\end{equation*}
Therefore,
\begin{equation}\label{eq:logistic-psi-one-eta-one-squared}
  \Psione{\vx}{\eta(\zeta)} \eta^{(1)}(\zeta)^2 = \eta(\zeta)^2 - 2 \eta(\zeta) \logregpi(\vx) + \logregpi(\vx).
\end{equation}
If $\abs{ \zeta - \mu_0(\vx) } \leq r$, then $\abs[\big]{ \eta(\zeta) - \eta \prn{ \mu_0(\vx) } } = \abs[\big]{ \eta(\zeta) - \logregpi(\vx) } \leq r$ (since the derivative of $\eta(\cdot)$ does not exceed 1). Since $0 < \min_{\vx \in \mathcal{X}} \logregpi(\vx) \prn{ 1 - \logregpi(\vx) } \leq \max_{\vx \in \mathcal{X}} \logregpi(\vx) \prn{ 1 - \logregpi(\vx) } < 1$, for small enough $r$ the right-hand side of~\eqref{eq:logistic-psi-one-eta-one-squared} for such $\zeta$ is also bounded away from zero and one.

Finally,
\begin{equation*}
  \Psitwo{\vx}{\eta} = \frac{\partial}{\partial \eta} \Psione{\vx}{\eta} = \frac{2 \prn[\big]{ \eta^3 - 3 \logregpi(\vx) \eta^2 + 3 \logregpi(\vx) \eta - \logregpi(\vx) }}{\eta^3 (1 - \eta)^3}.
\end{equation*}

Again since $0 < \min_{\vx \in \mathcal{X}} \logregpi(\vx) \prn{ 1 - \logregpi(\vx) } \leq \max_{\vx \in \mathcal{X}} \logregpi(\vx) \prn{ 1 - \logregpi(\vx) } < 1$, for such $\zeta$ that $\abs{ \zeta - \mu_0(\vx) } \leq r$ and $r$ small enough, the product $\eta(\zeta) \prn{ 1 - \eta(\zeta) }$ is bounded away from zero. So for such $\zeta$, $\abs{ \Psitwo{\vx}{\eta(\zeta)} }$ is uniformly bounded.

\begin{lemma}[Class \texorpdfstring{$G_{1}$}{G1}]\label{lem:logistic-entropy-one}
  The class
  \begin{equation*}
    \mathcal{G}_1 = \crl[\bigg]{ \mathcal{X} \times \mathbb{R} \ni (\vx, y) \mapsto \frac{\eta(\vp(\vx)\trans \vbeta) - y}{\eta^{(1)}(\vp(\vx)\trans \vbeta)} - \frac{\logregpi(\vx) - y}{\logregpi(\vx) \prn[\big]{ 1 - \logregpi(\vx) }}: \norm{ \vbeta - \betazero } \leq r }
  \end{equation*}
  with a large enough constant envelope satisfies the uniform entropy bound~\eqref{eq:uniform-entropy-bound} with $A \lesssim 1$ and $V \lesssim K$.
\end{lemma}

\begin{proof}[Proof of \cref{lem:logistic-entropy-one}]
The class
\begin{equation*}
  \mathcal{G}_{1 1} = \crl[\big]{ (\vx, y) \mapsto \vp(\vx)\trans \vbeta: \vbeta \in \mathbb{R}^K }
\end{equation*}
is VC with index not exceeding $K + 2$ by Lemma 2.6.15 in~\parencite{vaart96_weak_conver_empir_proces}. Since in a fixed bounded interval $\eta(\cdot)$ is Lipschitz and $(\vx, y) \mapsto y$ is one fixed function, the class
\begin{equation*}
  \mathcal{G}_{1 2} = \crl[\big]{ (\vx, y) \mapsto \eta(\vp(\vx)\trans \vbeta) - y: \norm{ \vbeta - \betazero }_{\infty} \leq r }
\end{equation*}
with a large enough constant envelope (recall that $\mathcal{Y}$ is a bounded set) satisfies the uniform entropy bound~\eqref{eq:uniform-entropy-bound} with $A \lesssim 1$ and $V \lesssim K$. Since $1 / \eta^{(1)}(\cdot)$ in a fixed bounded interval is Lipschitz, the same is true of
\begin{equation*}
  \mathcal{G}_{1 3} = \crl[\big]{ (\vx, y) \mapsto \eta^{(1)}(\vp(\vx)\trans \vbeta)^{-1}: \norm{ \vbeta - \betazero }_{\infty} \leq r },
\end{equation*}
where we used again that under these constraints $\eta^{(1)}\prn{\vp(\vx)\trans \vbeta}$ is bounded away from zero. This implies by \cref{lem:permanence-properties} that it is true of $\mathcal{G}_{1 2} \cdot \mathcal{G}_{1 3} - \riskder{y}{\eta\prn{\muzero{\vx}}}$ (since $\riskder{y}{\eta\prn{\muzero{\vx}}}$ is one fixed function), which is what we need.
\end{proof}

\begin{lemma}[Class \texorpdfstring{$\mathcal{G}_3$}{G3}]\label{lem:logistic-entropy-two}
  The class
  \begin{align*}
    \mathcal{G}_3 := \crl[\Big]{&\mathcal{X} \times \mathbb{R} \ni (\vx, y) \mapsto\\
                      &\brk[\bigg]{ \frac{\eta\prn{\vp(\vx)\trans \vbeta} - y}{\eta^{(1)}\prn{\vp(\vx)\trans \vbeta}} - \frac{\eta\prn{\vp(\vx)\trans \betazero} - y}{\eta^{(1)}\prn{\vp(\vx)\trans \betazero}} } \Indlet\crl{\vx \in \mathcal{N}_{[c \log n]}(\delta)} \colon\\
                      &\norm{ \vbeta - \betazero[q] }_{\infty} \leq r, \delta \in \Delta}
  \end{align*}
  with a large enough constant envelope satisfies the uniform entropy bound~\eqref{eq:uniform-entropy-bound} with $A \lesssim 1$ and $V \lesssim \log^d n$.
\end{lemma}

\begin{proof}[Proof of~\cref{lem:logistic-entropy-two}]
For a fixed vector space $\mathcal{B}$ of dimension $\dim \mathcal{B}$,
\begin{equation*}
  \mathcal{W}_{\mathcal{B}} := \crl[\bigg]{ (\vx, y) \mapsto \frac{\eta\prn{\vp(\vx)\trans \vbeta} - y}{\eta^{(1)}\prn{\vp(\vx)\trans \vbeta}}: \vbeta \in \mathcal{B}, \norm{ \vbeta - \betazero[q] }_{\infty} \leq r }
\end{equation*}
with a large enough constant envelope (recall that $\mathcal{Y}$ is a bounded set) satisfies the uniform entropy bound with $A \lesssim 1$ and $V \lesssim \dim \mathcal{B}$ by the same argument as in the proof of \cref{lem:logistic-entropy-one}. Therefore, for any fixed $c > 0$, $\delta \in \Delta$, the class
\begin{equation*}
  \mathcal{W}_{3, \delta} := \crl[\bigg]{ (\vx, y) \mapsto \frac{\eta\prn{\vp(\vx)\trans \vbeta} - y}{\eta^{(1)}\prn{\vp(\vx)\trans \vbeta}} \Indlet\crl{\vx \in \mathcal{N}_{[c \log n]}(\delta)}: \vbeta \in \mathbb{R}^K, \norm{ \vbeta - \betazero[q] }_{\infty} \leq r }
\end{equation*}
with a large enough constant envelope also satisfies the uniform entropy bound with $A \lesssim 1$ and $V \lesssim \log^d n$ because it is contained in the product of $\mathcal{W}_{\mathcal{B}_{3, \delta}}$ for some vector space $\mathcal{B}_{3, \delta}$ of dimension $\dim \mathcal{B}_{3, \delta} \lesssim \log^d n$ and a fixed function $\Indlet\crl{\vx \in \mathcal{N}_{[c \log n]}(\delta)}$. Subtracting a fixed bounded function does not change this fact, so the same is true of
\begin{align*}
  \mathcal{G}_{3, \delta} := \crl[\Big]{&(\vx, y) \mapsto\\
  &\brk[\bigg]{ \frac{\eta\prn{\vp(\vx)\trans \vbeta} - y}{\eta^{(1)}\prn{\vp(\vx)\trans \vbeta}} - \frac{\eta\prn{\vp(\vx)\trans \betazero} - y}{\eta^{(1)}\prn{\vp(\vx)\trans \betazero}}} \Indlet\prn{\vx \in \mathcal{N}_{[c \log n]}(\delta)}:\\
  &\vbeta \in \mathbb{R}^K, \norm{ \vbeta - \betazero[q] }_{\infty} \leq r}
\end{align*}
Since there are $O\prn{h^{-d}}$ such classes and $\log(1 / h) \lesssim \log n$, using the chain~\eqref{eq:log-d-entropy-chain} we obtain that $\mathcal{G}_3$ satisfies the uniform entropy bound~\eqref{eq:uniform-entropy-bound} with $A \lesssim 1$ and $V \lesssim \log^d n$.
\end{proof}

\paragraph{Verifying \cref{ass:complexity}}
Classes $\mathcal{G}_2$, $\mathcal{G}_4$, $\mathcal{G}_5$ are just singletons (and the existence of corresponding envelopes holds trivially). Classes $\mathcal{G}_1$ and $\mathcal{G}_3$ are tackled in \cref{lem:logistic-entropy-one} and \cref{lem:logistic-entropy-two}.

\paragraph{Verifying \cref{ass:complexity-addition}}
This is verified (in a more general setting) in \cref{sec:simplified-example-verifying-addition-for-nonconvex-bahadur}.

\paragraph{Verifying \eqref{eq:l2-norm-holder-type-restriction}}
This is obvious since $\mathcal{Q}$ is a singleton.

\paragraph{Verifying the condition in \cref{lem:consistency-nonconvex}\ref{cond:lower-bound-of-psi-one}}
Recall that in this case
\begin{equation*}
  \Psione{\vx}{\eta} = \frac{\eta^2 - 2 \eta \logregpi(\vx) + \logregpi(\vx)}{\eta^2 (1 - \eta)^2}.
\end{equation*}

The numerator is always positive since $0 < \logregpi(\vx) < 1$, and the denominator is also positive since $\eta \in (0, 1)$. Since $\Psione{\vx}{\eta}$ is continuous in both arguments and the image of a compact set under a continuous mapping is compact, we see that for any fixed compact subset of $(0, 1)$, $\Psione{\vx}{\eta}$ is bounded away from zero uniformly over $\vx \in \mathcal{X}$ and $\eta$ lying in this compact subset.

\paragraph{Verifying the condition in \cref{lem:consistency-nonconvex}\ref{cond:nonconvex-consistency-additional-complexity-assumption}}
In this verification, we will use $\theta_1 := \bp(\vx)\trans \bbeta$, $\tilde{\theta}_1 := \bp(\vx)\trans \widetilde{\bbeta}$ and $\theta_2 := \bp(\vx)\trans \betazero$ to simplify notations. Note that for $\theta$ lying in a fixed compact, both functions $\log \eta(\theta)$ and $\log(1 - \eta(\theta))$ are Lipschitz in $\theta$, so if $\norm{\bbeta}_{\infty} \leq R$ and $\norm{\widetilde{\bbeta}}_{\infty} \leq R$, we have
\begin{equation*}
  |\riskfunc{y}{\eta(\theta_1)} - \riskfunc{y}{\eta(\theta_2)} - \riskfunc{y}{\eta(\tilde{\theta}_1)} + \riskfunc{y}{\eta(\theta_2)}| \lesssim |\theta_1 - \tilde{\theta}_1| \lesssim |\bbeta - \widetilde{\bbeta}|,
\end{equation*}
where the constants in $\lesssim$ are allowed to depend on $R$ but not on $\bbeta, \widetilde{\bbeta}$, $\vx$ or $y$ (we used that $y$ and $1 - y$ are bounded by 1). The result follows.

This concludes the proof of \cref{prop:gl-holder-logistic-regression}.

\subsection{Examples of basis functions}\label{sec:examples-basis-functions}

Fix $0 \leq s < m$. Take any segment $[a, b]$ and partition it into $J$ sub-segments by taking $J - 1$ knots $a < \tau_1 < \ldots \tau_{J - 1} < b$. Consider a non-increasing tuple of real numbers $\prn{t_k}_{k = 1}^{2 m + (m - s) (J - 1)}$ defined as follows:
\begin{align*}
  &t_1 = \ldots = t_m = a,\\
  &t_{m + 1} = \ldots = t_{2 m - s} = \tau_1,\\
  &\cdots\\
  &t_{m + (m - s)(J - 2) + 1} = \ldots = t_{m + (m - s) (J - 1)} = \tau_{J - 1},\\
  &t_{m + (m - s) (J - 1) + 1} = \ldots = t_{2 m + (m - s) (J - 1)} = b,
\end{align*}
and put for $l \in \crl{1, 2, \ldots, m + (m - s) (J - 1)}$
\begin{gather*}
  p_l(x) = (-1)^m (t_{l + m} - t_l) [t_l, \ldots, t_{l + m}] (x - t)_+^{m - 1},\quad x \in [a, b),\\
  p_{m + (m - s) (J - 1)}(b) = \lim_{x \toup b} p_{m + (m - s) (J - 1)}(x),
\end{gather*}
where $[t_l, \ldots, t_{l + m}] g(t, x)$ denotes the divided difference \cite[Definition 2.49]{Schumaker2007} of $t \mapsto g(t, x)$ over points $t_l, \ldots, t_{l + m}$, and $\alpha_+ := \alpha \vee 0$. We will call the basis $\crl{p_l(\cdot)}_{l = 1}^{m + (m - s) (J - 1)}$ the (normalized) \textit{B-spline basis} of order $m$ with knots $\tau_1 < \ldots < \tau_{J - 1}$, each of multiplicity $m - s$. The functions $\crl{p_l(\cdot)}_{l = 1}^{m + (m - s) (J - 1)}$ form a basis in the function space
\begin{equation*}
\crl[\big]{f(\cdot) \in \smoothnoprn{s - 1}{[a, b]}: f(\cdot)\text{ is a polynomial of degree $m - 1$ on each $[\tau_j, \tau_{j + 1}]$, $j \in \crl{0, \ldots, J - 1}$}},
\end{equation*}
where we put $\tau_0 = a$, $\tau_J = b$ for simplicity (\cite[Corollary 4.10]{Schumaker2007})\footnote{The symbol $\smoothnoprn{-1}{[a, b]}$ here corresponds to the family of all functions $[a, b] \to \mathbb{R}$ (with no smoothness restrictions).}; moreover, each $p_l(\cdot)$ is strictly positive on $(t_l, t_{l + m})$, zero on the complement of $[t_l, t_{l + m}]$, and does not exceed $1$ (\cite[Theorem 4.9]{Schumaker2007}).

This definition can be extended to dimension $d$ by considering $\mathcal{X} = \bigotimes_{\ell = 1}^d [a_{\ell}, b_{\ell}]$ and the corresponding tensor products
\begin{equation*}
p_{l_1, \ldots, l_d}(\vx) = p_{l_1}(x_1) \ldots p_{l_d}(x_d)
\end{equation*}
for all tuples $\vl = (l_1, \ldots, l_d)$ such that $l_{\ell} \in \crl{1, 2, \ldots, m + (m - s) (J_{\ell} - 1)}$. The $d$ partitions of $[a_{\ell}, b_{\ell}]$ induce a tensor-product partition $\Delta$ of $\mathcal{X}$ (with $\kappa = J_1 \ldots J_d$ cells), and the basis $\vp(\vx)$ (arranged in a lexicographic order of $\vl$) is called the \textit{tensor-product B-spline basis} of order $m$ associated with the tensor-product partition $\Delta$ (\cite[Definition 12.3]{Schumaker2007}).

We refer to \cite{Schumaker2007} for details.

\begin{proposition}[Verifying \cref{ass:gl-local-basis,ass:gl-approximation-error}]
Suppose \cref{ass:dgp-holder,ass:gl-quasi-uniform-partition} hold with the tensor-product partition $\Delta$ as described above. Let $\vp(\vx)$ be a tensor-product B-spline basis of order $m$ associated with $\Delta$. Then $\vp(\vx)$ satisfies \cref{ass:gl-local-basis,ass:gl-approximation-error}.
\end{proposition}

\begin{proof}
It follows from the more general argument for Lemma SA-6.1 in \cite{Cattaneo-Farrell-Feng_2020_AOS}.
\end{proof}

Taking $s = m - 1$ recovers standard tensor-product B-splines with simple (multiplicity 1) partition knots, whereas $s = 0$ corresponds to piecewise polynomials, where, in particular, each basis function is only supported on one cell (making \cref{lem:holder-consistency-unconnected,lem:nonconvex-cons-unconnected} applicable). Additional examples are provided in \cite[Section SA-6]{Cattaneo-Farrell-Feng_2020_AOS}.

\section{Other parameters of interest}\label{sec:other-parameters-of-interest}

This section formalizes the discussion in Section \mainref{sec: Extensions} of the paper. The following theorem is now a simple corollary of the previous results presented in this supplemental appendix.

\begin{theorem}[Other parameters of interest]\label{th:other-parameters}\leavevmode
  \begin{assertions}[wide]
  \item\label{other-params-one} Suppose all the conditions of \cref{th:strong-approximation-binscatter-yurinski}\ref{th:yurinski-normal} hold with $\vv = \vzero$ and $\ruc = o(1)$. Then
  \begin{align*}
    &\sup_{\vq \in \mathcal{Q}} \sup_{\vx \in \mathcal{X}} \abs[\Bigg]{\frac{\eta \prn[\big]{\muhat[\vq]{\vx}} - \eta\prn[\big]{\muzero[\vq]{\vx}}}{\bigl|\eta^{(1)}(\muzero[\vq]{\vx})\bigr| \sqrt{\bigomegabar[\vq][\vzero] / n}} - \sign\crl[\big]{\eta^{(1)}(\muzero[\vq]{\vx})} \elltwobarqx[\vq][\vzero]\trans \zn(\vq)}\\
    &\quad = o_{\Problet}\prn{\rstr} + O_{\Problet}\prn[\Big]{\sqrt{n h^d} \prn[\big]{\ruc[2] + \rbah}}
  \end{align*}
  with $\zn(\vq)$ defined in \cref{th:strong-approximation-binscatter-yurinski}.

  \item\label{other-params-two} Fix any $k \in \crl{ 1, \ldots, d }$. Suppose all the conditions of \cref{th:strong-approximation-binscatter-yurinski}\ref{th:yurinski-normal} hold with $\vv = \ve_k$, where $\ve_k = (0, \ldots, 0, 1, 0, \ldots, 0)\trans \in \mathbb{R}^d$ with 1 at the $k$th place. Then
  \begin{align*}
    &\sup_{\vq, \vx} \abs[\Bigg]{\frac{\eta^{(1)} \prn[\big]{\muhat[\vq]{\vx}} \muhat[\vq][(\ve_k)]{\vx} - \eta^{(1)} \prn[\big]{\muzero[\vq]{\vx}} \muzero[\vq][(\ve_k)]{\vx}}{\bigl|\eta^{(1)} (\muzero[\vq]{\vx})\bigr| \sqrt{\bigomegabar[\vq][\ve_k]}}\\
    &\qquad - \sign\crl[\big]{\eta^{(1)}(\muzero[\vq]{\vx})} \elltwobarqx[\vq][\ve_k]\trans \zn(\vq)}\\
    &\quad = o_{\Problet}\prn{\rstr} + O_{\Problet}\prn[\Big]{\sqrt{n h^d} \prn[\big]{\ruc[2] + \rbah + h \ruc}}.
  \end{align*}

  \item\label{other-params-three} If $\riskderenvel{\vx_i}{y_i}$ is \subgauss\ conditionally on $\vx_i$, then in \cref{other-params-one,other-params-two} $\rstr$ can be replaced with $\rstrsubgauss$.
  \end{assertions}
\end{theorem}

\begin{proof}
\ref{other-params-one} We have uniformly over $\vq$ and $\vx$
\begin{align*}
  &\eta \prn{ \muhat[\vq]{\vx} } - \eta \prn{ \muzero[\vq]{\vx} }\\
  &\quad \overset{\text{(a)}}{=} \eta^{(1)} \prn{ \muzero[\vq]{\vx} } \prn{ \muhat[\vq]{\vx} - \muzero[\vq]{\vx} } + \frac{\eta^{(2)}(\xi)}{2} \prn{ \muhat[\vq]{\vx} - \muzero[\vq]{\vx} }^2\\
  &\quad \overset{\text{(b)}}{=} \eta^{(1)} \prn{ \muzero[\vq]{\vx} } \prn{ \muhat[\vq]{\vx} - \muzero[\vq]{\vx} } + O_{\Problet} \prn{ \ruc[2] }\\
  &\quad \overset{\text{(c)}}{=} \eta^{(1)} \prn{ \muzero[\vq]{\vx} } \bart[\vq][\vzero] \sqrt{\bigomegabar[\vq][\vzero] / n} + O_{\Problet} \prn{ \ruc[2] + \rbah}\\
  &\quad \overset{\text{(d)}}{=} \eta^{(1)} \prn{ \muzero[\vq]{\vx} } \sqrt{\bigomegabar[\vq][\vzero] / n} \brk[\Big]{ \bart[\vq][\vzero] + O_{\Problet} \prn[\big]{ \sqrt{n h^d} \prn{\ruc[2] + \rbah} } }.
\end{align*}
Here, (a) is by Taylor expansion, with some $\xi = \xi_{\vq, \vx}$ between $\muhat[\vq]{\vx}$ and $\muzero[\vq]{\vx}$. In (b), we used $\ruc = o(1)$ (giving that $\eta^{(2)}(\xi)$ does not exceed a fixed constant not depending on $\vq$ or $\vx$) and \eqref{eq:gRcohQ}. (c) is by \eqref{eq:TvfoLC} and since $\eta^{(1)} ( \muzero[\vq]{\vx} )$ is uniformly bounded. (d) is by $h^{-2 \multindnorm{\vv} - d} \lesssim_{\Problet} \inf_{\vq, \vx} |\bigomegabar[\vq]|$ (by \cref{lem:asympt-variance}) and since $\abs{\eta^{(1)}( \muzero[\vq]{\vx})}$ is bounded away from zero by \cref{ass:gl-holder}\ref{ass:psi-eta}.

Rewriting, we obtain
\begin{equation*}
  \sup_{\vq, \vx} \abs[\bigg]{ \frac{\eta (\muhat[\vq]{\vx}) - \eta( \muzero[\vq]{\vx})}{\eta^{(1)}( \muzero[\vq]{\vx}) \sqrt{\bigomegabar[\vq][\vzero] / n}} - \bart[\vq][\vzero] } \lesssim_{\Problet} \sqrt{n h^d} \big(\ruc[2] + \rbah\big).
\end{equation*}
It is left to combine this with \eqref{eq:yurinski-corollary-statement} and use the triangle inequality.

\ref{other-params-two} We have uniformly over $\vq$ and $\vx$
\begin{align*}
  &\eta^{(1)} \prn{ \muhat[\vq]{\vx} } \muhat[\vq][(\ve_k)]{\vx} - \eta^{(1)} \prn{ \muzero[\vq]{\vx} } \muzero[\vq][(\ve_k)]{\vx}\\
  &\quad \overset{\text{(a)}}{=} \eta^{(1)} \prn{ \muzero[\vq]{\vx} } \prn[\big]{ \muhat[\vq][(\ve_k)]{\vx} - \muzero[\vq][(\ve_k)]{\vx} }\\
  &\qquad + \prn[\big]{ \eta^{(1)} \prn{ \muhat[\vq]{\vx} } - \eta^{(1)} \prn{ \muzero[\vq]{\vx} } } \muhat[\vq][(\ve_k)]{\vx}\\
  &\quad \overset{\text{(b)}}{=} \eta^{(1)} \prn{ \muzero[\vq]{\vx} } \prn[\big]{ \muhat[\vq][(\ve_k)]{\vx} - \muzero[\vq][(\ve_k)]{\vx} } + \eta^{(2)}(\zeta) \prn{ \muhat[\vq]{\vx} - \muzero[\vq]{\vx} } \muhat[\vq][(\ve_k)]{\vx}\\
  &\quad \overset{\text{(c)}}{=} \eta^{(1)} \prn{ \muzero[\vq]{\vx} } \bart[\vq][\ve_k] \sqrt{\bigomegabar[\vq][\ve_k] / n}\\
  &\qquad + \eta^{(2)}(\zeta) \prn[\big]{ \muhat[\vq]{\vx} - \muzero[\vq]{\vx} } \muhat[\vq][(\ve_k)]{\vx} + O_{\Problet} \prn{ h^{-1} \rbah }\\
  &\quad \overset{\text{(d)}}{=} \eta^{(1)} \prn{ \muzero[\vq]{\vx} } \bart[\vq][\ve_k] \sqrt{\bigomegabar[\vq][\ve_k] / n}\\
  &\qquad + \eta^{(2)}(\zeta) \prn[\big]{ \muhat[\vq]{\vx} - \muzero[\vq]{\vx} } \prn[\big]{ \muhat[\vq][(\ve_k)]{\vx} - \muzero[\vq][(\ve_k)]{\vx} }\\
  &\qquad + \eta^{(2)}(\zeta) \prn[\big]{ \muhat[\vq]{\vx} - \muzero[\vq]{\vx} } \muzero[\vq][(\ve_k)]{\vx} + O_{\Problet} \prn{ h^{-1} \rbah }\\
  &\quad \overset{\text{(e)}}{=} \eta^{(1)} \prn{ \muzero[\vq]{\vx} } \bart[\vq][\ve_k] \sqrt{\bigomegabar[\vq][\ve_k] / n} + O_{\Problet} \prn[\big]{h^{-1} (\ruc[2] + \rbah) + \ruc}\\
  &\quad \overset{\text{(f)}}{=} \eta^{(1)} \prn{ \muzero[\vq]{\vx} } \sqrt{\bigomegabar[\vq][\ve_k] / n} \brk[\Big]{\bart[\vq][\ve_k] + O_{\Problet} \prn[\Big]{h \sqrt{n h^d} \prn[\big]{h^{-1} (\ruc[2] + \rbah) + \ruc}}}.
\end{align*}
Here, (a) is just rewriting. In (b) we used the mean-value theorem, with some $\zeta = \zeta_{\vq, \vx}$ between $\muhat[\vq]{\vx}$ and $\muzero[\vq]{\vx}$. In (c) we used \eqref{eq:TvfoLC} with $\vv = \ve_k$ and that $\eta^{(1)}(\muzero[\vq]{\vx})$ is bounded uniformly over $\vq, \vx$. (d) is just rewriting. In (e) we used $\ruc = o(1)$ (giving that $\eta^{(2)}(\zeta)$ does not exceed a fixed constant not depending on $\vq$ or $\vx$), \eqref{eq:gRcohQ} with $\vv = \vzero$ and $\vv = \ve_k$, and uniform boundedness of first partial derivatives of $\muzero[\vq]{\cdot}$. (f) is by $h^{-2 \multindnorm{\vv} - d} \lesssim_{\Problet} \inf_{\vq, \vx} |\bigomegabar[\vq]|$ (by \cref{lem:asympt-variance}) and since $\big| \eta^{(1)} (\muzero[\vq]{\vx}) \big|$ is bounded away from zero by \cref{ass:gl-holder}\ref{ass:psi-eta}.

Rewriting, we obtain
\begin{align*}
  &\sup_{\vq, \vx} \abs[\Bigg]{ \frac{\eta^{(1)} \prn{ \muhat[\vq]{\vx} } \muhat[\vq][(\ve_k)]{\vx} - \eta^{(1)} \prn{ \muzero[\vq]{\vx} } \muzero[\vq][(\ve_k)]{\vx}}{\eta^{(1)} \prn{ \muzero[\vq]{\vx} } \sqrt{\bigomegabar[\vq][\ve_k]}} - \bart[\vq][\ve_k] }\\
  &\quad \lesssim_{\Problet} \sqrt{n h^d} \prn[\big]{ \ruc[2] + \rbah + h \ruc }.
\end{align*}
It is left to combine this with \eqref{eq:yurinski-corollary-statement} and use the triangle inequality.

\ref{other-params-three} The argument is the same as for Parts (i) and (ii) with $\rstr$ replaced by $\rstrsubgauss$.
\end{proof}

\section{Simulation Evidence}\label{sec:simulation-evidence}

We conducted a small simulation experiment to demonstrate the finite sample properties of the partitioning-based $M$-estimation methodology. We studied pointwise and uniform (over $q \in \mathcal{Q}$, $\vx \in \mathcal{X}$) estimation and inference for the conditional distribution function (Example \mainref{example: Conditional Distribution Regression} in the paper, and discussed in Section \ref{sec:example-distribution-regression} above).

We considered four data generating processes. Models 1 through 3 were chosen to satisfy $y_i = m(\vx_i) + \epsilon_i$, where $\vx_i \thicksim \text{Uniform}(\mathcal{X})$ with $\mathcal{X} = [0,1]^d$, $\epsilon_i \thicksim \text{Normal}(0,1)$, $\vx_i \independent \epsilon_i$, $m(\vx)$ is defined in \cref{tab:m-def}. For Model 4, we took the treatment effect model from \cite{donald2014estimation} (with a slight change so that the support of $\vx_i$ is also $[0, 1]$):
\begin{equation*}
y_i = T_i \brk[\bigg]{\inde{V_i \leq 1 - x_i} \frac{V_i^2}{1 - x_i} + \inde{V_i > 1 - x_i} V_i} + (1 - T_i) \brk[\bigg]{\inde{U_i \leq x_i} \frac{U_i^2}{x_i} + \inde{U_i > x_i} U_i},
\end{equation*}
where $T_i = \inde{U^T_i < x_i}$, and $x_i, U_i^T, U_i, V_i \sim \text{Uniform}[0, 1]$ are independent. The estimand is $\muzero[q]{\vx}$, where
\begin{align*}
  \eta(\muzero[q]{\vx}) = F(q | \vx) = \Prob{y_i \leq q \given \vx_i = \vx}
\end{align*}
with $q \in \mathcal{Q} = [-0.2, 0.2]$ for Models 1 through 3 and $\mathcal{Q} = [0.2, 0.8]$ for Model 4; the link function $\eta(\cdot)$ is the complementary log-log link: $\eta(t) = 1 - e^{- e^t}$.

\begin{table}[h!]
  \caption{Definition of the $m(\vx)$ function.}
  \label{tab:m-def}
  \centering
  \begin{tabular}{lcc}
    \toprule
    & $d = 1$ & $d = 2$\\
    \midrule
    Model 1 & $\sin(2 \pi x) / 2$ & $\prn{\sin(2 \pi x_1) + \sin(2 \pi x_2)} / 4$ \\
    Model 2 & $\sin(5 x) \sin(10 x) / 2$ & $\sin(5 x_1) \sin(10 x_2) / 2$\\
    Model 3 & $\prn{1 - (4 x - 2)^2}^2 / 4$ & $\prn{1 - (4 x_1 - 2)^2}^2 \sin(5 x_2) / 4$\\
    \bottomrule
  \end{tabular}
\end{table}

We use the tensor-product B-spline basis $\vp(\vx)$ on $[0, 1]^d$ as defined in \cref{sec:examples-basis-functions} with equally spaced knots. (Recall that the number of knots on each segment is $J - 1$, and thus $\kappa = J^d$.) The estimator $\muhat[q]{\vx} = \vp(\vx)\trans \betahat[q]$ is constructed by solving the optimization problem
\begin{equation*}
\betahat[q] \in \argmin_{\bb \in \mathbb{R}^K} \sum_{i=1}^n \riskfunc[q]{y_i}{\eta(\bp(\vx_i)\trans \bb)}
\end{equation*}
for each $q$ of interest separately (with $\rho(\cdot)$ defined as in \cref{prop:gl-holder-distribution-regression}), using the L-BFGS-B algorithm initialized at $\boldsymbol{0}$. The matrices $\qbar[q]$ and $\sigmabar[q]$ were estimated using, respectively,
\begin{align*}
  \qhat[q] = 2 \Expectlet_n \brk{\vp(\vx_i) \vp(\vx_i) \trans \eta^{(1)}(\muhat[q]{\vx_i})^{2}}
\end{align*}
and
\begin{align*}
    \sigmahat[q] = 4 \Expectlet_n \brk{\vp(\vx_i) \vp(\vx_i) \trans \widehat{F}_{Y|X}(q|\vx_i) \prn{ 1 - \widehat{F}_{Y|X}(q|\vx_i) } \eta^{(1)}(\muhat[q]{\vx_i})^{2}},
\end{align*}
with $\widehat{F}_{Y|X}(q|\vx) = \eta(\muhat[q]{\vx})$. Thus, the variance estimator is $\bigomegahat[q][][\vx] = \vp(\vx)\trans \qhat[q][-1] \sigmahat[q] \qhat[q][-1] \vp(\vx)$.

The simulation experiment considered two sample sizes ($n = 5\,000$ and $n = 50\,000$), two dimensions ($d=1$ and $d=2$), and $1\,000$ replications. The pointwise estimation and inference results are presented in \cref{tab:dim1-point-small,tab:dim2-point-small} ($n = 5\,000$) and \cref{tab:dim1-point,tab:dim2-point} ($n = 50\,000$). For each model, we consider two choices of $J$ and three evaluation points for $(q, \vx)$. The root mean squared error (RMSE) for point estimators, the coverage rates and the average widths of pointwise 95\% nominal confidence
intervals (CIs) are reported.
The uniform estimation and inference results are presented in \cref{tab:dim1-uniform-small,tab:dim2-uniform-small} ($n = 5\,000$) and \cref{tab:dim1-uniform,tab:dim2-uniform} ($n = 50\,000$). We consider two choices of $J$, and use a discrete grid of $8$ equally spaced points in place of the continuous segment $\mathcal{Q}$ and $10$ (for $d = 1$) or $5 \times 5 = 25$ (for $d = 2$) points in place of $\mathcal{X}$. The uniform $95\%$ nominal confidence bands (CBs) are constructed using the plug-in approach as in \cref{sec:plug-in}, and the covariance structure of the discrete analogue of $\znhat(q)$ is obtained as in \cref{rem:feasibility-distr-reg}. The maximum estimation error, the uniform coverage rates and the average widths of confidence bands are reported.

\begin{table}[h!]
\caption{Pointwise simulation results for points $\crl{(q^{(k)}, x^{(k)})}_{k = 1}^3$, averaged across $1\,000$ replications with $n = 5\,000$, $d = 1$, where $(q^{(1)}, q^{(2)}, q^{(3)}) = (-0.2, 0, 0.2)$ for Models 1 through 3, $(q^{(1)}, q^{(2)}, q^{(3)}) = (0.45, 0.6, 0.75)$ for Model 4; $(x^{(1)}, x^{(2)}, x^{(3)}) = (0.3, 0.1, 0.2)$.}
\label{tab:dim1-point-small}
\centering

\py{merge_dicts([3, 5], False)}
\setlength{\tabcolsep}{5.75pt}
\begin{tabular}{@{}c c ccc ccc ccc@{}}
\toprule
& & \multicolumn{3}{c}{RMSE} & \multicolumn{3}{c}{Coverage} & \multicolumn{3}{c}{CI width} \\
\cmidrule(lr){3-5}\cmidrule(lr){6-8}\cmidrule(l){9-11}
Model & $J$ & Point 1 & Point 2 & Point 3 & Point 1 & Point 2 & Point 3 & Point 1 & Point 2 & Point 3 \\
\midrule
\addlinespace
\multirow{2}{*}{1} & 4 & 0.088 & 0.074 & 0.074 & 0.952 & 0.953 & 0.944 & 0.336 & 0.289 & 0.287 \\
 & 6 & 0.110 & 0.077 & 0.088 & 0.946 & 0.951 & 0.944 & 0.420 & 0.303 & 0.339 \\
\addlinespace
\multirow{2}{*}{2} & 4 & 0.069 & 0.071 & 0.070 & 0.953 & 0.947 & 0.946 & 0.270 & 0.277 & 0.274 \\
 & 6 & 0.088 & 0.076 & 0.085 & 0.950 & 0.942 & 0.948 & 0.343 & 0.291 & 0.321 \\
\addlinespace
\multirow{2}{*}{3} & 4 & 0.073 & 0.099 & 0.068 & 0.953 & 0.943 & 0.945 & 0.279 & 0.370 & 0.258 \\
 & 6 & 0.084 & 0.110 & 0.072 & 0.952 & 0.956 & 0.954 & 0.337 & 0.435 & 0.292 \\
\addlinespace
\multirow{2}{*}{4} & 4 & 0.065 & 0.058 & 0.060 & 0.952 & 0.961 & 0.942 & 0.252 & 0.237 & 0.225 \\
 & 6 & 0.078 & 0.063 & 0.068 & 0.963 & 0.953 & 0.952 & 0.308 & 0.246 & 0.263 \\
\bottomrule
\end{tabular}

\end{table}

\begin{table}[h!]
\caption{Uniform simulation results, averaged across $1\,000$ replications with $n = 5\,000$, $d = 1$.}
\label{tab:dim1-uniform-small}
\centering

\py{merge_dicts([3, 5], True)}
\begin{tabular}{c c c c c c}
\toprule
Model & $J$ & $\sup_{q, {x}} \abs{\muhat[q]{x} - \muzero[q]{x}}$ & Uniform coverage & Av. CB width\\
\midrule
\addlinespace
\multirow{2}{*}{1} & 4 & 0.151 & 0.950 & 0.399 \\
 & 6 & 0.181 & 0.948 & 0.475 \\
\addlinespace
\multirow{2}{*}{2} & 4 & 0.146 & 0.939 & 0.381 \\
 & 6 & 0.170 & 0.940 & 0.453 \\
\addlinespace
\multirow{2}{*}{3} & 4 & 0.187 & 0.922 & 0.441 \\
 & 6 & 0.223 & 0.952 & 0.550 \\
\addlinespace
\multirow{2}{*}{4} & 4 & 0.173 & 0.953 & 0.415 \\
 & 6 & 0.201 & 0.947 & 0.494 \\
\bottomrule
\end{tabular}

\end{table}

\begin{table}[h!]
\caption{Pointwise simulation results for points $\crl{(q^{(k)}, \vx^{(k)})}_{k = 1}^3$, averaged across $1\,000$ replications with $n = 5\,000$, $d = 2$, where $(q^{(1)}, q^{(2)}, q^{(3)}) = (-0.2, 0, 0.2)$, $\vx^{(1)} = (0.3, 0.1)$, $\vx^{(2)} = (0.1, 0.4)$ and $\vx^{(3)} = (0.2, 0.2)$.}
\label{tab:dim2-point-small}
\centering

\py{merge_dicts(4, False)}
\setlength{\tabcolsep}{5.75pt}
\begin{tabular}{@{}c c ccc ccc ccc@{}}
\toprule
& & \multicolumn{3}{c}{RMSE} & \multicolumn{3}{c}{Coverage} & \multicolumn{3}{c}{CI width} \\
\cmidrule(lr){3-5}\cmidrule(lr){6-8}\cmidrule(l){9-11}
Model & $J$ & Point 1 & Point 2 & Point 3 & Point 1 & Point 2 & Point 3 & Point 1 & Point 2 & Point 3 \\
\midrule
\addlinespace
\multirow{2}{*}{1} & 4 & 0.272 & 0.194 & 0.243 & 0.957 & 0.956 & 0.954 & 1.052 & 0.768 & 0.949 \\
 & 6 & 0.455 & 0.256 & 0.366 & 0.944 & 0.957 & 0.935 & 1.417 & 1.027 & 1.328 \\
\addlinespace
\multirow{2}{*}{2} & 4 & 0.273 & 0.182 & 0.237 & 0.957 & 0.934 & 0.941 & 1.059 & 0.662 & 0.904 \\
 & 6 & 0.425 & 0.231 & 0.318 & 0.950 & 0.939 & 0.955 & 1.422 & 0.855 & 1.248 \\
\addlinespace
\multirow{2}{*}{3} & 4 & 0.244 & 0.260 & 0.219 & 0.926 & 0.936 & 0.950 & 0.896 & 0.944 & 0.839 \\
 & 6 & 0.314 & 0.784 & 0.290 & 0.940 & 0.932 & 0.950 & 1.154 & 1.558 & 1.129 \\
\bottomrule
\end{tabular}

\end{table}

\begin{table}[h!]
\caption{Uniform simulation results, averaged across $1\,000$ replications with $n = 5\,000$, $d = 2$.}
\label{tab:dim2-uniform-small}
\centering

\py{merge_dicts(4, True)}
\begin{tabular}{c c c c c c}
\toprule
Model & $J$ & $\sup_{q, {\vx}} \abs{\muhat[q]{\vx} - \muzero[q]{\vx}}$ & Uniform coverage & Av. CB width\\
\midrule
\addlinespace
\multirow{2}{*}{1} & 4 & 0.487 & 0.940 & 1.163 \\
 & 6 & 0.730 & 0.950 & 1.593 \\
\addlinespace
\multirow{2}{*}{2} & 4 & 0.491 & 0.935 & 1.173 \\
 & 6 & 0.732 & 0.944 & 1.608 \\
\addlinespace
\multirow{2}{*}{3} & 4 & 0.478 & 0.950 & 1.172 \\
 & 6 & 0.806 & 0.921 & 1.656 \\
\bottomrule
\end{tabular}

\end{table}

\begin{table}[h!]
\caption{Pointwise simulation results for points $\crl{(q^{(k)}, x^{(k)})}_{k = 1}^3$, averaged across $1\,000$ replications with $n = 50\,000$, $d = 1$, where $(q^{(1)}, q^{(2)}, q^{(3)}) = (-0.2, 0, 0.2)$ for Models 1 through 3, $(q^{(1)}, q^{(2)}, q^{(3)}) = (0.45, 0.6, 0.75)$ for Model 4; $(x^{(1)}, x^{(2)}, x^{(3)}) = (0.3, 0.1, 0.2)$.}
\label{tab:dim1-point}
\centering

\py{merge_dicts([0, 2], False)}
\setlength{\tabcolsep}{5.75pt}
\begin{tabular}{@{}c c ccc ccc ccc@{}}
\toprule
& & \multicolumn{3}{c}{RMSE} & \multicolumn{3}{c}{Coverage} & \multicolumn{3}{c}{CI width} \\
\cmidrule(lr){3-5}\cmidrule(lr){6-8}\cmidrule(l){9-11}
Model & $J$ & Point 1 & Point 2 & Point 3 & Point 1 & Point 2 & Point 3 & Point 1 & Point 2 & Point 3 \\
\midrule
\addlinespace
\multirow{2}{*}{1} & 4 & 0.027 & 0.024 & 0.023 & 0.947 & 0.950 & 0.956 & 0.106 & 0.091 & 0.091 \\
 & 6 & 0.034 & 0.025 & 0.027 & 0.949 & 0.950 & 0.957 & 0.132 & 0.096 & 0.107 \\
\addlinespace
\multirow{2}{*}{2} & 4 & 0.026 & 0.023 & 0.022 & 0.906 & 0.940 & 0.943 & 0.085 & 0.088 & 0.086 \\
 & 6 & 0.028 & 0.023 & 0.026 & 0.942 & 0.963 & 0.947 & 0.108 & 0.092 & 0.101 \\
\addlinespace
\multirow{2}{*}{3} & 4 & 0.023 & 0.035 & 0.021 & 0.946 & 0.904 & 0.960 & 0.088 & 0.117 & 0.081 \\
 & 6 & 0.028 & 0.035 & 0.024 & 0.946 & 0.942 & 0.959 & 0.107 & 0.136 & 0.092 \\
\addlinespace
\multirow{2}{*}{4} & 4 & 0.020 & 0.019 & 0.019 & 0.955 & 0.947 & 0.943 & 0.080 & 0.075 & 0.071 \\
 & 6 & 0.024 & 0.019 & 0.022 & 0.954 & 0.951 & 0.950 & 0.097 & 0.078 & 0.083 \\
\bottomrule
\end{tabular}

\end{table}

\begin{table}[h!]
\caption{Uniform simulation results, averaged across $1\,000$ replications with $n = 50\,000$, $d = 1$.}
\label{tab:dim1-uniform}
\centering

\py{merge_dicts([0, 2], True)}
\begin{tabular}{c c c c c c}
\toprule
Model & $J$ & $\sup_{q, {x}} \abs{\muhat[q]{x} - \muzero[q]{x}}$ & Uniform coverage & Av. CB width\\
\midrule
\addlinespace
\multirow{2}{*}{1} & 4 & 0.047 & 0.943 & 0.126 \\
 & 6 & 0.057 & 0.953 & 0.150 \\
\addlinespace
\multirow{2}{*}{2} & 4 & 0.054 & 0.697 & 0.120 \\
 & 6 & 0.053 & 0.944 & 0.143 \\
\addlinespace
\multirow{2}{*}{3} & 4 & 0.066 & 0.871 & 0.139 \\
 & 6 & 0.070 & 0.941 & 0.173 \\
\addlinespace
\multirow{2}{*}{4} & 4 & 0.055 & 0.938 & 0.131 \\
 & 6 & 0.063 & 0.943 & 0.156 \\
\bottomrule
\end{tabular}

\end{table}

\begin{table}[h!]
\caption{Pointwise simulation results for points $\crl{(q^{(k)}, \vx^{(k)})}_{k = 1}^3$, averaged across $1\,000$ replications with $n = 50\,000$, $d = 2$, where $(q^{(1)}, q^{(2)}, q^{(3)}) = (-0.2, 0, 0.2)$, $\vx^{(1)} = (0.3, 0.1)$, $\vx^{(2)} = (0.1, 0.4)$ and $\vx^{(3)} = (0.2, 0.2)$.}
\label{tab:dim2-point}
\centering

\py{merge_dicts(1, False)}
\setlength{\tabcolsep}{5.75pt}
\begin{tabular}{@{}c c ccc ccc ccc@{}}
\toprule
& & \multicolumn{3}{c}{RMSE} & \multicolumn{3}{c}{Coverage} & \multicolumn{3}{c}{CI width} \\
\cmidrule(lr){3-5}\cmidrule(lr){6-8}\cmidrule(l){9-11}
Model & $J$ & Point 1 & Point 2 & Point 3 & Point 1 & Point 2 & Point 3 & Point 1 & Point 2 & Point 3 \\
\midrule
\addlinespace
\multirow{2}{*}{1} & 4 & 0.085 & 0.064 & 0.075 & 0.947 & 0.940 & 0.952 & 0.326 & 0.239 & 0.295 \\
 & 6 & 0.108 & 0.084 & 0.102 & 0.952 & 0.941 & 0.951 & 0.423 & 0.313 & 0.407 \\
\addlinespace
\multirow{2}{*}{2} & 4 & 0.083 & 0.053 & 0.074 & 0.952 & 0.945 & 0.943 & 0.328 & 0.206 & 0.281 \\
 & 6 & 0.107 & 0.067 & 0.095 & 0.951 & 0.949 & 0.959 & 0.424 & 0.263 & 0.385 \\
\addlinespace
\multirow{2}{*}{3} & 4 & 0.072 & 0.076 & 0.068 & 0.953 & 0.941 & 0.944 & 0.278 & 0.291 & 0.262 \\
 & 6 & 0.090 & 0.107 & 0.088 & 0.952 & 0.948 & 0.952 & 0.352 & 0.419 & 0.350 \\
\bottomrule
\end{tabular}

\end{table}

\begin{table}[h!]
\caption{Uniform simulation results, averaged across $1\,000$ replications with $n = 50\,000$, $d = 2$.}
\label{tab:dim2-uniform}
\centering

\py{merge_dicts(1, True)}
\begin{tabular}{c c c c c c}
\toprule
Model & $J$ & $\sup_{q, {\vx}} \abs{\muhat[q]{\vx} - \muzero[q]{\vx}}$ & Uniform coverage & Av. CB width\\
\midrule
\addlinespace
\multirow{2}{*}{1} & 4 & 0.146 & 0.953 & 0.363 \\
 & 6 & 0.193 & 0.960 & 0.484 \\
\addlinespace
\multirow{2}{*}{2} & 4 & 0.146 & 0.958 & 0.366 \\
 & 6 & 0.198 & 0.944 & 0.489 \\
\addlinespace
\multirow{2}{*}{3} & 4 & 0.145 & 0.930 & 0.365 \\
 & 6 & 0.198 & 0.950 & 0.494 \\
\bottomrule
\end{tabular}
\end{table}


\clearpage
\section*{References} 
\begingroup
\renewcommand{\section}[2]{}	
\bibliography{bib}
\bibliographystyle{jasa}
\endgroup